\def\real{{\tt I\kern-.2em{R}}}  %This appears to be the corect archive
\def\nat{{\tt I\kern-.2em{N}}}    %version. as of file date. 
\def\eps{\epsilon}
\def\realp#1{{\tt I\kern-.2em{R}}^#1}
\def\natp#1{{\tt I\kern-.2em{N}}^#1}
\def\hyper#1{\,^*\kern-.2em{#1}}
\def\monad#1{\mu (#1)}
\def\dmonad#1{\mu^\prime(#1)}
\def\nmonad#1{\mu^n(#1)}
\def\mmonad#1{\mu^m(#1)}
\def\St#1{{\tt st}#1}
\def\st#1{{\tt st}(#1)}
\def\hyperreal{{^*{\real}}}
\def\hyperrealp#1{{\tt ^*{I\kern-.2em{R}}}^#1} 
\def\hypernat{{^*{\nat }}}
\def\hypernatp#1{{{^*{{\tt I\kern-.2em{N}}}}}^#1} 
\def\eskip{\hskip.25em\relax}

\def\Hyper#1{\hyper {\eskip #1}}
\def\leaderfill{\leaders\hbox to 1em{\hss.\hss}\hfill}
\def\srealp#1{{\rm I\kern-.2em{R}}^#1}

\def\power#1{{{\cal P}(#1)}}
\def\iff{\leftrightarrow}

\def\pars{\par\smallskip}
\def\parm{\par\medskip}

\def\ref#1{$^{#1}$}

\def\m@th{\mathsurround=0pt}
\def\rightarrowfill{$\m@th \mathord- \mkern-6mu \cleaders\hbox{$\mkern-2mu 
\mathord- \mkern-2mu$}\hfil \mkern-6mu \mathord\rightarrow$}
\def\leftarrowfill{$\mathord\leftarrow
\mkern -6mu \m@th \mathord- \mkern-6mu \cleaders\hbox{$\mkern-2mu 
\mathord- \mkern-2mu$}\hfil $}
\def\noarrowfill{$\m@th \mathord- \mkern-6mu \cleaders\hbox{$\mkern-2mu 
\mathord- \mkern-2mu$}\hfil$}
\def\orgate{$\bigcirc \kern-.80em \lor$}
\def\andgate{$\bigcirc \kern-.80em \land$}
\def\inverter{$\bigcirc \kern-.80em \neg$}
 \magnification=1000
\tolerance 10000
\hoffset=0.25in
\hsize 6.00 true in
\vsize 8.75 true in
\font\eightsl=cmsl8
\pageno=1
\headline={\ifnum\pageno= 1\hfil {\quad}
\hfil\else\ifodd\pageno\rightheadline \else\leftheadline\fi\fi} 
\def\rightheadline{\eightsl \hfil Infinitesimal Modeling\hfil} 
\def\leftheadline{{} \eightsl \hfil  Infinitesimal Modeling
 \hfil} 
\voffset=1\baselineskip

\settabs 3 \columns
\baselineskip 14pt
\pageno=1
\centerline{}
\vskip 1.45in
\centerline{\bf Nonstandard Analysis Applied to}
\centerline{\bf Advanced Undergraduate Mathematics}
\centerline{$\diamondsuit$ Infinitesimal Modeling and Very Elementary Physics$\diamondsuit$}
\vskip .75in
\centerline{\bf Robert A. Herrmann}
\vskip 4.5in
\centerline{A (July - August) 1989, 1990, 1991 Instructional Development Project from the}
\centerline{Mathematics Department}
\centerline{United States Naval Academy}
\centerline{572C Holloway Road}
\centerline{Annapolis$,$ Maryland 21402-5002}
\eject
\ \ \par
\vskip 2.5in
\centerline{IMPORTANT NOTICE}\pars
\indent Since the writing of this book was financed entirely by a designated 
grant from the Federal Government that was specifically obtained for this 
sole purpose then a copyright for this specific book cannot be obtained by 
its author. 
Any portion of its contents can be copied and used without seeking permissions 
from the author. However$,$ when such copying or use is made of this material$,$ 
it is 
necessary that the author and the U. S. Naval Academy be indicated as 
the source of the  material being used. All typographical error have NOT be corrected.  
Further note that certain new results that appear in this book 
will be published under the author's name in scholarly 
journals.\footnote*{This publication process is now being instituted with 
certain results published as of January 1992$,$ accepted for publication or
in preparation. Further$,$ some typographical errors have been corrected as of this 17 June 1997 version.}  
\par\vfil\eject

\centerline{}
\centerline{\bf CONTENTS}
\bigskip
\indent Chapter 1\par
{\bf Introduction}\par
\line{\indent 1.1 \ \ A Brief History \leaderfill 4}
\line{\indent 1.2 \ \ Manual construction \leaderfill 8}
\bigskip
\indent Chapter 2\par
{\bf Infinitesimals$,$ Limited And Infinite Numbers}\par
\line{\indent 2.1 \ \ Some Notation and Definitions\leaderfill 10}
\line{\indent 2.2 \ \ Basic Algebra\leaderfill 12}
\line{\indent 2.3 \ \ Euclidean N-spaces\leaderfill 13}
\line{\indent 2.4 \ \ The Standard Part Operator\leaderfill 13}
\line{\indent 2.5 \ \ A Slight Quandary\leaderfill 15}
\bigskip
\indent Chapter 3\par
{\bf Some Set Theory$,$ Convergence And Leibniz's Principle}\par
\line{\indent 3.1 \ \ Some Set Theory\leaderfill 16} 
\line{\indent 3.2 \ \ Convergence\leaderfill 18}
\line{\indent 3.3 \ \ Constants and More Constants\leaderfill 19}
\line{\indent 3.4 \ \ The Leibniz Principle of *-transfer\leaderfill 21}
\line{\indent 3.5 \ \ A Few Simple Applications\leaderfill 23}
\bigskip
\indent Chapter 4\par
{\bf Some Modeling With The Infinite Numbers}\par
\line{\indent 4.1 \ \ Historical Confusion\leaderfill 26}
\line{\indent 4.2 \ \ The Internal Definition Principle\leaderfill 27}
\line{\indent 4.3 \ \ Hyperfinite Summation\leaderfill 30}
\line{\indent 4.4 \ \ Continuity and a Few Examples\leaderfill 32}
\bigskip 
\indent Chapter 5\par
{\bf Standard Rules For Integral Modeling}\par
\line{\indent 5.1 \ \ The Riemann Styled Integral\leaderfill 35}
\line{\indent 5.2 \ \ The Infinite Sum Theorems\leaderfill 37}
\line{\indent 5.3 \ \ Extensions\leaderfill 39}
\line{\indent 5.4 \ \ Applications and The Standard Modeling Rules\leaderfill 
54}
\line{\indent 5.5 \ \ Extensions of the Standard Rules\leaderfill 39}
\bigskip
\indent Chapter 6\par
{\bf Nonstandard Rules For Integral Modeling}\par
\line{\indent 6.1 \ \ Historical Examples\leaderfill 44}
\line{\indent 6.2 \ \ The Monadic Environment\leaderfill 46}
\line{\indent 6.3 \ \ Simple Applications\leaderfill 49}
\line{\indent 6.4 \ \ The Method of Constants\leaderfill 50}
\line{\indent 6.5 \ \ The Hyperfinite Method\leaderfill 54}
\line{\indent 6.6 \ \ Instruction\leaderfill 55}
\line{\indent 6.7 \ \ Realism\leaderfill 56}
\bigskip
\indent Chapter 7\par
{\bf Pure Infinitesimal Integral Modeling}\par
\line{\indent 7.1 \ \ Brief Discussion\leaderfill 58}
\line{\indent 7.2 \ \ Geometric Elements\leaderfill 59}
\bigskip
\indent Chapter 8\par
{\bf Refinements For Integral Modeling}\par
\line{\indent 8.1 \ \ A Very General Approach\leaderfill 65}
\line{\indent 8.2 \ \ The Line Integral\leaderfill 65}
\line{\indent 8.3 \ \ Order Ideals and Approximations\leaderfill 68}
\line{\indent 8.4 \ \ {\it n}th Order Increments\leaderfill 70}
\line{\indent 8.5 \ \ Microgeometry - Tangents to Curves\leaderfill 71}
\line{\indent 8.6 \ \ Microgeometry - Surface Elements\leaderfill 72}
\line{\indent 8.7 \ \ Microgeometry - Other Stuff\leaderfill 74}
\line{\indent 8.9 \ \ Gauge Integrals\leaderfill 75}
\bigskip    
\indent Chapter 9\par
{\bf The Beginnings Of Differential Equation Modeling}\par
\line{\indent 9.1 \ \ Brief Discussion\leaderfill 76}
\line{\indent 9.2 \ \ The Limit\leaderfill 77}
\line{\indent 9.3 \ \ Fluxions and Dynamic Geometry\leaderfill 77}
\line{\indent 9.4 \ \ Fluxions and Higher Order Infinitesimals\leaderfill 81}
\line{\indent 9.5 \ \ What is a Tangent?\leaderfill 82}
\line{\indent 9.6 \ \ What is an (k-surface) Osculating Plane?\leaderfill 84}
\line{\indent 9.7 \ \ Curvature\leaderfill 86}
\bigskip
\indent Chapter 10 \par
{\bf The Differential And Physical Modeling}\par
\line{\indent 10.1 \ Basic Properties\leaderfill 88}
\line{\indent 10.2 \ Some General Observations\leaderfill 91}
\line{\indent 10.3 \ Vibrating Membrane\leaderfill 92}
\line{\indent 10.4 \ Internal Heat Transfer\leaderfill 98}
\line{\indent 10.5 \ Concluding Remarks\leaderfill 101}
\bigskip \par
\indent {\bf Appendices}\par
\line{\indent For Chapter 2\leaderfill 103}
\line{\indent For Chapter 3\leaderfill 106}
\line{\indent For Chapter 4\leaderfill 115}
\line{\indent For Chapter 5\leaderfill 120}
\line{\indent For Chapter 6\leaderfill 127}
\line{\indent For Chapter 7\leaderfill 136}
\line{\indent For Chapter 8\leaderfill 137} 
\line{\indent For Chapter 10\leaderfill 133} 
\bigskip\par
\line{\indent{\bf References}\leaderfill 149}
\bigskip \par
\line{\indent{\bf Special Symbols}\leaderfill 151}
\bigskip\par 
\line{\indent{\bf Very Elementary Physics}\leaderfill 154}
\bigskip\par
\line{\indent{\bf Index for pages 1-153}\leaderfill 183}
\bigskip\par
\indent {\bf Applications}\par
\line{\indent Length Of A Continuous Curve (Ex. 4.4.4.A)\leaderfill 33}
\line{\indent Fractals (Ex. 4.4.2)\leaderfill 33}
\line{\indent The 2-dimensional Area Between Two Curves\leaderfill 40}
\line{\indent Volume Obtained By A 2-dimensional Integral\leaderfill 41}
\line{\indent Mass Obtained By 3-dimensional Integral\leaderfill 41}
\line{\indent Jordan-Volume Obtained By A 2-dimensional Integral\leaderfill 
42}
\line{\indent Incompressible Fluid\leaderfill 49}
\line{\indent Moment Of Inertia\leaderfill 50}
\line{\indent Incompressible Fluid Second Derivation\leaderfill 52}
\line{\indent Moment Of Inertia Second Derivation\leaderfill 53}
\line{\indent The 2-dimensional Area Between Two Curves (Single 
Integral)\leaderfill 59}
\line{\indent Volume Of Revolution (Single Integral)\leaderfill 61}
\line{\indent Surface Area of Revolution (Single Integral)\leaderfill 61}
\line{\indent Volume Obtained By 2-dimensional Integral\leaderfill 63}
\line{\indent Value Of An Electric Field Vector\leaderfill 64}
\line{\indent Energy Expended Within a Force Field \leaderfill 67}
\line{\indent Tangents To Curves\leaderfill 83}
\line{\indent The Osculating Plane\leaderfill 84}
\line{\indent Curvature\leaderfill 86}
\line{\indent Vibrating Membrane \leaderfill 93}
\line{\indent Internal Heat Transfer \leaderfill 98}
\vfil\eject
\centerline{Chapter 1.}
\medskip
\centerline{\bf INTRODUCTION}\par
\bigskip
\leftline{1.1 \underbar{A Brief History.}}\par             
\medskip 
\indent Scientists who use mathematical analysis as a tool have traditionally relied 
upon a vague process called ``infinitesimal reasoning'' - a process that from 
the time of {Archimedes} until 1961 had no fixed rules nor consistent language.  
However$,$ application of this intuitive process is the exact cause that has led 
to our great analytical successes both in scientific and engineering 
endeavors. Unfortunately$,$ it also led to great controversy.\par 
     Beginning in about 1600 a schism developed between some mathematicians 
and the foremost appliers of this analytical tool. {Leibniz} approved entirely 
of the concept of the infinitely small or infinitely large numbers but stated 
that they should be treated as ``{ideal}'' elements rather than real numbers. He 
also believed that they should be governed by the same laws that then 
controlled the behavior of the ordinary numbers. He claimed$,$ but could not 
justify the assertion$,$ that all arguments involving such ideal numbers could 
be replaced by arguing in terms of objects that are large enough or small 
enough to make error as small as one wished. {De l'Hospital} [1715] when he 
wrote the first Calculus textbook used the terminology exclusively and 
utilized a formal ``definition - axiom'' process supposedly delineating the 
notion of the infinitesimal. Unfortunately$,$ his first axiom is logically 
contradictory. {D'Alembert} insisted that the Leibniz concepts were without 
merit and only a process using a modified ``limit'' idea was appropriate. 
 {Euler} contended in opposition to D'Alembert that the Leibniz approach was the 
best that could be achieved and fought diligently for the acceptance of these 
ideal numbers.\par

 Due to what appeared to be logical 
 {inconsistencies} within the methods$,$ those mathematicians trained in classical 
logic began to demand that applied mathematicians produce ``proofs'' of their 
derivations. In answer to this criticism {Kepler} wrote$,$ {\it ``We could obtain 
absolute and in all respects perfect demonstrations from the books of 
Archimedes themselves$,$ were we not repelled by the thorny reading thereof.''} 
The successes of these vague methods and those scientists and mathematicians 
such as Leibniz$,$ Euler and Gauss who championed their continued use 
quieted the ``unbelievers.'' It should be noted that the concern of the 
critics was based upon the fact that they used the same vague processes and 
terminology in their assumed rigorous demonstrations.\par 
The major difficulty was the fact that mathematicians had not as yet developed 
a precise language for general mathematical discourse$,$ nor had they even 
decided upon accepted definitions for such things as the real numbers. Within 
their discussions they conjoined terms such as ``{infinitely small}'' with 
the term limit in the hopes of bringing some logical consistency to their 
discipline.\par  
   The situation changed abruptly in 1821. {Cauchy}$,$ the foremost 
mathematician of this period$,$ is believed by many to be the 
founder of the modern limit concept that was eventually formalized by 
 {Weierstrass} in the 1870's. A reading of Cauchy's Cours d'Analyse (Analyse 
Alg\'ebrique)[1821] yields the fact$,$ even to the causal observer$,$ that he 
relied heavily upon this amalgamation of terms and in numerous cases utilized
infinitesimal reasoning entirely for his ``rigorous'' demonstrations. He 
claimed to establish an important theorem using his methods - a theorem that 
 {Abel} [1826] showed by a counterexample to be in error. No matter how
mathematicians of that time period described their vague infinitesimal
methods they failed to produce the appropriately altered theorem - a modified
theorem that is essential to Fourier and Generalize Fourier Analysis.\par  
    Beginning in about 1870$,$ all of the language and methods of 
infinitesimal reasoning were replaced in the mathematical 
textbooks by the somewhat nonintuitive approximation methods we 
term the ``$\delta - \eps$'' approach. These previous difficulties are the 
direct causes that have led to the modern use of axiom techniques and the great
linguistic precision exhibited throughout modern mathematical literature.\par
   However$,$ scientists and engineers continued to use the old \underbar 
{incorrect} infinitesimal terminology. As an example$,$ Max {Planck} wrote in his 
books on theoretical mechanics that {\it ``a finite change in Nature always 
occurs in a finite time$,$ and hence resolves into a series of infinitely small 
changes which occur in successive infinitely small intervals of time.''} He 
then attempts to instruct the student in how to obtain mathematical models 
from this general description. Unfortunately$,$ at that time$,$ such terms as 
``infinitely small'' had no mathematical counterpart.\par 

     In many textbooks that claim to bridge the gape between 
abstract analysis and applications$,$ students often receive the 
impression that there is no consistent and fixed method to obtain 
applied mathematical expressions and indeed it takes some very special type of 
``intuition'' that they do not possess. In fact$,$ {Spiegel} in his present day 
textbook ``Applied Differential Equations'' writes the following when he 
discusses how certain partial differential equations should be ``derived.'' He 
states that rigorous methods should not be attempted by the student$,$ but {\it 
``it makes much more sense$,$ however$,$ to use plausible reasoning$,$ intuition$,$ 
ingenuity$,$ etc.$,$ to obtain such equations and then simply postulate the 
equations.''} \par  
         
     In 1961$,$ Abraham {Robinson} of Yale solved the infinitesimal 
problem of Leibniz and discovered how to correct the concept of 
the infinitesimal. This has enabled us to return to the more intuitive 
analytical approach of the originators of the Calculus. {Keisler} writes  
that this achievement  {\it ``will probably rank as one of the major 
mathematical advances of the twentieth century.''} Robinson$,$ who from 1944 - 
1954 developed much of the present supersonic aerofoil theory$,$ suggested that 
his discovery would be highly significant to the applied areas.  Such applied 
applications began in 1966$,$ but until 1981 were confined to such areas as 
Brownian motion$,$ stochastic analysis$,$ ultralogic cosmogonies$,$ quantum field 
theory and numerous other areas beyond the traditional experience of 
the student.\par    

   In 1980$,$ while teaching basic Differential Equations$,$ this author  was 
disturbed by the false impression given by Spiegel in the above quotation 
relative to the one dimensional wave equation. It was suggested that I apply my 
background in these new infinitesimal methods and find a more acceptable 
approach. The approach discovered not only gives the correct derivation for 
the n-dimensional general wave equation but actually solves the {d'Alembert - 
Euler} problem and gives a fixed derivation method to obtain the partial 
differential equations for mechanics$,$ hydrodynamics and the like. These 
rigorous derivation methods will bridge the gape between a student's 
laboratory$,$ and basic textbook descriptions for natural system behavior$,$ and 
the formal analytical expressions that mirror such behavior. Indeed$,$ slightly 
more refined procedures can even produce the relativistic alteration taught in 
modern physics. Moreover$,$ pure nonstandard models are now being used to 
investigate the properties of a substratum world that is believed to 
directly or indirectly 
effort our standard universe. These include pure nonstandard models for the
fractal behavior of a natural system$,$ nonstandard quantum fields$,$ a necessary 
and purely nonstandard model for a cosmogony (or pregeometry) that generates 
many different standard cosmologies as well as automatically yielding a 
theory of ultimate entities termed subparticles. \par  
         
  The major goal for writing this and subsequent manuals is to 
present to the faculty$,$ and through them to the student$,$ these 
rigorous alterations to the old infinitesimal terminology so that the 
student can once again benefit from the highly intuitive processes of 
infinitesimal reasoning - so that they can better grasp and understand exactly 
why infinitesimal models are or are not appropriate and when appropriate why 
they predict natural system behavior. Except for the basic calculus and the 
more advanced areas$,$ there are no textbooks nor any properly structured 
documentation available which presents this material at the undergraduate 
level. In my opinion it will be 15 to 20 years$,$ if not much longer$,$ before 
such material is available in the commercial market and instructors properly 
trained.  An immediate solution to this problem will give your students a 
substantial advantage over their contemporaries at other institutions and 
place your institution in the forefront of what will become a major worldwide 
trend in mathematical modeling.\par 
\medskip 
\leftline{1.2 \underbar{Manual Construction.}}\par
\medskip
\indent The basic construction of these manuals will be considerably different 
from the usual mathematical textbook. No proofs of any of the fundamental 
propositions will appear within the main body of these manuals. However$,$ all 
propositions that do not require certain special models to establish are 
proved within the various appendices. A large amount of attention is paid 
to the original intuitive approaches as envisioned by the creators of the 
Calculus and how these are modified in order to achieve a rigorous 
mathematical theory.\par 
   Another difference lies in the statements of the basic analytical 
definitions. Many definitions are formulated in terms of the original 
infinitesimal concepts and not in terms of those classical approximations 
developed after 1870. Each of these definitions is shown$,$ again in an 
appropriate appendix$,$ to be equivalent to some well-known  ``$\delta -\eps$'' 
expression. Moreover$,$ since these manuals are intended for individuals who 
have a good grasp of either undergraduate analysis or its application to 
models of natural system behavior then$,$ when appropriate$,$ each concept is 
extended immediately to Euclidean n-spaces.\par 
Nonstandard analysis is NOT a substitute for standard analysis$,$ it is a 
necessary rigorous extension. Correct and efficient infinitesimal modeling 
requires knowledge of both standard and nonstandard concepts and procedures. 
Indeed$,$ the nonstandard methods that are the most proficient utilize all of
known theories within standard mathematics in order to obtain the basic 
properties of these nonstandard extensions. It is the inner play between such 
notions as the standard$,$ internal and external objects that leads to a truly 
significant comprehension of how mathematical structures correlate to patterns 
of natural system behavior. Our basic approach employs simple techniques
relative to abstract model theory in order to take full advantage of all
aspects of standard mathematics. The introduction of these techniques is in 
accordance with this author's intent to present the simplest and 
direct approach to this subject.\par
    Since it is assumed that all readers of these manuals are well-versed 
in undergraduate Calculus$,$ then your author believes that is it unnecessary 
to follow the accepted ordering of a basic Calculus course; but$,$ rather$,$ 
he will$,$ now and then$,$ rearrange and add to the standard content. This will 
tend to bring the most noteworthy aspects of infinitesimal modeling to your 
attention at the earliest possible moment. I have this special remark for the 
mathematician. {\it These manuals are mostly intended for those who apply 
mathematics to other disciplines. For this reason$,$ many definitions$,$ proofs and 
discussions are presented in extended form. Many would not normally appear in 
a mathematicians book since they are common knowledge to his discipline. Some 
would even be considered as ``trivial.'' Please be patient with my 
exposition.}  
\par 
   It has taken 300 years to solve what has been termed ``The problem of 
 {Leibniz}'' and it should not be assumed that the solution is easily grasped  
or readily obtained. You will experience some startling new ideas and 
encounter procedures that may be foreign to you. Hopefully$,$ experience$,$ 
intuition and knowledge are not immutable. It is my firm belief that$,$ though 
proper training$,$ these three all important aspects of scientific progress  
can be expanded in order to reveal the true$,$ albeit considerably different,
mathematical world that underlies all aspects of rigorous scientific
modeling. It has been hoped for many years that individuals who have a vast 
and intuitive understanding of their respective disciplines would learn these 
concepts and correctly apply them to enhance their mathematical models.
It is through your willingness to discard the older less productive$,$ and
even incorrect$,$ modeling language that this goal will eventually be met. 
\vfil\eject
\vfil
\eject
\centerline{Chapter 2.}
\medskip
\centerline{\bf INFINITESIMALS$,$ LIMITED}
\centerline{\bf AND INFINITE NUMBERS}
\bigskip
\leftline{2.1 \underbar{Some Notation and Definitions.}}\par             
\medskip 
\indent When {\tt Robinson [1961]} first introduced his new concepts he used both 
abstract algebraic and logic notions. A few years later$,$ for simplicity in 
exposition$,$ the basic algebraic ideas where expanded and$,$ indeed$,$ many
of these algebraic results appear here for the first time. Utilizing some 
fundamental facts about algebraic structures the general properties of the 
infinitesimals can be readily obtained. Your author has successfully used this 
approach while instructing a series of 400 level courses. This algebraic
approach is based upon 
but two assumed requirements$,$ requirements that are established from 
fundamental set theory in a later section. First$,$ however$,$ we need to recall 
certain elementary definitions and identify notation.\par
Throughout these manuals the symbol $\real$  denotes the real numbers which we assume 
is a complete ordered field.\par
(I) Assume that there exists another ordered field
$\hyperreal$ with the properties that $\real$ is a ordered subfield of
$\hyperreal$ (i.e. $\real$'s basic order and field properties are those of 
$\hyperreal$ but restricted to members of $\real$) {\bf AND}
$\real \not= \hyperreal.$\par
 The set $\hyperreal$ is called 
by various names$,$ {\bf hyperreals$,$ star-reals} or {\bf extended reals}. Also 
let the symbol $\nat$ denote the natural numbers
(including zero).\par
The ordered field $\cal F$ is assumed to contain a copy of non-negative integers (natural numbers)  
$\nat$ formed by adjoining to the additive identity (the zero of $\cal F$) 
finite sums of the multiplicative identity (the 1 of $\cal F $). You can then 
simply consider $\nat \subset \cal F,$ in general. The concept {\bf complete}
is discussed in most elementary analysis books where it is almost always shown 
that {\sl if $\cal F$ is a complete ordered field and $r\in\cal F$$,$ then there 
exists some $n\in \nat$ such that $\vert r \vert < n.$} This property for a
field is called the {\bf Archimedean}  property for a field. 
As will be shown$,$ an important and 
general algebraic result states that {\bf whatever the hyperreals may be
they \underbar{cannot} be Archimedean.}\par
\indent   (II) The second assumption is the following: {\bf assume that
there exists some nonzero $\eps \in \hyperreal$ such that for each positive
$r \in \real$; it follows that}
$$0< \eps < r.$$
\indent Historically$,$ it is unfortunate that abstract algebra was not 
investigated at a much earlier date since much of the difficulty Leibniz and 
Euler faced in having the infinitesimals accepted as genuine entities would 
have been erased.\par
\vskip 18pt
\indent {\bf Theorem 2.1.1.} {\sl Assumption (II) holds for $\hyperreal$ if
and only if $\hyperreal$ is not Archimedean.}\par
\vskip 18pt
\indent What Theorem 2.1.1 means is that the more we know about ordered fields 
that are not Archimedean the better. Of course$,$ such fields also are not 
complete. It is now possible to define explicitly the ``infinitesimals.''
Note that we interpret the logical symbol $\land$ as the word ``and,'' the 
symbol $\forall$ by any one of the expressions ``for each$,$ for all$,$ for 
every'' and the symbol $\exists$ by one of the expressions ``there exists 
some$,$ there exists one$,$ there exists an."
Since such logical notation will become significant later in this manual we 
might as well practice its use as soon as possible. Let ${\real}^{+}$ denote 
the {positive real numbers}.\par
\vskip 18pt
\hrule
\smallskip
\indent {\bf Definition 2.1.1 (Infinitesimals).} Let
$\monad 0 = \{x\vert (x\in 
\hyperreal)\land \forall r(r\in {\real}^+\to \vert x \vert < r\}.$ Or in 
words$,$ $\monad 0$ is the set of all hyperreal numbers $x$ such that if $r$ is a 
positive real number$,$ then the absolute value of $x$ is less than $r.$ The set 
$\monad 0$ is called the set of {\bf INFINITESIMALS}. In the literature the 
set $\monad 0$ is also denoted by the sysmbols $M_1$ and $o.$\par
\smallskip
\hrule
\vskip 18pt

\hrule
\smallskip
\indent {\bf Definition 2.1.2. (Limited).}  Let
${\cal O} = \{ x \vert (x \in \hyperreal)\land\exists r(r \in {\real}^+ \land
\vert x \vert < r \}.$ Or in words$,$ $\cal O$ is the set of all hyperreal 
numbers 
$x$ such that there exists some positive real number $r$ such that $\vert x 
\vert < r.$ The set $\cal O$ is called the set of {\bf LIMITED} numbers. 
Robinson first called $\cal O$ the ``finite'' numbers and denoted them by the 
symbol $M_0.$\par
\smallskip
\hrule
\vskip 18pt
\indent The term ``limited'' is relatively new in the literature and many 
nonstandard analysts still employ the term ``finite.'' When this is done$,$ there
is some confusion when the term finite is used in the ordinary since of set 
theory. The concept of limited can also be interpreted as meaning that each of 
these hyperreal numbers is ``limited by'' or ``bounded by'' a real number. 
There are many reasons why the term ``bounded'' would be confused 
with the same term as employed in standard analysis and for this reason is not 
used for the concept of limited.\par
\vskip 18pt
\hrule
\smallskip
{\bf Definition 2.1.3. (Infinite).} The set $\hyperreal - \cal O= \hyperreal_\infty$ is the set 
of {\bf INFINITE} hyperreal numbers. Or in words$,$ those hyperreal numbers 
that are not limited are the infinite. The infinite numbers are also called 
the unlimited numbers.\par
\smallskip
\hrule
\vskip 18pt
I point out that the set of infinite numbers can also be characterized analytically
as follows: $x\in \hyperreal$ is infinite if and only if for every $r \in 
{\real}^+,\ \vert x \vert > r.$ This characterization is often very useful.  The next definition relates the infinitesimals to 
the original concept of when two numbers are ``infinitely close.'' This is 
probably the most significant concept for infinitesimal modeling and 
deserves attention.\par
\vskip 18pt
\hrule
\smallskip
\hrule
\smallskip
{\bf Definition 2.1.4. (Infinitely Close).} Two hyperreal numbers $x,y$ are
{\bf INFINITELY CLOSE} if $x - y \in \monad 0.$ Or in words$,$ if there 
difference is an infinitesimal. The symbol used for infinitely close is
$\approx .$\par
\smallskip
\hrule
\smallskip
\hrule
\vskip 18pt
How do we extend these definitions for the case of the Euclidean n-spaces,
$\realp n$ where we are using the Euclidean norm $\Vert (x_1,\ldots,x_n) \Vert 
= \sqrt {x_1^2 + \cdots + x_n^2}$?\par 
\vskip 18pt
\hrule
\smallskip
{\bf Definition 2.1.5. (Euclidean Extensions).} In Definitions 2.1.1 and 2.1.2 substitute for the absolute value 
symbol $\vert \cdot \vert$ the norm symbol $\Vert \cdot \Vert,$ as it would be defined relative to $\hyperreal^n.$ Of course$,$ we 
should also substitute the term ``vector'' or the symbol $\realp 
n$ when appropriate.\par
\smallskip
\hrule
\vskip 18pt
In the next section$,$ the basic algebraic properties for the infinitesimals$,$ 
limited numbers etc. are explored along with the relationships between the 
above hyperreal concepts and those of the hyperreal n-spaces.\par
\vfil\eject
\leftline{2.2 \underbar{Basic Algebra.}}\par
\medskip
You should expect that the basic properties for the above defined entities will 
be presented rapidly since the proofs only appear in the appendix. This has 
both advantages and disadvantages. These properties are couched in terms of 
some very well-known algebraic structures and we need a very brief refresher 
course relative to field theory.\par
   The ordered fields $\real,\ \hyperreal$ have the weaker structural 
property of being a {\bf ring.} As far as nonempty subsets of a 
field $\cal F$ are concerned {\sl rings may be characterized as any nonempty 
$A \subset \cal F$ that are closed under the operations of subtraction and
multiplication.} A ring is also closed under addition$,$ has a zero and 
additive inverses. Rings share with the field itself the right 
and left distributive laws$,$ and the usual associative and commutative 
properties used in parentheses manipulation. Rings can differ greatly from a 
field in that they need not contain an element that has a multiplicative 
inverse. In the case that $\cal F$ is ordered$,$ then rings$,$ in general$,$ only 
share the basic simple order property for the field $\cal F$ when it is 
considered to be restricted to the ring. Our last general ring concept it that 
of the ``{ideal}.'' Let the ring $\Re \subset \cal F$. Then nonempty
$\Im \subset \Re$ is an {\bf ideal} {\sl of (or in) $\Re$ if $\Im$ is a 
subring (i.e. a ring with respect to the ring structure of $\Re$ and a subset 
of $\Re$) and for each $r \in \Re$ and each $x \in \Im$ it follows that
$xr \in \Im.$} Thus an ideal ``absorbs'' the members of $\Re$ by 
multiplication. Before proceeding with our first list of properties notice 
that the infinitesimals $\monad 0 =\{x\vert x\in \hyperreal \land x\approx 
0\}.$\par
\vskip 18pt
\indent {\bf Theorem 2.2.1.} {\sl The set of limited numbers$,$ $\cal O,$ is a 
subring of $\hyperreal,\ \real \subset \cal O$ and $\cal O$ is not a 
field$,$ but if $x\in {\cal O} - \monad 0$$,$ then $x^{-1} \in \cal O$.}\par
\vskip 18pt
\indent {\bf Theorem 2.2.2.} {\sl The infinitesimals$,$ $\monad 0,$ form a 
subring of $\cal O$ and $\Gamma $ is infinite if and only if there is some
nonzero $\eps \in \monad 0$ such that $\Gamma = 1/{\eps}.$}\par
\vskip 18pt
\indent {\bf Theorem 2.2.3.} {\sl The set of infinitesimals$,$ $\monad  0,$ is an 
ideal of $\cal O.$}\par
\vskip 18pt
\hrule
\smallskip
\hrule
\smallskip
\indent The fact that $\monad 0$ forms an ideal of $\cal O$ is very 
significant for the correct theory of infinitesimals. Not only are the 
infinitesimals closed under finite addition and finite product$,$ but absorb$,$ 
under product$,$ all of the real numbers. It's unfortunate that Leibniz and 
others could not establish such results rigorously since$,$ if they could have$,$ 
undoubtedly much of the criticism of their concepts would not have 
developed.\par 
\smallskip
\hrule
\smallskip
\hrule
\vskip 18pt
\indent We now consider what happens when we take any real number and add to 
it the infinitesimals.\par
\vskip 18pt
\hrule
\smallskip
\indent {\bf Definition 2.2.1. (Monad).}  Let $r\in \real$. Then a {\bf MONAD OF 
(ABOUT)} $r\ $ is the set $\monad r =
\{x\vert x\in \hyperreal \land (x -r)\in 
\monad 0\} = \{x\vert x\in\hyperreal \land x\approx r\}.$\par
\smallskip
\hrule
\vskip 18pt
{\bf Theorem 2.2.4.} {\sl The binary relation $\approx$ is an equivalence 
relation on $\hyperreal .$}\par
\vskip 18pt
{\bf Theorem 2.2.5.} {\sl For each $x,y\in \real,$\hfil\break
\indent\indent (i) $\monad x \cap \monad y = \emptyset$ if and only if $ x 
\not= y,$\hfil\break
\indent\indent (ii) ${\cal O} = \bigcup \{\monad r\vert r\in \real\}.$}\par 
\medskip
{\bf Corollary 2.2.5.1.} {\sl If $x,y \in \real,\ x < y,\ z\in \monad x,\ 
w\in\monad y,$ then $z < w.$}\par
\medskip
{\bf Corollary 2.2.5.2.} {\sl If $x,y\in \monad r,\ z\in\hyperreal,\ x<z<y,$ 
then $z\in \monad r.$}\par
\vskip 18pt
The important Theorem 2.2.5 can be expressed in words by stating that the set 
of monads forms an ordered partition of the limited numbers.  It is beginning to 
appear as if the infinitesimals are indeed behaving in the manner first 
envisioned by the founders of the infinitesimal calculus. HOWEVER$,$ in at 
least one respect the monads as well as the set of limited numbers do not 
share an important real number property. From the previous discussion$,$ $\hyperreal$
is known not to be complete. But$,$ are there significant sets that have upper 
bounds in $\hyperreal$ and do not have a least upper bound?\par
\vskip 18pt
{\bf Theorem 2.2.6.} {\sl Each monad and the set of limited numbers are bounded 
above {\rm [}resp. below{\rm ]}$,$ but do not possess a least upper bound
{\rm [}resp. greatest lower bound{\rm ]}.}\par
\vskip 18pt
\leftline{2.3 \underbar{Euclidean n-spaces.}}
\bigskip
\indent The following are additional propositions that extend some of the 
above results to Euclidean n-spaces. Whenever possible the same notation is 
used and should be understood from the context. It should be obvious how the 
previous definitions extend to Euclidean n-spaces.\par
\indent {\bf Theorem 2.3.1.} {\sl For any $n\in\nat,\ n\geq 1,$  the set of limited 
vectors in $\hyperrealp n$ is equal to $${\cal O}^n = \overbrace{{\cal 
O}\times\cdots\times{\cal O}}^{n\rm\;factors}.$$}\par
\vskip 18pt
{\bf Theorem 2.3.2.} {\sl For any $n\in\nat,\ n\geq 1,$ then monad of $\vec v
=(x_1,\ldots ,x_n)\in \realp n$ is equal to $$\monad {\vec v} = 
\overbrace{\monad {x_1}\times\cdots\times\monad {x_n}}^{n\rm\;factors}.$$}\par 
\vskip 18pt  
{\bf Theorem 2.3.3.} {\sl For any $n\in\nat,\ n\geq 1,$ and for each $\vec 
v,\vec w\in 
\realp n,$\hfil\break
\indent\indent (i) $\monad {\vec v} \cap \monad {\vec w} = \emptyset$
if and only if $ \vec v 
\not= \vec w,$\hfil\break
\indent\indent (ii) ${\cal O}^n= \bigcup \{\monad {\vec v}\vert \vec v\in \realp n\}.$}\par 
\vskip 18pt
\leftline{2.4 \underbar{The Standard Part Operator}}\par
\medskip
In infinitesimal modeling$,$ various methods exist that allow us to 
investigate what might be termed as the {\sl micro-effects} that occur within
a world 
called the {\sl Nonstandard Physical World$,$ (i.e. {NSP-world})} or
{\sl Deductive World (i.e.{D-world}).} When these micro-effects are modeled by 
means of the infinitesimals$,$ certain mathematical operators applied to  
such infinitesimals yield standard mathematical objects that are often perceived to 
measure the corresponding {\sl natural world (i.e.{N-world})} effects that 
govern phenomenological behavior. The next operator is probably the most 
significant of these standardizing processes.\par
\vskip 18pt\vfil\eject
\hrule
\smallskip
\hrule
\smallskip
\indent {\bf Definition 2.4.1. (Standard Part).} Let $x\in \cal O$.
Define the function ${\rm st}\colon {\cal O} \to \real$ as follows:\par
\indent\indent (i) let $\st x = r,$ where $r \in \real$ is the unique real 
number such that $x\in \monad r.$\par
\indent\indent (ii) The function st is often extended to all of $\hyperreal$ 
by letting $\st x = +\infty,$ when $x\in \hyperreal - \cal O$ and $x 
> 0$ or $\st x = -\infty,$ when $x\in \hyperreal - \cal O$ and $x<0.$ The map 
st is called the {\bf standard part operator}.\par
\smallskip
\hrule
\smallskip
\hrule
\vskip 18pt
\indent Since $\real \subset \cal O$ then$,$ as will be seen$,$ the range of
st is $\real$ (i.e. st 
is a surjection.) Please note {\bf the basic properties of the standard part 
operator are highly important in elementary infinitesimal analysis.}\par
\vskip 18pt
{\bf Theorem 2.4.1.} {\sl Let $x,y\in \cal O$. Then\hfil\break
\indent\indent (i) $x \approx y$ if and only if $\st x = \st y,$\hfil\break
\indent\indent (ii) $x \approx \st x,$\hfil\break
\indent\indent (iii) if $x\in \real,$ then $\st x = x,$\hfil\break
\indent\indent (iv) if $x \leq y,$ then $\st x \leq \st y,$\hfil\break
\indent\indent (v) if $\st x \leq \st y,$ then either $x \leq y$ or $x-y \in
\monad 0$ with unknown order.}\par
\vskip 18pt
Part (iii) of Theorem 2.4.1 is what implies that the range of st = $\real.$ 
The next theorem is established from the basic definitions and the fact that 
$\monad 0$ is an ideal of $\cal O.$ The map st is purely algebraic in
character and its application is
a remarkable indication of how pure abstract algebra can be utilized to 
obtain classical results.  For algebraists$,$ these 
pure algebraic characterizations are discussed following the next theorem.\par
\vskip 18pt
{\bf Theorem 2.4.2.} {\sl Let $x,y\in \cal O$. Then \hfil\break
\indent\indent (i) $\st {x\pm y} = \st x \pm \st y,$ \hfil\break
\indent\indent (ii) $\st {xy} = (\st x)(\st y).$}\par
{\bf Corollary 2.4.2.1} {\sl Let $x,y\in \cal O$. Then \hfil\break
\indent\indent (i) if $\st y \not= 0,$ then $\st {x/y} = {\st x}/{\st y}.$
\hfil\break
\indent\indent (ii) if $y = \root n \of x,$ then $\st y = {(\st x)}^{(1/n)}$$,$ 
where it is always the case that if $x\in \cal O$$,$ then $\root n \of x \in 
\cal O$.}\par
\vskip 18pt
{\bf Theorem 2.4.3.} {\sl For $\vec v\in {\cal O}^n$$,$ let $\st {\vec v} =
(\st {x_1},\ldots,\st {x_n}).$ The set ${\cal O}^n$ forms a vector space with 
respect to the ring $\cal O$ (i.e. a module) and as such the map {\rm st} 
distributes over the vector space algebra as well as the component defined dot 
and for $n = 3$ the cross product operators.}
\vskip 18pt
Notice that even though $\monad 0$ is an ideal of $\cal O$ it is most 
definitely only a subring of $\hyperreal.$ Indeed$,$ products are somewhat 
indefinite. For example$,$ let $0 \not= \eps\in \monad 0.$ Then ${\eps}^2 \in 
\monad 0.$ 
But $\eps({\eps}^{-2}) \in \hyperreal - \cal O,$ while $\eps(1/\eps) \in \cal 
O$ and ${\eps}^2(1/\eps) \in \monad 0.$\par
\medskip
\indent [{\it \underbar{This is for the algebraist.} The above theorems show that the 
mapping st is a ring epimorphism with ker(st) = $\monad 0$.\par
\vskip 18pt\vfil\eject
{\bf Theorem 2.4.4.} The set $\monad 0$ is a maximual ideal in $\cal O$ and the 
quotient ring ${\cal O}/\monad 0$ is isomorphic to $\real.$}]
\vskip 18pt
\leftline {2.5 \underbar{A Slight Quandary.}}\par
In 1961 when Robinson first published his new theory$,$ it was restricted to 
subsets of $\hyperreal$ and functions defined on $\hyperreal.$ Shortly after 
this {\tt Luxemburg [1962]} expanded upon Robinson's paper but still retained the 
same general restrictions. It became obvious that any extensive applications 
of infinitesimal analysis to functional analysis would require an extension to
other set-theoretic 
objects. {\tt Robinson [1966]} did just this but at a great expense to 
elementary exposition. He used a purely logical approach and the simplified
theory 
of {types}. Indeed$,$ this author entered this field in 1968 and first used the 
type-theoretic approach. There is no doubt that the type-theoretic approach is 
beyond almost all mathematicians and appliers of mathematics that have not had 
specialized training. {\tt Machover and Hirschfeld [1969]} introduced a 
simplification to the Robinson type-theoretic approach by restricting their 
nonstandard analysis to set theory itself. Unfortunately$,$ their approach 
requires that individuals utilize two ``set theories'' - a pseudoset theory 
and the standard set theory. Your author wrote his doctoral dissertation in 
pseudoset theory. At approximately the same time {\tt Robinson and Zakon [1969]} 
published a paper that further simplified Robinson's approach. They were able 
to show how set theory itself could be used for all the necessary formations 
of nonstandard analysis and that pseudoset theory was no longer necessary. 
Since 1969 there have been other attempts at simplifying the foundations of 
nonstandard analysis but it is this author's belief that they are not 
appropriate to those mathematicians who are not accustomed to special logical 
procedures.\par The quandary I face is that there are no other appropriate 
simplifications of Robinson's theory except for the 1969 Robinson - Zakon 
approach. This approach definitely needs to be presented within this basic 
manual so that you can encounter the full power of infinitesimal analysis. The 
difficulty is my wish not to present numerous definitions and 
constructions from mathematical logic in the next chapter; but$,$ rather$,$ 
actually to do some infinitesimal analysis immediately. 
There is a procedure that does allow this if we postpone until a later chapter 
the reasons why this procedure functions.\par 
\vfil\eject
\centerline{Chapter 3.}
\medskip                
\centerline{\bf SOME SET THEORY$,$ CONVERGENCE}
\centerline{\bf AND LEIBNIZ'S PRINCIPLE}
\bigskip
\leftline{3.1 \underbar{Some Set Theory.}}\par    
\medskip 
Historically we have over 5000 years of standard mathematics available. It is 
certainly reasonable to utilize all of this standard mathematics in the 
development of what has become know as ``nonstandard analysis''. The term 
``nonstandard'' should not be taken to mean that the properties to be 
discussed are not relative to the known properties for standard mathematics. 
Indeed$,$ the opposite is the case. In the development of nonstandard 
mathematics$,$ it has always been required that nonstandard structures have 
all the same ordinary properties as the standard mathematical structures.
They are nonstandard in the sense that they have additional useful properties
not possessed by the standard structures. The simplification procedures 
mentioned at the conclusion of the last chapter have$,$ for the most part$,$ 
eliminated certain model theoretic difficulties associated with Robinson's 
application of the simplified theory of types. Depending upon the type of
set-theoretic entity used no special considerations need to be considered. 
Even though what is about to be presented seems very straightforward and is 
easy to grasp$,$ it should not be assumed that these
set-theoretic procedures hold in the exact form given for more extensive
set-theoretic processes that are not specifically mentioned. [Note: The formal set-theort used is found in {\bf Suppes [1960]}.]\par
  Why do we need to consider basic set theory at all? Are not the properties of 
the infinitesimals$,$ limited and infinite numbers sufficient for analysis? 
Didn't Leibniz stress the calculus of the infinitesimal and infinite numbers 
as the basis for analysis? 
Since mathematicians investigate the {\it relations between individual 
objects} and our simplification is set-theoretical then basic set theory
is required. Indeed$,$ even an {order pair} $(a,b)$ 
is set-theoretically defined as $(a,b) = \{\{a,b\},\{a\}\}.$\par
{\bf Theorem 3.1.1.} {\sl It is permissible to assume that\par
\indent\indent (i) if $A\subset \real,$ then $A\subset 
\hyperreal,$\hfil\break
\indent\indent (ii) if $A\subset \realp n,$ then $A\subset \hyperrealp n,$ 
\hfil\break
\indent\indent (iii) if $A\subset (\realp n)\times(\realp m),$ then 
$A \subset (\hyperrealp n)\times(\hyperrealp m).$}\par
\bigskip
Notice that the subsets considered in Theorem 3.1.1 exhaust almost all of the 
mathematical objects studied in an ordinary course in undergraduate analysis 
and much more. However$,$ at present$,$ the notion of a ``{constant}'' 
in a formal mathematical language needs to be discussed. From the viewpoint
of a formal language$,$ it is 
assumed that every subset of $\real,\ \realp n$ or $(\realp n)\times(\realp m)$
is associated with at least one constant symbol that ``names'' that subset. 
These constants are employed in two contexts. First$,$ as language constants and 
nothing more$,$ then as names for mathematical objects within our special set 
theory. Technically these are two different usages which are not usually 
mentioned in introductory texts. In infinitesimal analysis these technicalities
are significant$,$ however. In the next theorem$,$ the notation $\hyper A$ also 
has two meanings. $\hyper A$ is a constant symbol in a mathematical language 
and$,$ in most cases$,$ represents a mathematical entity different from $A.$\par
\vskip 18pt
{\bf Theorem 3.1.2.} {\sl If $A \subset \real$ {\rm [}resp. $\realp n,
\ (\realp n)\times(\realp m)${\rm ]}$,$ then there exists $\hyper A\subset \hyperreal$
 {\rm [}resp. $\hyperrealp n,\ (\hyperrealp n)\times(\hyperrealp m)${\rm ]}
such that \par
\indent\indent (i) $A \subset \hyper A$ and \hfil\break
\indent\indent (ii) $A = \hyper A$ if and only if $A$ is finite.}\par
\vskip 18pt
Thus Theorem 3.1.2. tells us that a collection of new sets exist - the
``star-sets'' - and obviously we need a certain amount of new terminology
in order 
to discuss their properties effectively. The basic operators that define the 
ordered field $\real$ may also be considered as subsets of $\real\times\real$ 
or $\realp 3.$ Do we then consider these operators extended to $\hyperreal$ 
and thus ``star'' the operator such as writing $a\ {\hyper +}\ b,$
where $a,b \in 
\hyperreal$? If one wants to be technical about the matter$,$ then the answer 
would be yes. However$,$ in practice these specific operators are not so 
denoted. On the contrary$,$ one may consider the *-field operators as basic 
and that $\real$ is defined by restricting them to the \underbar{set} $\real.$
In all cases$,$ the 
particular defining field operator is determinable from the context.\par
{\bf In all that follows let} ${\cal U} = \real\cup\realp n\cup (\realp 
n)\times(\realp m).$ Now a major consideration in the writing of these 
manuals relates to mathematical rigor. Even though all stated propositions 
are established in the appendixes$,$ how explicit should the main text be when 
confronted with the basic construction of a nonstandard model for analysis? 
Analysis of the fundamental objects in $\cal U$ often requires collections of 
members from $\cal U$ that are themselves not members of $\cal U.$ In order to 
satisfy the modern trend that rejects pronouncements that are not at least 
superficially established$,$ a slight acquaintance with the actual construction 
of the standard model for this analysis seems appropriate. It is definitely not 
essential that you have any in depth knowledge of the fine details of the 
construction of this standard structure. But$,$ be assured that every standard
set-theoretic object that will ever be consider in this analysis is a member
of the standard structure.\par
For any set $W$ in our set theory$,$ let ${\cal P}(W)$ denote the set of all 
subsets of $W,$ where the operator $\cal P$ is often called the 
{\bf {power set operator}.} Hence $A \in {\cal P}(W)$ if and
 only if $A \subset W.$
Starting with $\real = X_0$ an object called a {\bf {superstructure}} is 
constructed by induction as follows: assume that $X_n$ has been defined for 
$n \in \nat.$ Then define $X_{n+1} = {\cal P}(\bigcup \{X_i\vert 
 \leq i\leq n\}.$ 
Then the standard structure - {\bf the superstructure} - is the set ${\cal H} = 
\bigcup \{X_n\vert n \in \nat \}.$ Various properties relative to $\cal H$ are 
discussed in Appendix 3.
\vskip 18pt
\hrule
\smallskip
{\bf Definition 3.1.1. (Individual$,$ Entity$,$ Star-Sets).} Each member of $\real$
is called an {\bf {INDIVIDUAL}}$,$ while each \underbar{set} in $\cal H$ is
called an {\bf{ENTITY}} and for each entity $A$$,$ the set $\hyper A$ is
called an {\bf{EXTENDED STANDARD}} set. \par
\smallskip
\hrule
\vskip 18pt
It turns out that for any $W \in \cal H$ the set ${\cal P}(W) \in \cal H$ and
there exists the set $\Hyper {{\cal P}(W)}$.  $\{$What appears between these two braces is a refinement for what appeared in this place originally. It clears up a certain notational confusion. The star 
operator  ``*'' behaves like a mapping from ${\cal P}(W)$ into $\Hyper {{\cal P}(W)}$. As such$,$ there is a slight confusion in symbols that ``name'' the images. The usual notation for the range would be $\Hyper {({\cal P}(W))}.$ The same notation holds for $A \in {\cal P}(W),$ where $\Hyper (A)$ is the image notation. In order to relate this notation to our original when ``*'' is considered as an mapping$,$ we let $\hyper A$ be a name for the image $\Hyper (A)$. That is $\Hyper (A)= \hyper A.\}$\par
The next few propositions reveal the fundamental 
behavior of * and the first shows that * at least preserves all \underbar 
{finite} set-theoretic operations.\par
\vskip 18pt

{\bf Theorem 3.1.3.} {\sl Unless otherwise
stated all constants represent individuals or entities.\par
\indent\indent (i) $a\in A$ if and only if $\hyper a\in \hyper A;\ A\not= B$ 
if and only if $\hyper A \not= \hyper B.$\par
\indent\indent (ii) $A\subset B$ if and only if $\hyper A \subset \hyper 
B.$\par
\indent\indent (iii) $\hyper {\{A_1,\ldots,A_k\}} = \{\hyper {A_1},\ldots,
\hyper {A_k} \}.$\par
\indent\indent (iv) $\hyper {\eskip (A_1,\ldots,A_k)} = (\hyper {A_1},\ldots,
\hyper {A_k}).$\par
\indent\indent (v) $ (A_1,\ldots,A_k) \in A$ if and only 
$(\hyper {A_1},\ldots,\hyper {A_k}) \in \hyper A.$ \par
\indent\indent (vi) Let $A,B \subset \cal U.$ Then $\hyper {\eskip (A \cup B)} = \hyper A
\cup \hyper B,\ \hyper {\eskip (A \cap B)} = \hyper A \cup \hyper B,
\ \hyper {\eskip (A - 
B)} = \hyper A - \hyper B,\ \hyper {\eskip (A\times B)} = \hyper A \times \hyper 
B,\ \hyper {\eskip \emptyset} = \emptyset.$\par
\indent\indent (vii) If $A\in \real,$ then $\hyper A = A.$}\par
Recall that if $R$ is any n-ary relation $(n > 1),$ then the {\bf {\it i}th 
projection} is the set $P_i(R) = \{x_i\vert (x_1,\ldots,x_i,\ldots,x_n)\in R\},$ 
where $1 \leq i \leq n.$ Also since  $R\subset A^n$ if and only if $\hyper 
{R} \subset (\hyper A)^n$ then $\hyper {R}$ is any n-ary 
relation in $(\hyper A)^n$ and the set-theoretic projections are defined for 
$\hyper {R}.$ If $R$ is an binary relation$,$ then $P_1$ is called the 
 {\bf domain} and $P_2$ the  {\bf range} of $R.$ \par
\vskip 18pt
{\bf Theorem 3.1.4.} {\sl Let $\emptyset \not= R\subset 
 A^n,\ (n>1),$\par
\indent\indent (i) $\hyper {P_i(R)} = P_i(\hyper R).$\par
\indent\indent (ii) If $R$ is a binary relation (i.e. n = 2)$,$ then $\hyper 
{R^{-1}} = (\hyper R)^{-1}.$}\par 
\vskip 18pt
The next theorem presents two technical results and is stated here for 
completeness.\par
\vskip  24pt
{\bf Theorem 3.1.5.} {\sl Let $A$ be an entity. Then\par
\indent\indent (i) $\hyper {\{(x,x)\vert x\in A\}} = \{(x,x)\vert x \in \hyper 
A \}.$\par
\indent\indent (ii) $\hyper {\{(x,y)\vert x\in y\in A\}}= \{(z,w)\vert z\in 
w\in \hyper A\}.$}\par
\vskip 18pt
Theorems 3.1.3$,$ 3.1.4 can be paraphrased by simply considering the * map to 
be a natural map which appears to distribute over finite ``everything.''\par
\vskip 18pt 
[{\it \underbar{This is {for the algebraist}.} The map * when restricted to
${\cal P}(W)$ is a homomorphism on the Boolean algebra
${\cal P}(W).$ The map * is a monomorphism and it is a significant fact 
that it is not an epimorphism.}]\par
\bigskip
\leftline{3.2 \underbar{Convergence.}}\par
Since the time of {Zeno} [350 BC] the concept of the {non-finite} has been somewhat 
controversial. In modern axiomatic set theory$,$ except for considerations of 
consistency$,$ such controversies have little meaning. However$,$ in physical 
modeling$,$ accepting the reality of the {non-finite} has numerous philosophical 
ramifications. Unfortunately$,$ these manuals are not the appropriate form to 
discuss these side issues. We have mentioned the concept from an intuitive
set-theoretic 
point of view and this will suffice since infinitesimal analysis eliminates,
to a great degree$,$ much of the more esoteric and philosophic concerns that 
plagued early mathematical discourse.\par
{\tt Leibniz [1701]} discussed what he claimed to be a relation between his
calculus of infinitesimal and infinite numbers and the ``method of Archimedes'' 
(the so-called method of exhaustion). He$,$ of course$,$ never established this 
claim and did not seem to use any definable approximation process such as the 
``$\delta - \eps.$'' He mentions such things as fractions with infinite numbers 
as denominators$,$ infinite sequences and the like. In 
particular$,$ Leibniz$,$ de l'Hospital and even Cauchy specifically required$,$ as 
axiomized by de l'Hospital$,$ that ``a curved line be considered an infinite
assemblage of straight line segments each infinitely small.'' This concept 
should be distinguished from the completely different idea of a curve as 
composed of an infinite number of points$,$ with no magnitude$,$ - the  
 {\it {indivisibles}} - as put forth by Cavalieri [1635]. Leibniz also instituted a
notion he called the ``{sovereign principle}.'' This principle served as his 
justification that the infinitesimals and his infinite numbers obey the same 
rules as the ordinary real numbers. As will be seen$,$ many of these early ideas 
were not correct from the rigorous view point and$,$ moreover$,$ Robinson's methods are very 
broad in context since they may be applied to all standard mathematical 
structures not merely to $\realp n.$ \par 

Even though such statements as {\it a sequence 
converges to a real number ``as n goes to infinity or as n grows without 
bound or as n gets infinitely large''} and other such intuitive expressions
may be eliminated entirely from our language of convergence and replaced with 
statements involving the behavior of extended standard functions that are 
property defined on $\hypernat,$ it is a misconception that the formal
limit definition with its approximating techniques is no longer 
needed. From the view point of mathematical modeling$,$ nonstandard and standard 
analysis complement each other. The effects that appear in the natural world 
are modeled by standard analysis. The limit approach$,$ with all of 
its approximation techniques$,$ is often considered to be a natural world
approximation for events modeled within the NSP-world by nonstandard objects. 
Our basic definition for sequential convergence is stated in a language that 
parallels some of Leibniz thoughts; but$,$ as established in appendix 3 it is 
equivalent to the Bolzano-Weierstrass concept.
{\bf In all that follows,
let ${\nat}_\infty$ denote the infinite natural numbers.}\par                                                    
\vskip 18pt  
\hrule
\smallskip
{\bf Definition 3.2.1. {(Sequential Convergence)}.} Let the sequence $S\colon 
\nat 
\to \realp n,\ (n\geq 1).$ Then $S$ converges to $\vec r\in \realp n$ (notation 
$S\to \vec r$)$,$ if for
each $\Gamma \in {\nat}_\infty$ it follows that $\hyper S(\Gamma) \in 
 \monad {\vec r}.$
Thus $S$ converges to $\vec r \in \realp n$ if $\hyper S$ maps all of the
infinite natural 
numbers into a single monad (i.e. $\hyper S[{\nat}_\infty] \subset \monad 
{\vec r}$).
Recall that the value $S(n)$ of a sequence is often written as $S_n.$ \par
\smallskip
\hrule
\vskip 18pt
{\bf Theorem 3.2.1.} {\sl Definition 3.2.1 (i.e. $S\to \vec r$) is equivalent
to the limit definition for convergence of a sequence 
(i.e. $\lim_{n\to \infty} S(n) = \vec r$).} \par
\vskip 18pt 
The next result simply recasts definition  3.2.1 in terms of the algebraic 
standard part operator. However$,$ this seemingly ``trivial'' fact yields 
immediately all of the basic ``{limit theorems}'' that tend to be a first 
stumbling block for the first year calculus student. 
\par
\vskip 18pt
{\bf Theorem 3.2.2.} {\sl Let the sequence $S\colon \nat \to \realp n,
\ (n\geq 1).$ Then $S\to \vec r \in \realp n$ if and only if for each 
$\Gamma \in \nat_\infty,\ \st {\hyper S_\Gamma} = \vec r.$}\par
\vskip 18pt
Notice that when definition 3.2.1 is stated without the symbolism then it closely 
parallels the conceptional process outlined by Leibniz. {\bf A sequence 
converges to a real number $r$ if its value for each infinite number is 
infinitely close to $r.$} Thus Leibniz and others explained the result that 
the sequence $S(n) = n^{-2}$ converges to 0 by asserting that if
$\Gamma$ is an 
infinite number$,$ then ${\Gamma}^2$ is an infinite number. This implies that
${\Gamma}^{-2}$ is an infinitesimal and thus $\hyper S(\Gamma)$ is
infinitely close to 0. (They did not use the * idea$,$ however. And$,$ when no 
confusion results$,$ the * is often omitted from the extended standard function 
notation.) Unfortunately$,$ 
in order to apply such intuitive procedures one must know how the *-function
behaves. In the next section$,$ such an investigation is initiated.\par
\bigskip
\leftline{3.3 \underbar{Constants and More Constants.}}
\medskip
In section 3.2$,$ it was pointed out that Leibniz believed that the 
infinitesimal and infinite numbers ``obey the same rules as the ordinary real 
numbers.'' It was necessary to develop abstract model theory before Leibniz's 
vague comment could be rigorously justified. As previously discussed$,$ every
object in ${\Re} ={\cal U} \cup {\cal P}({\cal U})$ is denoted by a ``constant''
symbol as is every object in $\cal H.$ 
When theorems about the set ${\cal U} \cup {\cal P}({\cal U})$
are written in terms of these constants$,$ variables and the symbols $\eps$ and 
$=$$,$ an intuitive ``{first-order language}'' is often used.  Indeed$,$ every
theorem and definition that appears in 
this manual has been stated or can be reformulated in such a ``language.'' 
Here are a few examples of how these intuitive expressions appear$,$ where
all the constants are assumed to represent 
members of $\cal H.$ \par
\vskip 18pt
{\bf {Example 3.3.1.}} In our set theory you have the basic definition for 
the union of two sets $A,B \subset \Re.$ This might be expressed as follows:
for each $x \in \Re,\ x \in A \cup B= C$ if and only if $x\in A$ or $x\in 
B.$ This can also be expressed in a formalized first-order statement as 
$\forall x(x\in \Re \rightarrow (x\in A \cup B= C \leftrightarrow (x\in A)\lor
(x\in B))).$\par
\vskip 18pt
In the above first-order statement$,$ it is required that the 
variable symbol $x$ be restricted to members of $ \Re$. {\it The absolute 
requirement is that all quantified language variables must be restricted to specific sets that are
entities (i.e. to sets that are elements of various $X_p$) and that are 
represented by constants within our language.} Mostly$,$ the
variables may simply be restricted to the set $\Re,$ with certain exceptions 
mention in future chapters. However$,$ better 
comprehension is often produced when the language variables are restricted 
to sets of immediate interest.  
\vskip 18pt
{\bf {Example 3.3.2.}} Suppose that you have a sequence that is strictly 
increasing. How might this be stated? For each $x$ and $y$ if $x,y \in \nat$ 
and $x < y,$ then $S(x) < S(y).$\par
\vskip 18pt
{\bf {Example 3.3.3.}} Since the usual operator and relation symbols have been
retained as constant symbols$,$ it is a simple matter to give explicit 
functional statements. For each $n \in \nat,\ S(n) = \sin {(n^2)}.$ Notice that
the symbol $n$ has now been used as a variable and that the symbol $\sin$ is 
a constant that represents the sine function.\par
\vskip 18pt
\hrule
\smallskip
\hrule
\smallskip
Obviously$,$ in order to formulate properly these first-order expressions it is 
essential that individuals gain facility with the basic language of 
mathematics. In particular$,$ great precision must be maintained in both written 
and oral mathematical exposition. This is also one of the goals of the core 
mathematics program at the Academy.\par
\smallskip
\hrule
\smallskip
\hrule
\vskip 18pt
Are there significant mathematical concepts that might not be expressible by 
such first-order statements? Suppose you wanted to express the general concept
of mathematical {induction} as put forth by {Peano} - a concept considered to be 
one of the most basic in all of mathematics. {\sl If Q is any property which may 
or may not hold for the natural numbers$,$ and if (1) the natural number $0$ 
satisfies property Q and (2) whenever a natural number $n$ satisfies property 
Q$,$ then $n+1$ has property Q$,$ then all natural numbers satisfy property Q.} The 
difficulty in expressing the induction axiom lies in the requirement that we 
express the phrase ``Q is any property.''\par
 There are two ways of expressing a property within a first-order language. 
First$,$ a property is expressed by other first-order statements or secondly,
 it may be 
claimed that a property can be represented by a member of ${\nat}^m.$ In the 
first case$,$ we need a variable that$,$ at least$,$ represents predicates in our
language and this is not allowed in a first-order language. 
A third approach would be to replace the single induction axiom by 
infinitely many axioms 
where each axiom refers to one specific property expressed by one specific 
first-order expression. Unfortunately$,$ it can be shown that you would not 
obtain \underbar{all} of the possible properties by this process.\par
 For the second set-theoretic case$,$ it turns out that all of the possible 
 {properties for the natural numbers} still cannot be obtained since we do not 
have a set in our slightly restricted set theory that contains \underbar{all} 
of the objects that represent \underbar{all} of the properties that can be 
expressed by our first-order set-theoretic language. Further$,$ using the basic 
language of set theory one can conceive of a property determined by the 
predicate P(x) = ``$\emptyset \subset x,$'' where $\emptyset$ (the empty set) is the 
object that satisfies the set-theoretic expression: there exists a set $x$ 
such that for all sets $y,\ y\notin x.$  However$,$ for set theory to be 
consistent it must be assumed usually that if there does exist something that 
satisfies this property$,$ it is not to be called a set. Consequently$,$ one must 
be very sure that an informally stated
notion is indeed expressible by 
means of our simplistic first-order language$,$ even though it may not be done 
within an informal discussion or argument.\par   
   Every object in $\cal H$ is denoted by a logical constant. Some of these 
constants are in the customary form such as $\nat,\ \real,\ 2,\ \sqrt 2,\ +,\ 
\sin$ and all the very well-known symbols used in standard analysis.
As an abbreviation$,$ 
let $C({\cal H})$ denote the set all such constants. As previously indicated$,$ in nonstandard analysis the set 
$\Hyper {\cal H} = \bigcup \{\hyper 
 X_n\vert n \in \nat\}$ is constructed and it is assumed that each object in 
$\Hyper {\cal H}$ is denoted by a constant. The set of all such 
constants is denoted by $C(\Hyper {\cal H}).$ In particular$,$ it follows that    
 $\Hyper {\Re} = \Hyper {\cal U} \cup \hyper {\cal P}({\cal U})$ and each 
object in $\Hyper {\Re}$ is also denotable by a constant. (Note: 
$\hyper {({\cal P}({\cal U}))}$ is 
denoted by $\hyper {\cal P}({\cal U})$ since the symbol ${\cal P}({\cal U}),$
in this form$,$ is considered as one constant.)  Now if $A \in C({\cal H}),$ it
denotes an object in $\cal H$ and the notation  $\hyper A$ is the constant
that denotes the object 
in $\Hyper {\cal H}$ obtained by means of the * map. {\bf Also note that
Theorems 3.1.2 
and 3.1.3 indicate that some members of $\Hyper {\cal H}$ are denotable
by  
both a starred and unstarred constant and by a previous convention the basic 
field 
operators defined on $\hyperreal$ are not starred. In these cases$,$ we tend to 
use only the unstarred notation.} Theorems 3.1.2 and 3.1.3 imply that 
there are infinitely many members of $\Hyper {\cal H}$ that are not 
named by starred members of $C({\cal H}).$ But$,$ nevertheless$,$ they do have constant
names in an \underbar{extended language}. \par
\vskip 18pt
{\bf {Example 3.3.4.}} We know that ${\nat}_\infty \not= \emptyset.$ Moreover,
$\hypernat \subset \hyper {\eskip {\cal U}}$ implies that each member of $\hypernat$ is a member of
$\Hyper \Re$ and $\Re$ is a member of some $X_p.$ We have used capital Greek letters as the names for some of 
the members of ${\hypernat}_\infty.$ By the * process we only obtain the stars
of the natural numbers$,$ where by convention we drop the *. Since the * process 
is one-to-one into $\Hyper {\cal H}$ by Theorem 3.1.3 then no member of
${\hypernat}_\infty$ is the * of any of the constants in $C({\cal H}).$\par
  
\vskip 18pt
\leftline{3.4 \underbar{The {Leibniz' Principle} of {*-transfer.}}}\par 
\medskip
Thus far it may not appear that the previous rules and conventions associated
with infinitesimal analysis yield a simplification. However$,$ I invite you to 
compare the second chapter of Robinson's 1966 book with above procedures and judge 
for yourself. But$,$ why do we need to consider these specialized first-order 
languages in the first place? The reason lies within certain very powerful 
results from the subject area of mathematical logic. These fundamental 
propositions coupled 
with our conventions lead to a completely correct formulation of the Leibniz 
Principle$,$ which now becomes a theorem. But$,$ first$,$ one final procedure needs
to be discussed prior to formalizing this highly useful principle.\par
\vskip 18pt
\hrule 
\smallskip
{\bf Definition 3.4.1. {(*-transfer)}.} Let $\Phi$ be an intuitive first-order 
sentence written with respect to the above rules and only containing constants 
from $C({\cal H})$ and the basic set-theoretic binary operators $\in,\ =,
\ \cup,\ \cap,\ \subset,\ \times,\ -$ etc. Then $\Hyper {\Phi}$ is the same 
sentence as $\Phi$ except \underbar{every} constant that appears in $\Phi$ is now 
proceeded\break
by a *. In this process$,$ the conventions as to when the * may be dropped are 
also followed.\par
\smallskip
\hrule
\vskip 18pt
The sentence $\Hyper {\Phi}$ is called the {\bf *-transform} of the sentence
$\Phi.$ The mathematical object $\Hyper {\cal H}$ is actually obtained 
by application of the {compactness theorem} for a formal first-order language
or by an algebraic construction
called the {ultraproduct} construction [{\tt Stroyan and Luxemburg [1976]}] 
and a process known as the 
 {Mostowski collapse} [{\tt Barwise [l977]}$,$ {\tt Herrmann [1986]}]. 
 The Leibniz Principle is an immediate consequence of the 
fundamental definition of what it means to say that  ${\cal H}$ and $\Hyper  
{\cal H}$ are the universes for first-order models.
You can learn about first-order models by referring to 
any good book in mathematical logic.\par
\vskip 18pt
{\bf Theorem 3.4.1 (The Leibniz' Principle).} {\sl A sentence $\Phi$ holds true
for members of ${\cal H}$ if and only if the sentence $\Hyper {\Phi}$ holds true for 
members of $\Hyper {\cal H}.$}\par                                       
\vskip 18pt
What Theorem 3.4.1 says is that if we let $K$ be the set of all of those 
specially written sentences that 
hold true for members of ${\cal H}$$,$ then the members of $\Hyper {\cal H}$
along with the basic set-theoretic operators form a model for the set of
all sentences obtained from $K$ by *-transfer. 
 {\it Probably much more significant is the fact that infinitely many other 
statements not obtainable by *-transfer hold true for members of 
$\Hyper {\cal H}.$} Each member of $\Hyper {\cal H}$ is called an 
{\bf {internal}} entity
(or internal individual if it is a member of $\hyperreal$). Observe that 
an object $A$ is internal if and only if there exists some $p \in  \nat$ such 
that $A \in \hyper X_p.$ Consequently$,$  
the extended standard sets are internal entities; but$,$ there are many 
internal entities that are not extended standard sets and internal 
individuals that are not individuals. Shortly the additional terms {\it extended
standard} and {\it internal} are more fully exploited. But$,$ first a few 
examples.\par  
\vskip 18pt
{\bf {Example 3.4.1.}} Let $S_1,\ S_2$ be two sequences. The definition of the 
addition of two such functions can be stated as follows: for every $x$ if
$x\in \nat,$ then the function $C = S_1 + S_2$ if and only if $C(x) = S_1(x) +
S_2(x).$ The *-transform becomes: for every $x$ if
$x\in \hypernat,$ then the function $\Hyper C = \hyper S_1 + \hyper S_2$ 
if and only if $\Hyper C(x) = \hyper S_1(x) + \hyper S_2(x).$ Notice that such 
forms as $\Hyper C(x)$ mean that the * is applied to the $C$ only. Thus *-
transfer extends to all of the ordinary definitions of this type.\par

\vskip 18pt  
{\bf {Example 3.4.2.}} Let's look at a few of the parts of Theorem 3.1.3. 
(i) Using just constants themselves is legal. Thus if $a,A\in C({\cal H}),$ and
the expression $a\in A$ holds for ${\cal H}$$,$ then the *-transform is simply the 
expression $\hyper a\in \hyper A.$ (iii) Notice that if $\{A_1,\ldots,A_n\}$ 
is a finite set of members of $\cal H$$,$ then there is a constant $C$ in 
$C({\cal H})$ that represents this set. Clearly$,$ it is
 unnecessary to mention continually that a symbol is a constant on one 
hand and a name for an object in ${\cal H}$ on the other. The context of a 
statement will usually serve to indicate a constants usage. This first-order 
statement is: for every $x,\ x\in C$ if and only if $x = A_1 \lor 
x=A_2\lor\cdots \lor x = A_k.$ This is a collection of finitely many symbols and 
is thus allowed. The *-transfer becomes: for every $x,\ x\in \hyper C$ if and only
if $x = \hyper A_1 \lor x=\hyper A_2\lor\cdots \lor x = \hyper A_k.$ [Note: 
translate $\lor$ by the word  ``or.''] \par
[Note: Our defined operators $+,\ \cdot,\ \vert \cdot\vert,\ \Vert \cdot\Vert,$ etc. are usually considered but the *-transfer of the operators (the nonstandard extension to $\hyperreal,\ \hyperreal^n$) as
they are defined on $\real$ and $\real^n.$]\par
\vskip 18pt
\hrule
\smallskip
\hrule
\smallskip
{\bf {Example 3.4.3. An argument}.} *-transfer is certainly 
important$,$ but it is our ability to argue by using statements that contain 
symbols that are \underbar{not} obtainable by *-transfer which is the key to
 nonstandard 
analysis. Indeed$,$ as will be established$,$ there are three different categories
of constant mathematical symbols used in the following argument. Let
$\Gamma \in \hypernat-\nat.$ [$\Gamma \in C(\Hyper {\cal H})$ and there 
is no symbol for $\hypernat - \nat$ in $C({\cal H}).$]
Since $\Gamma \notin \nat$ [it turns out that $\nat \notin \Hyper {\cal H})$]
then for each $x\in 
\nat,\ x < \Gamma$ for if not$,$ then there exists some $n\in \nat$ such
that $0<\Gamma < n$ and for each $x\in \nat,\ x \not= \Gamma.$  However$,$  
we also know that the set $\{x\vert 0<x<n\}$ is a finite set and thus each 
member is an element of $\nat$ (even under *-transfer by our conventions). 
This would yield a contradiction. Now since $\real$ [$\real \in C({\cal H})$]
is Archimedean then for each $r \in \real$ there exists some $n\in \nat$ such 
that $\vert r \vert < n.$ Hence$,$ $\vert r \vert < \Gamma$ implies that
$\Gamma \in {\nat}_\infty.$ Therefore$,$ $\hypernat - \nat \subset 
{\nat}_\infty.$ Since it is obvious that ${\nat}_\infty \subset \hypernat -
\nat$ then this implies that ${\nat}_\infty = \hypernat -\nat.$\par
\smallskip
\hrule
\smallskip
\hrule
\vskip 18pt
{\bf Theorem 3.4.2.} {\sl The set of infinite natural numbers ${\nat}_\infty =
\hypernat - \nat.$}
\vskip 18pt
Mathematicians are so practiced in arguments such as illustrated in example 
3.4.3 that it is often not realized that there are at least five categories of 
constants as symbols being used. (1) There are the {\bf {standard (unstarred) 
constants}} that appear in $C({\cal H})$. (2) The {\bf {extended standard 
constants}} 
that appear in $C(\Hyper {\cal H}).$ (3) The {\bf {internal constants}} 
that appear in $C(\Hyper {\cal H})$ but are not obtained by *-transfer. 
(4). Then {\bf {external constants}}$,$ such as ${\nat}_\infty$ that represent 
mathematical objects but do not fit categories (1) (2) (3). Last$,$ there are 
the constants the comprise that word forms of the  ``metalanguage'' that 
discusses and establishes things about the mathematical objects$,$ such as the 
phrase ``there is not a symbol for.'' Theorem 3.1.2 and part (vii) of Theorem 
3.1.3 indicate the only instances where two categories - (1) and (2) - overlap 
in the sense that there are definitely two distinct constants representing the 
same mathematical object. With respect to the next definition$,$ these basic 
facts about the usage of the constants - especially their correspondence 
to various mathematical objects - should always be kept in mind.\par
\vskip 18pt
\hrule
\smallskip
{\bf Definition 3.4.2. (Internal$,$ External.)}  Any member of $\Hyper {\cal H}$
will be termed$,$ in 
general$,$ an {\bf {INTERNAL ENTITY}} or {\bf {INTERNAL INDIVIDUAL}}. Any
$A \subset \hyper X_p $ that IS NOT INTERNAL is termed
an {\bf {EXTERNAL ENTITY}.}\par
\smallskip
\hrule
\vskip 18pt
Definition 3.4.2 now allows for a somewhat better but equivalent statement of
the Leibniz Principle for objects such as $\Re.$ \par
                              
{\bf Theorem 3.4.3. {(Leibniz Principle Restated)}.} {\sl  A sentence $\Phi$
holds true
for members of $\Re = {\cal U} \cup {\cal P}({\cal U})$ if and only if the
sentence $\Hyper {\Phi}$ holds true for the members of $\hyper {\eskip \cal U}$
or internal subsets of $\hyper {\eskip \cal U}.$}\par    
\vskip 18pt                                   
\leftline{3.5 \underbar{A Few {Simple Applications}.}}
\medskip
Using infinitesimal procedures$,$ there are simple and specifically describable
methods for determining the limit of a sequence$,$ where Theorem 
3.2.1 allows us to use the expression ``limit of a sequence'' as an
abbreviation for the phrase ``real number to which the a sequence converges.''
The next applications amply illustrate the most significant of these 
procedures.\par
\medskip
{\bf Application 3.5.1.} {\it This is an example of the direct application of 
the algebraic properties of the infinite and infinitesimal numbers.}\par
Let $p\in \nat,\ p > 0$ and assume that for each 
$0 < n\in \nat,\ S_n = (1/n)^p.$ Then $S \to 0.$\par
Proof. Let $\Gamma \in \hypernat - \nat.$ Then $(1/\Gamma)\in \monad 0$ 
implies that $(1/\Gamma)^p \in \monad 0.$ Thus result follows from Definition 
3.2.1.\par
\medskip
{\bf Application 3.5.2.} {\it A {bounding method}.}\par
Let $a\in \real,\ 0 < \vert a \vert < 1.$ For each $0 < n\in \nat,$ let $S_n =
a^n.$ Then $S \to 0.$\par
Proof. By induction it is not difficult to show that if $0\leq b\in \real,$ 
then for each $n\in \nat,\ (1 + b)^n \geq (1+nb).$ Let $b = (1/\vert a \vert) 
- 1.$ Then $b > 0$ and $\vert a \vert = 1/(1 + b).$ Thus for $n\in\nat,$
$$0 <(\vert a\vert)^n ={1\over{(1 + b)^n}} \leq {1\over{1 + nb}}< {1\over{bn}}
,\ (n >0).$$
Thus by *-transfer for each $n \in \hypernat $
$$0<(\vert a \vert)^n < (1/b)(n)^{-1},\ (n > 0).$$
In particular$,$ for each $\Gamma \in \hypernat - \nat,$
$$0<(\vert a \vert)^\Gamma < (1/b)(\Gamma)^{-1}.$$
The result now follows since $a^{\Gamma}\in \monad 0.$ \par
\medskip
{\bf Application 3.5.3} {\it A {bounding method$,$ redefinition} and the 
standard part operator.}\par
Let $1 < a \in \real.$ For each $0 < n\in \nat,$ let $S_n = (a)^{1/n}.$ Then 
$S \to 1.$\par
Proof. Define the sequence $Q_n = a^{1/n}-1,\ 0 < n\in \nat.$ Then 
$$a=(1 + Q_n)^n = 1 + nQ_n + other\  positive\  terms,$$ where $0<n\in \nat.$
Hence $a > nQ_n,\ \forall n\in \nat,\ (n >0).$ By *-transfer $a > nQ_n,\ \forall
n\in \hypernat\ (n>0).$ In particular$,$ $ a > \Gamma \hyper Q_\Gamma,\ \forall \Gamma 
\in \hypernat - \nat.$ Hence $$\forall \Gamma \in \hypernat - \nat,\ 0
< \hyper Q_\Gamma < a(1/\Gamma) \in \monad 0.$$
Thus $Q \to 0$ implies that $\st {\hyper Q_\Gamma} = 0 = \st {\hyper S_\Gamma} 
- \st 1 = \st {\hyper S_\Gamma} - 1,\ \forall \Gamma \in \hypernat - \nat.$
The result follows from Theorem 3.2.2.\par
\medskip
{\bf Application 3.5.4.} {\it A bounding method$,$ redefinition and the 
standard part operator.}\par
Let $S_n = {\root n \of n},\ \forall n\in \nat,\ (n>0).$ Then $S \to 1.$\par
Proof. Let $Q_n = S_n -1 = {\root n \of n} -1,\ \forall n\in \nat,\ (n>0).$ Then
$$ n = (1 + Q_n)^n \geq (n(n-1)/2)Q_n^2,\ \forall n\geq 2,$$
by the binomial expansion. Hence $$ 0\leq Q_n\leq \sqrt {2\over{n-1}},\ \forall 
n\geq 2.$$                                                                    
By *-transfer$,$ $$ 0\leq Q_\Gamma\leq \sqrt {2\over{\Gamma-1}},\ \forall \Gamma \in 
\hypernat - \nat.$$                       
But$,$ $ \sqrt {2\over{\Gamma-1}} \in \monad 0,\ \forall \Gamma \in 
\hypernat - \nat$ implies that $Q \to 0$; which implies that $S \to 1.$\par
The bounding techniques and algebraic manipulations illustrated by above
examples are actually of the same type that were used originally with old 
style infinitesimal analysis to argue for these results. As previously 
mentioned the standard part operator establishes the classical limit theorems 
which$,$ of course$,$ can be applied always. However$,$ the standard part operator 
can be applied directly without referring to the limit theorem at all.\par
\medskip
{\bf Application 3.5.5} {\it {Using the standard part operator}.}\par
Find the limit of the sequence $S_n =((1/n)^{10}) {\root n \of n},\ n > 0.$ 
Let arbitrary $\Gamma \in \hypernat -\nat.$ Then $\st {1/\Gamma} = 0$ implies 
that $\st {(1/\Gamma)^{10}} = 0.$ From application 3.5.4$,$ $\st {(\Gamma)^{1/\Gamma}}=
1.$ Thus $\st {\hyper {S_\Gamma}} = 0 \cdot 1 = 0.$ Hence$,$ $S\to 0.$\par
One of the most significance statements made by a researcher is ``What 
if...?'' The basic propositions of infinitesimal analysis tend to lend 
themselves to many 
such ``{What if...?}'' type questions. Consider$,$ for example$,$ Theorem 3.2.2. 
One is prone to ask; what if $\forall \Gamma \in {\nat}_\infty,\ 
\st {\hyper S_\Gamma} \in {\cal O}^n$?  Or$,$ what if there exists some $\Gamma 
\in {\nat}_\infty$ and $\st {\hyper S_\Gamma} \in {\cal O}^n$? We state a theorem 
that represents an interesting result relative to the last ``What if 
...?'' - a result established in appendix 3 solely by means of the standard 
part operator.\par
\bigskip
{\bf Theorem 3.5.1.} {\sl Let $S\colon \nat \to \real.$ If $S$ is an 
increasing {\rm [}resp. decreasing{\rm ]} sequence and there exists some
$\Gamma \in {\nat}_\infty$ such that $\hyper S_\Gamma \in \cal O$$,$ then $S
\to \st {\hyper S_\Gamma.$}}
\vfil\eject
\centerline{Chapter 4.}
\medskip                
\centerline{\bf SOME MODELING WITH}
\centerline{\bf THE INFINITE NUMBERS}
\bigskip
\leftline{4.1 \underbar{Historical Confusion.}}\par    
\medskip 
First$,$ I point out that the infinite numbers are completely different from the 
idea of extending the real numbers by adjoining the new objects $\pm \infty$ 
to $\real$ and impressing upon these objects certain topological and algebraic 
properties. This can be readily shown by considering the concept of when
$\lim_{n \to \infty} = +\infty.$ Let ${\real}_\infty^+$ denote the positive 
infinite hyperreal numbers.\par
\vskip 18pt
\hrule
\smallskip
{\bf Definition 4.1.1. ($S \to +\infty$).} Let $S\colon \nat \to \real.$ Then 
$S \to +\infty$ if for each $\Gamma \in {\nat}_\infty,\ \hyper S(\Gamma) \in
{\real}_\infty^+$.\par
\smallskip 
\hrule
\vskip 18pt
{\bf Theorem 4.1.1.} {\sl Definition 4.1.1 is equivalent to the limit 
definition for a sequence $S \to +\infty.$ }\par
\vskip 18pt
Theorem 4.1.1 indicates that the values of $\hyper S$ are scattered throughout 
the set ${\real}_\infty^+.$\par
\vskip 18pt
{\bf Theorem 4.1.2.} {\sl Let $r \in {\hyperreal}^+,\ S\colon \nat \to \real$
 and $S\to +\infty.$ Then for each
$\Omega \in {\nat}_\infty$ there exists some $\Delta \in {\nat}_\infty$ such 
that $\hyper S(\Omega) + r \leq \hyper S(\Delta).$} \par
\vskip 18pt
Thus the sequential property that $S \to +\infty$ does not correspond to the
idea that a sequence might ``converge'' to an infinite number. Indeed$,$ 
we have not even mentioned the idea of a monad about an infinite number even 
though it is possible to extend monad theory to cover such cases. This is 
particular significant for applied modeling when the infinite series is 
discussed relative to the sequence of partial sums each infinitely small.\par
   Recall that except for modern times most mathematics was 
almost {used exclusively for applied physical or geometric modeling} and was 
not consider as a study of abstract entities. 
 {De l'Hospital} apparently believed in the {objective reality of both the
infinitesimal and infinite} quantities as well as the {existence of sets that
contained infinitely many members}. To him$,$ they existed in the natural world. 
When de l'Hospital wrote that it was a requirement that one must regard a 
curve as a totality of an infinity of straight line segments$,$ each infinitely 
small: or {\it``(which is the same) as a polygon with an infinite number of 
sides$,$ each infinitely small$,$  
which determine by the angle at which they meet$,$ the curvature of the 
curve...''} he apparently meant that all of these intuitively expressed objects 
exist in reality. Leibniz did not approve of this interpretation. Over and 
over again$,$ Leibniz proclaimed that such objects as the infinitesimal or 
infinite numbers were ``ideal'' or imaginary. He claimed that they are 
theoretically useful but that they did not correspond to real natural things. 
Moreover$,$ Leibniz apparently accepted the concept of a \underbar{potentially} 
infinite set as well$,$ rather than the possibility of an objectively real infinite 
set of objects. Almost all of the ancient controversy as to the reality of such 
notions apparently came about as a direct result of the basic philosophical 
and theological predilections of the investigators. Their general philosophic 
belief systems were reflected in both their scientific and mathematical 
views.\par 
   With respect to the physical possibility of there being objectively real
objects that are characterized by infinitesimals or even infinite numbers$,$ 
Robinson has replace these abstruse 
philosophical considerations with the following observation. In the first 
fundamental paper delineating his theory$,$ he wrote: ``{\it For phenomena on a 
different scale$,$ such as are considered in Modern Physics$,$ the dimensions of a 
particle may not be observable directly. Accordingly$,$ the question whether 
or not a scale of non-standard analysis is appropriate to the physical world 
really amounts to asking whether or not such a system provides a better 
explanation of certain observable phenomena than the standard system of real 
numbers. The possibility that this is the case should be borne in mind.''} 
{\rm [Fine Hall$,$ Princeton University]} {\tt Robinson [1961]} \par   
As previously mentioned Zeno proposed his famous physical paradox of Achilles' 
and the Tortoise in their never ending (potentially infinite time) foot race.
It is claimed$,$ that this paradox is resolved by modern mathematics through 
application of the infinite series. Such a series can supposedly model this 
physical foot race - a statement that is obviously false. Such a series would 
require the non-mythlogical object$,$ the Tortoise$,$ and the mythological 
Archilles to have variable ``sandal'' sizes - sizes that decreased to the 
point of being unmeasurably small.\par
    De l'Hospital accepted a curve as being
identifiably the same as a polygon with a fixed infinite number of line
segments comprising its sides.  On the other hand$,$ {Eudoxus [370 BC]}
devised the method of exhaustion which assumes the true existence of a
finite sequence of inscribed and circumscribed polygons. In general$,$ for a
closed non-polygonal curve none of these Eudoxus polygons were
considered to be 
the curve under investigation; but$,$ rather$,$ 
by the ``continuity process'' they would continually squeeze the curve between 
these two types of polygons and ``exhaust'' the space in between.
By this process the length of a curved segment was conceived of as an
intuitive sequence composed of portions of the polygon's perimeters.
Thus developed the idea of
a partial sum that 
represented the sum of the lengths of the sides of an n-gon - a finite sum 
that remained finite but acquired more and more terms.  Those that employed 
this method often guessed at a specific formula then justified their guess
by indirect and not direct argument. In modern times$,$ Planck described the 
accepted procedure for modeling the behavior of a natural system when he wrote
that: {\it ``a finite change in Nature always occurs in a finite time$,$ and 
hence resolves into a \underbar{series} of infinitely small changes which 
occur in successive infinitely small intervals of time.''} What might the 
term ``series'' mean in Planck's statement? \par                         

Leaving aside the ontological question associated with the notion of the
infinite$,$ it will be demonstrated that Robinson's theory of the infinitesimal
and infinite hyperreal numbers 
brings a concrete and rigorous language to the above vague methods of physical 
and geometric modeling. Unfortunately$,$ in doing so$,$ many of these previous ideas
will require modification.\par 
\bigskip
\leftline{4.2. \underbar{The Internal Definition Principle.}} \par
\medskip
Clearly there is a need to acquire a better understanding of the relationship
between the concepts of the 
infinite series$,$ the sequence of partial sums$,$ the internal process of 
partial summing with its arithmetic and the notion expounded by Planck. 
But$,$ prior to examining these concepts$,$ yet another technical procedure needs 
to be discussed.\par
In definition 3.4.1$,$ our first-order language was extended to include 
the basic set-theoretic binary operators. As explained in the Appendix to 
Chapter 3 immediately after the proof of Theorem 3.1.3$,$ our first-order 
language may also include the symbol $(\cdot,\ldots,\cdot)$ for n-tuple 
formation where the coordinates are either variables or constants.
Under *-transfer these operators and the n-tuple formation symbol are not 
starred in any $\Phi$ that has been properly formulated with the variables 
restricted to entities in $\cal H.$ \par
\medskip
{\bf {Example 4.2.1.}} Assume that you are given some relation $R$$,$ a nonempty 
set $A$ and the next statement. For each $x$ if $x \in {\cal 
P}({\cal U}),$ then $(x,A) \in R$ and $x \subset A.$ Or$,$ as a formal 
statement $\forall x(x \in {\cal P}({\cal U}) \to (x,A) \in R \land 
x \subset A).$ Then the *-transfer would read: For each $x$ if $x \in \hyper 
{\cal P}({\cal U}),$ then $(x,\hyper A) \in \hyper R$ and $x \subset \hyper 
A.$\par
\bigskip
Please note that we needed to star the symbol ${\cal P}({\cal U})$. This would 
also be the case if a variable appeared where the $\cal U$ appears.\par
 When a collection $C$ of mathematical expressions are written they include
``mathematical'' variables and constants.
Within our standard set theory the constants represent individuals or 
entities. From the view point of the mathematical structure$,$ the entities 
represent such things as basic sets$,$ operators$,$ relations$,$ specific n-tuples 
and other definable objects. In general$,$ the variables either vary over 
every member of some set or represent some one unknown element. It is common 
practice for mathematicians to use different variable 
symbols to represent elements of distinct sets or distinct positions in
n-tuples and the like. Within $C$ certain of these 
variables may also be ``quantified.'' What this means is that if $v$ is one 
of the variables in $C$$,$ than there also appears in $C$ in the customary 
location the phrase ``for all $v$'' (i.e. $\forall v$) or ``there exists
some $v$'' (i.e. $\exists v$). Any variable in $C$ that is not associated with
some quantifier is termed a {\bf {free variable}.} Now in formal logic the 
concept of the free and quantifier bounded variable must be more 
carefully described since formally the same variable can appear both in a free 
and not free position. Since we are trying to be as non-technical as possible 
and since it is customary to use many different variable symbols in 
mathematical prose$,$ this somewhat vague free variable definition should 
suffice. If it still is confusing$,$ then most elementary logic books have a 
more formal presentation.\par
\vskip 18pt
{\bf {Example 4.2.2.}}  When the 
definition of continuity is expressed one usually states the following: The 
function $f\colon A \to \real$ is continuous at $p \in A$ if for each $\eps 
\in {\real}^+$ there exists some $\delta \in {\real}^+$ such that whenever 
$0\leq \vert x - p \vert < \delta$ and $x \in A,$ then $0\leq \vert f(x) -f(p) 
\vert < \eps.$ The symbols $f,\ A,\ \real,\ 
 {\real}^+,\ p,\ -,\ \vert,\ \leq,\ <,\ 0\ $ are all considered as constants.
The 
symbols $\eps,\ \delta,\ x\ $ are variables. This entire definition may be 
formally expressed by 
$\Phi =\forall x(x \in {\real}^+ \to \exists y(y \in {\real}^+\land \forall z(
z \in A \land 0\leq \vert z - p \vert < \delta \to \vert f(z) - f(p) \vert < 
 \eps))).$ 
\vskip 18pt
Suppose that you are interested in the subset of $A$ in example 4.2.2 that 
comprises all of the points of continuity. In this case$,$ the constant $p$ is 
considered a variable and one writes such a set as $\{p\vert p\in A \land 
 {\Phi}(p) \}.$ Built into our set builder notation is the quantifier 
$\forall.$ ``The set \underbar{of all} $p \in A$ such that ${\Phi}(p)$ holds 
true.'' The expression ${\Phi}(p)$ is now considered a formula in one 
variable $p.$ Note that $p$ is a free variable. Of course$,$ these set builder 
formulas should present no difficulties since these language constructions 
are the ordinary and customary ones used by the mathematical community.
Two more examples should sufficiently illustrate this easily grasped relation
between our first-order language and common mathematical usage. \par
\vskip 18pt
{\bf {Example 4.2.3.}} Most of the 
time when we write mathematical formula in 
variables they are considered as free. Such an expression as $y = 3x$ has been 
defined for many different structures. Even though to obtain the graph of such 
an expression one needs to know the domain and codomain$,$ this is not necessary 
until the structure itself is considered. Thus the graph may be the 
$\{(x,y)\vert x \in [0,3] \subset \real \land y \in \real \land y = 3x\}.$ On 
the other hand$,$ we might also have $\{(x,y)\vert x \in \realp 7 \land y \in 
 \realp 7 \land y = 3x \}.$ Thus $y = 3x$ is a mathematical formula$,$ 
${\Phi}(x,y),$ in two variables with a constant 3 and a constant operator
(multiplication) that can be used to generate many different sets in our
set theory.\par
\vskip 18pt
 Formulas that appear in set builder notation can also contain quantifiers 
as the next example indicates.
\vskip 18pt
{\bf {Example 4.2.4.}} Let the sets
$A,B \in X_p$ and let
$ B^A$ denote the set of all functions with domain $A$ and codomain $B.$ 
Then $B^A \in X_{p+3}.$ One might want to consider a special subset of $B^A$ 
defined by
$\{x\vert x \in B^A\ and\ there\ exists\ y \in A \ such\ that\ x(y) > 3 \}$
 (i.e. $\{x\vert x \in B^A\land \exists y(y \in A \land x(y) > 3 \})$ or
 $\{x\vert x \in B^A\land \exists y \exists z(y \in A \land z \in B \land 
(y,z) \in x \land z > 3 \}).$  This set exists by the axiom of comprehension 
and every though the defining expression contains more than one variable,
all but the $x$ are bounded by a quantifier that appears immediately to the
left.\par
\vskip 18pt  
The basic reason for discussing and presenting these simple examples lies in  
two most interesting results that allow for the set-theoretic generation of all 
extended standard or internal sets and n-ary relations by means of such set 
builder formulas. Of course$,$ n-ary relations are sets but are singled out 
specifically due to their obvious usefulness.
Recall that $C({\cal H})$ denotes the set of all constants that signify members 
of the set $\cal H$ while $C(\Hyper {\cal H})$ the names for members of $\Hyper 
{\cal H}.$ A formula in 
our first-order language is called standard [resp. internal] {\bf bound} if each 
quantified 
variable is restricted to an object represented by a constant in $C({\cal H})$ [resp. $C(\Hyper {\cal H})$]. What this means is that, for a $\cal H$ [resp. $\Hyper {\cal H}$] interpretation of the formula, each bounded variable must be interpretation as varying over a set contained in $\cal H$ [resp. $\Hyper {\cal H}$]. The formula in Example 4.2.2 is bound. A formula like $\forall x \exists y(x \in y)$ is not bound. \par
\vskip 18pt
{\bf Theorem 4.2.1. (The Extended Standard Definition Principle)}\par 
 {\sl \indent\indent (i) A set $A$ in our set theory is an extended standard 
set (i.e. there exists some $B \in C({\cal H})$ such that $A = \hyper B$) if 
and only if there exists some standard set $D$ and a standard bound formula 
$\Phi (x)$ in one free variable where each constant in $\Phi (x)$ is a member
of $C({\cal  H})$ and
$$ A = \{x\vert x \in \hyper D \land \Hyper {\Phi}(x) \}.$$ \pars
\indent\indent (ii) A set $A$ is an extended standard n-ary relation ($n > 1$)
if and only if there exist n standard sets $D_1,\ldots,D_n$ and a standard bound 
formula $\Phi (x_1,\ldots,x_n)$ in n free variables where each constant in 
$\Phi (x_1,\ldots,x_n)$ is a member of $C({\cal H})$ and 
$$A = \{(x_1,\ldots,x_n)\vert
 x_1 \in \hyper D_1 \land \cdots \land x_n \in \hyper D_n \land \Hyper 
{\Phi (x_1,\ldots,x_n)}\}.$$}\par
\vskip 18pt
{\bf Theorem 4.2.2. (The Internal Definition Principle)}\par 
 {\sl \indent\indent (i) A set $A$ in our set theory is an internal 
set if and 
only if there exists some internal set $D$ and a internal bound formula $\Phi (x)$ in 
one free variable where each constant in $\Phi (x)$ is a member of 
$C(\Hyper {\cal H})$ and
$$ A = \{x\vert x \in D \land \Phi (x) \}.$$ \par
\indent\indent (ii) A set $A$ is an internal n-ary relation ($n > 1$)
if and only if there exist n internal sets $D_1,\ldots,D_n$ and an internal bound 
formula $\Phi (x_1,\ldots,x_n)$ in n free variables where each constant in 
$\Phi (x_1,\ldots,x_n)$ is a member of $C(\Hyper {\cal H})$ and
 $$A = \{(x_1,\ldots,x_n)\vert
 x_1 \in D_1 \land \cdots \land x_n \in D_n \land \Phi (x_1,\ldots,x_n) \}.$$}\par
\vskip 18pt
Internal objects are basic to nonstandard analysis since it is only internal 
properties that hold for the model $\Hyper {\cal H}.$ As far as physical modeling 
 is concerned {\bf internal objects represent NSP-world effects that 
directly or indirectly yield the observed natural world behavior being
modeled approximately by a standard mathematical structure.}
\bigskip
\leftline{4.3. \underbar{Hyperfinite Summation.}} \par
\medskip
    In example 4.2.4 the set-theoretic notation $B^A$ is utilized to denote 
the set of all functions with domain $A$ and codomain $B$. As previously 
mentioned there are procedures that require objects in $\cal H$ that may not 
be members of $\Re.$ Assume that ${\cal A},\ {\cal B}$ are 
sets of subsets and ${\cal A},\ {\cal B} \in X_p.$ Now let the sets $A\in {\cal 
A},\ B \in  \cal B.$ Then $A,B \in X_{p-1},\ (p \geq 2),\ A \cup B \subset
 X_0 \cup X_{p-2}$ implies that $A \cup B \in X_{p-1}$ and if $a \in A,\ 
b\in B,$ then $\{a\},\{a,b\} \in X_{p-1}.$ Hence $(a,b) \in X_p.$ 
Thus if $f \in B^A,$ then $f \in X_{p+1}.$ Indeed$,$ $A \times B \in X_{p+1}.$
Therefore$,$ $B^A \in X_{p+2}.$ It is possible to  
consider the formation of each set $B^A$ as an operator $F(x,y)$ 
where $x \in \cal A$ and $y \in \cal B.$ 
The *-transfer process can be extended to this operator in the same manner as 
is done with the $\cup,\ \cap,\ \times,\ -,\ (\cdot,\ldots,\cdot)$ operators in 
the sense that is it not starred when written entirely in variable form. For 
specific members of $C(\Hyper {\cal H}),$ if 
$D \in \hyper {\cal A}$ and $E \in \hyper {\cal B},$ then $\hyper F(D,E) 
 \in \hyper X_{p+2}$ and $\hyper F(D,E)$ 
 is the internal set of all internal functions
with domain $D$ and codomain $E.$ This last fact comes from the *-transfer 
of a general characterization for the set-theoretic concept of $x^y.$\par
\vskip 18pt
{\bf {Example 4.3.1.}} Let $p \in \hypernat.$ Then the set 
$\{x\vert x \in \hypernat \land 0 \leq x \leq p\} = [0,p]$ is an internal 
subset of $\hypernat.$ Notice that the symbol $[0,p]$ can be used as an 
abbreviation for this set's defining property (i.e. $x \in [0,p]$ if and only if  $x \in 
\hypernat \land 0\leq x \leq p$) which can be restated in an 
appropriate first-order expression and substituted for the notation $[0,p].$
With this in mind it is clearly possible to now consider $p$ as a 
variable.\break
\hrule
\smallskip
Hence it follows that $BHF =\{y\vert y \in \Hyper {\cal P}(\nat) \land
\exists p(p \in \hypernat \land y = [0,p])\}$ is an internal (indeed$,$ an 
extended standard set) of subsets
of $\hypernat.$ The set $BHF$ is called the set of all {\bf {basic 
hyperfinite subsets} of} $\hypernat.$ \par
\smallskip
\hrule
\vskip 18pt
A nonempty $A\subset B$ is 
{\bf {finite}} if there exists some $n \in \nat$ and a function 
$f\colon [0,n] \to B$ such that the range of $f$ (i.e. $P_2(f)$) $= A.$ 
The intuitive idea of a finite set appears to be equivalent to this 
functional definition if you are willing to accept such things as the 
``finite''
axiom of choice and elementary procedures of recognition. The intuitive idea 
of the finite is based upon the human recognition of a distinction between
symbols written on paper as they are considered geometric forms and 
consequently it is related to the most basic aspects
of concrete geometry. It is this accepted recognition of the differences 
between geometric forms that than allows one to give a concrete meaning to a
correspondence between these forms and $\nat.$ One does not really establish 
that such a correspondence exists but its existence is accepted as part of 
the metamathematical methods.\par
 Theorem 3.1.4 states 
that for any standard function $f,\ \Hyper (P_i(f)) = P_i(\hyper f),\ i =1,2.$
The projections $P_i$ can be considered as maps from the set of all nonempty 
subsets of $A_1 \times \cdots \times A_n$
 into the sets $A_i$$,$ where $1\leq i \leq n.$ A modification of proof of 
Theorem 3.1.4 yields\par
\vskip 18pt
{\bf Theorem 4.3.1.} {\sl Let $n > 1.$ Then for each internal $R \subset
 \hyper A_1 \times \cdots \times \hyper A_n$ and for each 
$i,\ 1\leq i \leq n$ it follows that $\hyper P_i(R) = P_i(R),$ where 
$P_i(R)$ is an internal subset of $\hyper A_i.$}\par
\vskip 18pt
\hrule
\smallskip
\hrule
\smallskip
{\bf Definition 4.3.1 (Hyperfinite).} An 
internal subset $A$ of a set $\hyper B$ is {\bf hyperfinite} if it is empty 
or there exists 
some $[0,p] \in BHF$ and $f \in \hyper F([0,p],\hyper B)$ such 
that $P_2(f) = A.$ \par
\smallskip
\hrule 
\smallskip
\hrule
\bigskip
{\bf Theorem 4.3.2.} {\sl Let $F(B)$ be the set of all finite subsets of $B$. 
Then nonempty $A \subset \hyper B$ is hyperfinite if and only if 
$A \in (\Hyper {F})(B)),$ where $F$ is considered as an operator that generates all of the finite subsets of a set.}\par
\bigskip
{\bf Theorem 4.3.3.} {\sl Any nonempty finite set of internal individuals or
entities is internal and hyperfinite.}\par
\vskip 18pt
With respect to Definition 4.3.1 the maps in each
$\hyper F([0,p],\hyper B)$ are internal and behave like internal (partial) sequences. 
The ordinary finite manipulation we do with finite sets of real numbers can 
be extended to the hyperreals by means of *-transfer and by describing these 
processes by means of (partial) sequences defined on various $[0,p].$ As far 
as the NSP-world is concerned hyperfinite sets have all the same set-theoretic 
first-order properties as the finite sets and we do not usually establish 
these basic hyperfinite properties each time they are first employed. Thus$,$ not 
only is the union of finitely many hyperfinite subsets of a set $\hyper B$
a hyperfinite subset of $\hyper B$ but the union of a hyperfinite collection
of hyperfinite subsets of $\hyper B$ is a hyperfinite subset of $\hyper 
B.$ {\it However$,$ from the external or metamathematical point of view most 
hyperfinite sets are not finite as the next result indicates.}\par
\vskip 24 pt
{\bf Theorem 4.3.4.} {\sl Let $A \in \Re$ and assume that $A$ is infinite.
Then there exists a hyperfinite set $F$ such that $F \not = A,\ F \not= 
\hyper A$  and  $A \subset F \subset \hyper A.$}\par     
\vskip 18pt
\hrule
\smallskip
\hrule
\smallskip
It is precisely the concept of the hyperfinite that leads not only
to a clear 
understanding of the processes that underlie a convergent infinite series$,$ but 
also leads to the basic notion of the integral and Planck's meaning of the 
term ``series.''\par  
\smallskip
\hrule
\smallskip
\hrule
\vskip 18pt

{Finite summation} of elements of $\real$ can be
consider a function  $\Sigma$
 defined on $\{{\real}^{[0,n]}\vert n \in \nat \}$ with values in $\real$ and 
a corresponding function for the finite summation of members of $\realp n.$
It is customary to express the values as $\sum_{i=0}^n a_i.$ By *-transfer if 
$\Gamma \in {\nat}_\infty,$ then the {\bf {hyperfinite sum}}
 $\Hyper \sum_{i=0}^\Gamma a_i \in \real.$ By convention$,$ the symbol
 $\Hyper \sum_{i=0}^\Gamma$ is written 
as $\sum_{i=0}^\Gamma.$ It is a simple matter to translate Definition 3.2.1 for 
sequential convergence into the following theorem for convergence of a 
infinite series. \par
\vskip 18pt
{\bf Theorem 4.3.5.} {\sl An infinite series$,$ $\sum_{i=0}^\infty a_i,$ converges 
to $r \in \real$ if and only if for each $\Gamma \in {\nat}_\infty$ it follows 
that $\sum_{i=0}^\Gamma a_i \in \monad r.$}\par
\vskip 18pt
My experience indicates that Theorem 4.3.5 is not the most paramount 
application 
of the concept of hyperfinite summation for either physical or geometric 
modeling and$,$ indeed$,$ does not correspond to Planck's description for a series 
of infinitely small changes. In the next section$,$ examples are discussed that 
substantially indicate the true character of this concept when applied to 
geometry or natural system behavior.\par
\bigskip
\leftline{4.4. \underbar{Continuity and a Few Examples.}}\par
\medskip
De l'Hospital's concept of what constitutes a {\bf curve} is not satisfactory for 
Robinson's theory. Recall that one acceptable analytical definition for
the notion
of a {curve} in $\realp n$ is the following: a curve is a continuous map 
$c\colon [0,1] \to \realp n.$ This is equivalent to considering $c$ as 
determined by n continuous coordinate
functions 
$x_i = f_i(t),\ 1 \leq i \leq n$ each defined on $[0,1] \subset \real.$ 
Of course$,$ the geometric curve $C$ 
determined by these functions is usually considered as the set 
$\{(x_1,\ldots,x_n)\vert t \in [0,1]\}.$ The *-transform of these 
defining functions leads to the functions $x_i = \hyper f_i(t),\ 1 \leq i 
\leq n$ each defined on $\Hyper [0,1] \subset \hyperreal$ and they generate the 
``hypercurve" $\Hyper C \subset \hyperrealp n.$ The analytic geometry of 
$\hyperrealp n$ is similar to the customary geometry except that it must be
considered 
 {non-Archimedean} in character. From the viewpoint of the geometry of the NSP-world$,$ if 
the $C$ is not linear$,$ then $\Hyper C$ is not linear and this would entail a 
necessary rejection de l'Hospital's infinitesimal description. As to the 
definition of continuity the following captures the envisioned belief that 
continuous functions preserve the infinitely close.\par 
\vskip 18pt
\hrule
\smallskip
{\bf Definition 4.4.1. ({Continuity}).} For any nonzero $n,m \in \nat$ and any 
nonempty $A \subset \realp n$ a function $f\colon A \to \realp m$
is {\bf CONTINUOUS} at $p \in A$ if $\hyper f[\monad p \cap \hyper A] \subset
\monad {f(p)}.$ Also$,$ $f$ would be {\bf {UNIFORMLY CONTINUOUS}} on $A$ if for 
each $p,q \in \hyper A$ such that $p \approx q,$ then $\hyper f(p) \approx 
\hyper f(q).$\par
\smallskip
\hrule
\vskip 18pt
Observe that uniform continuity seems to preserve the infinitely close in 
the most satisfactory manner. One the other hand$,$ pointwise continuity is a 
monad preserving property.\par
\vskip 18pt
{\bf Theorem 4.4.1.} {\sl Definition 4.4.1 for continuity and uniform 
continuity is equivalent to the classical $\delta - \eps$ definition.}\par 
\vskip 18pt
One of the most powerful ideas in elementary analysis is that of the compact 
subset [or subspace if you wish] of the space $\realp n.$ Rather than dwell 
upon the many equivalent standard definitions for this notion a direct 
nonstandard assault is very enlightening since it reveals immediately the 
relationship between compactness and continuity as well as an intuitive 
comprehension of what is being compressed or compacted.\par
\vskip 18pt
\hrule
\smallskip
{\bf Definition 4.4.2. ({Compactness}).} For any nonzero $n \in \nat$ a 
nonempty set $A \subset \realp n$ is {\bf COMPACT} if 
$\hyper A \subset \bigcup \{\monad r\vert r \in A\}.$\par
\smallskip
\hrule
\vskip 18pt
{\bf Theorem 4.4.2.} {\sl Definition 4.4.2 for compactness is equivalent to 
the standard definition utilizing open covers.}\par
\vskip 18pt
The reason that some 17'th century  geometers considered non-linear curves 
to be 
collections of infinitesimal line segments was in their desire to use 
infinitesimal analysis to measure a curves length by corresponding this 
measure to the polygons of Eudoxus. Even though de l'Hospital's original 
description is inadequate$,$ a modification does secure the accepted analytical 
results. \par
\vskip 18pt
{\bf {Example 4.4.1.A.}} Let $\Gamma \in {\nat}_\infty.$ Then $F =\{t_i\vert t_i 
= i/\Gamma \land 0 \leq i \leq \Gamma \}$ is an internal and hyperfinite 
subset of $\Hyper [0,1].$ By *-transfer, $F$ behaves like an ordered partition 
of the interval $[0,1]$ as defined in the standard sense. Such a set is 
termed a {\bf {fine partition}} (i.e. hyperfinitely many members of $[0,1]$ 
generating subintervals that are infinitesimal in length). The internal set 
$F$ generates the 
internal set of ``points'' $P = \{(\hyper f_1(t_i),\ldots,\hyper 
f_n(t_i))\vert t_i \in F\}$ that are members of the hypercurve $\Hyper C.$
Now for each $i = 0,\ldots,\Gamma - 1,$ and each $j,\ 0\leq j \leq n$ 
let $\hyper f_j(t_{i+1}) - \hyper 
f_j(t_i) = d(j,i).$ (If $c$ is continuous$,$ then each $d(j,i) \in \monad 0.$) 
 For each $i \in \hypernat$ such that $0 \leq i \leq 
\Gamma -1,$ the internal set $\ell_i = \{(x_1,\ldots,x_n)\vert
\forall j \in \hypernat,\ 0\leq j \leq n,\  x_j =
\hyper f_j(t_i) + t(d(j,i)) \land t \in 
\Hyper [0,1] \}$ is a hyperline segment connecting the two points 
$(\hyper f_1(t_i),\ldots,\hyper f_n(t_i)),\ (\hyper f_1(t_{i+1}),
\ldots,\hyper f_n(t_{i+1}))$ on the curve $\Hyper C.$ From this one obtains 
the internal hyperpolygonal curve ${\cal P}_\Gamma = \bigcup \{\ell_i \vert 
0\leq i \leq \Gamma - 1\}.$ As to the length of ${\cal P}_\Gamma $ simply 
extend the concept of length in the classical sense by defining for each 
$i = 0,\ldots,\Gamma -1$ the vector $\vec v_i = (d(1,i),\ldots,d(n,i)) \in 
\hyperrealp n.$ Then let the hyperfinite sum $\sum_{i=0}^{\Gamma - 
1} \Vert \vec v_i \Vert = \vert {\cal P}_\Gamma \vert \in \hyperreal.$ Even though$,$ 
in general$,$ you would have a different hyperpolygon with a different hyperreal 
length for $\forall \Gamma \in {\nat}_\infty$ {\tt Robinson [1966$,$ 84-86]} showed 
that if 
$c$ is continuously differentiable$,$ then for all $\Gamma \in {\nat}_\infty, \
\vert {\cal P}_\Gamma \vert \in \monad r$ and the real number $r$ was the 
length of the curve obtained in the classical sense by means of the 
integral.\par
\vskip 18pt
{\bf {Example 4.4.1.B.}} Under the same 
criterion as stated in example 4.4.1.A that $c$ is continuously 
differentiable$,$ the length of a curve 
is actually closer to the limit concept then it is to the de l'Hospital 
description. This is seen by simply following the same process but replacing 
$\Gamma$ with an arbitrary nonzero $n \in \nat.$ As $n$ increases this yields 
an increasing sequence $\vert {\cal P}_n \vert.$ Application of Theorem 3.5.1
implies that if there exists but one $\Gamma \in {\nat}_\infty$ such that 
$\vert \hyper {\cal P}_\Gamma \vert \in \cal O,$ then from example 4.4.1.A 
this sequences converges to $\st {\vert \hyper {\cal P}_\Gamma \vert}$ and
has the same standard part for all $\Gamma \in {\nat}_\infty.$ \par
\vskip 18pt
The process used to obtain the length of a curve in example 4.4.1.A shows that 
for most ordinary curves there are infinity many hyperpolygons that have the 
same standard part generated length. Thus each can be used as a NSP-world
representative for 
the curve itself$,$ at least as far as length is concerned. The fact that there 
does not$,$ in general$,$ exist a unique hyperpolygon is a disadvantage from the 
viewpoint of the founders of the infinitesimal method and forces a rejection 
of the de l'Hospital description. However$,$ if in the natural world a curve is 
envisioned to be a path of motion produced by physical processes$,$ then the 
lack of uniqueness could be a advantage. Under this interpretation$,$ it would 
indicate that different and possibly interesting NSP-world ultranatural 
processes yield the same effect when they are restricted to the natural 
world.\par
 There has arisen recently a significant application of the seeming 
esoteric idea that there may exist a multitude of distinct NSP-effects that 
yield the same natural world effects. Next is an example of how this might 
occur. \par
\vskip 18pt
{\bf {Example 4.4.2. (Fractals)}} Science
has become interested in order and design as reflected in what has become 
known as ``fractal'' behavior. For this example$,$ let $c\colon [0,1] 
\to \realp n$ be a curve; but$,$ assume that the sequence of polygon 
approximations$,$ $\vert {\cal P}_i\vert,$ discussed in example 4.4.1.B has the 
property that $\vert {\cal P}_i\vert \to +\infty.$ This is apparently one of 
the salient features of a fractal curve.
In {\tt Herrmann [1989]} it is shown that for any nonempty compact 
$K \subset \real$ and for any continuous $c\colon K \to \realp n$ there exists 
an internal $G \colon \Hyper K \to \hyperrealp n$ such that $G$ is *-
differentiable of any order $m \in \hypernat,\ G$ has a well-defined 
hyperreal length and $\st G = c.$ From the viewpoint of the NSP-world of 
processes and paths of motion$,$ this $G$ represents the same effects as does
the function $c$ except that $G$ has an associated length concept and is 
ultrasmooth. 
Moreover$,$ the internal object $G$ is somewhat less arbitrary in character 
than are those in example 4.4.1 since it is selectable from a specific 
algebra of functions. This may be significant since when it has been
rigorously shown that certain physical attributes are representable by 
fractal curves then the theory that models such physical attributes is usually 
associated with some specific algebra of functions. Thus a standard fractal 
curve may be replaced by the standard part of an ultrasmooth curve with a 
well-defined length.\par
\vskip 18pt
The idea of hyperfinite summation as representing the geometric length of a 
curve is$,$ of course$,$ closely associated with the elementary integral. Indeed$,$ 
when we quoted Planck's fundamental description for physical modeling 
a question was asked$,$ {``What might the term `series' mean in Planck's
statement?''} Apparently$,$ what Planck meant by this term is the type of hyperfinite 
summation that$,$ as seen in the next chapter$,$ produces the integral.\par                           
\vfil\eject
\centerline{Chapter 5.}
\medskip                
\centerline{\bf STANDARD RULES}
\centerline{\bf FOR INTEGRAL MODELING}
\bigskip
\leftline{5.1 \underbar{The Riemann Styled Integral.}}\par    
\medskip                                                  
The mathematical concept variously termed ``integration'' (i.e. to bring 
together the parts or to make whole) was$,$ until after the time of Cauchy$,$ always 
considered to be 
a specifically defined summation process. In 1823 Cauchy wrote the following
description$,$ where $f\colon [a,b] \to \real$ is assumed to be continuous. ``...if one divides $X-x_0$ into 
infinitesimally small elements $x_1 - x_0, x_2 - x_1,\ldots, X- x_{n-1}$ the sum
$$S = (x_1 - x_0)f(x_0) + (x_2  - x_1)f(x_1) + \cdots + (X -x_{n-1})
f(x_{n-1})$$ converges to a limit represented by the definite integral
$\int_{x_0}^X f(x)\,dx.$''
\indent From our new nonstandard point of view such a sum is produced by
a {hyperfinite partition} [see Example 4.4.1.A] generated 
by some $\Gamma \in {\nat}_\infty$ and the internal hyperfinite sequence of
values $f(x_i),\ 0 \leq i \leq {\Gamma -1}.$  Obviously$,$ the Riemann sum$,$ 
where $f$ may be evaluated at any member of each 
subinterval $[x_i,x_{i+1}],$ is styled after this Cauchy definition.         
However$,$ does the standard part 
of such a hyperfinite sum exist and is it independent of the partition?
Moreover$,$ can the concept be extended to bounded not necessarily continuous 
functions? \par
In the following very brief discussion, neither the Stieltjes nor Lebesgue 
generalization is considered$,$ even though these have been extensively 
investigated by nonstandard means. Indeed$,$ research indicates that the 
use of arbitrary partitions along with the 
 {Darboux} concept of the upper and lower sums and the upper and lower integrals 
which he proved to be equivalent to the {Riemann integral} are note worthy in  
that they more easily yield the rigorous proofs that establish the 
properties of the highly applicable converging Riemann sum notion. Since our
paramount concern is modeling with the Riemann styled integral$,$ our stated 
results are in terms of such easily conceived hyperfinite sums. A simple
nonstandard definition for an integral of bounded functions - the 
 {H-integral} - is given in   
{\tt Herrmann [1985]}. In this paper$,$ it is shown that the H-integral is 
equivalent to the {Darboux integral}. [Note: in Herrmann [1985] the Darboux 
integral is called the {Riemann integral}. Further$,$ Theorem 3.3 and Corollaries 
3.3.1$,$ 3.3.2$,$ 3.3.3$,$ 3.3.4 in this paper are incorrect as stated. However$,$ our 
use of the results from this paper are not related to these few erroneous 
conclusions.] As is well-known the Darboux integral is equivalent to the Riemann 
integral conceived of as approximated by the {Riemann Sums}. It seems expedient$,$ 
however$,$ to consider all such integral concepts extended to reasonable subsets 
of $\realp n$ using the simplest possible procedures [{\tt Apostal [1957]}$,$ 
 {\tt Spivak [1965]}].\par
For an n-dimensional space $n \geq 1$$,$ the closed set $R = [a_1,b_1] \times 
\cdots \times [a_n,b_n],\  a_i < b_i,\ 1\leq i \leq n$ is called 
 {\bf a {rectangle}.} Of course$,$ if $n = 1,$ then a ``rectangle'' is but a 
closed interval. As usual$,$ consider for each $[a_i,b_i]$ a  {\bf {partition}} 
 $P_i$ as a finite set of members of $[a_i,b_i]$ such that $a_i,b_i \in P_i$ 
and where $P_i$ is considered as ordered. This is often explicitly written as 
$P_i = \{x_{i0},\ldots, x_{ik}\},\  a_i = x_{i0} < x_{i1} <\cdots<  x_{ik} = b_i.$ 
This determines the closed one-dimensional subintervals 
 $[x_{i(p-1)}, x_{i{p}}],\  1\leq p \leq k.$ In brief$,$ this 
process obtains a partition $P = P_1 \times \cdots \times P_n$ of $R$ and 
a finite collection of closed n-dimensional subrectangles  $R_q$ obtained by 
considering $([x_{10}, x_{11}] \cup \cdots \cup [x_{1{k-1}},x_{1k}]) \times 
\cdots \times ([x_{n0}, x_{n1}] \cup \cdots \cup [x_{n{m-1}},x_{nm}]).$ Each 
$R_q$ has a measure$,$ $v(R_q)\in \real$$,$ assigned to it which  
is intuitively the product of the lengths of the sides. For simplicity of 
notion the definition of the measure $v(R_q)$ is left intuitively 
understood.\par
Probably the simplest partition to consider would be the one termed a 
{\bf {simple}} partition. These are formed by selecting $n$ nonzero natural 
numbers $m_1,\ldots, m_n$ and dividing each interval $[a_i,b_i]$ into an equal 
length partition by adding to each successive partition point the number
$(b_i -a_i)/(m_i).$ This concept is extended to the nonstandard world by 
selecting $n$ infinite natural numbers $\Gamma_1,\cdots,\Gamma_n$ and 
generating for each interval $[a_i,b_i]$ an internal hyperfinite partition$,$ 
$P_i,$ each subinterval of which has positive infinitesimal length 
$(b_i -a_i)/(\Gamma_i)= dx_i.$ Then the partition $P=P_1 \times \cdots \times 
P_n$ is a {\bf {simple fine 
partition}} of $R.$ Such a partition yields an internal set of hyperrectangles 
$R_i$ such that $\hyper v(R_i) = dx_1 \cdots dx_n=dX  \in \monad 0^+.$ 
You could be much more general and consider the {\bf {fine}} partitions which 
are internal collections of hyperfinitely many members of $\Hyper [a_i,b_i]$
such that the length of any subinterval is an infinitesimal. In Cauchy's 
definition he evaluated a function at specific endpoints of each subinterval. 
For a bounded function $f\colon R \to \real$ this evaluation concept can 
also be applied in the case of the subrectangles 
into which $R$ is partitioned by evaluating the function $f$ at$,$ say$,$ 
the corner nearest to the origin. However$,$ it has become customary to be 
somewhat more general and include the concept of the intermediate partition.
Let $P$ be a partition of the rectangle $R$ and assume that $P$ determines 
the set of subrectangles $\{R_q \vert 1\leq q \leq m\}.$  An {\bf {intermediate 
partition}}$,$ $Q$$,$ is any finite sequence of vectors $\{\vec v_q \},$ where 
$\vec v_q
\in R_q$ for each $q$ such that $1 \leq q \leq m.$ \par
\vskip 18pt
\hrule
\smallskip
{\bf Definition 5.1.1. ({The Integral})}. Let $f\colon R \to \real$ 
be bounded and $\cal P$ the set of simple partitions of $R.$ 
Then $f$ is said to be {\bf INTEGABLE} if there exists 
some $r \in \real $ and a \underbar{simple fine partition,} $P \in \hyper 
{\cal P}$ such that for each of its internal intermediate partitions 
$\{\vec v_q\},$ where $1 \leq q \leq \Gamma \in {\nat}_\infty,$ 
$$ \sum_{k=1}^\Gamma \hyper f(\vec v_q)\hyper v(R_q) \in \monad r.$$\par
\smallskip
\hrule
\vskip 18pt
{\bf Theorem 5.1.1.} {\sl A bounded function $f\colon R \to \real$ is 
integrable if and only if it is integrable in the sense of Darboux and 
Riemann Sums.}\par
\bigskip
{\bf Theorem 5.1.2.} {\sl If bounded $f\colon R \to \real$ is integrable$,$ then
there exists a unique $r \in \real$ such that  
for every fine partition $P = \{\vec x_0,\ldots \vec x_\Omega \},\ \Omega \in 
{\nat}_\infty$ and every internal intermediate partition $Q = \{\vec v_q\},\ 
1 \leq q \leq \Gamma \in {\nat}_\infty$ it follows that
$$\sum_{k=1}^\Gamma \hyper f(\vec v_q)\hyper v(R_q) \in \monad r.$$}\par
\vskip 18pt
The unique real number that exists by Theorem 5.1.2 is$,$ from theorem 5.1.1$,$ the 
classical value of the definite integral and hence if bounded $f\colon R \to \real$
 is integrable in the sense of Definition 5.1.1$,$ then we may write
$$\int \cdots \int_R f(\vec x)\, dx_1\cdots dx_n = \int _R f(\vec x)\, dX=
 \st {\sum_{k=1}^\Gamma \hyper f(\vec v_q)\hyper v(R_q)}.$$\par
Observe that Theorem 5.1.2 indicates that {Cavalieri's} notion of the 
 {indivisible} line segment as being the foundation for the definite integral is 
untenable. Even though the above hyperfinite sums correlate directly 
to the 
intuitive concept of the definite integral$,$ even for possibly discontinuous 
functions$,$ this fact alone does not lead to the appropriate selection of 
specific integrands that will produce meaningful geometric or physical 
measures. {\it Indeed$,$ one of the most significant aspects of modern infinitesimal 
analysis is that there does exist describable modeling procedures that 
rigorously establish that a specific integrand does provide the requisite 
value for a specific geometric or physical quantity.}\par
\bigskip  
\leftline{5.2 \underbar{The Infinite Sum Theorems.}}\par    
\medskip                                                  
Throughout applied mathematics numerous linear functionals are utilized to 
discuss and predict geometric or physical qualities. Such functionals are 
essential to the indirect 
verification of many physical theories for within the laboratory 
environment it is the predicted values displayed by elaborate machinery 
that often yield the only indications that unobserved events may be occurring.
As indicated in section 5.1$,$ there are
now \underbar{rigorous} {rules} that lead to an immediate adoption of the 
integral as the appropriate modeling structure
when one mentally conceives of such events and 
applies experience to determine  
the geometric or physical properties that might establish that it is  
likely that certain hypothesized behavior is actually occurring.\par
Prior to Robinson's discovery certain vaguely described rules did appear in the 
mathematical literature; but none was consistently defined in a rigorous 
language nor did they have any 
particular relation to the successfully applied intuitive notions we term 
infinitesimal reasoning. As an intermediate step in establishing a consistent 
and complete approach to this problem$,$ it is now possible to describe 
explicitly sufficient infinitesimal conditions that establish 
the integral as the correct modeling structure. These intermediate 
rules have 
become known as {\it {The Infinite Sum Theorems.}} For any bounded function 
$f\colon [a,b]\to \real$ a generalizations of these rules 
can be found in Herrmann [1985]. In appendix 5$,$ this rule is further 
generalized and applied to the case that bounded $f\colon R \to \real.$ 
Obviously$,$ a {\bf {subrectangle},} $R_S,$ is a rectangle that is a subset of 
$R.$ Let $\cal C$ 
be the set of all simple partitions of $R,\ {\cal S}(P)$ the set of all 
subrectangles generated by $P \in {\cal C}$ and ${\cal C}_{PSR} 
= \{R_S \vert \exists P(P \in {\cal C} \land R_S \in {\cal S}(P)\}$ 
the set of all subrectangles contained in any simple partition of $R.$ \par
\vskip 18pt
\hrule
\smallskip
{\bf Definition 5.2.1. ({Simply Additive}).} Let $\{R_q\vert 1 \leq q \leq 
m\}$                                                 
be any simple partition of $R$ and let $B$ be any map defined on the 
collection ${\cal C}_{PSR}$ and having real number values. Then 
$B$ is said to be {\bf SIMPLY ADDITIVE} if for each $\{R_q\vert 1 \leq q \leq 
m\} = {\cal S}(P),\ P \in {\cal C}$ it follows that
$$B(R) = \sum_{k=1}^m B(R_k).$$ 
\smallskip
\hrule
\vskip 18pt
Obviously$,$ $\hyper B$ is defined on the set $\Hyper {\cal C}_{PSR}$ which 
contains all of the infinitesimal subrectangles contained in any simple fine 
partition of $R.$  I point out that our first infinitesimal sum 
theorem actually holds for a slightly more general partition than a simple 
partition of $R$ - the special partition. However$,$ a simple partition is a 
special partition and in most applied cases the simple partition suffices.\par
\vskip 18pt
{\bf Theorem 5.2.1. (An Infinite Sum Theorem.)} {\sl  Let bounded 
$f\colon R \to \real$ 
and simply additive $B\colon {\cal C}_{PSR} \to \real.$ If there exists a 
simple fine partition $\{R_q \vert 1 \leq q \leq \Gamma \}$ and for each $R_q$ 
there exists some $\vec p \in R_q$ such that 
$$ \hyper B(R_q)/dX \approx \hyper f(\vec p),\eqno(*)$$
then $f$ is integrable and 
$$B(R) = \int_R f(\vec x)\, dX.$$}\par
\vskip 18pt
Two observations about Theorem 5.2.1. The infinitesimal $dX$ need not be 
considered the finite product of coordinate measures but may also take on the 
character of such physical quantities as the finite product of 
infinitesimal 
momenta$,$ the finite product of infinitesimal probabilities and even the finite 
product of infinitesimal charges or infinitesimal numbers of elementary 
particles if such things can be conceived of in objective reality.  Further$,$ 
it is somewhat unfortunate for applied mathematics that the converse of 
Theorem 5.2.1 does not hold. N.J. {\tt Cutland [1986]} has supplied your author 
with an example of a function defined on $[0,1]$ that is Darboux integrable 
(hence integrable) but if you define $B(R_S) = \int_{R_S} f(x) \, dx,$ then $\hyper B$ does not satisfy 
property (*).  On the other hand$,$ if $f$ is continuous on $R,$ then (*) does 
hold for such integrally defined functionals and$,$ indeed$,$ a much stronger property 
called {\it supernearness} holds as well. Let ${\cal C}_{SR}$ be the set of all 
subrectanghles contained in $R.$ {\bf For simplicity of notation$,$ throughout 
this manual$,$ maps such as $B$ are denoted as being$,$ at least$,$ defined on sets 
such as ${\cal C}_{PSR},\ {\cal C}_{SR}$ etc. Simple additivity and additivity 
will greatly enlarge their domains of definition.}\par
\vskip 18pt
\hrule
\smallskip
{\bf Definition 5.2.2. ({Supernearness}).}\hfil\break Let $(x_1,\ldots,x_n),\ (y_1,\ldots, 
y_n) \in \hyper R$ and $y_i -x_i \in {\monad 0}^+,\ 1\leq i\leq n$; and let 
$R_S = \{(z_1,\ldots,z_n)\vert \forall i(\ 1\leq i\leq n\to x_i \leq z_i 
\leq y_i)\land (z_i \in \hyperreal)\}$ denote an infinitesimal subrectangle of $R.$     
A map $B\colon {\cal C}_{SR} \to \real$ is 
{\bf {SUPERNEAR}} to bounded $f\colon R \to \real$ if for every infinitesimal 
subrectangle $R_S$ of $\hyper R$ and for every $\vec p \in R_S$ it follows that
$$ \hyper B(R_S)/dX \approx \hyper f(\vec p),\eqno (**)$$
where $dX = \prod_{i=1}^n (y_i - x_i)\in {\monad 0}^+.$ \par
\smallskip
\hrule
\vskip 18pt
{\bf Theorem 5.2.2.} {\sl A bounded function $f\colon R \to \real$ is 
continuous if and only if there exists a map $B\colon {\cal C}_{SR} \to \real$ 
that is supernear to $f.$} \par
\vskip 18pt
In the proof of Theorem 5.2.2 the following interesting integral property is 
established.\par
\vskip 18pt
{\bf Corollary 5.2.2.} {\sl Suppose that $f\colon R \to \real$ is continuous.
For each $R_S \in {\cal C}_{SR}$ define $B(R_S) = \int_{R_S} f(\vec x)\, dX.$
Then $B$ is supernear to $f.$} \par
\vskip 18pt

{\bf Theorem 5.2.3.} {\sl Let bounded $f\colon R \to \real.$ If $B\colon {\cal 
C}_{SR} \to \real$ is supernear to $f$ and simply additive on each simple 
partition of each $R_S,$ then $f$ is continuous on $R$ and 
$$B(R_S) = \int_{R_S} f(\vec x)\, dX$$ for each $R_S \in {\cal C}_{SR}.$}\par
\vskip 18pt                                                                   
In applications of the integral to geometric and physical problems it is 
       usually assumed that the map $B\colon {\cal C}_{SR}\to \real$ is$,$ at least,
 {\bf {additive}} on ${\cal C}_{SR}.$ Recall that this means that if 
nonempty $\{R_i\bigm \vert 1 \leq i\leq n\} \subset {\cal C}_{SR}$ is 
pairwise disjoint or pairwise has only  boundary 
points in common$,$ then $B(\cup\{R_i\}) = \sum_{i=1}^nB(R_i).$ Note that if 
$B$ is additive on ${\cal C}_{SR},$ then $B$ is simply additive on each member 
of $\cal C$ and for each simple partition of $R_S\in {\cal C}_{SR}.$ \par
\vskip 18pt
{\bf Corollary 5.2.3.1} {\sl Let bounded $f\colon R \to \real.$ If $B\colon 
{\cal C}_{SR} \to \real$ is supernear to $f$ and additive on ${\cal C}_{SR},$ 
then $f$ is continuous on $R$ and 
$$B(R_S) = \int_{R_S} f(\vec x)\, dX$$ for each $R_S \in {\cal C}_{SR}.$}\par
\vskip 18pt          
{\bf Corollary 5.2.3.2} {\sl Let bounded $f\colon R \to \real.$ There exists one 
and only one map $B\colon {\cal C}_{SR} \to \real$ that is supernear to $f$ 
and either simply additive on each simple partition of each $R_S$ 
or additive on ${\cal C}_{SR}$}.\par
\bigskip  
\leftline{5.3 \underbar{Extensions.}}\par    
\medskip                                                  
In general$,$ the bounded real valued function $f$ need not be defined on such a 
convenient set as $R.$ If $f\colon D \to \real$ is defined on a bounded set 
$D,$ then the most expedient procedure to follow is to define a function
$\hat f\colon R \to \real,$ where $D \subset R,$ by $\hat f(\vec x) = f(\vec x)$ for each 
$\vec x \in D$ and $\hat f(\vec x) = 0$ for each $\vec x \in R - D.$ With this
case then$,$ as is customary$,$ let $\int_D f(\vec x)\, dX = \int_R \hat f(\vec x),
\ dX.$\par
As far as a map such as $B\colon {\cal C}_{SR} \to \real$ is concerned$,$ the 
additivity of $B$ may be extended to all Jordan-measurable subsets of $R.$ The 
fact that $B$ may be additive on a lesser collection of subsets of $R$ will 
suffice for the basic modeling rules described in the next section. These 
modeling rules are very specific in character and if the proper simplistic 
assumptions for $B$ are utilized$,$ then they lead directly to the appropriate 
infinite sum theorem and its associated integral equivalence.\par 
%==========TABLE STRUTS

%24 POINT DEEP ROWS
\newbox\medstrutbox
\setbox\medstrutbox=\hbox{\vrule height14.5pt depth9.5pt width0pt}
\def\medstrut{\relax\ifmmode\copy\medstrutbox\else\unhcopy\medstrutbox\fi}

\bigskip

\newdimen\leftmargin
\leftmargin=0.0truein
%for textbook use 1.1truein
\newdimen\widesize
\widesize=3.0truein
\advance\widesize by \leftmargin       
%\moveleft\leftmargin
\hfil\vbox{\tabskip=0pt\offinterlineskip
\def\tablerule{\noalign{\hrule}}

\halign to \widesize{\medstrut\vrule#\tabskip=0pt plus2truein

&\hfil\quad#\quad\hfil&\vrule#

\tabskip=0pt\cr\tablerule

%&\multispan3 \hfil {\tablerm TABLE TITLE} \hfil&\cr\tablerule
%\biggerstrut&ENTRY A&&ENTRY B&\cr\tablerule

&$\Rightarrow$ IMPORTANT $\Rightarrow$&\cr\tablerule
}} \par
$\Rightarrow$ In the following applications$,$ the standard requirements are 
stated in 
terms of what we perceive to be global behavior of well-known ordinary 
functionals 
and their relation to standard characterizing properties. {\bf Many of these 
observation are not obvious. The following applications are actually intended 
to foster an appreciation for the nonstandard modeling rules and procedures 
that appear in Chapter 6 - rules that lead more directly to the 
appropriate conclusions.} The reason we present the following applications is 
that this global approach is used in the more elementary textbooks. $\Leftarrow$          
\medskip 

\leftline{5.4 \underbar{Applications and the Standard Modeling Rules.}}\par                                                         
\medskip      
One  of the unusual aspects of the Infinite Sum Theorem 5.2.1 is that it 
does not view 
the functional $B$ directly but$,$ rather$,$ a {``mean value''} must be considered if 
the attention is directed toward the concept of being infinitely close. As the 
derivations in the following elementary applications indicate this problem is 
submerged within the derivation itself and does not usually occur when the 
properties of the basic functional are proposed. Our applications are mostly 
geometric and elementary in character$,$ while the major applications to the 
physical sciences will appear in their respect manuals. {\it Further$,$ we 
concentrate upon those applications that traditionally appear in the customary 
core calculus$,$  elementary differential equations and physical science 
courses.}\par
 \medskip
 {\bf Application 5.4.1} {\it The 2-dimensional area between two continuous 
curves.}\par 
 First assume that all of the following 
functions are continuous on their indicated domains. Give $h\colon [c,d] \to 
\real$ and $k\colon [c,d] \to \real,$ where we denote by the symbol $ h \leq 
k$ the condition that $h(x) \leq k(x)$ for each $x\in [c,d].$ Let's look at 
the original idea behind an area function $A([c,d],h,k)$ that measures the 
intuitive area between these two curves and over the interval $[c,d].$ Define 
a {\bf {basic region}} $R$ for any pair of functions $h,k,\ h \leq k$ and over 
any interval $[c,d]$ contained in their common domain by $R = \{(x,y)\vert x 
\in [c,d] \land h(x) \leq y \leq k(x) \}.$ Assume that we have two fixed 
functions $f,g,\ f \leq g$ defined on $[a,b].$ The follow {rules (axioms) 
appear to model} our intuitive notion of an area function.\par
 (i) The area 
function $A$ is$,$ at least$,$ defined on all basic regions determined by $[c,d] 
\subset [a,b].$ \par
 (ii) For the above two fixed functions $f,g$ defined on 
$[a,b]$ the area function $A$ is$,$ at least$,$ additive on the set ${\cal 
C}_{SR}$ of $[a,b].$\par 
(iii) If $D,E$ are two basic regions in the domain of 
$A$ and $D \subset E,$ then $A(D) \leq A(E).$\par 
(iv) If two functions $h,k$ 
are constant over any $[c,d] \subset [a,b],$ then $A([c,d],h,k) = (d-c)(k-
h).$\par 
These four properties for the area function $A$ are certainly 
reasonable and seem to model the intuitive notions from elementary plane 
geometry. We now formally establish that $$ A([a,b],f,g) = \int_a^b (g(x) - 
f(x)) \, dx.\eqno (I_1)$$ \par
\indent {\sl Derivation.} Let $[x_1,x_1 
+ dx]$ be any hyperinterval generated by any simple fine partition of $\Hyper 
[a,b].$ Note in this case $dx = (b-a)/\Gamma,$ where $\Gamma \in 
{\nat}_\infty.$ By considering the *-transfer of the standard extreme value 
theorem for continuous functions defined on closed intervals it follows that 
$\hyper f$ and $\hyper g$ attain their maximum and minimum values $\hyper 
f_M,\ \hyper g_M$ and $\hyper f_m,\ \hyper g_m$ respectively on    $[x_1,x_1 + 
dx].$ It is not difficult to model statements (i)$,$ (ii)$,$ (iii) (iv) set-
theoretically and extend these properties to the nonstandard world. We need 
only consider statements (i)$,$ (ii)$,$ (iii)$,$ (iv) as intuitively *-transformed 
by changing the terminology to ``hyper'' or ``*'' terminology. When this is 
done statements (i)$,$ (iii) and (iv) yield the result that $$(\hyper g_m - 
\hyper f_M)dx \leq \hyper A([x_1,x_1 +dx],\hyper f,\hyper g]) \leq (\hyper g_M 
- \hyper f_m)dx  \Rightarrow\eqno (1)$$ $$(\hyper g_m - \hyper f_M) \leq  (\hyper 
A([x_1,x_1 +dx],\hyper f,\hyper g]))/dx \leq (\hyper g_M - \hyper f_m).\eqno 
(2)$$ 
Since $f$ and $g$ are uniformly continuous on $[a,b]$ and $\hyper f,\hyper g$ 
attain their respective maximum and minimum value at members of $[x_1,x_1 
+dx]$ then definition 4.4.1 implies that $$\hyper f_m \approx \hyper f(x_1) 
\approx \hyper f_M,\ \hyper g_m \approx \hyper g(x_1) \approx \hyper 
g_M.\eqno (3)$$ 
Consequently$,$ $$(\hyper g_m - \hyper f_M) \approx (g(x_1) - f(x_1)) \approx 
(\hyper g_M - \hyper f_m) \Rightarrow\eqno (4)$$ $$(A([x_1,x_1 +dx],\hyper f,\hyper 
g]))/dx \approx  (g(x_1) - f(x_1))\eqno (5)$$ from the fact that $f,g$ are bounded and 
Corollary 2.2.5.2. Statement (i) allows application of the Infinite Sum 
Theorem  and integral equation $I_1$ is the consequence. \par 
\medskip 
Once 
equation $I_1$ is obtained then it may be checked against the standard area 
measures for the ordinary Euclidean plane figures in order to insure that it 
is indeed an extension. I point out that throughout many of these applications 
similar modeling rules such as (i)$,$ (ii)$,$ (iii) and (iv) are essential                        
if one wishes to achieve a formal derivation.  As will be illustrated there 
are notable exceptions to this general approach where one of the standard 
conditions (i) - (iv) may fail. However$,$ here are two more applications where 
the standard functional characterizations can be formulated \par 
\medskip 
{\bf Application 5.4.2. {\it Volume obtained by a 2-dimensional integral.}}\par 

Assume$,$ as in the previous application$,$ that all functions are continuous on 
their indicated domains. Suppose that two functions $h,k$ are defined on a 
rectangle $R_0\subset \realp 2$ and have the property that the $h\leq k$ on 
$R_0$. Generalizing the definition in application 5.4.1$,$ define a basic region 
$R_B$ in  $\realp 3$ for each a pair $h,k$ by $R_B =\{(x,y,z)\vert (x,y) 
\in R_0 \land h(x,y) \leq z \leq k(x,y)\}.$ As before$,$ we attempt to model the 
concept of a volume function $V(R_0,h,k)$ over any rectangle $R_0\subset 
\realp 2.$ Assume that we have two real valued fixed functions $f,g,\ f\leq g$ 
defined on a rectangle $R.$\par 
(i) The volume function $V$ is defined$,$ at 
least$,$ for all basic regions determined by rectangles that are subsets of 
$R.$\par 
(ii) For the above two functions $f,g,$ the function $V$ is$,$ at 
least$,$ additive on the set ${\cal C}_{SR}$ of $R.$\par 
(iii) If $D,E$ are two 
basic regions in the domain of $V$ and $D\subset E$$,$ then $V(D) \leq 
V(E).$\par 
(iv) If the functions $h,k$ are constant over any rectangle $R_0 
\subset R,$ then $V(R_0,h,k) = ({\rm area}\ R)(k - h).$ \par 
If $V$ satisfies 
these rules$,$ then $$V(R,f,g) = \int_R (g(\vec x) - f(\vec x))\, dX.$$ 
\indent {\sl Derivation.} Except for a very slight modification$,$ this is 
exactly the same as the derivation for application 5.4.1. Simply let $[x_1,x_1 
+ dx_1] \times [x_2,x_2 + dx_2]$ be a hyperrectangle generated by a simple fine 
partition of $\hyper R.$ In this case$,$ $dx_1 = (b-a)/\Gamma,\ dx_2 = (d-
c)/\Omega,\ \Gamma,\Omega \in {\nat}_\infty.$ Since $R$ is compact the 
remainder of this derivation is as in application 5.4.1.\par 
\medskip       
{\bf Application 5.4.3. {\it Mass obtained by a 3-dimensional integral.}}\par 
Assume that $\rho (\vec x)$ is a continuous point density function defined on  
a rectangle $R \subset \realp 3.$ General physical experience leads to the 
following characterizations for the elementary mass$,$ $M(R,\rho ),$ of $R$.\par 
(i) The mass is defined for$,$ at least$,$ the set of all 
subrectangles$,$ ${\cal C}_{SR},$ of $R.$ \par 
(ii) The mass is additive on the 
set ${\cal C}_{SR}.$ \par 
(iii) If $R_S$ is a subrectangle of $R$ and 
${\rho}_1(\vec x)$ is a continuous density function defined on $R_S$ with the 
property that ${\rho}_1(\vec x) \leq \rho (\vec x)$ for each $x \in R_S,$ then 
$M(R_S,{\rho}_1) \leq M(R_S,\rho ).$\par 
(iv) If the function $\rho$ is constant 
over any rectangle $R_0 \subset R,$ then $M(R_0,\rho) = \rho(v(R_0)).$\par 
If 
$M$ satisfies rules (i) -(iv)$,$ then $$M(R,\rho) = \int_R \rho (\vec x) \, 
dX.$$ \indent {\sl Derivation.} As in the previous cases$,$ one selects a 
simple fine partition of $\hyper R$ and lets $R_q$ be some hyperrectangle 
determined by such a simple fine partition. The above characterizations are 
extended by *-transfer to the NSP-world. From continuity$,$ $\hyper \rho$ 
attains its minimum value $\hyper \rho_m$ and maximum value $\hyper \rho_M$ at 
members of $R_q.$ From (i)$,$ (iii) and (iv) it follows that $$(\hyper \rho_m)dX 
\leq \hyper M(R_q,\hyper \rho)\leq (\hyper \rho_M)dX \Rightarrow\eqno (1)$$ $$\hyper 
\rho_m \leq (\hyper M(R_q,\hyper \rho))/dX \leq \hyper \rho_M.\eqno (2)$$ 
Let $x_1$ be 
any member of $R_q.$ From the uniform continuity of $\rho$ it follows that 
$$\hyper \rho_m \approx \hyper \rho (x_1)\approx \hyper \rho_M.\eqno (3)$$ 
Consequently$,$ $$(\hyper M(R_q,\hyper \rho))/dX \approx \hyper \rho (x_1)\eqno 
(4)$$ and 
the derivation follows from the Infinite Sum Theorem.\par 
\vskip 18pt
\hrule 
\smallskip 
\hrule 
\smallskip 
{\sl In the above applications$,$ the rules (i) -- (v) obviously depict  
these functionals from a standard point of view as described completely in terms of 
the standard world. Except within  the formal derivations$,$ there may appear to 
be no consuming need for any insight into the infinitesimal NSP-world. As is 
be amply illustrated in the next chapter$,$ one or more of these 
rules may not be self-evident when the integral is 
applied for both geometric and physical modeling. When this occurs$,$ then it is 
often the case that certain simplistic and local aspects of the standard world 
are axiomatically impressed upon the pure NSP-world. This then 
leads to rules that do include descriptions for pure NSP-world behavior.} 
\smallskip
\hrule
\smallskip
\hrule
\smallskip
\vskip 24 pt 
\leftline{5.5 \underbar{Extensions of the Standard Rules.}}\par
In practice$,$ the rectangle $R$ is too restrictive to be of much significance 
in applications. Fortunately$,$ there are techniques that will allow us to 
remove this restriction - techniques that lead to a straightforward modification 
of a few of the terms that appear in such standard rules as (i) -- (iv). \par
The appropriate alteration of these rules begins with the concepts briefly 
mentioned in section 5.3. \par
(1) Assume that $f\colon D \to \real$ is 
continuous on $D,$ where $D$ is a compact Jordan-measurable subset of 
$\realp n.$ [Apostal [1957]$,$ {\tt De Lillo [1982]}] The {Jordan-measurable} subsets of $R$ 
include those that commonly appear throughout basic applications. \par
 (2) Let 
${\cal C}_{JR}$ denote the set of all Jordan-measurable subsets of $R.$ Clearly$,$ 
${\cal C}_{SR} \subset {\cal C}_{JR}$ as is well-known. \par 
(3) Assume that 
$B\colon {\cal C}_{JR} \to \real.$ \par                                          

(4) Now extend $f$ to $\hat f$ and assume 
that $\hat f$ is integrable on $R.$ It is easy to show that the 
value of $\int_D f(\vec x)\, dX = \int_R \hat f(\vec x)\, dX$  is independent 
of the choice of $R.$\par
Our next task is to see how a simple modification of the standard rules leads 
to 
a refined derivation that establishes the same integral expression. {\bf In the 
following applications$,$ the rule modifications are written in italics and  $v(J)$ 
denotes the real Jordan content for any Jordan-measurable set $J\subset R.$}
\medskip
{\bf Application 5.5.1. {\it Volume obtained by a 2-dimensional integral.}}\par 

Assume$,$ as in application 5.4.2$,$ that all functions are integrable on 
their indicated {\it compact Jordan-measurable} domains. Suppose that two 
functions $h,k$ are defined on a 
{\it compact and Jordan-measurable} $J_0\subset R \subset \realp 2$ and 
have the 
property that the $h\leq k$ on 
$J_0$. Generalizing the definition in application 5.4.1$,$ define a
 basic region 
$J_B$ in  $\realp 3$ for each a pair $h,k$ by $J_B =\{(x,y,z)\vert (x,y) 
\in J_0 \land h(x,y) \leq z \leq k(x,y)\}.$ As before$,$ we attempt to model the 
concept of a volume function $V(J_0,h,k)$ over any $J_0 \in {\cal C}_{JR}.$ 
Assume that we have two real valued fixed and continuous functions 
$f,g,\ f\leq g$ defined on {\it compact and Jordan-measurable} $J\subset R.$\par 
(i) The volume function $V$ is defined$,$ at 
least$,$ for all basic regions determined by all {\it compact and 
Jordan-measurable sets} that are subsets of $R.$\par 
(ii) For the {\it integrable extensions} $\hat f,\hat g,$ {\it of} the above two functions  
the function $V$ is$,$ at 
least$,$ additive on the set ${\cal C}_{JR}$ of $R.$\par 
(iii) If $D,E$ are two 
basic regions in the domain of $V$ and $D\subset E$$,$ then $V(D) \leq 
V(E).$\par 
(iv) If the functions $h,k$ are constant over any $J_0 \in {\cal C}_{JR},$
then $V(J_0,h,k) = (v(J_0))(k - h).$
If $V$ satisfies 
these rules$,$ then $$V(J,f,g) = \int_R (\hat g(\vec x) - \hat f(\vec x))\, 
dX=\int_J (g(\vec x) - f(\vec x))\, dX.$$
\indent {\sl Derivation.}  Let $I= [x_1,x_1 + dx_1] \times [x_2,x_2 + dx_2]$ 
be a hyperrectangle generated by a simple fine 
partition of $\hyper R.$ In this case$,$ $dx_1 = (b-a)/\Gamma,\ dx_2 = (d-
c)/\Omega,\ \Gamma,\Omega \in {\nat}_\infty.$ By *-transfer$,$ we transfer the 
general results concerning Jordan-measurable sets to the NSP-world. Assume that 
$\hyper J \cap I = K\neq \emptyset.$ Since $I \in \hyper {\cal C}_{JR}$ then
$K \in \hyper {\cal C}_{JR}.$ Further since $I$ is *-compact then $K$ is 
*-compact. Noting that $\hyper f = \hyper {\hat f},\ \hyper g = \hyper {\hat 
g}$ on $K$$,$ then the *-extreme value theorem implies that $ \hyper {\hat f},\ 
\hyper {\hat g}$ attain their maximum and minimum values $\hyper {\hat f_M,},\ 
\hyper {\hat g}_M$ and $\hyper {\hat f}_m,\ \hyper {\hat g}_m$ respectively on 
$K.$ Now by application of *-additivity and (iv) and the fact that $K,I,I-K$ 
are *-Jordan-measurable (if $I-K = \emptyset,$ then we still let it be 
measurable with content equal to zero) we have that 
$\Hyper V(I,\hyper {\hat f},\hyper {\hat g}) = \Hyper V(K,\hyper {\hat 
f},\hyper {\hat g}) 
+ \Hyper V(I-K,\hyper {\hat f},\hyper {\hat g}) = \Hyper V(K,\hyper {\hat f},\hyper 
{\hat g}).$\par
Next we also apply (iii) and obtain 
$$(\hyper {\hat g}_m - \hyper {\hat f}_M )\hyper v(K) \leq (\hyper {\hat g}_m - 
\hyper {\hat f}_M )dX \leq \Hyper V(I,\hyper {\hat f},\hyper {\hat g})=$$ 
$$\Hyper V(K,\hyper {\hat f},\hyper 
{\hat g)} \leq(\hyper {\hat g}_M - \hyper {\hat f}_m )\hyper v(K) \leq (\hyper {\hat 
g}_M - \hyper {\hat f}_m )dX.\eqno (1)$$\par
\noindent Since $f =  \hat f,\ g = \hat g$ are uniformly continuous on $K$  and $K \neq 
\emptyset$ then for $x_1 \in K$ we have that
$\hyper {\hat f}_m \approx \hyper {\hat f}(x_1) \approx \hyper {\hat f}_M,\ 
\hyper {\hat g}_m \approx \hyper {\hat g}(x_1) \approx \hyper {\hat 
g}_M.$
Consequently,
$$(\hyper {\hat g}_m - \hyper {\hat f}_M )\approx (\hyper {\hat g}(x_1) - 
\hyper {\hat f}(x_1))\approx (\hyper {\hat g}_M - \hyper {\hat f}_m ).\eqno 
(2)$$
Application of expressions (1) and (2) yields
$$\Hyper V(I,\hyper {\hat f}, \hyper {\hat g})/dX \approx  (\hyper {\hat g}(x_1) - 
\hyper {\hat f}(x_1)). \eqno (3)$$\par
\indent For the case that $K = \emptyset$ it is obvious that $\Hyper V(I,\hyper 
{\hat f}, \hyper {\hat g}) = 0 = (\hyper {\hat g}(x_1) - 
\hyper {\hat f}(x_1))dX\, x_1 \in I.$ In this case we also have that 
expression (3) holds. 
Application of the Infinite Sum Theorem completes the derivation.\par
\medskip
The above example suffices to show how all of the previous standard modeling 
rules and applications 
can be extended to the case of the Jordan-measurable subsets and integrable 
functions.  However$,$ for many applications of integral modeling to 
geometric theories and natural system behavior such listed axioms for the 
behavior of the conjectured functionals are often not evident.  This is 
particularly so for standard axioms such as (iii) and (iv). To eradicate 
this difficulty$,$ a 
direct appeal is made to the NSP-world$,$ either to the infinitesimal terms of 
the hyperfinite 
sum that appears in definition 5.1.1 or to statement (*) of Theorem 5.2.1.\par
\vfil\eject
\centerline{Chapter 6.}
\medskip                
\centerline{\bf NONSTANDARD RULES}
\centerline{\bf FOR INTEGRAL MODELING}
\bigskip
\leftline{6.1 \underbar{Historical Examples.}}\par  
In 1855$,$ {\tt Maxwell [1890]} presented his fluid flow analogue model for Faraday's 
concept of both magnetic and electric lines of force. An analysis of Maxwell's 
imagery relative to our present understanding of the behavior of infinitesimal
quantities is very enlightening. Maxwell considers {``tubes'' of moving 
points} (not particles) of fluid and their paths of motion as a pure imaginary 
picture of what one might conceive of as line a force. Obviously$,$ the idea was not to consider the 
concept of ``force'' as an independent entity but$,$ rather$,$ 
to first picture ``something'' - the points of fluid material and their 
paths of motion - as representing the effects of unknown forces. 
``{\it The direction of motion of the fluid will in general be different at 
different points of the space which it occupies$,$ but since the direction is 
determined for every such point$,$ we may conceive a line to begin at any point 
and to continue so that every element of the line indicates by it direction 
the direction of motion at that point in space. Lines drawn in such a manner 
that their direction always indicates the direction of fluid motion are called 
\underbar{{lines of fluid motion}.}}'' [Maxwell$,$ 1890:160] \par
Maxwell then imagines a closed curve on a surface - a surface that ``cuts''
the lines of fluid motion - and the lines of fluid motion that intersect this 
surface curve. These curve generated flow lines then produce$,$ in his mind$,$ a 
tubular surface which he calls {\it {a tube of fluid motion}.} He then fills 
the interior of these tubular surface with the flow lines that intersect 
that portion of the surface which would have the curve as its boundary. He 
also assumes that the fluid is incompressible. This yields another {postulated 
property of these flow lines}. ``{\it The quantity of fluid which in a unit of 
time crosses any fixed section of the tube is the same at whatever part of the 
tube the section is taken....and no part runs through the sides of the tube$,$ 
therefore the quantity which escapes from the second section is equal to that 
which enters through the first.}'' [Maxwell 1890:161]\par
Maxwell then supplies a paramount nonstandard rule to the methods of 
infinitesimal model 
- a rule that has recently been called the concept of the infinitesimal   
microscope. ``{\it An infinite number of lines would have to be drawn at 
indefinitely small intervals; but since the description of such a system of 
lines would involve continual reference to the theory of limits$,$ it has been 
thought better to suppose the lines drawn at intervals depending on the 
assumed unit$,$ and afterwards to assume the unit as small of we please by 
taking a small submultiple of the standard unit.}'' [Maxwell 1890:161] Notice 
that Maxwell's statement about the necessity of limit theory is now known to 
be false. He may have made such a statement$,$ as did Kepler before him$,$ to 
placate those who might be more attuned to rigorous derivations.  Further$,$ 
following general scientific practices$,$ Maxwell 
does not establish his limit theory conclusions but$,$ rather$,$  
ascribes to an infinitesimalizing approach he claims is equivalent to the 
physical limit theory. More importantly we have our first vague nonstandard 
rule for infinitesimal modeling\par
\vskip 18pt
\hrule
\smallskip
\hrule
\smallskip
\centerline{\bf {VR1.}}
 These infinitesimalizing ideas are equivalent to a 
type of {infinite magnification of a infinitely small portion of the fluid} - a
 magnification that yields finitely many lines in our field of view that 
appear to be drawn at real finite distances apart.\par
\smallskip
\hrule
\smallskip
\hrule
\vskip 18pt
 I mention that such 
geometric notions as expressed in VR1 can indeed be formalized within the geometric theory of 
$\hyperrealp 3.$\par
Within Maxwell's research reports he states numerous times that the 
reasons for his derivations are ``{evident}.'' Thus$,$ he often gives no specific 
causes for his logical conclusions and leaves them axiomatic in character. It 
is$,$ however$,$ the vague methods of infinitesimalizing that continue to 
interests us - 
methods that are also often assumed to be ``evident'' to Maxwell's audience.
Intuitively$,$ as it will be illustrated$,$ physical and geometric infinitesimal 
integral modeling also displays yet another vague nonstandard rule.                                                            
\vskip 18pt
\hrule
\smallskip
\hrule
\smallskip
\centerline{\bf {VR2.}}
Infinitesimal integral modeling often makes a direct appeal to the 
Infinite Sum Theorems as well as a 
simplified interpretation of the equation (Theorem 5.1.2)  
$$\int \cdots \int_R f(\vec x)\, dx_1\cdots dx_n = \int _R f(\vec x)\, dX=
 \st {\sum_{k=1}^\Gamma \hyper f(\vec v_q)\hyper v(R_q)}.$$
In the infinite sum theorems the equivalence relation $\approx$ is 
replaced by an equality while the standard part operator is ignored and the 
integral is made equal to some type of summation process which is often  
conceived of as finite in character.\par
\smallskip
\hrule
\smallskip
\hrule
\vskip 18pt
How does Maxwell apply his fluid motion analogy to assumed continuously varying 
magnetic properties? Surprising$,$ he views them as discrete and constant with relation 
to his magnified portion the of fluid material. ``{\it The quantity of magnetism in 
any section of a magnetic body is measured by the number of lines of magnetic 
force that pass through it.}'' Maxwell [1890:182] This ``number'' is assumed 
to be a standard natural number. ``{\it If $i$ be the quantity of the 
magnetization at any point$,$ or the number of lines of force passing through
unit of area in the section of the solenoid$,$ then the total quantity of 
magnetization in the circuit is the number of lines which pass through a 
section$,$ $I = \sum i\, dy\, dz,$ where $dy\, dz$ is the element of the section$,$ and the 
summation is performed over the whole section.} [Maxwell 1890:183] In the 
magnified view$,$ the section is to be conceived of as a rectangle with actual 
real number area that is then made arbitrary small by Maxwell's small unit 
convention. The constant numbers $i$ are not assumed to be altered as a 
physical quantity 
by the small unit convention but are fixed constants. The cardinality (i.e. 
intuitively the number of terms in the summation) is ignored and this   
``summation'' is equated to the integral value $I.$ \par
{\bf Maxwell appears to have arrived at his conclusions by considering the simplest of laboratory experiences for assumed constant quantities and 
geometric configurations. He has then impressed these experiences upon the 
infinitesimal nonstandard world. By assuming that the outcome is somehow additive 
in character this leads directly to the integral model.} \par
 With respect to pure 
geometric nonstandard modeling the same general process$,$ with certain 
exceptions$,$ is also applied. It was seen in Example 4.4.1.A on page 33 that 
the length of an n-dimensional curve is viewed globally as approximated by a 
hyperpolygonal curve and following this the length of the curve is defined to be 
the hyperfinite sum
$\sum_{i=0}^{\Gamma - 1} \Vert \vec v_i \Vert 
= \vert {\cal P}_\Gamma \vert \in \hyperreal.$ However$,$ for curves that are 
continuously differentiable Robinson formally showed that each term of this 
summation 
could be replaced by the term $$\sqrt{\sum_{j=1}^n \hyper f_j^\prime(t_i)^2}\, 
dt$$ where this replacement term can be conceived of as the infinitesimal length of an 
infinitesimal line segment infinitely close (in the sense of the infinitesimal  
sum theorem) to the original line segments that comprise the hyperpolygonal 
curve. {\bf A proof of this can be found in appendix 6.} On the other hand$,$ 
one could proceed backwards$,$ as some geometers have done$,$ and discuss 
the vague notion of the ``element of length,'' $ds,$ while forcing the length of a curve 
to be the ``sum'' of such elements $ds.$ Indeed$,$ one of the greatest of all 
infinitesimal geometers does all of his analysis in terms of the still vague 
``element'' concept and uses$,$ what would be today$,$ totally unacceptable 
non-rigorous derivations. ``{\it Let us now examine the integral curvature of this 
triangle$,$ which is equal to $\int k\, d\sigma,\ d\sigma$ denoting a surface 
element of the triangle. Wherefore$,$ since this element is expressed by
$mdp\cdot dq,$ we must extend the integral $\int\! \int mdp\cdot dq$ over the 
whole surface of the triangle.}''[{\tt Gauss 1827; Art 20}] Gauss in his papers 
also states that various quantities  are equal when in reality they are but
infinitely close. Fortunately$,$ the algebraic manipulations of infinitesimal 
quantities as employed by Gauss and Maxwell were restricted to their ring 
properties. Unfortunately$,$ Gauss' derivations are highly non-rigorous in 
character.\par
 By comparing Robinson's techniques with many of significant mathematical 
models developed  over the past few centuries by application of intuitive 
infinitesimal modeling$,$ it is now possible to gain a little better 
insight into the methods used$,$ to make these methods somewhat more 
acceptable in rigor and to improve considerably upon the notions outlined 
in the vague 
nonstandard rules VR1 and VR2.\par
\bigskip
\leftline{6.2 \underbar{The Monadic Environment.}}\par
\medskip
Vague rule 1 can now be made essentially rigorous in character. When we 
model the most elementary geometric or simplistic physical behavior it is 
usually conceived of within bounded portions of $\realp n.$ Even though one 
may later remove the {boundedness concept}$,$ such modeling often begins with 
such a restriction. Furthermore$,$ this holds true whether or not one is 
concerned with real or complex variables. What is {conceived of as
``elementary geometric'' or ``simplistic physical''} behavior is most often 
fostered by individual experiences and a personal development of discipline 
intuition. Simplistic diagrams and {sketching} develop an intuition for 
geometric concepts$,$ while {basic laboratory experimentation} tends to yield 
to the conscientious investigator an intuitive understanding of basic 
natural system behavior. It is this {phenomenological approach} to simplistic
physical behavior that has led to the development of most of our present day 
intricate mathematical models that$,$ even though they may tend to predict 
observed behavior$,$ may not correspond in their entirety to physical reality. \par
\vskip 18pt
\hrule
\smallskip
\hrule
\smallskip
\hrule
\smallskip
\centerline{IR1}
Elementary geometric or simplistic physical behavior takes place within  
an m-dimensional monad$,$ $\monad {\vec p}\subset \hyperrealp m,$ where standard 
$\vec p\in \realp m.$
Such behavior may be intuitively *-transferred to similar behavior within 
$\monad {\vec p}.$\par
\smallskip
\hrule
\smallskip
\hrule
\smallskip
\hrule
\vskip 18pt
It might be argued that IR1 is too restrictive for many of our modern 
applied mathematical structures. It is interesting to note that many monadic 
properties generalize not only to general topological spaces (see many of 
papers published by Robert A. Herrmann from 1975 - 1984) but even to the more 
general pre and pseudo topological spaces [{\tt Herrmann 1980]}. Can we identify 
for $\hyperrealp m$ what$,$ at least partially$,$ constitutes basic simplistic 
behavior? The first step in this process is to study in $\realp n$ or the 
laboratory geometric or physical behavior restricted to entities 
termed {``elements.''} Further$,$ one is often only concerned with specific 
functionals associated with such objects. As previous mentioned$,$ in applied 
infinitesimal analysis$,$ infinitesimals that are denoted by such symbols as 
$dx$ need not correspond to geometric concept of length. Indeed$,$ they can 
be negative or correspond to such notions as velocity$,$ acceleration or even 
infinitesimal ``numbers of molecules.'' Nevertheless$,$ it is convenient to 
regard the basic elements as geometric in character and not to be concerned
with their specific functional or physical interpretation.\par
\vskip 18pt
\hrule
\smallskip
\hrule
\smallskip
\hrule
\smallskip
\centerline{IR2}
Let internal $L\colon \hyperrealp n \to \hyperrealp m$ be a designated linear 
transformation and $S\subset \hyperrealp n$ an infinitesimal 
subrectangle of $\hyper R,$ where rectangle $R \subset \realp n.$ 
{\bf {An m-dimensional element} is the configuration $L[S].$} The {\bf 
 {basic element}} is the infinitesimal subrectangle itself.\par
\smallskip
\hrule
\smallskip
\hrule
\smallskip
\hrule
\vskip 18pt
For all such designated $L$ considered in IR2 the elements can be 
characterized as 
m-dimensional {parallelepipeds}$,$ where the 1-dimensional parallelepiped is 
to be conceived of as a line segment. Further$,$ we have the following little 
theorem that shows the relationship between infinitesimal subrectangles and 
monads.\par
\vskip 18pt
{\bf Theorem 6.2.1.} {\sl Let $A$ be a compact subset of $\realp n$ and the 
infinitesimal subrectangle $R \subset \hyper A.$ Then there exists some $p \in 
A$ such that $R \subset \monad p.$}
\vskip 18pt
%==========TABLE STRUTS

%24 POINT DEEP ROWS
\newbox\medstrutbox
\setbox\medstrutbox=\hbox{\vrule height14.5pt depth9.5pt width0pt}
\def\medstrut{\relax\ifmmode\copy\medstrutbox\else\unhcopy\medstrutbox\fi}

\bigskip

\newdimen\leftmargin
\leftmargin=0.0truein
%for textbook use 1.1truein
\newdimen\widesize
\widesize=3.0truein
\advance\widesize by \leftmargin       
%\moveleft\leftmargin
\hfil\vbox{\tabskip=0pt\offinterlineskip
\def\tablerule{\noalign{\hrule}}

\halign to \widesize{\medstrut\vrule#\tabskip=0pt plus2truein

&\hfil\quad#\quad\hfil&\vrule#

\tabskip=0pt\cr\tablerule

%&\multispan3 \hfil {\tablerm TABLE TITLE} \hfil&\cr\tablerule
%\biggerstrut&ENTRY A&&ENTRY B&\cr\tablerule

&$\Downarrow$ IMPORTANT $\Downarrow$&\cr\tablerule
}} \par
\bigskip
$\Rightarrow$ The hypotheses of many of the following rules and theorems are stated in 
terms of a Jordan-measurable set $J$. Such premises may be weakened by 
restricting the functional to ${\cal C}_{SR}$ and assuming that $J$ is but a 
subrectangle. {\bf I mention the important but well-known fact that all of the 
usual geometric configurations utilized in the basic calculus and differential 
equation undergraduate courses are {Jordan-measurable}. Further$,$ the {elementary 
prototype} used in place of the general Jordan-measurable set is the 
subrectangle itself. Most experimental investigations do not go beyond the 
subrectangle. After the hypotheses that appear in the following rules and 
theorems are restricted to subrectangles and their conclusions are verified$,$ 
then their consequences are often extended to $J$ without further 
confirmation.
\rightline{$\Uparrow$}}\par
\vskip 18pt
In order to use the infinite sum theorem 
directly for the {\bf {basic elementary integral}} it is assumed that 
the internal linear map $L$ is the extended identity may $\Hyper I.$
Then the following are gleaned from the basic definitions or the intuitive 
methods of the geometer as well as from simple experiments on natural system 
behavior whether it be objectively real or imaginary.  
\par 
\bigskip 
\hrule
\smallskip
\hrule
\smallskip
\hrule
\smallskip
\centerline{IR3 -- {Infinitesimal Max. and Min. Rule}}
(1) We wish to measure a quantity $M$ for a compact Jordan-measurable set $J 
\subset R \subset \realp n,$ where $M$ is defined on and$,$ at least$,$ additive 
over members of the set ${\{\cal C}_{SR},R-J,J\}.$ Further$,$ if subrectangle 
$S \subset R - J,\ M(S) = 0,$ and $M(R-J) = 0.$ 
Let $v(J)$ denote the Jordan content.\par
(2) There is a generating function $f(\vec x)$ that is related to the 
functional 
$M$ in the following manner:\par
\indent\indent (i) The functions $f$ is continuous on $J.$\par
\indent\indent (ii) Let $P$ be some simple fine partition$,$ 
$S \in \hyper {\cal S}(P)$ and $K = \hyper J \cap S \not= \emptyset.$ 
Then there exist $\vec x_m \in K$ and $\vec x_M 
\in K$ such that $\hyper {f_m}=\hyper {f(\vec x_m)} = 
\Hyper {\inf} \{\hyper f(\vec 
x)\vert \vec x \in K\} = \inf \{\hyper f(\vec x)\vert \vec x \in K\}$ and
$\hyper {f_M} =\hyper {f(\vec x_M)} = \Hyper {\sup} \{\hyper f(\vec 
x)\vert \vec x \in K\} = \sup \{\hyper f(\vec x)\vert \vec x \in K\}$ and\par 
\indent\indent (iii) $(\hyper {f_m)}\, \hyper v(S) \leq \hyper M(S) 
\leq (\hyper f_M)\, \hyper v(S).$ [Note: this is the case where $\hyper L$ 
is the identity map.]
\smallskip
\hrule
\smallskip
\hrule
\smallskip
\hrule
\vskip 18pt
Obviously the rules in IR3 are closely related to those that appear in chapter 
5$,$ section 5.4. Except for 2(iii) these rules are usually tacitly 
assumed. It is useful to repeat the observation made in section 5.4. 
 {\it Certain aspects of rule IR3 could be relaxed if statement 2(iii) is 
formulated in terms of an infinitesimal ``mean value'' property for the 
functional $M.$ However$,$ it is the customary linear form in 2(iii) that appears 
throughout all of the traditional definitions - a form that we are trained to 
seek and experimentally justify.} Notice that the terms of a Riemann sum are 
represented in the inequality in 2(iii) by either $\hyper {f_m}$ or 
$\hyper {f_M}.$ The following theorem is almost obvious.
\vskip 18pt
{\bf Theorem 6.2.2.} {\sl If IR3  holds$,$ then $$M(J) = \int_Jf(\vec x)\, 
dX.$$}\par     
\vskip 18pt 
{\bf Theorem 6.2.3.} {\sl Let compact Jordan-measurable $J \subset R \subset 
\realp n.$ If continuous $f\colon J \to \real,$ then for any partition $P$ of 
$R$ and any $S \in {\cal S}(P),$ where $K = J \cap S \not= \emptyset$ 
there exist $\vec x_m \in K$ and $\vec x_M 
\in K$ such that $f_m=f(\vec x_m) = \inf \{f(\vec 
x)\vert \vec x \in K\}$ and
$f_M = f(\vec x_M) = \sup \{f(\vec 
x)\vert \vec x \in K\}.$} \par         
\vskip 18pt           
The title of this chapter is actually somewhat misleading. Even though all of 
the previous rules relative to the monadic environment are stated in terms of 
this monadic world$,$ generally for the integral$,$ individuals still rely upon 
standard world observations. It is the generation of {differential 
equation models} that utilize this infinitesimal locale exclusively. 
What mental or experimental procedures does a researcher employ in order to 
arrive at the conclusion expressed in IR3 part 2(iii)? {\bf Evidently$,$ these 
standard procedures must be closely related to the rule IR1 - IR3 even if 
they are expressed in the language of a standard mathematical structure.}
In derivations written 
prior to 1981$,$ we read that such conclusions are {``self-evident''} or 
``obvious.'' These self-evident features of informal infinitesimal modeling 
are$,$ of course$,$ some type of unmentioned {infinitesimal reasoning process}. 
Formally$,$ the unmentioned {infinitesimalizing procedure} can now be 
characterized as formal *-transfer and what needs to be determined in the 
large scale or macroscopic environment prior to such infinitesimalizing is 
contained in the premises of the next somewhat obvious proposition. We call a 
set$,$ ${\cal P},$ of simple partitions of $R$ {\bf acceptable} if there exists a simple fine partition $P  \in \hyper {\cal P}.$
\vskip 18pt
{\bf Theorem 6.2.4. ({Self-evident Max. and Min.})} {\sl Let the rectangle $R\subset \realp 
n$ and suppose that compact Jordan-measurable $J \subset R.$ Let $M$ be 
defined as in (1) of IR3$,$ continuous 
$f\colon J\to \real,\ {\cal P}$ an acceptable set of partitions of $R$ and any 
$P \in {\cal P}.$ If for any $S \in {\cal S}(P)$ such 
that $J \cap S \not= \emptyset$ it follows that 
$(f_m)v(S) \leq M(S) \leq (f_M)v(S),$ 
then the infinitesimalizing process IR3 holds.}\par
In order to better apply the self-evident theorem an intuitive discussion
of its content is in order. The values $f_m$ and $f_M$ are  
values for the original defining function restricted to $J \cap S.$ Thus 
essentially for the boundary type subrectangles $S$ (i.e. $S \not \subset
J$ but $S \cap J  \not= \emptyset$) these values have been extended to the 
entire rectangle $S.$ What has been done should be viewed as an application of 
these bounding values of $f$ to the extended configurations $\bigcup
\{S\vert S\in {\cal S}(P)\land J\cap S\not= \emptyset\}.$ {\bf [Important. 
See the simplest rules IR5$,$ IR6 in Appendix 6 where a continuity concept 
eliminates this difficulty.]}\par
\medskip   
\leftline{6.3 \underbar{Simple Applications.}}\par
\medskip
Assuming in the macroscopic world a {continuous distribution of matter} 
throughout a space 
region Synge and Griffith define the {moment of inertia} as $I = \int r^2\, dm$
and then state that ``$dm$ is the mass of an infinitesimal element...''
[{\tt Synge and Griffith [1959:173]}]. For a {uniform rod} of total mass $m$ and 
length $2a$ positioned on the $x$-axis$,$ they state that $dm= (m\, dx)/(2a).$
It is often the case that this concept is extended to the case the mass is not 
uniform but is rather determined by the continuous point density function
$\rho (x).$ In this particular case it is claimed that $dm = \rho (x)\, dx.$
However$,$ from the corrected notions of infinitesimal modeling this last 
statement is simply incorrect. Based upon Corollaries 5.2.3.1 and 5.2.3.2 and 
the definition of supernearness we can only be assured that for an 
infinitesimal subrectangle $R$ and for any $x \in R,\ dm = \hyper M(R,\hyper
\rho )$ is infinitely close to $\hyper \rho (x)\, dx.$ Indeed$,$ the notion  of 
 {``infinitely close''} in this context is not sufficient for a proper 
understanding of the relationship between $dm$ and $\hyper \rho (x)\, dx.$ As 
well be discussed later these two objects must$,$ for a given $dx,$ be 
``closer'' then indicated by the general infinitely close concept. It is$,$ 
therefore$,$ necessary to establish various  elementary applications of 
these infinitesimal rules in a manner distinctly different from the customary 
ones. I envision that many of our present day texts that claim  to teach 
the infinitesimal modeling of natural system behavior will need to be 
completely revised if rigor is to be incorporated.\par
\medskip
{\bf Application 6.3.1.} {\it The quantity of an incompressible fluid passing 
through a 2-dimensional rectangle {\rm [}resp. Jordan-measurable plane region 
$J.${\rm ]}}\par 
Suppose that we have for a macroscopic environment a function $V(x,y)$ that 
represents the point velocity$,$ in a normal direction$,$ of the incompressible 
fluid passing through a 2-dimensional rectangle $R.$ Then the amount of fluid  
passing through $R$ per unit time$,$ $Q(R),$ is $$\int_R v(x,y)\, dX.$$\par
{\sl Derivation.} Before we start this derivation observe that it must be 
considered slightly less rigorous than those that appear in chapter 5. What 
does experience indicates about such a value $Q(R)$? First$,$ since the fluid 
is incompressible then $Q$ is additive on ${\cal C}_{SR}.$ Observation also 
indicates that for a subrectangle $R_S$ of any simple partition
$(V_m)\, v(R_S) \leq Q(R_S) \leq (V_M)\, v(R_S).$ If you agree to this 
experiential argument$,$ then by Theorems 6.2.2 and 6.2.3 the result 
follows.\par
\bigskip
{\bf Application 6.3.2.} {\it Moment of inertia of a circular disc of radius $a$ 
about a line
$\ell$ through the center perpendicular to the plane of the disc.}\par
[Point mass method.] Consider the disc $J$ centered at the origin of our 2-dimensional Cartesian 
coordinate system. Let $\rho (x,y)$ represent a continuous point density 
function for $J.$ The moment of inertia$,$ $I_1,$ for a finite system of point 
masses$,$ $m_i,$ located on the disc each at a distance of $r_i$ from $\ell$
is $I_1 = \sum r_i^2\, m_i.$ This may be rewritten as $I_1 = \sum (x_i^2 + 
y_i^2)\, m_i.$ The first requirement for application of IR3 is that this idea 
be extended to a continuous density function$,$ $\rho (x,y),$ of point masses. 
This implies that we consider the continuous moment of inertia generating function
$f(x,y) =\rho (x,y)\, (x^2 + y^2).$ Using this assumed function it follows that
$$Q(J) = \int_J f(x,y)\, dX.$$
If $\rho$ is a constant$,$ then 
$$Q(J) = \int_J \rho (x,y)\, (x^2 + y^2)\, dX = \int_0^a 2\pi \rho r^3\, dr =
(\pi \rho a^4)/2) = (m/2)a^2.$$\par
{\sl Derivation.} It's clear that we have selected the basic aspects of this 
application to force it to have the properties expressed in IR3. However$,$ what 
experiences do we have with moments of inertia that will allow us to conclude 
that the inequality $(f_m)\, v(S) \leq I(S) \leq (f_M)\, 
v(S)$ holds? If you have such experiences$,$ then the result is 
immediate. If you don't or you cannot adequately explain you intuition$,$ then 
another derivation method would be required.\par
\bigskip
Application 6.3.2 and the like seem to have a very weak derivations. There 
are two notations that will aid in eradicating these derivation difficulties.
Indeed Maxwell explicitly states one of these procedures and there is a  
statement within the premises of application 6.3.2 that is significant and 
often appears when extensions are considered. { \bf This leads to two 
often used$,$ significant and powerful infinitesimal modeling procedures. 
The first is the 
 {extension of constants values} and the second for,
 {point definable quantities,} is the {extension of the finite} to the 
hyperfinite.}\par
\bigskip
\leftline{6.4. \underbar{The Method of Constants.}}\par
\medskip
The method of constants makes a direct appeal to expression (*) of the Infinite 
Sum Theorem (5.2.1) as it appears on page 37. Individuals observe simple 
properties about a functional relative to the assumption that the generating 
function can have constant values. These observations and a simplified 
 physical theory are then transferred to the NSP-world. They then assume that 
these simple properties hold for the NSP-
world and then without justification state that for a continuous generating 
function the value of the functional is but a ``sum '' of the appropriate (infinitesimal) 
quantities. A rigorous procedure is now possible.\par
\medskip
\hrule
\smallskip
\hrule
\smallskip
\hrule
\smallskip
\centerline{IR4-- {Method of Constants}}
(1) In what follows$,$ let for any $A \subset \realp n$ ``int'' denote the 
interior of $A.\ $ Let ${\cal A}=\{{\cal C}_{SR},\{{\rm int}(J\cap S)
\not= \emptyset\bigm\vert S \in {\cal C}_{SR}\}\}.$
We wish to measure a quantity $M$ for a Jordan-measurable set $J 
\subset R \subset \realp n,$ where $M$ is$,$ at least$,$ defined on and   
additive over the members of the set $\{{\cal A}, R-J,J\}$ and
for a subrectangle $S,\ S \subset R -{\rm int}(J)$ it follows that
$\ M(S) = 0$ and $M(R-J) = 0.$ Let $v(J)$ denote the Jordan content.\par
(2) There is a generating function $f(\vec x)$ that is related to the 
functional 
$M$ in the following manner:\par
\indent\indent (i) The function $f$ is bounded on $J.$\par
\indent\indent (ii) Let $P$ be any arbitrary simple fine partition$,$ 
$S \in \hyper {\cal S}(P)$ an arbitrary infinitesimal subrectangle and 
$\emptyset \not=  K = \Hyper {\rm int}(\hyper J \cap S).$ \par
\indent\indent (iii) There exists some $\vec x \in K$ such that 
$\hyper M(S)= \hyper f(\vec x)\hyper v(S)$ or  $\hyper M(S)/\hyper v(S) 
\approx \hyper f(\vec x).$  
\par
\smallskip
\hrule
\smallskip
\hrule
\smallskip
\hrule
\vskip 18pt

Please note carefully where IR4 differs from IR3. First$,$ $f$ is only assumed 
to be bounded and $J$ need not be compact. 
Also note that for $\emptyset \not= K,\ \hyper v(K) \not= 0$ since in the 
standard case nonempty and Jordan-measurable ${\rm int}(J \cap S)$ contains a 
rectangle.\par 
\vskip 18pt
{\bf Theorem 6.4.1.} {\sl If IR4 holds$,$ then $$M(J) = \int_Jf(\vec x)\, dX.$$}
\vskip 18pt
Once again we need some sort of {infinitesimal reasoning process} that leads to 
application of IR4. Within the laboratory or an imaginary mind experiment 
individuals often tacitly observe that the hypotheses of the next  
``self-evident'' theorem hold and$,$ without stating it$,$ take an intuitive *-
transform and obtain IR4. It is IR4$,$ or something akin to it$,$ that appears in 
the literature. However$,$ certain aspects of the next  result are necessary due 
to technical difficulties and one useful requirement is some what unexpected. 
\vskip 18pt
{\bf Theorem 6.4.2. ( {Self-evident Method of Constants })} 
{\sl Let the rectangle $R\subset \realp 
n$ and suppose that Jordan-measurable $J \subset R.$ Let $M$ be defined as in 
(1) of IR4$,$ continuous 
 $f\colon R\to \real,\ {\cal P}$ an acceptable set of partitions of $R$ and any
$P \in {\cal P}.$ If for any $S \in {\cal S}(P)$ such that $\emptyset \not= K= 
{\rm int}(J \cap S)$ there exists some $\vec x \in K$ and some $\vec y \in S$ 
such that (i) $M(K)= f(\vec x)\, v(K)$ and (ii) $M(S)= f(\vec y)\, v(S),$  
then the infinitesimalizing process IR4 holds for $f$ restricted to $J.$}\par
\vskip 18pt
As with IR3 the intuitive notion behind infinitesimal reasoning Theorem 6.4.2 
is the extension of the generating function values to the boundary type 
subrectangles S. Both of our self-evident theorems display an interesting 
phenomenon. The rules IR3 and IR4 require the existence of but one simple fine 
partition for application; but$,$ the self-evident theorems need an infinite 
collection of such partitions with the requisite properties in order to 
guarantee that such a partition exists in the NSP-world. In practice the 
requirement that $f$ be continuous on $R$ may be relaxed to piecewise 
continuity on a set of $J$ covering subrectangles of $R$ that at most overlap 
only on their boundaries. 
 {\bf Of interest is the necessity for premise (ii). In applications it is 
simply assumed from observation that there is a partition of 
subrectangles$,$ $S,$ that are  ``small
enough'' in size that there exists such a constant $f(\vec x),\ \vec x\in
{\rm int}(J\cap S)$ such that (i) holds. The same methodology should also 
indicate that there exists a constant $f(\vec y),\ \vec y \in S$ such that 
(ii) holds. However$,$ \underbar{the (ii) premise is never mentioned as a 
requirement.}}\par
\bigskip
{\bf Application 6.4.1. (Application 6.3.1 revisited.)} 
{\sl The quantity of an incompressible fluid passing 
through a 2-dimensional rectangle {\rm [}resp. Jordan-measurable plane region 
$J.${\rm ]}}\par 
Suppose that we have for a macroscopic environment a  continuous function  
$V(x,y)$ that 
represents the point velocity$,$ in a normal direction$,$ of the incompressible 
fluid passing through a 2-dimensional rectangle $R.$ Let $J$ be any Jordan-
measurable subset of $R.$ Then the amount of fluid  
passing through $J$ per unit time$,$ $Q(J),$ is $$\int_J v(x,y)\, dX.$$\par
{\sl Derivation.} It appears that the method of constants yields a more 
satisfactory derivation of this integral formula. Consider \underbar{any} 
simple 
partition $P$ and \underbar{any} $S \in {\cal S}(P).$ Let $\emptyset \not= 
K = {\rm int}(J \cap S).$ Then experience indicates that there is some 
$\vec x \in K$ such that $Q(K) = V(\vec x)\, v(K).$ Indeed$,$ this can be most 
easily argued by the intermediate value theorem if $J$ is connected. For the 
same reason there is some $\vec y \in S$ such that $Q(S) = V(\vec y)\, v(S).$
If you agree to these seeming innocuous statements relating scalar constant 
velocity normal to $J,$ then the result follows from the Method of Constants.\par
\medskip
Originally the basic infinitesimal reasoning behind Theorem 6.4.2 was simply 
that in the NSP-world $\hyper M(K) \approx \hyper M(S),$ where $\emptyset 
\not= \Hyper {\rm int}(\hyper J \cap S)$ and $S$ is a infinitesimal subrectangle of $\hyper 
R.$ However$,$ technically$,$ there seems to be no way to avoid that fact that 
this idea is incorrect. It may be a good starting point but$,$ infinitesimally,
it is necessary that $\hyper M(K)/\hyper v(K) \approx \hyper M(S)/\hyper v(S).$
\vskip 18pt
{\bf Theorem 6.4.3. (Extended Self-evident Method of Constants)}
{\sl Let the rectangle $R\subset \realp n$ and suppose that 
 Jordan-measurable $J \subset R.$ Let $M$ be defined as in (1) of IR4$,$ continuous 
 $f\colon R\to \real,$ continuous $g\colon R\to \real,\ {\cal P}$ an acceptable set of 
partitions of $R$ and any
$P \in {\cal P}.$ If for any $S \in {\cal S}(P)$ and $\emptyset \not= 
K= {\rm int}(J \cap S)$ there exists some 
$\vec x_1,\ \vec x_2 \in K$ and some $\vec y_1,\ \vec y_2 \in S$ such that 
(i) $M(K)= f(\vec x_1)\, g(\vec x_2)\, v(K)$ and (ii) $M(S)= 
f(\vec y_1)\, g(\vec y_2)\, v(S),$ then the infinitesimalizing process IR4 
holds for $fg$ restricted to $J.$}\par
\vskip 18pt
Why do we need the Extended Self-evident theorem? First$,$ many  
physical and geometric quantities are defined as the product of other 
previously defined  
generating functions. Moreover$,$ it is often the case that these quantities are 
actually defined for the sole purpose of 
applying the differential and integral calculus. One example of this should 
suffice. Consider the following supposedly non-calculus approach to the moment 
of inertia of a circular disc of radius $a$ about a line $\ell$ through the 
center perpendicular to the plane of the disc. In the book {\it Mechanics of 
Engineers} {\tt [Morley [1942]]} the following approach is used. \par

Assume that the density is unity. . . $,$ consider the disc divided 
into $n$ ``ring-shaped strips such as $PQ,$ each of width $a/n.$'' Morley then 
takes the distance of the $p$th strip from the center $O$ as $p \times a/n.$
He then lets the area be the same as that of the rectangle $2\pi \times pa/n 
\times a/n = 2\pi p\, a^2/n^2.$ He then states that the moment of inertia 
about $\ell$ would be this area times the distance of the outer edge $p \times 
a/n$ which yields $$2 \pi a^4p^3/n^4.\eqno(1)$$
If now we add these $n$ strips this would yield$,$ after simplification$,$  
the expression  $$\bigl(\pi a^4/2\bigr) \bigl(1 + {2 \over n} + {1\over 
{n^2}}\bigr).\eqno(2)$$\par                                       
However$,$ on the other hand$,$ if we choose $PO$ to be the distance of the strip 
from the line $\ell$ then assuming the everything else holds the expression 
for the sum all the $n$ strips would be 
$$\bigl(\pi a^4/2\bigr) \bigl(1 - {2 \over n} + {1\over{n^2}}\bigr).\eqno(3)$$\par 
Observe that for either (2) or (3) the limit as $ n \to \infty$ is 
$\pi a^4/2.$ Morley's definition requires that such a limit be taken. If these 
ideas are now applied to
parallel lines interior to the strips$,$ then the same conclusion would 
follow. I discuss aspects of this example more fully at the beginning of the next 
section. \par
{\bf {Remark 6.4.1}} Closer examination of the Self-evident Theorems 6.4.2 is useful.
For a particular linear function $M$ defined as in IR4$,$ the relation between 
the values $M(S)=0,$ where the subrectangle $S \subset R -{\rm int}(J),$ is of 
a special nature.
These values$,$ in general$,$ are assumed to be unaffected by the function 
$f.$ It is as if the function $f$ has degenerated to the zero function when it 
is observed by subrectangles exterior to ${\rm int}(J).$ 
In most practical cases 
this is exactly how it should be. Since we are not interested in the behavior of 
$f$ on such exterior objects. Two things often happen in practice. 
The function $f\colon J \to (a,b) \subset \real$ is only assumed to be 
continuous on $J.$ It 
does not matter whether or not $f$ is the restriction of a function that is 
continuous on $R.$ The other case is that the function $f$ is explicitly given 
and it is indeed continuous on some bounding set $R.$ If $J$ is compact and 
$f\colon J \to (a,b)$ is continuous and not explicitly expressed$,$ $M$ is 
defined as in IR4 and the remaining hypotheses hold for any continuous 
extension $F\colon R\to (a,b)$ of $f$ to $R,$ then by the 
Teitze Extension Theorem $M(J) = \int_Jf(\vec x)\, dX.$ The same modifications 
could be made$,$ if necessary$,$ to the Extended Self-evident Theorem 6.4.3.\par
\bigskip
{\bf Application 6.4.2.(Application 6.3.2 revisited.)} {\it Moment of inertia of a circular disc of radius $a$ 
about a line $\ell$ through the center perpendicular to the plane of the disc.}\par
 Consider the disc $J$ centered at the origin of our 2-dimensional Cartesian 
coordinate system. Let $\rho (x,y)$ represent a continuous point density 
function for $J$ which is assumed to be continuous on some $R \subset 
\realp 2$ such that $J\subset R.$ Then the moment of inertia$,$ $Q(J),$ of $J$ 
about $\ell$ a line perpendicular to the disc at its center is 
$$Q(J) = \int_J f(x,y)\, dX,$$
where $f(x,y) =\rho (x,y)\, (x^2 + y^2).$ Using this assumed 
function it follows that if $\rho$ is a constant$,$ then 
$$Q(J) = \int_J \rho (x,y)\, (x^2 + y^2)\, dX = \int_0^a 2\pi \rho r^3\, dr =
(\pi \rho a^4)/2) = (m/2)a^2.$$\par
{\sl Derivation.} Kinetic energy experimentation with a flywheel with movable weights 
attached to its surface indicates that for actually physical regions $J,$   
there does exist a ``small enough'' simple partition rectangle $S$ 
such that there are two points $ \vec x,\ \vec y = (y_1,y_2) \in 
{\rm int}(J \cap S)\not= \emptyset$ and  
$Q(J \cap S) = 
\rho (\vec x)\, (y_1^2 + y_2^2).$ For this particular subrectangle $S$ 
it is also observed that 
their exist two such points in $S$ with the same property for $Q(S).$ The 
same conclusion apparently would hold for all other simple partition 
rectangles with a ``smaller diagonal length.'' [This generalization to 
``all'' such rectangles appears reasonable.] Let 
$(m_1,m_2) \in {\nat}^2$ be the generator of the simple partition $P$ such 
that $S \in {\cal S}(P).$ Let the set of partitions ${\cal P}$ generated by 
the pair $\{(x,y)\bigm \vert (x,y) \in {\nat}^2 \land (x \geq m_1) \land 
(y \geq m_2)\}.$ Then ${\cal P}$ is an acceptable partition. Letting $f(x,y)= 
(x^2 + y^2),$ which is continuous on $R,$ then the result follows from 
the Method of Constant Theorem 6.4.3.
\bigskip
\leftline{6.5. \underbar{The Hyperfinite Method.}}\par 
\medskip
Except when the operational approach to physical quantities is used$,$ most 
authors when discussing the properties of rigid body motion immediately force 
upon the reader the imaginary notion of the point ``particle.'' 
Synge and Griffith$,$ in the text mentioned$,$ use exclusively this technique. Also$,$ even though 
it may not be apparent from his example in the previous section$,$ Morley 
motivates all of his derivations similarly. It is now possible to establish in 
a somewhat rigorous fashion that the technique of particles is adequate. \par
As a prototype$,$ we again concentration upon a  nonempty plane Jordan-measurable 
region $J.$ Letting 
$\rho$ be an appropriate density function and $d$ and appropriate distance 
function$,$ from line $\ell$ of rotation one establishes that the moment of 
inertia is $$I(J) = \int_J \rho(\vec x)\, d^2(\vec x)\, dX.\eqno (1)$$
By *-transfer of the hypotheses of Theorem 6.4.3 it follows that for a simple 
fine partition $P= \{S_1,\dots,s_\Gamma\}$ of $\hyper R$ and two hyperfinite 
sequences$,$ $Q_1^\prime= \{\vec x_1,\ldots,\vec x_\Gamma\},\ Q_2^\prime
= \{\vec y_1,\ldots,\vec y_\Gamma\},\ \vec x_i,\vec y_i$  members of  
$\emptyset \not= K_i = \Hyper {\rm int}(\hyper J \cap S_i)$ such that 
$\hyper I(K_i) = \hyper \rho (\vec x_i)\, \hyper d(\vec y_i)\hyper v(K_i).$ 
Observe that properties of  
Jordan-measurable sets and the measure  $v$ imply that ${\rm int}(\hyper J 
\cap S_i)$ is *-Jordan-measurable and that $\hyper v(K_i) = \hyper v(\hyper J \cap 
S_i);$ which leads to $\hyper I(\hyper J \cap S_i) = \hyper \rho (\vec x_i)
\, \hyper d(\vec y_i)\hyper v(K_i) =\hyper \rho (\vec x_i) \, \hyper d(\vec y_i)\hyper 
v(\hyper J \cap S_i) = \hyper I(J\cap S_i).$ There is a nonempty internal
${\cal A} \subset \hypernat$ such that ${\cal A} = \{n\bigm\vert n \in 
\hypernat \land 1 \leq n \leq \Gamma \land \hyper J \cap S_n \not= 
\emptyset\}.$ From the additivity of $I$ it follows that
$$\hyper I(J) = \sum_{j \in {\cal A}} \!\!\! \hyper \rho (\vec x_j) \, 
\hyper d(\vec y_j)\hyper v(\hyper J \cap S_j).\eqno (2)$$ \par                   
In the standard case for every $P \in {\cal P}$ and $\{S_1,\ldots,S_n\} = 
{\cal S}(P)$ there are$,$ of course$,$ two finite sequences  
$Q_1= \{\vec x_1,\ldots,\vec x_n\},\ Q_2
= \{\vec y_1,\ldots,\vec y_n\},\ \vec x_i,\vec y_i$  members of $\emptyset 
\not= {\rm int}(J\cap S_i)$ such that
$$I(J) = \sum_{j=1}^n\rho (\vec x_j) \, d(\vec y_j) v(J \cap S_j)= \st 
{\hyper I(J)}.\eqno (3)$$ 
It's equation (3) above that yields the concept of the particle point masses 
by defining $m_j = \rho (\vec x_j)\, v(J \cap S_j)$ and$,$ hence$,$ $$ I(J) =
 \sum_{j=1}^nm_j \, d(\vec y_j) v(J \cap S_j)\eqno (4)$$
If now one assumes the particle point mass equation (4) holds for each member 
of ${\cal P}$ then *-transfer yields
$$\hyper I(J) = \sum_{j \in {\cal A}} \!\!\! \hyper m_i\, 
\hyper d(\vec y_j)\hyper v(\hyper J \cap S_j).\eqno (5)$$ \par       
The process of introducing (i) the additivity of $I$$,$ (ii) simplifying the 
terms of the sum$,$ if possible$,$ by reducing to other geometrical quantities$,$ 
and (iii)  extending to the hyperfinite$,$ is called the 
{\bf {hyperfinite method}}. However$,$ \underbar{this method is unnecessary}
since it is 
but a simple extension of the Method of Constants. Note that for non-discrete
quantities it is$,$ 
technically$,$ not correct to say that we ``extend the finite sum of quantities 
to an infinite sum of such quantities'' or some similar expression as is 
often done by Maxwell and others.
The sum is not$,$ in general$,$ an (external) infinite sum$,$ but$,$ is an internal
hyperfinite sum. The conceptual and formal differences between these two 
concepts are considerable.\par 
But looking at our definition 5.1.1 for the integral$,$ we are lead to the$,$  
usually not appreciated$,$ approximation methods. It follows that for any 
position real $r$ there is a simple partition $P$ and a finite set of 
subrectangles $\{S_1,\ldots,S_2\} = {\cal S}(P)$ such that for any 
set of intermediate partition points $Q = \{\vec x_1,\ldots \vec x_n\},$
$$\Big|\sum_{i=1}^n \hat {\rho}(\vec x_i)\, \hat d^2(\vec x_i)\, v(S_i) - I(J) 
\Big| < r.\eqno (6)$$\par                                                  
Now each $S_i$ is of one of three types. (a) $S_i \subset R -{\rm int}(J),$ 
(b) $K = {\rm int}(J \cap S_i) =  {\rm int}(J) \cap {\rm int}(S_i) \not= 
\emptyset$ and ${\rm int}(S_i) \not\subset {\rm int}(J),$
(c) ${\rm int}(S_i)\subset {\rm int}(J).$ Let $\partial A,\ A \subset \realp 
n$ denote the {\bf {boundary points}} of $A.$ In case (a) since $v(\partial 
J) = 0,\ v(J) = v({\rm int}(J)),\ \partial J \subset R - {\rm int}(J)$ and 
$v(S) \not= 0$ then there exists some $\vec x \in S$ such that $\vec x \notin J.$
In case (b)and (c) there is some $\vec x \in {\rm int}(J) \cap S.$ 
Consequently$,$ there is a finite sequence of points $S =\{\vec y_1,\ldots,\vec 
y_n\}$ that contains a subsequence $S_0 = \{\vec y_{k(1)},\ldots,\vec y_{k(m)}\}$ such 
that
$$\Big|\sum_{i=1}^m {\rho}(\vec y_{k(i)})\, d^2(\vec y_{k(i)})\, v(S_i) - I(J) 
\Big| < r,\eqno (7)$$
where the $S_i$ are all of type (b) or (c) and each $\vec y_{k(i)} \in {\rm 
int}(J).$ Yet there exists a simple partition and another sequence $S^\prime =\{\vec z_1,\ldots,
\vec z_p\}$ that contains a subsequence $S_0^\prime = \{\vec 
z_{k(1)},\ldots,\vec z_{k(q)}\}$ such that 
$$\Big|\sum_{i=1}^q {\rho}(\vec z_{k(i)})\, d^2(\vec z_{k(i)})\, v(S_i) - I(J) 
\Big| < r,\eqno (8)$$
and all of the rectangles $S_i$ are interior rectangles of type (c) and each 
$\vec z_{k(i)} \in {\rm int}(J).$ \par
%\end%
Of course these special selections of members in ${\rm int}(J)$ extended to
the NSP-world and equation (8) reads as
$$\sum_{i=1}^\Gamma \hyper {\rho}(\vec z_{k(i)})\, \hyper 
d^2(\vec y_{k(i)})\, \hyper v(S_i) \approx I(J) ,$$
and all of the infinitesimal subrectangles $S_i$ are interior subrectangles 
of type (c) and each $\vec z_{k(i)} \in \Hyper {\rm int}(J).$  \par
Thus there are many different NSP-world configurations and expressions that 
have the same N-world effects. These effects are obtain by application of the 
standard part operator. However$,$ even 
though it is obvious that all of the procedures discussed in this section and 
the above numbered expressions give a very detailed and analytically correct 
approach to infinitesimal modeling - an approach that reveals  much about the 
nature of the NSP-world model - an early introduction of these infinitesimal 
concepts into an elementary exposition for physical modeling is often 
unnecessary and,
indeed$,$ they will tend to submerge the new discoveries expressed by the 
self-evident theorems. On the other hand$,$ under certain circumstances$,$ it  
appears necessary to include these infinitesimal notions. In the 
next section$,$ I diverge briefly into certain instructional aspects of these 
discoveries.\par
\bigskip
\leftline{6.6. \underbar{Instruction.}}\par
\medskip
The basic methods of infinitesimal modeling revolve about the application of 
simplified physical or geometric theories to ``simple configurations'' within 
the N-world environment. For the physical sciences$,$ natural system behavior 
is viewed locally with many of the requirements of the self-evident theorem 
assumed. Physical modeling is viewed as objective while geometric modeling is 
subjective. Physical modeling deals with observation and experience with the 
behavior of natural systems. Geometric modeling$,$ while originally motivated by 
physical concerns$,$ has become a subject of abstract definitions. This is 
obvious from the applications made in chapter 5 for there$,$ the properties of the 
geometric measures are paramount and are global properties 
obtained from mathematical experience. The local self-evident rules and even 
the infinitesimal rules in this chapter are more relevant to natural system 
behavior. However$,$ both the geometric and physical may be more closely 
associated with IR3 and IR4 then first assumed.\par
Experience dictates that the simplest and first level of comprehension for 
physical concerns is obtained from the two theorems on the 
Self-evident Method of Constants. They would be the easiest to apply for the 
neophyte. There is a reason$,$ 
however$,$ why in certain cases the actual infinitesimal rules IR3 and IR4 need 
to be applied. This is especially the case for geometric measures. These rules 
would be a second level of comprehension and this entails a certain basic 
familiarity with infinitesimal concepts. Moreover$,$ the generating functions 
for a perceived linear functional need not be the integrand utilized. 
Indeed$,$ many of the rules and integral notions within chapter 5 and chapter 6 
can be expressed by hyperfinite summation that yields internal functionals 
not just the standard extension of a standard concept. An example of this is 
our proof of the ``length of a curve integrand'' in appendix 6 on pages 
209--210.\par
In the proof for the length of a continuously differentiable curve 
$c\colon [0,1] \to \realp n$ an internal functional is defined for the length 
of a$,$ possibly broken$,$ hyperpolygonal line determined by a hypercurve $\hyper c
\colon \Hyper [0,1] \to \hyperrealp n$ with coordinates functions 
$\hyper f_j \colon \hyper [0,1] \to \hyperreal.$ This internal function may be 
consider defined as follows: let $\{t_i\}$ be any hyperfinite partition of 
$\Hyper [0,1].$ Let $T$ be any nonempty internal$,$ hence hyperfinite$,$ subset of
$\{t_i\}$ containing two or more members. For each pair of elements in $T$ 
generate the hyperline segments determined by the curve $\hyper c$ and 
consider $\Vert \vec v_i \Vert$ from Example 4.4.1.A.  Finally$,$ 
consider the internal hyperfinite sum operator over the internal set of the
$\Vert \vec v_i \Vert$ determined by $T.$ Let $L$ denote this internal 
functional. Now for each internal subinterval $S$ determined by $T$ define 
$L(S) = L(\Hyper {\rm int}(S)).$ This is our *-additive functional defined 
over the set of all such internal functions $c.$ Under the given hypotheses 
it is shown that for any internal subrectangle $S$$,$ taken from a simple fine 
partition of $\Hyper [0,1]$ there exists some $t_i \in S$ such that 
$ \Vert \vec v_i \Vert/dt \approx \sqrt {\sum_{j=1}^n (\hyper f_j^\prime 
(t_i))^2}.$                                                              
The expression on the right is the standard extension of the 
standard function $\sqrt {\sum_{j=1}^n (f_j^\prime (t_i))^2}$ while the 
expression on the left is $L(S)/dt.$ Thus our 
\underbar{basic definition} for the length of a curve - one that with a 
slight generalization is also shown to be equivalent to rectifiability -  
almost satisfies IR4 with the exception that IR4 is written in terms of a 
standard additive functional not an internal *-additive nonstandard entity. 
It is obvious how IR3 and IR4 would be modified to include such  internal 
functionals.\par
The fact that geometric definitions are subjective in character would allow 
us to define many such concepts entirely in terms of infinitesimal concepts.   
This would be an aid in developing many integral statements about 
n-dimensional geometric properties in terms of integrals defined on 
k-dimensional spaces$,$ where $k < n.$ \par
\bigskip
\leftline{6.7 \underbar{Realism.}}\par
I will not$,$ at this point$,$ dwell upon the philosophical modeling concept 
termed {``realism.''} This concept assumes that all mathematical objects within 
the mathematical formalism correspond to objects within objective (physical) 
reality. There is$,$ however$,$ an interesting historical fact relative to the 
rejection of realism. Even though {Bohr} rejected complete realism for his 
concept of quantum mechanics in order to force a type of {physical completeness} 
upon this theory - a completeness that is now known to be in error [{\tt Aerts 
[l984]}] - he was not the first to do so. In his 1909 treatise on his theory  
of electrons$,$ {Lorentz} accepted the notion of an ether but {rejected 
complete realism.} ``{\it I should add that$,$ while thus denying the real 
existence of ether stresses$,$ we can still avail ourselves of all the 
mathematical transformations by which application of the formula (43) may be 
made easier. We need not refrain from reducing the force to a 
surface-integral$,$ and for convenience's sake we may continue to apply to the 
quantities occurring in this integral the name of stresses. Only we must be 
aware that they are only imaginary ones$,$ nothing else than auxiliary 
mathematical quantities. Perhaps all this that has been said about the 
absolute immobility of the ether and the non-existence of the stresses$,$ may 
seem somewhat startling.}'' {\tt Lorentz [1952:31]} \par
As far as geometric modeling is concerned the concept of realism is not as 
significant as when it is applied to physical problems. Infinitesimal
modeling may be considered as a partial realism. The notion of a continuum 
model in an ultimately discrete world seems non-realistic. However$,$ what is 
being modeled is macroscopic and large scale behavior$,$ macroscopic and large 
scale effects upon observers and other objects. It is not the actual physical 
entities that are being modeled; but$,$ rather$,$ the effects these entities  
produce relative to a specific set of scenarios.\par
However$,$ on another level$,$ modern physical theories often deal with entities 
that are not directly observable$,$ only their indirect effects emerge within the 
laboratory setting. Technically$,$ these entities are speculations that may in 
time be replaced by yet other speculative objects. Indeed$,$ using special 
techniques$,$ the language that actually expresses physical 
theories is mathematically embedded into the natural numbers by a G\"odel 
coding; which is then embedded into a nonstandard structure. The linguistics 
of the physical theory produces a new collection of statements that tells us 
about the behavior of a new world$,$ called in general$,$ the NSP-world. {\tt [Herrmann 
1987]} Within this world we are not speculating ad hoc about infinitesimal 
objects$,$ it is the standard theory that generates their logical existence and 
even predicts some of their extraordinary properties. Whether or not such 
entities are accepted as ``real'' or not depends upon their usefulness. 
In case you may have missed it previous, I again quote the following from the first paper Robinson published  relative to his formal 
theory of infinitesimals. ``{\it For phenomena on a different scale$,$ such as 
considered in Modern Physics$,$ the dimensions of a particular body or process 
may not be observable directly. Accordingly the question whether or not a 
scale of non-standard analysis is appropriate to the physical world really 
amounts to asking whether or not such a system provides a better explanation 
of certain observable phenomena than the standard system of real numbers. 
The possibility that this is the case should be borne in mind.''} Fine Hall$,$ 
Princeton University. {\tt [Robinson l961]}
 
\vfil\eject
\centerline{Chapter 7.}
\medskip                
\centerline{\bf PURE INFINITESIMAL}
\centerline{\bf INTEGRAL MODELING}
\bigskip
\leftline{7.1 \underbar{Brief Discussion.}}\par  
This basic manual is intended to present brief accounts of various approaches 
to infinitesimal modeling so that individuals may select the method the is 
most appropriate for their discipline.  
In Chapter 5$,$ the standard and global type rules for integral modeling are 
introduced in a piecemeal fashion. Each of the applications in that chapter 
depends upon an in depth knowledge of the properties of a specific functional
that measures a specific geometric or physical quantity. Actually$,$ as is 
clear from the derivations$,$ it is more 
likely that the properties of the specific functional are selected in concert  
with the requirements of the Infinite Sum Theorem. As for Chapter 6$,$ although the 
Self-evident theorems are more general in character$,$ they still rely upon the 
Infinite Sum Theorem as the proofs of Theorems 6.2.2 and 6.4.1 indicate. 
However$,$ in general$,$ neither the approach of Chapter 5 nor the approach of 
Chapter 6 appears within the pre-1960 literature. Previous appliers of the 
concepts of infinitesimal reasoning relied heavily upon a {pure infinitesimal 
approach} that is highly discipline orientated and which makes a direct appeal 
to the {Riemann styled sum} and Definition 5.1.1. This pure infinitesimal 
approach relies upon expedient definitions for elusive geometric or physical 
entities within the NSP-world - entities called by the general term 
``{elements}.'' The collection for geometric modeling includes the basic 
elements of IR2 in Chapter 6 and various Euclidean compositions and 
decompositions of these configurations.\par
{\bf For each discipline$,$ the collection of elements is 
sequentially defined from the most basic through the more complex in a$,$ 
hopefully$,$ consistent manner.} It is almost always the case that the defined 
objects are taken from those in the N-world that behave in the simplest 
possible discipline defined manner and this simplistic behavior is then impressed by 
*-transfer upon the NSP-world. For geometry$,$ the selection of 
an acceptable set of elements is somewhat more arbitrary than for physical  
disciplines. The experimental or observational disciplines rely upon the 
concept of what might be termed as \underbar{simple or idealized behavior}
and well-grounded physical methods of approximation.  However$,$ after the 
infinitesimal elements have been selected then$,$ in all cases of which I am 
aware$,$ the rigorous derivations that these elements lead to an integral model 
are very similar. Thus$,$ there may be an unconscious interplay between 
the appropriate element$,$ as axiomatically selected$,$ and the rigorous method.
Examination of the literature leads to the following general observations.\par
(i) The infinitesimal elements are axiomatically selected$,$ per discipline$,$ and 
the concept of (hyperfinite) summation is applied. Whether or not this 
summation can be calculated by means of a Riemann styled integral is somewhat 
secondary.\par
(ii) In anticipation of an analytical approach$,$ it has become customary to 
employ certain symbols and terms that tend to describe inaccurately 
the actual situation. As a prototype consider the symbol $ds.$ In 
infinitesimal modeling this is a {\bf general symbol} that represents any of 
the lengths of hyperline segments that comprise any of the hyperpolygonal 
curves that determine the length of a geometric curve. With respect to our 
notation $ds = \Vert \vec v_i \Vert.$ Only when actual attempts are made to 
calculate the length of the entire curve are
additional analytical constraints placed upon the representing functions. For 
modeling purposes$,$ the 
geometric or physical elements are considered more basic in character in order 
to avoid$,$ if possible$,$ all of the well-known analytical difficulties.\par
(iii) Restriction to the monadic environment takes place only when an 
n-dimensional quantity is being measured with respect to an 
n-dimensional partition. Otherwise$,$ the infinitesimal elements span both the 
N-world and the NSP-world; portions are monadic and portions are not.\par
(iv) The selected infinitesimal elements are not unique$,$ even in their general 
character. \par
In the next section$,$ certain 
elements will be selected for problems in elementary infinitesimal geometry. 
Following this$,$ analytical constraints will be applied
so that elementary integral expressions can be rigorously derived.
These rigorous derivations display common derivation features that 
should be apparent. Indeed$,$ the common 
features already appear in the proof on page 127 that gives the length of a 
continuously differentiable curve relative to the hyperpolygonal 
approximating curves discussed in examples 4.4.1.A and 4.4.1.B.\par
\bigskip
\leftline{7.2 \underbar{Geometric Elements.}}
\medskip
The first and most basic measure that appears in Euclid's geometry book is the 
length of a line segment. I repeat$,$ with slight modifications$,$ the information in example 4.4.1.A for the 
generalization of this to curves.\par
Recall that a curve is a continuous map 
$c\colon [0,1] \to \realp n.$ This is equivalent to considering $c$ as 
determined by n continuous coordinate
functions 
$x_i = f_i(t),\ 1 \leq i \leq n$ each defined on $[0,1] \subset \real.$ 
The geometric curve $C$ 
determined by these functions is usually considered as the set 
$\{(x_1,\ldots,x_n)\vert t \in [0,1]\}.$ The *-transform of these 
defining functions leads to the functions $x_i = \hyper f_i(t),\ 1 \leq i 
\leq n,$ each defined on $\Hyper [0,1] \subset \hyperreal$ and they generate the 
``hypercurve" $\Hyper C \subset \hyperrealp n.$\par
Let $Q$ be any fine partition of $\Hyper [0,1].$  By *-transfer, $Q$ behaves 
like an ordered partition 
of the interval $[0,1]$ as defined in the standard sense and we write
$Q = \{t_i \bigm\vert 0\leq i \leq \Gamma\}.$ The internal 
hyperfinite set $Q$ generates the 
internal set of ``points'' $P = \{(\hyper f_1(t_i),\ldots,\hyper 
f_n(t_i))\vert t_i \in Q\}$ that are members of the hypercurve $\Hyper C.$
Now for each $i = 0,\ldots,\Gamma - 1,$ and each $j,\ 0\leq j \leq n$ 
let $\hyper f_j(t_{i+1}) - \hyper 
f_j(t_i) = d(j,i).$ Since $c$ is continuous then each $d(j,i) \in \monad 0.$ 
 For each $i \in \hypernat$ such that $0 \leq i \leq 
\Gamma -1,$ the internal set $\ell_i = \{(x_1,\ldots,x_n)\vert
\forall j \in \hypernat,\ 0\leq j \leq n,\  x_j =
\hyper f_j(t_i) + t(d(j,i)) \land t \in 
\Hyper [0,1] \}$ is a hyperline segment connecting the two points 
$(\hyper f_1(t_i),\ldots,\hyper f_n(t_i)),\ (\hyper f_1(t_{i+1}),
\ldots,\hyper f_n(t_{i+1}))$ on the curve $\Hyper C.$ From this one obtains 
the internal hyperpolygonal curve ${\cal P}_Q = \bigcup \{\ell_i \vert 
0\leq i \leq \Gamma - 1\}.$ As to the length of ${\cal P}_Q $ simply 
extend the concept of length in the classical sense by defining for each 
$i = 0,\ldots,\Gamma -1$ the vector $\vec v_i = (d(1,i),\ldots,d(n,i)) \in 
\hyperrealp n.$ Then let the hyperfinite sum $\sum_{i=0}^{\Gamma - 
1} \Vert \vec v_i \Vert = \vert {\cal P}_\Gamma \vert \in \hyperreal.$ {\bf 
For the curve$,$ the 
geometric element 
is the hyperline segment $\ell_i.$} Even though$,$ 
in general$,$ you would have a different hyperpolygon with a different hyperreal 
length for $\forall \Gamma \in {\nat}_\infty$ the following Proposition is 
proved in appendix 6. \par
\vskip 18pt
{\bf Theorem 7.2.1.} {\sl Consider continuous $c\colon [0,1] \to \realp n.$ Then 
$c$  is rectifiable if and only if there exists some $L \in \real$ such that 
for every fine partition $Q$ of $\Hyper [0,1]$
$$\st {\Hyper \vert {\cal P}_Q \vert} = L.$$}\par
\vskip 18pt
{\bf Application 7.2.1.} {\it The 2-dimensional area between two continuous 
curves using the 1-dimensional integral.}\par 
Give continuous $f\colon [a,b] \to \real$ and continuous $g\colon [a,b] \to \real,$
where $f(x) \leq g(x)$ for each $x\in [a,b].$ Then the area between the curves
$f$ and $g$ is
$$\int_a^b(g(x) - f(x))\, dx.$$\par
{\sl Derivation.} From the NSP-world view point each curve is viewed as any 
hyperpolygonal line created by any simple fine partition of $\Hyper [a,b].$
Let $Q$ be such a simple fine partition. For a given $I_i=[x_i,x_{i+1}] \in 
\hyper {\cal S}(Q),$ let $\ell(f)_i$ be the corresponding hyperline segment 
determined by the curve $f$ and $\ell(g)_i$ that determined by $g$ for the 
interval $I_i.$ Select as the infinitesimal element the {\bf hypertrapezoid}$,$ $T_i,$ 
with sides $\ell(f)_i,\ \ell(g)_i$ and the hyperline segments with end points
$\{(x_i,\hyper g(x_i),\ (x_i,\hyper f(x_i))\}$ and $\{(x_{i+1},
\hyper g(x_{i+1}),\ (x_{i+1},\hyper f(x_{i+1}))\}.$ 
Observe that this hypertrapezoid selection 
is consistent with the NSP-world view of a curve and the hypertrapezoid is 
composed of *-Euclidean composition or decomposition of our previously 
defined basic elements. By *-transfer the hypertrapezoid's hyperarea $\hyper A(T_i) =
(1/2)((\hyper g(x_i) - \hyper f(x_i)) + (\hyper g(x_{i+1}) - \hyper 
f(x_{i+1})))\, dx.$ As was done for the length of a curve$,$ let
$T(Q) = \bigcup\{T_i\bigm\vert 0\leq i\leq \Gamma -1\}$ and consider the 
hyperfinite sum  $$\sum_{I_i \in \hyper {\cal S}(Q)} \!\!\!\!\hyper A(T_i)= 
\hyper A(T(Q)).\eqno 
(1)$$\par
From this definition it follows that 
$$\st{\hyper A(T(Q))} = \int_a^b (g(x) - f(x))\, dx.\eqno(2)$$ \par
$\bigtriangleup$ Rather than relegate the proof that there exists a unique $r \in \real$ such 
that for any simple fine partition $Q$ of $\Hyper [a,b],\ \hyper A(T(Q)) \in \monad 
r$ to the appendix we present it as follows by investigating  the 
behavior of the expression $((\hyper g(x_i) - \hyper f(x_i)) + (\hyper g(x_{i+1}) - \hyper 
f(x_{i+1}))).$ From the uniform continuity of $g,\ f$ it follows that
$(\hyper g(x_i) - \hyper f(x_i)) \approx (\hyper g(x_{i+1}) - \hyper 
f(x_{i+1})).$ Hence$,$ $\hyper A(T_i) = ((\hyper g(x_i) - \hyper f(x_i))+ \delta_i)\, 
dx,$ where $\delta_i \in \monad 0.$  Consequently,
$$\sum_{I_i \in \hyper {\cal S}(Q)} \!\!\!\!\hyper A(T_i)= \sum_{I_i \in \hyper {\cal S}(Q)}\!\!\!\!
((\hyper g(x_i) - \hyper f(x_i))\, dx+ \sum_{I_i \in \hyper {\cal S}(Q)}\!\!\!\!
\delta_i\, dx.\eqno (3)$$                                                      
By considering the finite case it follows that there exists some
$\delta \in \{\vert\delta_i\vert\bigm\vert 0\leq i\leq \Gamma -1\}$ such that 
$\delta = \max\{\vert\delta_i\vert\bigm\vert 0\leq i\leq \Gamma -1\}$ and$,$ hence$,$ $\delta \in 
\monad 0.$ By *-transfer we have that 
$$\big|\!\!\!\!\sum_{I_i \in \hyper {\cal S}(Q)} \!\!\!\! \delta_i\, dx \big| \leq \!\!\!\!
\sum_{I_i \in \hyper {\cal S}(Q)} \!\!\!\!\vert \delta_i\vert\, dx \leq 
\delta\!\!\!\!\sum_{I_i \in \hyper {\cal S}(Q)} \!\!\!\!\, dx = \delta (b-a) \in \monad 0. 
\eqno(4)$$
Thus
$$ \sum_{I_i \in \hyper {\cal S}(Q)} \!\!\!\! \delta_i\, dx= \lambda \in \monad 
0.\eqno (5)$$
Therefore$,$ (2) can now be written as
$$\sum_{I_i \in \hyper {\cal S}(Q)} \!\!\!\!\hyper A(T_i)= \sum_{I_i \in \hyper {\cal S}(Q)}\!\!\!\!
((\hyper g(x_i) - \hyper f(x_i))\, dx+\lambda \eqno (6)$$                                                      
Since the function $(g(x) - f(x))$ is bounded and integrable and the 
partition $Q$ can also be considered as an internal immediate partition 
then Theorem 5.2.2  yields
$$\st{\!\!\!\!\sum_{I_i \in \hyper {\cal S}(Q)} \!\!\!\!\hyper A(T_i)}= \st 
{\!\!\!\!\sum_{I_i \in \hyper {\cal S}(Q)}\!\!\!\!((\hyper g(x_i) - \hyper 
f(x_i))\, dx} = \int_a^b(g(x) - f(x))\, dx.\eqno (7)$$
and this completes the proof. $\bigtriangleup$\par
\bigskip
The above derivation certainly appears interesting except for the obvious fact 
that the {\bf hypertrapezoid configuration is not a unique geometric element.}
The key to this derivation is the process that begins with the statement that 
the values are infinitely close ($\approx$) and the expressions 
(3)(4)(5)(6)(7). The derivation would hold for hyperrectangles$,$ hypertriangles and the like.
Hyperrectangles yield the same 1-dimensional integral expressions. 
Hypertriangles yield 1/2 the indicated integral expression$,$ and so forth. This 
non-uniqueness for elements is sometimes considered one of the basic difficulties  
with infinitesimal modeling both for the geometric and physical. However$,$ 
since the concept of area obtained by integrals is supposed to be an 
extension of the Euclidean area measure then any element that does not yield 
such an extension can be rejected. {\it For elementary calculus$,$ the use of the 
hypertrapezoid is very expedient when  the usual 1-dimensional integral 
expressions for the basic geometric measures are considered. 
It will be our element of 
choice in this case.}\par
\bigskip
{\bf Application 7.2.2.} {\it Volume of revolution using the 1-dimensional 
integral.}\par  
Given continuous $f\colon [a,b] \to \real$ where $f(x) \geq 0$ for each 
$x \in [a,b].$ Then the volume obtained by rotating this curve about the 
$x$-axis is
$$\int_a^b \pi f^2(x)\, dx.$$\par
{\sl Derivation.} Let $Q$ be a simple fine partition of $\Hyper [a,b].$
For a given $I_i = [x_i,x_{i+1}] \in \hyper {\cal S}(Q),$ let $\ell(f)_i$ be 
the corresponding hyperline segment determined by $f.$ Consider the 
hypertrapezoid composed of hyperline segments $\ell(f)_i,\ I_i$ and
hyperline segments with end point $\{(x_i,0),(x_i,\hyper f(x_i)\},\ 
\{(x_{i+1},0),(x_{i+1},\hyper f(x_{i+1})\}.$ The rotation of this 
hypertrapezoid about the $x$-axis is the frustum of a right cone$,$ $F_i.$
 By *-transfer$,$ the volume of such a *-geometric configuration is
$\Hyper V(F_i)) = (1/3)\pi (f^2(x_i) + f(x_i)f(x_{i+1}) + f^2(x_{i+1})).$ Assuming the 
usual *-additivity$,$ the volume of the entire configuration $F(Q) = \bigcup\{F_i \bigm 
\vert 0 \leq i \leq \Gamma -1\}$ is                                         
$$\sum_{I_i \in \hyper {\cal S}(Q)} \!\!\!\!\Hyper V(F_i)= \Hyper V(F(Q)).\eqno 
(8)$$\par                                                          
From this definition it will follow that
$$\st {\Hyper V(F(Q))} = \int_a^b \pi \, f^2(x) \, dx.\eqno(9)$$ \par
Noting that $f^2(x_i) + f(x_i)f(x_{i+1}) + f^2(x_{i+1}) \approx  
f^2(x_i) + f(x_i)f(x_i) + f^2(x_i) = 3f^2(x_i)$ then (9) follows in the 
same manner as in the derivation for application 7.2.1.\par
\bigskip
{\bf Application 7.2.3.} {\it The surface area of revolution using the 
1-dimensional integral.}\par
Given continuous differentiable $f\colon [a,b] \to \real$ where $f(x) \geq 0$ for each 
$x \in [a,b].$ Then the surface area obtained by rotating this curve about the 
$x$-axis is
$$\int_a^b 2\pi f(x)\sqrt {1 + (f^\prime (x))^2}\, dx.$$\par
{\sl Derivation.} Let $Q$ be a simple fine partition of $\Hyper [a,b].$
For a given $I_i = [x_i,x_{i+1}] \in \hyper {\cal S}(Q),$ let $\ell(f)_i$ be 
the corresponding hyperline segment determined by $f.$ Consider the 
hypertrapezoid composed of hyperline segments $\ell(f)_i,\ I_i$ and
hyperline segments with end point $\{(x_i,0),(x_i,\hyper f(x_i)\},\ 
\{(x_{i+1},0),(x_{i+1},\hyper f(x_{i+1})\}.$ The rotation of this 
hypertrapezoid about the $x$-axis is the frustum of a right cone$,$ $F_i.$
 By *-transfer$,$ the surface area of such a *-geometric configuration is
$S(F_i)) = \pi (\hyper f(x_i) + \hyper f(x_{i+1}))\Hyper \vert \ell(f)_i\vert,$ where 
$\Hyper \vert \ell(f)_i\vert$ is the *-length of the hyperline segment.
In appendix 6$,$ $\Hyper \vert \ell(f)_i\vert= \Vert \vec v_i \Vert.$
Assuming the usual *-additivity$,$ the surface area of the entire configuration 
$F(Q) =\bigcup\{F_i \bigm \vert 0 \leq i \leq \Gamma -1\}$ is                                         
$$\sum_{I_i \in \hyper {\cal S}(Q)} \!\!\!\!\hyper S(F_i)= \hyper S(F(Q)).\eqno 
(10)$$\par      
From this definition it will follow that
$$\st {\hyper S(F(Q))} = \int_a^b 2\pi f(x)\sqrt {1 + (f^\prime (x))^2}\, 
dx.\eqno(11)$$ \par
Now the complete derivation uses the proof of the integral length formula 
that appears in appendix 6. There it is shown that
$\Vert \vec v_i \Vert = \sqrt {1 + (\hyper f^\prime (x_i))^2}\, dx + 
\delta_i\, dx,$ where $\delta_i \in \monad 0.$ But$,$ $\hyper f(x_i) = \hyper 
f(x_{i+1}) + \lambda_i,\ \lambda_i \in \monad 0.$ Therefore$,$ $\hyper f(x_i) + 
\hyper f(x_{i+1}) = 2\hyper f(x_i) + \lambda_i.$ Since $f$ is bounded then   
$\hyper S(F_i)) = \pi (\hyper f(x_i) + \hyper f(x_{i+1}))\Hyper \vert \ell(f)_i\vert=
2\pi f(x_i) \sqrt {1 + (\hyper f^\prime (x_i))^2}\, dx + \gamma_i\, dx,\                                               
\gamma_i \in \monad 0.$ The derivation is completed by application of steps
(3)(4)(5)(6)(7) as demonstrated in application 7.2.1.\par
\vskip 18pt
%==========TABLE STRUTS

%24 POINT DEEP ROWS
\newbox\medstrutbox
\setbox\medstrutbox=\hbox{\vrule height14.5pt depth9.5pt width0pt}
\def\medstrut{\relax\ifmmode\copy\medstrutbox\else\unhcopy\medstrutbox\fi}

\bigskip

\newdimen\leftmargin
\leftmargin=0.0truein
%for textbook use 1.1truein
\newdimen\widesize
\widesize=3.0truein
\advance\widesize by \leftmargin       
%\moveleft\leftmargin
\hfil\vbox{\tabskip=0pt\offinterlineskip
\def\tablerule{\noalign{\hrule}}

\halign to \widesize{\medstrut\vrule#\tabskip=0pt plus2truein

&\hfil\quad#\quad\hfil&\vrule#

\tabskip=0pt\cr\tablerule

%&\multispan3 \hfil {\tablerm TABLE TITLE} \hfil&\cr\tablerule
%\biggerstrut&ENTRY A&&ENTRY B&\cr\tablerule

&$\Downarrow$ IMPORTANT $\Downarrow$&\cr\tablerule
}} \par
\bigskip
$\Rightarrow$ Please note that for applications 7.2.1$,$ 7.2.2$,$ and 7.2.3 the 
actual defining geometric quantities are given by equations (1)$,$ (8) and (10),
respectively. The geometric configurations$,$ $T(Q)$ and $F(Q)$ composed of the hyperfinite 
union of the respective infinitesimal elements may be considered as members 
of the 
nonstandard extension of the set of all ordinary Euclidean configurations. 
Configurations $T(Q),\ F(Q)$ are internal subsets of appropriate 
*-Euclidean entities that do exhibit a standard area or volume measure. By 
*-transfer of the standard properties of geometric measures relative to subsets 
it follows that the standard part of each of these expressions exists as a 
real number. What is established in these applications under the analytical 
constraints given is that the standard part is expressible by the indicated 
integral. $\Leftarrow$   
\vskip 18pt

Relative to the infinitesimal geometric elements$,$ I have not altered the 
requirement that for the n-dimensional integral these elements be the 
infinitesimal subrectangles. In the latter sections of this manual$,$ since it 
is written for individuals with a strong undergraduate mathematics background$,$ 
the {Jordan-measurable} generalization for the basic rectangular region was 
used. 
As far as a Jordan-measurable $J$ is concerned do we ever 
need to consider any subrectangle that is not in the interior of $J$? \par
 Let
Jordan-measurable $J \subset R\subset \realp n$ and let $Q$ be any fine 
partition of $R.$ It is sometimes useful to assume that $J$ is a closed subset 
of $R$ (hence$,$ compact). For$,$ if $J$ is not closed$,$ then recall that the Jordan-content of $J,\ v(J) 
= v(\overline{J}).$ {\bf Intuitively$,$ think of $J$ as any of the ordinary 
regions studied in elementary calculus with their intuitive boundary$,$ inner 
and exterior portions.} Using these ideas of boundary$,$ inner and exterior 
portions of $J,$ 
the set 
of all subrectangles determined by $Q$$,$ which is denoted by $\hyper {\cal S}(Q),$ 
may be separated formally into three disjoint hyperfinite collections. Let 
the set of {\bf {boundary subrectangles}} be
$\partial (Q) = \{S\bigm \vert S \in \hyper {\cal S}(Q) \land S \cap \hyper J 
\not= \emptyset \land S \cap (\hyperrealp n - \hyper J) \not= \emptyset\}.$ 
Observe that if $S \in \partial (Q),$ then $S \subset \monad p$ and 
$\monad p \cap \hyper J \not= \emptyset $ and 
$ \monad p \cap (\hyperrealp n - \hyper J) \not= \emptyset$ imply that 
$p \in \partial (J).$ Now consider the set of {\bf {exterior subrectangles}} 
${\rm ext}(Q) = \{S\bigm \vert S \in \hyper 
{\cal S}(Q) \land S \subset (\hyperrealp n - \hyper J) \}.$ Finally$,$ the set 
of {\bf {inner subrectangles}} is ${\rm inn}(Q) = 
{\cal S}(Q) - (\partial (Q) 
\cup {\rm ext}(Q)).$ Observe that the boundary subrectangles generate 
boundary points of $J.$ However$,$ by considering $J$ to be a rectangle or  
the interior of a rectangle then it is clear that there may exist some 
$S \in {\rm inn}(Q)$ or $S\in {\rm ext}(Q),$ respectively$,$ such that $S 
\subset 
\monad p$ and $p \in \partial (J).$ It is definitely the case$,$ however$,$ that
$S \subset \hyper J$ if and only if  $S \in {\rm inn}(Q).$ Of course the sets
$\partial (Q),\ {\rm inn}(Q)$ and ${\rm ext}(Q)$ are mutually disjoint.
This leads to the very useful\par 
{\bf Theorem 7.2.2.} {\sl Let Jordan-measurable $\ J \subset R 
\subset \realp n,$ bounded $f\colon J \to \real,$ and $Q$ be any fine 
partition of $\hyper R.$ Let hyperfinite $\hyper {\cal S}(Q) = \{S_i 
\bigm \vert 0\leq i \leq \Gamma -1\}.$  Assume that there exists a 
hyperfinite sequence $U_i,$ where $U_i \in S_i$ for each $S_i \in {\cal S}(Q)$.
Then 
$$\st {\!\!\!\!\!\!\sum_{S(i) \in {\rm inn}(Q)}\!\!\!\!\!\!\hyper f(U_i) 
\hyper v(S_i)}= 
\st {\!\!\!\!\!\!\sum_{S(i) \in {\cal S}(Q)} \!\!\!\!\!\!\hyper 
{\hat f}(U_i) \hyper v(S_i)}.$$}\par
\vskip 18pt
Theorem 7.2.2 indicates$,$ as expected$,$ that only the interior infinitesimal subrectangles are 
significant when modeling with respect to the Jordan-content of any set.
\bigskip
{\bf Application 7.2.4.} {\it Volume obtained by 2-dimensional integral.}\par
Suppose that compact Jordan-measurable $J \subset R \subset \realp n.$ Let 
continuous $f\colon J \to \real$ For each 
$\vec x \in J$ let $f(x) \geq 0.$ Then the volume between the surface 
determined by $f$ and the $xy$-plane is
$$\int_Jf(x)\, dX.$$
{\sl Derivation.} Let $Q$ be a simple fine partition of $R$ and consider some
$S_i \in {\cal S}(Q).$ From the NSP-world viewpoint$,$ the surface curves obtained 
by intersecting the surfaces with planes parallel to the coordinate planes are 
hyperpolygonal lines. Since the choice of the particular hyperline segments 
is arbitrary and $S_i = [x_i,x_{i+1}] \times [y_i,y_{i+1}]$ then consider the
*-Euclidean configuration $TR(S_i)$ composed of the {truncated hyperrectangular 
solid} with $S_i$ as its base and its top a parallelogram with adjacent sides
the hyperline segments $\{(x_i,y_i,\hyper {\hat f}(x_i,y_i),(x_i,y_{i+1},\hyper 
{\hat f}(x_i,y_{i+1})\}$ and 
$\{(x_i,y_i,\hyper {\hat f}(x_i,y_i)),(x_{i+1},y_i,\hyper 
{\hat f}(x_{i+1},y_i)\}$ 
Using the *-Euclidean measure for the volume $\Hyper V(T(S_i)$ then once again 
*-additivity yields for the configuration $TR(Q) = \bigcup \{TR(S_i)\bigm 
\vert 0\leq i\leq \Gamma - 1\}$ 
$$\sum_{S(i) \in {\cal S}(Q)} \Hyper V(TR(S_i)) = \Hyper V(TR(Q)).\eqno(12)$$
Investigating the various configurations $TR(S_i),$ and using the *-Euclidean 
measure of these configurations then Theorem 7.2.2 yields 
that 
$$\sum_{S(i) \in {\rm inn}(Q)}\!\!\!\!\!\!\Hyper V(T(S_i)) \approx 
\!\!\!\!\!\sum_{S(i) \in {\cal S}(Q)} \!\!\!\!\!\!\Hyper V(S_i).\eqno(13)$$
Since $J$ is compact then
$$\hyper f(x_i,y_i)\approx  \hyper f(x_i,y_{i+1}) \approx \hyper 
f(x_{i+1},y_i).\eqno (14)$$ \par
Applying the *-Euclidean measure formula yields that 
$\Hyper V(T(S_i) = (\hyper f(x_i,y_i) + \delta_i)\, dX$  and the method of 
Application 7.2.1 yields 
$$\st {\!\!\!\!\!\sum_{S(i) \in {\rm inn}(Q)}\!\!\!\!\!\!\Hyper V(T(S_i)} = 
\st {\!\!\!\!\!\sum_{S(i) \in {\cal S}(Q)} \!\!\!\!\!\!\Hyper V(S_i)} =
\int_Jf(x)\, dX.\eqno (15)$$ \par
\bigskip
For integral modeling$,$ {physical infinitesimal elements} are based upon the 
Method of Constants. It assumes that the physical quantity being considered 
has the same effect as if it were concentrated at some point within 
subrectangle or another similar object. 
The method is essentially outlined by Maxwell in his previously quoted 
descriptions. Clearly$,$ such infinitesimal physical modeling is  highly 
discipline orientated and is closely associated with simple mind experiments.
\bigskip
{\bf Application 7.2.5.} {\it The value of the electric field vector at a 
point $P$ exterior to the plane of a charged 2-dimensional closed and 
bounded Jordan-measurable region.}\par
Consider a Jordan-measurable region $J,$ a point $P = (a,b,c),\ 
(c \not= 0)$ in 
space exterior to the plane. Let $\rho (x,y,0)$ be a continuous charge density 
function defined on $J.$ The scalar value of the electric field at $P$ is
$$\vert \vec E(a,b,c) \vert = \int_J{{\rho (x,y,0)}\over{(x-a)^2 +(y-
b)^2 + c^2}}\, dX.$$                                      
\vskip 18 pt
{\sl Derivation.} (The method of point charges.) 
For a single point charge $q$ at a distance $r$ from $P,$ 
the definition of the scalar vaule of the electric field is $\vert \vec E\vert
 = q/r^2.$
Let $J \subset R \subset \realp 2,$  $(x_i,y_i,0)\in J,$ and consider a simple fine partition 
$Q$ of $R.$ If one considers a point charge with value $\hyper \rho (x,y,0) 
\hyper v(S_i);\ (x,y,0) \in S_i\in {\rm inn}(Q),$ then 
$\vert \vec E(a,b,c)\vert = (\hyper \rho (x^\prime,y^\prime,0)\, 
\hyper v(S_i))/((x-a)^2 + (y - b)^2 + c^2).$ But the function $h(x,y)= \rho 
(x,y,0)/((x-a)^2 + (y - b)^2 + c^2)$ is uniformly continuous on $J.$ Thus
$$\hyper \rho (x^\prime,y^\prime,0)\hyper v(S_i)=(\hyper \rho (x_i,y_i,0)+ \delta_i) 
\hyper v(S_i)=$$ $$\hyper \rho (x_i,y_i,0) \hyper v(S_i) + \delta_i \hyper 
v(S_i),$$ where $\delta_i \in  \monad 0.$ 
Assuming that the scalar value of the electric field for any nonempty finite 
set of point charges is the sum of the individual values and that within any 
2-dimensional rectangle$,$ $S,$ with charge density $\rho (x,y,0)$ there is some 
point where the entire charge can be considered as concentrated$,$ or that for a 
constant charge density $\rho$ the total charge is $\rho v(S),$
then *-transfer and the method used in applications 7.2.1 and 7.2.4 imply that 
the scalar value of the electric field at $P$ due 
to the charge on $J$ is
$$\vert\vec E_J(P)\vert =\st {\vert \hyper {\vec E_J(P)} \vert} =$$
$$\st {\!\!\!\!\!\!\!
\sum_{S(i) \in {\rm inn}(Q)}\!\!\!\!\!\!\!(\hyper \rho (x^\prime,y^\prime,0)\, 
\hyper v(S_i))/((x-a)^2 + (y - b)^2 + c^2)} =$$
$$ \int_J{{\rho (x,y,0)}\over{(x-a)^2 +(y-
b)^2 + c^2}}\, dX.$$\par
\bigskip
{\bf Remark:} In 7.2.5$,$ it is not correct to simply consider the points 
$(x^\prime,y^\prime,0) \in S_i.$ It is always necessary that the chosen 
intermediate partition be \underbar{internal}. Further$,$ extending the concept 
of the constant charge density is a major approach to infinitesimal modeling
since if $\rho (x,y,0)$ is continuous on compact $J,$ then $\rho$ behaves 
in a constant-like manner in that if infinitesimal $S \subset J$ and 
$p,q \in S,$ then $p \approx q$ implies that $\st {\rho(p)} = \st {\rho(q)}.$                                    
Using the physical element concept$,$ this is called the {\bf {elemental 
method of constants}.}
\vfil\eject
\centerline{Chapter 8.}\medskip
 \centerline{\bf REFINEMENTS FOR}
\centerline{\bf INTEGRAL MODELING}
\bigskip
\leftline{8.1 \underbar{A Very General Approach.}}\par  
{\tt Hurd and Loeb [1985]} construct a very general integral concept.  
For example$,$  consider any hyperfinite 
$\{x_1,\ldots,x_\Gamma \} \subset X,$ let $B$ be any set of internal
hyperreal-valued functions defined on $X,$ and $\{a_1,\ldots,a_\Gamma \}$ a fixed set of
hyperreal nonnegative numbers. They then consider the hyperfinite sum operator
$\sum_\Gamma$ defined on each $f  \in B$ by $\sum_\Gamma f = \sum_1^\Gamma
a_i\, f(x_i).$ This and other examples are generalized and an entire theory 
of integration is developed that 
incorporates various classical generalizations of the Riemann integral$,$ 
especially the Lebesque. They apply their theory to stochastic processes such 
as the Poisson process and Brownian motion. Since the background necessary to 
study their generalization is beyond the scope of these manuals$,$ their theory 
will not be presented. Indeed$,$ except for concept of hyperfinite summation,
  their general approach is probably unsuited for elementary modeling. Our 
goal in this last chapter on integral modeling is to examine more closely the 
specific contents of the hyperfinite sum as defined in definition 5.1.1 
relative to modifications of its geometric or physical meaning.\par
Recall definition 5.1.1. Let $f\colon R \to \real$ 
be bounded and $\cal P$ the set of simple partitions of $R.$ 
Then $f$ is said to be integrable if there exists 
some $r \in \real $ and a simple fine partition$,$ $P \in \hyper 
{\cal P}$ such that for each of its internal intermediate partitions 
$Q=\{\vec v_q\},$ where $1 \leq q \leq \Gamma \in {\nat}_\infty,$ it follows 
that $ \sum_{k=1}^\Gamma \hyper f(\vec v_q)\hyper v(R_q) \in \monad r.$ \par
Obviously$,$ modifications can be made in the concept of the basic partition $P$ 
or the intermediate partition $Q.$ Such modifications are discussed at the 
conclusion of this chapter. It is important to stress at this point something 
that is not apparent about this definition. It is immediate from examination 
of the proofs in appendix 5$,$ that$,$ basically$,$ there are two reasons why this 
definition works. First$,$ $\hyper v(R_q) \in \monad 0$. And$,$ secondly$,$ due to 
behavior of the measure $v,$ the integral value $r$ is independent of the fine 
partition chosen. This must be taken into consideration if the partition 
concepts are not to be altered.  For the basic modifications that follow$,$ 
we extend the partitioning requirement to \underbar{any} arbitrary fine partition and 
any arbitrary intermediate partition and modify the terms of the hyperfinite 
sum with a view towards applications as they appear in elementary calculus 
courses.\par
\bigskip
\leftline{8.2. \underbar{The Line Integral.}}
\medskip  
Let bounded $\phi\colon [a,b] \to \real.$ Our first consideration is to modify 
the value of the *-measure $\hyper v([t_{i-1}i,t_i]),$ where $[t_{i-1},t_i] 
\subset \Hyper [a,b].$ We know the importance of the ``increment'' in the 
calculus$,$ hence$,$ our first modification is to replace $\hyper v([t_{i-
1},t_i])$
with the standard extension of the increment operator. In place of 
$\hyper v([t_{i-1},t_i])$ write $\hyper \phi (t_i) - \hyper \phi (t_{i-1}).$
What happends if $\phi$ is continuous? Well$,$ in that case$,$ if $t_i - t_{i-1} 
\in \monad 0,$ then $\hyper \phi (t_i) - \hyper \phi (t_{i-1})\in \monad 0$ and one of the 
most basic requirements for that factor of the term of the hyperfinite sum is 
met.\par
\vskip 18pt
\hrule
\smallskip
{\bf Definition 8.2.1. ({Riemann - Stieltjes Integral}).} Let bounded
$f\colon [a,b]\to \real$ and bounded $\phi \colon [a,b] \to \real.$ 
Then $f$ is 
{\bf RIEMANN STIELTJES Integrable with respect to $\phi$} if there exists a real
$r$ such that for any fine partition $P= \{a = t_0,\ldots,t_\Gamma=b\}$ of 
$[a,b]$ and any internal intermediate 
partition $Q=\{t_1^\prime,\ldots,t_\Gamma^\prime\}$ of $P$ it follows that     
$$\sum_{i=1}^\Gamma \hyper f(t_i^\prime)\, (\hyper \phi (t_i) - \hyper \phi (t_{i-
1})) \in \monad r.$$\par
\hrule
\vskip 18pt
Rather than investigate the Riemann-Stieltjes integral as defined in 8.2.1,
I pass directly to the line integral of elementary calculus. Referring to 
example 4.4.1A$,$ let bounded $c\colon [a,b] \to \realp n, \ c(t) = 
(c_1(t),\ldots,c_n(t))$ be considered a curve with graph $C.$ Assume that 
$E\subset \realp n,\ C \subset E$ and bounded $F\colon E \to \realp n.$   
Write $F$ as $F(\vec x) = (f_1(\vec x),\ldots, f_n(\vec x)).$ Consider the 
composite function $(Fc)\colon [a,b] \to \realp n$ defined by
$(Fc)(t) = (f_1(c(t)),\ldots,f_n(c(t))),\ t \in [a,b].$ 
Notice that if $t_j^\prime \in [t_{j-1},t_j] \subset \Hyper [a,b],$ then$,$ letting
$\bullet$ denote the ``dot'' (inner) product,
$$\hyper (Fc)(t_j^\prime) \bullet \hyper {\vec v_j} 
= \sum_{i=1}^n \hyper f_i(\hyper c(t_j^\prime))\, (\hyper c_i(t_{j-1}) -
\hyper c_i(t_j)),$$ 
 where $\hyper {\vec v_j} =\hyper c(t_j) - \hyper c(t_{j-1})$ represents the directed 
hyperline segment portion of some hyperpolygonal curve representation for $C.$
Now taking a hyperfinite sum yields
$$\sum_{j=1}^\Gamma \hyper (Fc)(t_j^\prime) \bullet \hyper {\vec v_j} =
\sum_{j=1}^\Gamma \bigl(\sum_{n=1}^n \hyper f_i(\hyper c(t_j^\prime))\, 
(\hyper c_i(t_{j-1}) -
\hyper c_i(t_j))\bigr)=$$
$$\sum_{i=1}^n\bigl(\sum_{j=1}^\Gamma  \hyper f_i(\hyper c(t_j^\prime))\, 
(\hyper c_i(t_{j-1}) -
\hyper c_i(t_j))\bigr).$$\par
Hence$,$ if for each $i=1,\ldots,n,\ \sum_{j=1}^\Gamma  \hyper f_i(\hyper 
c(t_j^\prime))\, (\hyper c_i(t_{j-1}) -
\hyper c_i(t_j)) \in \monad {r_i},$ then $\sum_{j=1}^\Gamma \hyper 
(Fc)(t_j^\prime) \bullet \hyper {\vec v_j} \in \monad {r_1 + \cdots + r_n}.$ On the 
other hand$,$ if $\sum_{j=1}^\Gamma \hyper 
(Fc)(t_j^\prime) \bullet \hyper {\vec v_j} \in \monad r$ and for each $i=1,\ldots,n,\ 
 \sum_{j=1}^\Gamma  \hyper f_i(\hyper c(t_j^\prime))\, (\hyper c_i(t_{j-1}) -
\hyper c_i(t_j)) \in {\cal O},$ then $\st {\sum_{j=1}^\Gamma  \hyper f_i
(\hyper c(t_j^\prime))\, (\hyper c_i(t_{j-1}) -
\hyper c_i(t_j)}= r_i$ implies that $r_1 + \cdots r_n = r.$ This leads to the 
notion of the line integral as an extension of the Riemann-Stieltjes 
integral. \par
\vskip 18pt
\hrule
\smallskip
{\bf Definition 8.2.2. ({Line Integral}).} 
Let bounded $c\colon [a,b] \to \realp n, \ c(t) = 
(c_1(t),\ldots,c_n(t))$ be considered a curve with graph $C.$ Assume that 
bounded $F\colon C \to \realp n,$
where  $F(\vec x) = (f_1(\vec x),\ldots, f_n(\vec x)).$ Consider the 
composite function $(Fc)\colon [a,b] \to \realp n$ defined by
$(Fc)(t) = (f_1(c(t)),\ldots,f_n(c(t))),\ t \in [a,b].$ Then $F$ is {\bf LINE 
INTEGRABLE with respect to C} if there exists some real $r$ such that for any 
fine partition $P= \{a = t_0,\ldots,t_\Gamma=b\}$ of 
$[a,b]$ and any internal intermediate 
partition $Q=\{t_1^\prime,\ldots,t_\Gamma^\prime\}$ of $P$ it follows that     
$$\sum_{j=1}^\Gamma \hyper (Fc)(t_j^\prime) \bullet \hyper {\vec v_j} 
\in \monad r.$$          
In which case we write  
$$r = \st {\sum_{j=1}^\Gamma \hyper (Fc)(t_j^\prime) \bullet \hyper {\vec 
v_j}} =\int_C 
F\bullet
d\vec R.$$\par
\smallskip
\hrule
\vfil
\eject
%==========TABLE STRUTS

%24 POINT DEEP ROWS
\newbox\medstrutbox
\setbox\medstrutbox=\hbox{\vrule height14.5pt depth9.5pt width0pt}
\def\medstrut{\relax\ifmmode\copy\medstrutbox\else\unhcopy\medstrutbox\fi}

\bigskip

\newdimen\leftmargin
\leftmargin=0.0truein
%for textbook use 1.1truein
\newdimen\widesize
\widesize=3.0truein
\advance\widesize by \leftmargin       
%\moveleft\leftmargin
\hfil\vbox{\tabskip=0pt\offinterlineskip
\def\tablerule{\noalign{\hrule}}

\halign to \widesize{\medstrut\vrule#\tabskip=0pt plus2truein

&\hfil\quad#\quad\hfil&\vrule#

\tabskip=0pt\cr\tablerule

%&\multispan3 \hfil {\tablerm TABLE TITLE} \hfil&\cr\tablerule
%\biggerstrut&ENTRY A&&ENTRY B&\cr\tablerule

&$\Downarrow$ IMPORTANT $\Downarrow$&\cr\tablerule
}} \par
\bigskip
$\Rightarrow$ In the next derivation$,$ the significant 
method of the maximum and minimum is applied to physical elements. This approach 
is different than that used in integral rule IR3. Moreover$,$ a special 
relationship between the work done by a force field over a hyperline segment 
and points on a rectifiable curve is advanced. This relationship aids 
in our comprehension of energy related NSP-world properties.$\Leftarrow$\par
\medskip
{\bf Application 8.2.1.} {\sl Energy expended within a force field while 
moving along a curve.}\par
Let $c\colon [a,b] \to \realp n$ be a continuous differentiable curve with 
graph $C.$ Assume that continuous $F\colon E \to \realp n,$ open $E \supset 
C.$ The work done 
in moving 
through the force field on the path $C$ is
$$W(C) = \int_C F\bullet d\vec R = 
\int_a^b\bigl(\sum_{i=1}^nf_i(c_1(t),\ldots,c_n(t))c_i^\prime(t)\bigr)\, dt.$$\par
{\sl Derivation}.  In experimental physics$,$ the concept of ``{work}''  ({energy expended}) 
is introduce. {\bf All 
one needs to do is to establish its properties for a polygonal curve.}  
Suppose we have continuous force field $F\colon E \to 
\realp n.$ Let ${\cal P}_k\subset \realp n$ be any finite polygonal curve$,$ 
$\ell_j$ one of 
the line segment portions of ${\cal P}_k$ with $\vec v_j= (c_1(t_{j}) - 
c_1(t_{j-1}), \ldots,c_n(t_{j}) - c_n(t_{j-1})),$ denoting this line 
segment considered as a directed line segment in the direction of motion 
through the field $F.$ If $F$ is constant on $\ell_j,$ then the work done 
moving along $\ell_j$ is defined as $W(\ell_j) = F \bullet (\vec v_j/\Vert 
\vec v_j \Vert)\Vert \vec v_j \Vert,$ where length of $\ell_j= \Vert \vec v_j 
\Vert.$
 What if the force field is not constant? Consider ${\cal P}_k$ as 
represented by a continuous $\ell\colon [a,b] \to \realp n$ and assume that 
$F$ is defined on $\ell.$ Then for a given 
$\ell_j= \{(x_1(t),\ldots,x_n(t))\bigm \vert t \in [t_{j-1},t_j]\}$ there exists 
some $t_m,\ t_M$ such that $W_m(\ell_j) = F(\ell(t_m))\bullet \vec 
v_j\leq W(\ell_j)=F(\ell(t))\bullet\vec v_j= W_M(\ell_j)=F(\ell(t_M))\bullet 
\vec v_j$ for each $t \in [t_{j-1},t_j].$ 
Let's make the one assumption that the actual amount of work expended moving 
along the line segment $\ell_j$ is $W(\ell_j)$ and that 
$W_m(\ell_j) \leq W(\ell_j) \leq W_M(\ell_j).$
Then from continuity there exists some $h_j^\prime \in [t_{j-1},t_j]$ such that 
$W(\ell_j) = F(\ell(h_j^\prime))\bullet \vec v_j.$ 
The idea of the nonconstant force field over a line segment is embedded into the 
the NSP-world by *-transfer assuming that what has been established above 
holds for all such polygonal curves. Hence$,$ let ${\cal P}_\Omega$ be a 
hyperpolygonal representation for the curve generated by a fine partition$,$ 
$\ell_j$ an hyperline segment in ${\cal P}_\Omega.$ Since $F$ is continuous on 
$E$ then $\hyper F$ is defined on ${\cal P}_\Omega.$ It follows 
that $\Hyper W(\ell_j) = \hyper F(\ell_j (h_j^\prime))\bullet \hyper {\vec 
v_j}.$ 
For polygonal curves$,$ in general$,$ the work done is an additive function. 
Thus for the hyperpolygonal curve ${\cal P}_\Omega$ 
$$\Hyper W({\cal P}_\Omega) = \sum_{j=1}^\Omega \hyper F(\ell_j (h_j^\prime))
\bullet \hyper {\vec v_j}.\eqno (1)$$ \par
We now show that there exists a real number $W(C)$ such that for any fine 
partition ${\cal P}_\Gamma$ and \underbar{any} intermediate partition $Q = 
\{t_1^\prime,\ldots,t_\Gamma^\prime\},\ \st {\Hyper W({\cal P}_\Gamma)} = 
W(C).$\par 
Consider
$$\hyper F(\hyper c(t_j^\prime))\bullet \hyper {\vec 
v_j}= \sum_{i=1}^n \hyper f_i(\hyper 
c_1(t_j^\prime),\ldots,\hyper c_n(t_j^\prime))\,(\hyper c_i(t_j) - 
\hyper c_i(t_{j-1})).\eqno (2)$$
The curve $c$ being continuously differentiable on $[a,b]$ implies that 
$\hyper c_i(t_j) - \hyper c_i(t_{j-1}) = \big(\hyper 
c_i^\prime(t_j)+\delta_{ij}\big)\,(t_j -t_{j-1}),$
where $\delta_{ij} \in \monad 0.$
Hence$,$ $\hyper F(\hyper c(t_j^\prime))\bullet \hyper {\vec 
v_j}=$
$$\sum_{i=1}^n \hyper f_i(\hyper 
c_1(t_j^\prime),\ldots,\hyper c_n(t_j^\prime))\,\hyper c_i^\prime(t_j)\, (t_j -t_{j-1}) + 
\left(\sum_{i=1}^n \delta_{ij}\right)\, (t_j -t_{j-1}).\eqno (3)$$
 Consequently$,$ $\hyper F(\hyper c(t_j^\prime))\bullet \hyper {\vec 
v_j}=$
$$ \sum_{i=1}^n \hyper f_i(\hyper 
c_1(t_j^\prime),\ldots,\hyper c_n(t_j^\prime))\,\hyper c_i^\prime(t_j)\, (t_j -t_{j-1}) + 
\delta_j\,(t_j - t_{j-1}).\eqno (4)$$                                                     
Uniform continuity of the $f_i$ yields that $\hyper f_i(\hyper 
c_1(t_j^\prime),\ldots,\hyper c_n(t_j^\prime)) = \hyper f_i(\hyper 
c_1(t_j),\ldots,\hyper c_n(t_j)) +\lambda_{ij},$ where $\lambda_{ij} \in 
\monad 0.$ Once again  this yields $\hyper F(\hyper c(t_j^\prime))\bullet \hyper {\vec 
v_j}=$
$$ \sum_{i=1}^n \hyper f_i(\hyper 
c_1(t_j),\ldots,\hyper c_n(t_j))\,\hyper c_i^\prime(t_j)\, (t_j -t_{j-1}) + $$ 
$$ \lambda_j\,(t_j - t_{j-1})+\delta_j\,(t_j - t_{j-1})\eqno (5)$$\par  
Now continuing the basic {\bf elemental derivation process} leads to
$\sum_{j=1}^\Gamma\bigl(\hyper F(\hyper c(t_j^\prime))\bullet \hyper {\vec 
v_j}\bigr) =$
$$\sum_{j=1}^\Gamma\bigl(\sum_{i=1}^n \hyper f_i(\hyper 
c_1(t_j),\ldots,\hyper c_n(t_j))\,\hyper c_i^\prime(t_j)\, (t_j -t_{j-1})\bigr) + 
\delta,\ \delta \in \monad 0.\eqno (6)$$
Therefore$,$ from Theorem 5.1.2
$$W(C)= \st {\sum_{j=1}^\Gamma\bigl(\hyper F(\hyper c(t_j^\prime))\bullet 
\hyper {\vec v_j}\bigr)} = $$ $$
\int_a^b\bigl(\sum_{i=1}^nf_i(c_1(t),\ldots,c_n(t))c_i^\prime(t)\bigr)\, 
dt.\eqno (7)$$
The fact that this is a line integral follows from Definition 8.2.1. Finally,
it will almost always be the case that the special energy property $\hyper 
F(\ell_j(t_j^\prime)) \bullet \hyper {\vec v_j} =  \hyper F(\hyper 
c(t_j^\prime))
\bullet \hyper {\vec v_j} + \eps_j \Vert \hyper {\vec v_j} \Vert,\ \eps_j \in 
\monad 0,\ t_j^\prime= h_j^\prime$ holds. [See note [2] on page 148.]
Since $C$ is rectifiable then 
$\sum_{j=1}^\Gamma \hyper F(\ell_j(t_j^\prime)) \bullet \hyper {\vec v_j} 
\approx \sum_{j=1}^\Gamma \hyper F(\hyper c(t_j^\prime))\bullet \hyper {\vec 
v_j}.$\  We are 
using the hyperpolygonal representations for $C$ as the basic NSP-world entity 
to determine the N-world physical effects. It is clear that the 
appropriate measure for the work expended moving along the path $C$ should be 
the unique value $\st {\sum_{j=1}^\Gamma \hyper F(\ell_j(t_j^\prime)) 
\bullet \hyper {\vec v_j}} $ obtained in (7). This completes the derivation.\par
\bigskip
{\bf Obviously$,$ derivation 8.2.1 also establishes the elementary method for the 
calculation of a line integral.}\par
\bigskip
\leftline{8.3. \underbar{Order Ideals and Approximations.}}
\medskip  
One often reads in the literature that such and such an expression is a ``first- 
order approximation'' or some such phrase. These vague approximation concepts 
can be discussed from the infinitesimal viewpoint and$,$ indeed$,$ lead to 
the notions of the ``{microconstruction}'' and ``{microeffects}.''\par
\vskip 18pt
{\bf Theorem 8.3.1.} {\sl For each $\eps \in \monad 0$ the set
$o(\eps) = \{\eps\, h \bigm \vert h \in \monad 0 \}$ is an ideal in $\monad 
0.$}\par
\vskip 18pt
\hrule
\smallskip
{\bf Definition 8.3.1. ({Order Ideals}).} For a given $\eps \in \monad 0$ 
the set $o(\eps)$ is called an {\bf ORDER IDEAL (of infinitesimals)}.\par
\smallskip
\hrule
\vfil
\eject 
We briefly investigate some of the basic properties of the order 
ideals. (I note that these order ideal properties appear for the first time 
in this manual.) First$,$ it is obvious that $o(\eps) = o(-\eps).$ In all that 
follows$,$ let ${\monad 0}^+ = \{x\bigm\vert x \in \monad 0 \land x \geq 0\}$ be 
the set of all {\bf {nonnegative infinitesimals}}. The next 
theorem seems to be one of the more significant ones relative to order 
ideals.\par
\vskip 18pt
{\bf Theorem 8.3.2.} {\sl Let $\eps \in {\monad 0}^+.$ Suppose that
$w \in \hyperreal$ and $0\leq w\leq \eps\, h \in o(\eps).$ Then 
$w \in o(\eps).$}\par
\vskip 18pt
{\bf Theorem 8.3.3.} {\sl Let $\eps,\, \delta \in {\monad 0}^+.$ If 
$0 \leq \delta \leq \eps,$ then $o(\delta) \subset o(\eps)$ and $o(\delta)$ is a 
ideal in $o(\eps).$}\par
\vskip 18pt
Are there order ideals such that $o(\delta) \subset o(\eps)$ and $o(\delta) \not= 
o(\eps)$?\par
(1) Let $0 < \delta \leq \eps,\ \delta,\eps \in {\monad 0}^+.$ Then 
$o(\delta\eps) \subset o(\eps)$ and $o(\delta\eps) \not= o(\eps).$\par
(2) Let $n \in {\nat}^+ = \nat - \{0\},\ \eps \in {\monad 0}^+.$ Then 
$o(\eps^n) \subset o(\eps^{n-1}) \subset \cdots \subset o(\eps)$ and
$o(\eps^i) \not= o(\eps^j);\ 1 \leq i,j \leq n;\ i \not= j.$\par

Suppose that you have the set $\{\eps_1,\ldots,\eps_n\}.$ Then consider the 
chain $C_1:\ 
o(\eps_1\times\cdots\times\eps_n) \subset o(\eps_1\times\cdots\times\eps_{n-
1}) \subset \cdots \subset o(\eps_1).$ The chain $C_1$ is just one possible 
chain of order ideals leading from $o(\eps_1\times\cdots\times\eps_n)$ to 
$o(\eps_1).$ Such chains are used for comparison purposes and$,$ in this case,
the ideal $o(\eps_1\times\cdots\times\eps_n)$ is called an {\bf {n'th order 
ideal}} where$,$ in general$,$ such ideals as $o(\eps)$ are called {\bf {first-order
ideals.}} Part of this chapter will deal with the relation between n'th order 
ideals and the infinitesimal concept of n'th order approximations. 
It is useful to consider 
other operational methods that might generate different order ideals rather than 
simply restricting their generation to products of infinitesimals.\par
Let $\{\eps_1,\ldots,\eps_k\} \subset \monad 0.$ Define
$o(\eps_1,\ldots,\eps_k) = \{\eps_1\, h_1 + \cdots + \eps_k\, h_k \bigm \vert
h_i \in \monad 0 \land 1\leq i\leq k\}.$\par
\vskip 18pt
{\bf Theorem 8.3.4.} {\sl Let $\eps =\max{\{\vert \eps_1\vert,\ldots,\vert 
\eps_k \vert\}}.$ Then $o(\eps_1,\ldots,\eps_k) = o(\eps).$}\par
\vskip 18pt
{\bf Theorem 8.3.5.} {\sl Let $\eps =\max{\{\vert \eps_1\vert,\ldots,\vert 
\eps_k \vert\}}.$ Then
 $$o(\eps_1,\ldots,\eps_k) = o(\sqrt {\eps_1^2 + \cdots
+\eps_k^2}) = o(\eps).$$}\par
\vskip 18pt 
Thus neither the difference of infinitesimals nor the *-Euclidean norm is a 
useful process for the generation of higher order ideals. As will be 
illustrated through out the remainder of this manual$,$ the n'th order ideals are 
related to the notion of the n'th order approximation. Refer back to all of 
our derivations where the basic elemental derivation process is used. In each 
case$,$ a single term of the required hyperfinite sum of elemental measures$,$ say
$\hyper M(\cdot),$  
is investigated for an 
arbitrary simple partition of infinitesimal volume $dX.$ The derivation 
shows that this term is equal to the value of a standard extension$,$ say 
$\hyper F(\vec x_i),$ plus $h\, dX,$ where $h \in \monad 0.$  Consequently$,$ 
$$\hyper M(\cdot) - \hyper F(\vec x_i) \in o(dX).$$  This expression also 
yields an equivalence relation which is often denoted by the algebraic 
notation
$$\hyper M(\cdot) \approx \hyper F(\vec x_i) \pmod {o(dX)}$$ and
is expressed by stating  that  $\hyper M(\cdot)$ and 
$\hyper F(\vec x_i)$ are {\bf infinitely close of order $dX$.} Hence$,$ to apply 
this derivation process the values $\hyper M(\cdot)$ are not just infinitely 
close 
to  $\hyper F(\vec x_i)$$,$ in a general sense; but$,$ they are infinitely close
of order $dX.$ Noting that for an n-dimensional integral $dX = dx_1 \times 
\cdots \times x_n$ it follows that$,$ from a comparative viewpoint$,$ $\hyper 
M(\cdot)$ and $\hyper F(\vec x_i)$
 can be considered as {\bf {infinitely close of order n.}}
 In the literature$,$ you will also find the less 
descriptive expression ``$\hyper M(\cdot)$ is infinitely close to 
$\hyper F(\vec x_i)$ compared to $dX$'' as a synonym for ``infinitely close 
of order $dX$.''
The reason it appears necessary that such quantities need to
 be infinitely close of order $dX$ is that the steps in the derivation process 
that proceed from this step require a hyperfinite sum to be extracted and the 
results must remain infinitely close.\par 
One little observation about the  order ideals that are created by products.
If $0 \not= r \in \Hyper [-1,1]-\monad 0,$ then for any $\eps \in \monad 0,$ it follows that 
$o(\eps\, r) = o(\eps).$ To see this$,$ Theorem 8.3.3 yields that
$o(\vert \eps\, r \vert) = o(\eps\, r) \subset o(\vert \eps \vert) = o(\eps).$
Let $\eps\, h \in o(\eps).$ Since $h/r \in \monad 0,$ then $(\eps\, r)(h/r) 
\in o(r\eps)$ implies that $\eps\, h \in o(\eps\, r).$ Further$,$ note that
if $0\not= \delta \in \monad 0,$ then $\eps\, r \in o(\eps) - o(\delta\, 
\eps.)$\par    
\bigskip
\leftline{8.4. \underbar{nth Order Increments.}}\par
\medskip
In the next section$,$ we investigate exactly what one means by a tangent line to a 
curve $c$ at a point $p,$ where $c$ is differentiable at $p.$ First$,$ however$,$ it 
is useful for this and the future sections on modeling by means of the 
derivative or differential to formally consider the NSP-world view of the {\it n}th 
order increment ({nth difference}). For $n \in {\nat}^+$ and bounded 
$f\colon [a,nb] \to \real,$ recall 
that the {\bf {nth order increment}}$,$ $\Delta^nf(x,b),$ is defined by 
induction$,$ where for $x \in [a,nb],$ by $\Delta f(x,b) = f(x + b) - f(x).$ This 
leads to the general expression
$$\Delta^nf(x,b) = \sum_{k=0}^n(-1)^k{n \choose k}f(x + (n-k)b) = \sum_{k=0}^n
(-1)^k{n \choose k}f(x + kb).$$ 
Notice that as an operator the *-nth order increment$,$ $\hyper 
(\Delta^n)\hyper f(x,b) =
\Delta^n\hyper f(x,b).$\par
\vskip 18pt
{\bf Theorem 8.4.1.} {\sl Let $1 \leq n \in \nat.$ Suppose that $f^{(n-
1)}\colon [a, nb] \to \real$ and that $f^{(n)}\colon (a,nb) \to \real,$ where
$f^{(k)}$ denotes the {\it k}th derivative of $f.$ Then there exists some $t \in
(a,nb)$ such that $\Delta^n f(a,b) = f^{(n)}(t)\, b^n.$}\par
{\bf Corollary 8.4.1.1} {\sl  Let $1 \leq n \in \nat.$ Suppose that $f^{(n-
1)}\colon [a, b] \to \real$ and that $f^{(n)}\colon (a,b) \to \real,$ then for 
each $dx\in {\monad 0}^+$ and $c \in \Hyper [a,b),$ there exists some $t  \in 
(c,c +ndx)$ such that $\Delta^n\hyper f(c,c+dx) = \hyper f^{(n)}(t)\,(dx)^n.$}\par   
\vskip 18pt
 Theorem 8.4.1 holds if the hypotheses are appropriately altered 
to $f^{(n-1)}\colon [a-na,b] \to \real$ and that $f^{(n)}\colon (a-na,b) \to 
\real.$ In this case$,$ Corollary 8.4.1.1 may be altered to $dx \in \monad 0,\ 
dx < 0,\ c \in \hyper (a,b],\ t \in (c+ndx,c)$ and $\Delta^n\hyper f(c+dx,c) 
= \hyper f^{(n)}(t)\,(dx)^n.$  \par
\vskip 18pt
{\bf Theorem 8.4.2.} {\sl Let $1 \leq n \in \nat.$ Suppose that $f^{(n-
1)}\colon [a, b] \to \real$ and that $f^{(n)}\colon (a,b) \to \real.$ 
If $c \in (a,b),$ then for each $dx\in \monad 0,\ dx \geq 0 $ {\rm [}resp. $dx 
<0${\rm ]}
$$f^n(c)\, (dx)^n \approx \Delta^n \hyper f(c,c+dx), {\rm [resp.}f(c+dx,c){\rm 
]} \pmod {o((dx)^n)}.$$}\par
\vskip 18pt
One important aspect of Theorem 8.4.2 is that the quantities $f^n(c)\, (dx)^n 
$ and $\Delta^n \hyper f(c,c+dx)$ are not simply infinitely close; but$,$ 
rather$,$ are {infinitely closed of order $n.$}  For an infinitesimal$,$ 
$dx,$ Robinson and those that founded the infinitesimal 
calculus consider $\Delta^n \hyper f(c,c+dx)$ to be the {\it n}th order 
differential of $f$ at $c.$ However$,$ most other authors still retain 
the notion that $d^nf(c) = f^{(n)}(c)\,(dx)^n = f^{(n)}(c)\,dx^n;$ 
which I shall retain as well. Under the hypotheses of Theorem 8.4.2 it  
follows for nonzero positive [resp. negative] infinitesimal$,$ $dx,$ that 
$d^nf(c)/dx^n \approx \Delta^n \hyper f(c,c+dx), {\rm [resp.}f(c+dx,c){\rm 
]} = \hyper f^{(n)}(t)$ and$,$ hence$,$ $\st {d^nf(c)/dx^n} =
\st {\Delta^n \hyper f(c,c+dx), {\rm [resp.}f(c+dx,c){\rm ]}} =  f^{(n)}(c).$ 
\par
Finally$,$ all that has been said about {\it n}th order ideals is extended to 
m-dimensional objects of the form $ o^n(\eps_1,\ldots,\eps_m) 
= o(\eps_1) \times \cdots \times o(\eps_m).$
\bigskip
\leftline{8.5. \underbar{Microgeometry - Tangents to Curves.}}\par
\medskip  
A {tangent to a curve} $c\colon [a,b] \to \realp p$ is usually defined 
as the intuitive limit of a set of secants. How might this be viewed within 
the NSP-world? In order to analytically answer this question$,$ let $c$ be 
differentiable at $t \in (a,b).$ Next$,$ let any $r \in \Hyper [-1,1].$ Fix 
$\eps \in {\monad 0}^+$ and let $f_i$ be the coordinate functions of $c.$ The 
Fundamental Theorem  of differential calculus states$,$ in infinitesimal form,
that there exists some $h_i \in \monad 0$ such that
$$\hyper f_i(t + r\, \eps) =f_i(t) +f_i^\prime (t)(r\, \eps) + h_i(r\, 
\eps).\eqno (1)$$ 
For any $r_1 \in \Hyper [-1,1] - \monad 0$$,$ equation (1) then yields$,$ since 
$o(r_1\, \eps) = o(\eps)$ 
$$\hyper c(t + r_1\, \eps) = c(t) + c^\prime(t)(r_1\, \eps) + (h_1,\ldots,h_p)(r\, 
\eps),\eqno (2)$$
$$\hyper c(t + r_1\, \eps) \approx c(t) + c^\prime(t)(r_1\, \eps) \pmod 
{o^n(\eps)}.\eqno (3)$$\par
Each component of $o^n(\eps)$ is a first-order ideal and what follows next is 
a direct result of this fact and the concept of the {\bf resolving power of a 
microscope.} Expression (1) is first transformed into an external 
relation determined by
$$g_i(t + r\, \eps) =f_i(t) +f_i^\prime (t)(r\, \eps),\eqno (4)$$ 
$$\vec g(t + r\, \eps) =c(t) + c^\prime(t)(r\, \eps).\eqno (5)$$\par
Assume that we are ``looking at'' the geometric NSP-world situation with an 
infinite powered microscope with ``{first-order $\eps$-resolving power.}'' 
Physically$,$ this corresponds to the idea that if the distance between two 
objects is a member of a first-order ideal $o(\eps)$$,$ then the objects 
cannot be resolved (i.e. cannot be distinguished one from the other.)  A simple 
proof shows that if $r\in \Hyper [-1,1],$ then the best we can say is 
that
$$\Vert \hyper c(t + r\,\eps) - \vec g(t + r\,\eps)\Vert \in o(\eps).\eqno (6)$$ 
Equation (6) does not mean that for selected $r\,\eps$ the value 
$\Vert c(t + r\,\eps) - \vec g(t + r\, \eps)\Vert$ may not be in an {\it n}th order 
ideal. Indeed$,$ if $r = \eps^{n-1},$ then this would be the case. 
Thus (6) means that$,$ in general$,$ $o(\eps)$ is the ``smallest'' order ideal
that can be guaranteed to contain this value in all possible cases.  
For an ``{$\eps$-infinitesimal microscope}'' (abbr: $\eps$-IM) with first-order 
$\eps$-resolving power$,$ (6) implies that within the field of view 
$$\hyper c(t + r\,\eps) = \vec g(t + r\,\eps)=c(t) + c^\prime(t)(r\, \eps).
 \eqno (7)$$  \par                                                      
A microscope is suppose to magnify$,$ however. How is this feat accomplished. 
The point$,$ $c(t),$ is being considered as the center of view of the 
$\eps$-IM and this point is translated to the origin. This gives us the 
expression
$\hyper c_1(t + r\,\eps) = \hyper c(t + r\,\eps) - c(t) = c^\prime(t)(r\,\eps).$                                                        
The magnification and resolving power of ordinary microscopes are related. 
Hence$,$ for consistency$,$ consider for any $x_j \in \hyperreal$ and $\delta_j \in 
\monad 0$ the general {\bf {infinite magnification 
operator}} $m(x_1\,\delta_1,\ldots,x_p\delta_p) = (x_1,\ldots,x_p),$ which is 
also assumed to be \underbar{linear}. If you wish 
to specify a specific shaped field of view$,$ say a p-dimensional closed 
sphere$,$ you can also restrict this magnification to $\sum_{j=1}^px_j^2 \leq 
2.$ Obviously$,$ both the translation and magnification can be combined into 
one operator. Letting$,$ as before$,$ $r \in\Hyper [-1,1]$ the final view in 
the $\eps$-IM only shows the hyperline segment
$\{m(f_i^\prime(t)r\, \eps,\ldots,f_p^\prime(t)r\, \eps) \vert r\in 
\Hyper [-1,1]\} = \{(f_1^\prime(t)r,\ldots, f_p^\prime(t))\vert r\in 
\Hyper [-1,1]\} = \{c^\prime(t) r\vert r \in  \Hyper [-1,1]\}.$\par
\bigskip
\leftline{8.6 \underbar{Microgeometry - Surface Elements.}}\par   
\medskip
The well-known difficulties of determining a single Euclidean configuration as 
an appropriate approximation for the surface of a 3-dimensional object will 
not be discussed in this manual. [See {\tt Cesari [1956]}] Instead$,$ I pass 
directly to analytical considerations with the appropriate constraints. 
Let open $G \subset \realp 2.$ Assume that $\vec r\colon G \to \realp 3$ and 
that for $(u,v) \in G$ the continuous partial derivatives 
$\vec r_u(u,v),\ \vec r_v(u,v)$ 
exist. Let $\vec r_u(u_0,v_0) = \vec a,\ \vec r_v(u_0,v_0) =\vec b$ and 
$\vert \vec a \times \vec b\vert \not= 0.$ In order to be consistent with 
requirement IR2$,$ define$,$ in matrix notation on column vectors$,$ the linear 
transformation $L\colon \realp 3 \to \realp 3$ as follows:
$$L= \pmatrix{r_1(u_0,v_0)&a_1&b_1\cr
              r_2(u_0,v_0)&a_2&b_2\cr
              r_3(u_0,v_0)&a_3&b_3\cr}.\eqno (1)$$ \par
In the previous section$,$ a hyperline segment portion of the tangent line was 
viewed within an IM. The same procedures are now applied to the tangent plane 
to the surface generated by $\vec r.$ Let $(s,t) \in \Hyper [-1,1] \times
\Hyper [-1,1]= I;\ \delta,\eps \in {\monad 0}^+.$ Then the set
$R_S = \{(u,v)\bigm\vert (u,v) = (u_0 + s\,\delta,v_0 + t\,\eps)
\land (s,t) \in I\}$ is an 
internal infinitesimal rectangle and $R_S \subset \monad {(u_0,v_0)} \in \Hyper 
G.$ The linear transformation $\hyper L$ transforms the internal set of vectors
$\{(1,u-u_0,v-v_0)^T\bigm\vert (u,v) \in R_S\}= D$ onto a configuration 
$\hyper L(D)$ which is a *-Euclidean hyperparallelogram containing 
$\vec r(u_0,v_0).$ When considered restricted to $D,$ this linear 
transformation can be expressed by 
$$\vec k(u,v) = \vec r(u_0,v_0) + (u - u_0)\vec a + (v - v_0)\vec b.\eqno (2)$$  Using 
the *-Euclidean concepts$,$ this hyperparallelogram is a NSP-world portion of the 
standard tangent plane to the surface at $\vec r(u_0,v_0).$ Two adjacent 
sides of this hyperparallelogram have end points 
$\vec k(A) = \vec r(u_0,v_0) + (-\delta) \vec a + (-\eps)\vec 
b;\ \vec k(B) =\vec r(u_0,v_0) + (-\delta)\vec a + (\eps)\vec b;\ 
\vec k(C) = \vec r(u_0,v_0) + (\delta)\vec a + (-\eps)\vec b.$ This 
yields two hyperline segments of hyperlength $\vert \vec {AB}\vert =
\Vert 2\delta\, \vec a\Vert$ and 
$\vert \vec {AC} \vert = \Vert 2\eps\, \vec b\Vert,$ with the hyperarea being $\Vert (4\delta\, \eps)
\vec a \times \vec b \Vert = 4\delta\, \eps \Vert \vec a \times \vec b \Vert=
\Vert \vec a \times \vec b \Vert\, dX.$ \par
Now to analyze the order ideal relationship between members of this 
hyperparallelogram and the surface itself$,$ assume that $\vec r_u$ and $\vec 
r_v$ are continuous at $(u_0,v_0).$ The Fundamental Theorem of 
Differential Calculus in infinitesimal form$,$ implies that for any 
$s\,\delta$ and any $t\,\eps,\ (s,t)\in I$ there exists  
$ \vec \eta \in \mu^3(\vec 0)$ such that 
$$\Hyper {\vec r}(u_0 + s\,\delta,v_0 + 
t\,\eps) = 
\vec r(u_0,v_0) + s\,\delta\,\vec a + t\,\eps\,\vec b + \Vert(s\,\delta,
t\,\eps)\Vert\,\vec\eta.\eqno (3)$$
Considering any $(s,t) \in \Hyper [-1,1] \times \Hyper [-1,1]$ a 
simple proof yields that 
$$\Vert \Hyper {\vec r}(u_0 + s\,\delta,v_0 + t\,\eps) - 
\Hyper {\vec k}(u_0 + s\,\delta,v_0 + t\,\eps)\Vert \in o(\max 
\{\delta,\eps\}) = o(\lambda).\eqno (4)$$\par
The magnification operator for our $(\delta,\eps)$-IM is the mapping
$m(x_1\,\delta+ y_1\,\eps,\ldots,x_3\,\delta + y_3\,\eps) =
m(x_1+ y_1,\ldots,x_3 + y_3),$ where the $x$'s and $y$'s are hyperreal numbers.
Translating and magnifying $\hyper k(u_0 + s\,\delta,v_0 +t\,\eps)$ 
yields as the final view in the $(\delta,\eps)$-IM the hyperparallelogram 
$\{s\vec a + t\vec b\bigm\vert (s,t) \in I.\}$ {\bf BUT,} these results are 
more significant than a simple exercise in analyzing the IM view of 
the tangent plane.\par
Thus far we have decided upon three types of {geometric elements}.\par
(A) Rectifiable curves $\Leftrightarrow$ {\bf {hyperline segments}} and {\bf 
 {hyperpolygonal curves}.}\par
(B) For 1 - 3 dimensional geometric measures by means of the 
1-dimensional integral $\Leftrightarrow$ {\bf {hypertrapezoids}.}\par
(C) For n-dimensional integrals $\Leftrightarrow$ {\bf 
 {infinitesimal rectangles}.}\par
These elements are intuitive in character and only by means of a restrictive 
analytical description are  physical quantities relative to them$,$ such as mass$,$ 
infinitesimal energy and the other applications given in this manual$,$  
actually calculable by means of the integral. \underbar{However}$,$ as 
evident from all that has preceded$,$ {\bf most individuals consider the integral as 
but a hyperfinite sum of entities that are intuitively defined} and do not$,$ generally$,$ 
concerned themselves with the difficulties in calculation. I have been 
slightly  
restrictive in some of the basic definitions by requiring  that functions that 
generate hyperfinite sums$,$ at least$,$ be bounded. Of course$,$ in certain cases 
this restriction might be relaxed. Indeed$,$ in the older literature$,$ geometers 
utilized intuitive infinitesimal geometry and these notions were  not expressed 
originally in terms of any such analytical constraints.\par 
If ${\cal S}$ denotes the geometric point-set called a surface$,$ then $\hyper 
{\cal S}$ is the 
hypersurface. If$,$ intuitively$,$ $T \subset{\cal S}$ is the set of 
surface points 
at which tangent planes ${\cal T}$ to ${\cal S}$ exist$,$ then $\Hyper T$ is the set of 
points in $\hyper {\cal S}$ at which the {hypertangent planes} $\Hyper {\cal T}$ exist.    
For the surface integral
the geometric element -  {\bf the {surface element$,$ $\sigma$}} - is an infinitesimal 
parallelogram containing a point from the {hypersurface}.  This surface 
element is considered to be contained in a hypertangent plane. As such this 
element has an infinitesimal area $d\sigma.$ {\bf From the viewpoint of 
infinitesimal modeling this description of the {geometric surface element} is adequate.}
From the view point of surface integral calculation$,$ since there are 
infinitely many geometric surface elements of different area$,$ such a 
description is not sufficient. \par  
The only question that remains is which collection of surface elements should 
be required for an \underbar{analytical} definition? For this elementary manual$,$ 
the *-Euclidean area notion will be maintained along with a fixed set of 
vectors normal to the hypertangent planes. Let $E \subset \realp 2$ and 
bounded $\vec r\colon E \to \realp 3.$ Suppose that ${\cal S} = \{\vec 
r(u,v)\bigm \vert (u,v) \in E\}$ and that for nonempty $T^\prime \subset E$ 
the set of points $\vec r\,[T^\prime]=T$ is called a set of {\bf {tangent points 
to the surface} ${\cal S}.$} Further$,$ there exists a mapping 
$\vec{\nu} \colon T \to (\realp 3-\{{\vec 0}\})$$,$ where each $\nu(\vec 
t\,\,)$ 
is called a 
{\bf {normal vector}} to the surface ${\cal S}.$ Then each $\vec t\in T$ 
defines a 
unique {\bf {tangent plane}} ${\cal T}(\vec t,\vec{\nu}\,(\vec t\,\,)) 
= \{(x_1,\ldots,x_3)\bigm\vert
(x_1,\ldots,x_3) \in \realp 3 \land \vec{\nu}\,(\vec t\,\,)
 \bullet ((x_1,\ldots,x_3) -\vec t\,\,)=0\} \subset \realp 3.$ Assume that 
$E \subset R \subset \realp 2.$ If $P^\prime$ is any partition of $R,$ 
then there 
exists a nonempty finite set of subrectangles $R_i^\prime \subset R$ 
such that $T^\prime \cap R_i^\prime \not= \emptyset,\ 1\leq i\leq k$ 
and for each such $R_i^\prime$ there exist the intermediate 
partitions $Q^\prime = \{t_1^\prime,\ldots,t_k^\prime\}$ such that 
$t_i^\prime \in R_i^\prime\cap T^\prime.$ Hence$,$ for any fine partition 
$P$ of $\hyper R$ there exists 
a hyperfinite set of infinitesimal subrectangles $\tau(P) = \{S\bigm \vert
S\cap \Hyper T^\prime \not= \emptyset \land S \in P\}$ and a corresponding set 
of internal intermediate partitions $\eta(P).$ These ideas and 
notations are used in the next definition.  
\vskip 18pt\vfil\eject
\hrule
\smallskip
{\bf Definition 8.6.1. ({Surface Integral}.)} Let the surface with its tangent 
planes be defined as in the above paragraph and assume that bounded
$F\colon {\cal S} \to \real.$ Let $E \subset R \subset \realp 2.$ Then $F$ is 
said to be {\bf SURFACE INTEGRABLE with respect to ${\cal S}$} if there exists 
some simple fine partition $P$ of $R$ and some $r \in \real$ such that for 
each intermediate partition $Q \in \eta(P)$
$$\sum_{S \in \tau(P),\ t^\prime \in S\cap Q}\!\!\!\!\!\!\!\!\!\!\!\hyper F( 
\Hyper {\vec r}\,(t^\prime))\Vert \hyper {\vec {\nu}}\,(\Hyper {\vec r}\,
(t^\prime))\Vert \hyper v(S)\in \monad r.$$ \par
\smallskip
\hrule 
\bigskip 
[{\bf Remark.} In the definition of the line integral$,$ it is required that all 
fine partitions be considered. This was done so that consideration could be 
given to rectifiable curves that need not be continuously differentiable. It 
is obvious that Definition 8.6.1. is styled solely for the integral as defined 
by 5.1.1.]\par
Since it is clear from Theorem 7.2.2 and Definition 5.1.1 that for 
$F$ to be integrable with respect to $\cal S$ the values $\hyper F( 
\Hyper {\vec r}\,(t^\prime))\Vert \hyper {\vec {\nu}}\,(\Hyper {\vec r}\,
(t^\prime))\Vert$ must be obtained from an integrable function defined on  
an appropriate Jordan-measurable $J\subset R\subset \realp 2,$ I see no need to 
state the various well-known functions that lead to this conclusion - with one 
exception. Let compact $J \subset R \subset \realp 2$ and bounded $\vec r\colon J \to 
\realp 3$ generate a surface $\cal S.$ How do we generalize the concept of 
Euclidean area to such a surface? \par
Suppose that $\vec r\,(u,v) = \vec p + u\,\vec a + v\, \vec b,$ where 
$\vec p,\ \vec 
a\not= \vec 0,\  \vec b\not=\vec 0$ are fixed members of $\realp 3$ and
$(u,v) \in [a,b] \times [c,d]= J.$ Then $\vec r\,[J]$ is a paralelogram with 
area equal to $(d-c)(b-a)\Vert \vec a \times \vec b\Vert.$ But$,$ $\vec r_u = \vec a$ 
and $\vec r_v = \vec b.$ Assume that $K = {\rm int}(J) - W,$ where 
$W$ is a set of Jordan-content zero. Then taking other simple Euclidean 
surfaces generated by 
corresponding simple sets of surface equations leads to the conclusion that 
it should be required for functions $\vec r\colon J\to \realp 3$ that $\vec r_u(u,v)$ and 
$\vec r_v(u,v)$ exist and $\Vert \vec r_u(u,v)\times 
\vec r_v(u,v)\Vert \not = 0$ for $(u,v) \in K,$ if 
area preserving surface integration is desired. This is precisely the starting 
point for elementary analysis. Notice that the linear transformations 
that correspond to IR2 are pointwise defined and each is considered to 
generate only the specific hyperparallelogram associated with a specific term 
in the hyperfinite sums used in Definition 8.6.1.  
I leave it to the reader to combine the 
elemental method of constants on surface elements with Definition 8.6.1 
to obtain the concept of flux of a vector field across a surface. \par
\bigskip
\leftline{8.7 \underbar{Microgeometry - Other Stuff.}}\par   
\medskip
What about the geometry of the coordinate transformations? In general$,$ it 
appears necessary to consider stronger constraints then for surface 
integration. If a map $f\colon R \to \realp m,$ where $R \subset \realp m,$ 
is to be 
considered a coordinate transformation$,$ then it  needs to be considered 
locally infinitely close to a linear transformation  determined by a nonzero Jocobian. 
What this signifies is that each m-dimensional infinitesimal subrectangle is 
mapped onto an {infinitesimal parallelepiped - a hyperparallelepiped}.  
Thus the infinitesimally partitioned space is distorted into a space of 
hyperparallelepipeds.
 All the details as to why this appears necessary can be found
in {\tt Stroyan and Luxemburg} beginning in section 5.8 [110]. \par
In the above mentioned section of the book by Stroyan and Luxemburg$,$ you will 
also find in section 5.8 the infinitesimal approach to orientated partitions 
via differential forms. In section 5.9$,$ the authors investigate the 
infinitesimal calculus on manifolds. In all of these sections$,$ the authors 
extend the integral concept to internal functions that usually preserve monads.
Thus their integral concept$,$ although it is the same as has been presented here 
for nonstandard extensions of standard functions$,$ is actually defined on a 
wider class of objects.\par
\bigskip
\leftline{8.8 \underbar{Gauge Integrals}}\par   
\medskip
Recently$,$ additional emphases has been given to a generalization of the 
Riemann sum  called the {\bf {gauge integral}} [{\tt Swartz and Thomson
 [1988]}]. This integral concept is also called the {Riemann-complete 
integral}$,$ [{\tt Henstock [1961]}]; or {S-integral}$,$ [{\tt Mawhin [1985]}]; among 
other names. The gauge integral can also be 
generalized to the {\tt Jarnik$,$ Kurzweil$,$ Schwabik [1983]} integral 
(the {{\bf M-integral}}$,$ Mawhin [1985].) Indeed$,$ a Riemann sum styled generalization 
leads to the classical Lebesgue integral [{\tt McShane [1973]}]. From the 
infinitesimal viewpoint$,$ all of these generalizations are very similar.\par
For the S and M-integrals$,$ a positive real valued function$,$ $\delta,$  called a {\bf 
 {gauge}} 
is utilized to select not only a special partition but a special intermediate 
partition as well.\par
 Within the NSP-world a set of gauges determines an internal 
mapping $\hat {\delta}$ from $\hyper R\subset \hyperrealp m$ into 
$\monad 0$ called a {\bf {microguage}.} However$,$ 
associated with every gauge$,$ $\delta,$ is an object $(Q_\delta,P_\delta)= 
\{(t_1,S_1),\ldots,(t_k,S_k)\}$ where $P_\delta = \{S_1,\ldots,S_k\}$ is a 
subdivision of $R$ by nonoverlapping subsets  
and  $Q_\delta = \{t_1,\ldots,t_k\}$ is a partial sequence where each 
$t_i \in R.$ Such an object  $(Q_\delta,P_\delta)$ called a 
{\bf {$\delta$-fine 
partition}}. Thus for a microguage $\hat {\delta}$ there exists a hyperfinite 
subdivision  
$P_{\hat{\delta}}= \{S_1,\ldots,S_\Gamma\}$ of 
$\hyper R$ and an hyperfinite intermediate partition $Q_{\hat{\delta}}= 
\{t_1,\ldots,t_\Gamma\}.$ The pair $(Q_{\hat{\delta}},P_{\hat{\delta}})$ is called 
a {\bf {micropartition.}} I point out that each member of $P_{\hat{\delta}}$ is 
an infinitesimal subrectangle in this case. 
For a function $f\colon R \to \realp k$$,$ 
Mawhin [1986] shows that {\sl $f$ is S-integrable if and only if there exists an 
$\vec r \in 
\realp k$ such that for each micropartition $(Q_{\hat{\delta}},P_{\hat{\delta}})$ of 
$\hyper R$
$$\sum_{i=1}^\Gamma \hyper f(t_i)\, \hyper v(S_i) \in \monad {\vec r}.$$}\par
Mawhin gives a completely similar characterization for the M-integral but 
replaces the general micropartition  with the more specialized {regular 
micropartition.}\par
Looking at another type gauge$,$ which I shall call the {\bf $L$-gauge}$,$ 
$\lambda,$ and a corresponding partition pair 
$(Q_{\lambda},P_{\lambda}),$ where it is \underbar{not} assumed that
if $t_i \in Q_{\lambda},$ then $t_i \in S_i \in P_{\lambda},$ McShane [1973]
defines his integral by the same Riemann sum technique. The {\bf {McShane 
integral}} 
is equivalent to the Lebesgue integral. Using the same technique employed by
Mawhin$,$ the following is established in appendix 8.
\vskip 18pt 
{\bf Theorem 8.8.1.} {\sl A bounded function $f\colon R \to \real$ is 
Lebesgue integrable with value $r\in  \real$ if and only if  
for each $L$-micropartition $(Q_{\hat{\lambda}},P_{\hat{\lambda}})$ of 
$\hyper R$
$$\sum_{i=1}^\Gamma \hyper f(t_i)\, \hyper v(S_i) \in \monad  r.$$}\par
\vskip 18pt
Thus from the infinitesimal viewpoint there is no basic difference between any 
of these integral concepts discussed in this section. {\bf They are all hyperfinite sums$,$ where each term is 
the scalar product of a limited object (the value of a bounded 
function) by the volume of an infinitesimal subrectangle.}  
\vfil\eject
\centerline{Chapter 9.}
\medskip    
\centerline{\bf THE BEGINNINGS OF}              
\centerline{\bf DIFFERENTIAL EQUATION MODELING}  
\bigskip 
\leftline{9.1. \underbar{Brief Discussion.}}\par
In the previous sections$,$ we presented the standard rules$,$ nonstandard rules 
and the highly successful method of elements in order to establish integral 
models. The strengths of the derivations depended first$,$ upon the acceptance of an 
appropriate correspondence between the terms within a specific discipline and 
objects within the mathematical structure$,$ and then upon three specific methods.
 The presence of the Infinite Sum 
Theorem and the two Self-Evident Theorems led to what appear to be 
strong derivation arguments. As far as the very popular ``{method of 
elements}'' 
is concerned$,$ for the geometric case the nonuniqueness of the chosen elements  
introduced a definite weakness into the method. For the physical case$,$ the 
point charge 
method of application 7.2.5 and the elements of energy used in application 
8.2.1 may not be considered particularly realistic in character. Interestingly 
enough$,$ however$,$ the {\bf {elemental method of constants}} which allows function 
evaluation at ``convenient points'' is often an easily accepted facet of 
infinitesimal modeling. The reason for this is that experience has shown that$,$ 
for the physical applications$,$ quantities that are infinitely close within the 
NSP-world produce the same natural world effects. \par
Differential equation models are somewhat more difficult to justify for 
numerous reasons. Primarily$,$ this is due to the often  vague 
``{rules of correspondence}'' which must be introduced. 
These rules$,$ stated in a secondary 
metalanguage$,$ create a correspondence$,$ $\Phi,$ between a list of terms or relations 
selected from a discipline theory and a list of terms or relations selected from 
a mathematical structure. For differential equation modeling$,$ the relation
$\Phi$ is defined in a more piecemeal fashion than is the case 
for integral modeling. Whereas in integral modeling$,$ the hyperfinite sum 
concept$,$  the Self-Evident Theorems$,$ and the like require that $\Phi$ maps 
certain terms and corresponding relations to specific mathematical objects$,$ 
this is not the process dictated for differential equation modeling.\par
Within a specific discipline what 
constitutes {simplistic} or {idealized behavior} is not just determined by 
nonmathematical experience and intuition.  Using a back-and-forth 
correspondence technique$,$ specific simplistic behavior that approximates the 
more complex is selected in concert with an available mathematical structure. 
This simplistic behavior is then imposed$,$ by intuitive *-transfer$,$ upon an 
infinitesimal portion of the NSP-world.\par
The historical development of the differential calculus does not follow 
the customary 
patterns incorporated into our modern textbooks. Indeed$,$ the originators of this 
subject relied heavily upon visualization and the experiential aspects of 
objective reality$,$ not only for their basic modeling techniques$,$ but also for 
their methods of algebraic manipulation. The scholastic training of the 1600's 
emphasized certain acceptable approaches that did not clash with 
philosophical dictum. Their methods became controversial when 
abstractions were made that led  away from what some considered to be 
perceivable reality.\par
What are the intuitive notions that the originators of the infinitesimal 
concept attempted to model and can we learn the foundations of such modeling 
from their writings? First$,$ it was assumed that geometry was the servant of 
{natural philosophy}. In the 1686 preface of the first edition of Newton's
``Mathematical Principles of Natural Philosophy''  he writes
{\it``Therefore {geometry} is founded in mechanical practice$,$ and is nothing but 
that part of universal mechanics which accurately proposes and demonstrates 
the art of measuring.''} [{\tt Newton [1686]}] Berkeley tells us that 
``{\sl Lines are supposed to be generated by the {motion of points}$,$ planes by 
the motion of lines$,$ and solids by the motion of plans.''} [{\tt Berkeley [1734]}]
Consequently$,$ during the period of its formation the infinitesimal calculus
and its application to natural philosophy was based upon {point-motion} that 
produced geometric configurations. Newton took the concept of the {locus}$,$ made 
additional observations and$,$ with respect to natural philosophy$,$ claimed 
that such motion was caused by ``forces.'' These forces coupled with 
experiential evidence lead to other motions that produced other geometric 
configurations that can predict other motions and observable patterns. 
Thus {physical observation} of {moving objects} came first; followed by the 
concept of {motion}; which 
leads to geometry; which$,$ coupled with {force} relations$,$ leads to other geometric 
patterns; which were interpreted as paths of motion. Finally$,$ these paths of 
motion predicted the location of observable moving objects. As is seen in 
the next section$,$ it is the coalescing of motion$,$ distance traveled and 
absolute time that produced the infinitesimal calculus and
is the fundamental methodology that leads to geometric differential 
modeling.\par
\bigskip
\leftline{9.2. \underbar{The Limit.}}
\medskip
For elementary limit theory all that one needs is a simple alteration of 
Definition 4.4.1. For $\vec p \in \realp n,$ let the {\bf {deleted monad}} be 
$\dmonad {\vec p} = \monad 
{\vec p} - \{\vec p\}.$ \par 
\vskip 18pt
\hrule
\smallskip
{\bf Definition 9.1.1. ({Limit of a Function}).} For any nonzero $n,m \in \nat$ and any 
nonempty $A \subset \realp n$ a function $f\colon A \to \realp m$
has the {\bf LIMIT = $\vec L$ as  $\vec x \to \vec p \in \realp n$} if
$\dmonad {\vec p} \cap \hyper A \not= \emptyset$ and  
$\hyper f[\dmonad {\vec p} \cap \hyper A] \subset \monad {\vec L}.$ \par
\smallskip
\hrule
\vskip 18pt
The following theorem is obtained in the same manner as the corresponding one for 
continuity.\par
\vskip 18pt
{\bf Theorem 9.1.1.} {\sl Definition 9.1.1 for the limit of a function is
equivalent to the classical $\delta - \eps$ definition.}\par 
\vskip 18pt
The basic properties of the limit operator follow from those of the 
standard part operator. And$,$ for this reason$,$ the next theorem$,$ which  
follows from Chapter 10$,$ section 1$,$ corresponds
for extended standard functions to the classical definition of the derivative.
\vskip 18pt
{\bf Theorem 9.1.2. } {\sl Suppose that nonempty open $A \subset \realp n$ and 
$\vec a \in A.$ Then $f\colon A \to \realp m$ is differentiable at $\vec a$ if and only if  there 
exists a unique linear mapping $T_a\colon {\realp n} \to \realp m$ such that for 
each $\vec {\eps} \in \monad {\vec 0},$
$$\hyper f(\vec a + \vec {\eps}\,) - f(\vec a) = \Hyper T_a(\vec {\eps}\,) 
+ (\Vert \vec {\eps}\, \Vert)\lambda(\vec {\eps}\,),$$
where $\lambda(\vec {\eps}\,) \in \monad {\vec 0}.$}
\vskip 18pt
\leftline{9.3. \underbar{Fluxions and Dynamic Geometry.}}\par
\medskip
In the middle 1600's Newton utilized a purely dynamic method to arrive at his 
geometry. He introduced a new type of dynamics that for some natural 
philosophers was not related to the material world and would need to be 
rejected. I will not go into this controversy in this manual. However$,$ Newton's 
{dynamical (loci) methods}$,$ restated from the view point of Robinson's theory,
may be a better 
approach to the {foundations of infinitesimal modeling} than presently presented 
within the standard calculus course. Unfortunately$,$ some of Newton's actual 
demonstrations of the more complex geometric concepts$,$ such as curvature$,$ are 
not valid from the viewpoint of Robinson's theory and require slight 
alterations. When these alterations are conjoined with Robinson's 
theory$,$ then elementary demonstrations are easily obtained and comprehended.\par
In Newton's paper [{\tt Summer 1665}]$,$  an {algorithm} is given that yields the 
relations between the ``fluxions'' $p,\ q$ associated with the variables 
$x,\ y.$ These variables are related by an algebraic expression that is 
assumed to generated a geometric configuration. In Newton [{\tt Oct. 1665 - May 
1666: 383}] the algorithm is specifically described. How Newton$,$ by 
observation$,$ arrived at this algorithm and what exactly  $p,\ q$ represent is
discussed later in his paper. A better explanation of how he formulated his 
algorithm and the meaning of the term fluxion appears in his Oct. 1666
tract [{\tt Newton, 1666}]. Under proposition 7 [{\tt Newton 1666:402}]$,$ he explains his 
algorithm$,$ step by step. After some examples$,$ he discusses how he arrived at 
this algorithm and what {fluxions} signify [{\tt Newton 1666:414}]. First$,$ he 
considers two ``bodies $A$$,$ $B$ moving uniformly.'' He lets an algebraic 
expression $f(x,y) = 0$ represent a relation between the distance traveled 
by these two bodies. Then Newton introduces the concept of the  
distance traveled by a body having uniform velocity $p,$ usually$,$ over a ``moment'' of ``infinitely 
small'' time $o.$ Newton represents the distance each body travels by the sum 
of line segment lengths. Body $A$ first travels along $\overline{ac}$ and at the 
same time body $B$ travels along $\overline{bg}.$ Now in an ``infinitely 
small'' period of time$,$ $o$$,$ body $A$ travels along $\overline{cd}$ and 
during the same time interval body  
$B$ travels along the segment $\overline{gh}.$ He states that the motion is not$,$ 
in general uniform$,$ but it is {\it ``as if the body $A$ with its velocity $p$ 
describe the infinitely little line $\overline{cd}= p \times o$ in one moment$,$ 
in that moment the body $B$ with the velocity $q$ will describe the line
$\overline{gh} = q \times o.$ So that if the described lines be 
$\overline{ac} = x,$ and $\overline{bg}= y,$ in one moment$,$ they will be 
$\overline{ad} = x + p\,o,$ and $\overline{bh} = y + q\,o$ in the next.''}\par
Newton claims that the fluxions  $p$ and $q$ are a type of velocity (what 
type?) and he proceeds to 
demonstrate how relations between these fluxions$,$ in particular the relations 
relative to the quotient $q/p,$ are obtained. It's within 
this demonstration that contradictions occur. Newton writes 
{\it``Now if the equation expressing the relation between the lines $x $ and 
$y$ be $x^3 -abx + a^3 -dy^2 = 0.$ I may substitute $x+ p\,o$ and $y + q\,o$ 
into the place of $x$ and $y$; because (by the above) they as well as $x$ and 
$y$ do signify the lines described by the bodies $A$ and $B.$} [Of course$,$ this 
statement would only be true if fluxions or the motion of the bodies is 
uniform over a \underbar{standard} time interval$,$ $o$$,$ and the ordinary 
Galilean physics is applied.] {\it By doing so there results
$$x^3 + 3pox^2 + 3p^2o^2x+ p^3o^3-dy^2-2dpoy-dq^2o^2$$
$$-abx - abpo$$
$$+a^3 = 0. \eqno (1)$$
But $x^3 -abx + a^3 -dy^2 = 0$ (by supp). Therefore there remains only
$$3pox^2 + 3p^2o^2x+p^3o^3 -2dqoy -dq^2o^2-abpo = 0.\eqno (2)$$
Or dividing it by $o$ it is 
$$3px^2 + 3p^2ox + p^3o^2 -2dqy -dq^2o -abp=0. \eqno (3)$$}
[Thus for the algebraic processes of the 1600's $o$ is a nonzero real number.
Newton goes on to write:] {\it Also those terms are infinitely little in which 
$o$ is. Therefore omitting them there results
$$3px^2 -abp - 2dqy = 0. \eqno (4)$$
The like may be done in all other equations.} Newton would then continue and 
express his important ratio 
$${q \over p} = {{3x^2 -ab}\over{2dy}}.\eqno (5)$$ 
Obviously step (4) is not justified and to some$,$ such as Berkeley$,$ 
contradicts the nature of the infinitely small \underbar{nonzero} quantity $o.$
After this$,$ and in other 
demonstrations$,$ Newton indicates that he guessed at 
portions of his fluxion creating algorithm by applying this steps (1) - (5) 
to numerous algebraic expressions and making certain observations as to the 
physical appearance of such equations (4) and (5). This algorithm$,$ it is claimed$,$ simply 
eliminated the need to
apply continually the above$,$ often criticized$,$ process. Newton repeats similar
derivations in his Winter 1670 - 1671 tract as well as suggesting that the 
delineated process may be applied to relations between three or more 
variables. \par
Assuming that the paths of two objects can be so linearly displayed$,$ then 
modern infinitesimal analysis seems to reveal Newton's intentions and correct 
his vague logic. 
First$,$ note that Newton assumes without justification that the $q$ and $p$ 
that appear in equations (1) - (3) are the same as those that appear in (4). 
Since this should not be so assumed then in all cases suppose that when the 
process arrives at an expression such as (4) that the $q$ is replaced by 
$Q$ and the $p$ by $P.$ Newton seems to describe his notion of the 
``infinitely small'' in his Scholum following Lemma XI of Book one of his 1686 
treatise. In particular$,$ he discusses the ``ultimate velocities'' $Q$ and $P$ 
and his ratios $Q/P$ that he calls ``ultimate ratios.'' {\it Perhaps it may be 
objected$,$ that there is no ultimate proportion of evanescent quantities; 
because the proportion$,$ before the quantities have vanished$,$ is not ultimate$,$ 
and when they are vanished is none. But by the same argument it may be alleged 
that a body arriving at a certain place$,$ and there stopping$,$ has no ultimate 
velocity; because the velocity$,$ before the body comes to the place$,$ is not its 
ultimate velocity; when it has arrived$,$ there is none. But the answer is easy; 
for by the ultimate velocity is meant that with which the body is moved$,$ 
neither before it arrives at its last place and the motion ceases$,$ nor after$,$ 
but at the very instant it arrives; that is$,$ the velocity with which the body 
arrives at  its last place$,$ and with which the motion ceases.  And in like 
manner$,$ by the ultimate ratio of evanescent quantities is to be understood the 
ratio of the quantities not before they vanish$,$ nor afterwards$,$ but with which 
they vanish....For  those ultimate ratios with which quantities vanish are not 
truly the ratios of ultimate quantities$,$ but limits towards which the ratios 
of quantities decreasing without limit do always converge; and to which they 
approach nearer than by any given difference$,$ but never go beyond$,$ nor in 
effect attain to$,$ till the quantities are diminished \underbar{in 
infinitum.}}  What has been described by Newton is very close to the classical 
limit notion.  However$,$ from his applications and arguments this does not seem 
to be what Newton truly believed 
but only a popular exposition that would not offend the geometers of his day.
As is well-known Newton was very fearful of criticism and$,$ even though he 
would use his fluxion methods in private to model physical behavior$,$ he 
did not perform fluxion computations directly within this all important 
research document. 
\par
In equations (1) - (3) suppose that the $o$ is a 
standard real number and that it is claimed that (3) holds in $\cal M$ for any 
such $o$ that is an ``infinitely small '' real number. Berkeley [{\tt Berkeley [1734]}] 
indicates that the Leibniz school would include quantities that are 
``smaller than'' those 
real numbers previously termed as infinitely small. By *-transfer of 
Berkeley's 
description it follows that expressions such as (3) would hold for 
any infinitesimal $\eps.$ Further$,$ Newton's notion of motion and velocities can 
be embedded into an appropriate $\hyperrealp n.$ All other objects that 
appear in such an equation as (3) are considered fixed standard constants or 
relations in  $o.$
 In order to determine what the quantities $p$ or $q$ are 
measuring$,$ let's simplify our investigation and 
assume that the relationship between the distances traveled by body 
$A$ and $B$ is implicitly functional. Indeed$,$ let $f\colon [a,b] \to \real$ 
and$,$ as with Newton$,$ denote $y = f(x),$ where fixed $x \in (a,b).$ 
Replicating 
Newton's process down to step (3) obtains 
$$q ={{f(x + p\, o) - f(x)}\over {o}}.\eqno (6)$$ 
Equation (6) gives a relation between $q$ and $p$ \underbar{and} the {\bf 
arbitrary}
infinitely small \underbar{real number} $o.$ Newton requires $q$ and $p$ to 
be some type of velocity (fluent) and it is not unreasonable that whatever 
kind of velocity these quantities might be measuring that they be bounded 
over his set of infinitely small times. Further$,$ from the above descriptions 
the ``infinitely small''$,$ we suppose that if $o_1$ is infinitely 
small$,$ then any quantity $o$ such that $\vert o \vert < \vert o_1 \vert$ is 
also infinitely small. Thus even if we consider $o$ as a real number$,$ then 
*-transfer yields
$$\hyper q ={{\hyper f(x + \hyper p\, \eps) - f(x)}\over {\eps}},\eqno (7)$$ 
where $\eps \in \monad 0.$\par
If $\hyper p=0,$ then $\hyper q = 0.$  Suppose$,$ therefore$,$ that $\hyper p
\not= 0.$ Then 
$$\hyper q = \hyper p\left({{f(x + \hyper p\, \eps) - f(x)}\over {\hyper 
p\,\eps}}\right).\eqno (8)$$
Assuming that $f$ is differentiable at $x$ then application of definition 
9.1.1 implies that
$$Q = \st {\hyper q} = \st {\hyper p}\St {\left({{f(x + \hyper p\, \eps) 
- f(x)}\over {\hyper p\,\eps}}\right)} = P\,f^\prime(x).\eqno (9)$$ \par
Can we be certain that Newton's ratio $q/p$ may be formally written as 
$Q/P = f^\prime$? In his 1666 tract$,$ Newton claims to obtain the fundamental 
ratio $q/p$ by first substituting for $x$ the expression $x + o$ and  
letting $y = y+ (q/p)\, o,$ where $o$ is an arbitrary infinitely small real 
number. Developing the above formal derivation leads to 
$${q \over p} ={{f(x + o) - f(x)}\over {o}}.\eqno (10)$$   
$$\hyper {\left({q\over p}\right)} ={{\hyper f(x + \eps) - f(x)}\over 
{\eps}}, \eqno (11)$$
$$\St {\hyper {\left({q \over p}\right)}}=\St {\left({{\hyper f(x + \eps) 
- f(x)}\over {\eps}}\right)} = f^\prime(x)= {Q\over P}.\eqno (12)$$ \par
As Newton did in his applications$,$ where $q$ and $p$ where characterized 
merely as  
velocities$,$ he would also arrive at an expression for $q$ only$,$ by assuming 
that since $q$ and $p$ are introduced as arbitrary parameters then 
one could simply let $p = 1.$ The final result he would obtain$,$ prior to 
omitting the higher order $o$'s$,$ would be as expressed by equation (10). Thus 
his two methods are not quite equivalent; but$,$ rather$,$ equation (12) should 
probably be expressed  as 
$$\st {\hyper q}= \St {\left({{\hyper f(x + \eps) 
- f(x)}\over {\eps}}\right)} =f^\prime(x)= Q.\eqno (13)$$\par
Equation (13) is significant for$,$ at least$,$ two reasons in his applications.
First$,$ for {\it unit constant} velocity $p,$ the velocity measured by $q$
is the {\it average velocity} and as such $q$ would behave in the 
infinitesimal world in accordance {Galilean physics}$,$ where infinitesimal time 
$\eps$ is considered as a parameter. Secondly$,$ Newton often used the 
expression
$${q \over p} ={{f(x + o) - f(x)}\over {o}}, p = 1,\eqno (14)$$ 
to obtain what we now know to be the second derivation of $f$ at $x$. This he 
did by writing $z = q$ and then expressing the relation between the velocity 
$r$ (in our notation not his) of $z$ as 
$$z + R\, o = {{f((x+o) + o) - f(x+o)}\over {o}}. \eqno (15)$$
This obviously leads to 
$$r = {{f(x + 2\, o) - 2f(x + o) + f(x)}\over {o^2}} \Rightarrow\eqno (16)$$
$$R = f^{(2)}(x),\eqno (17)$$
assuming that the second derivative of $f$ exists at $x$ and applying Theorem 
8.4.2. 
In his writings 
Newton used his famous  ``dot'' notations which$,$ unfortunately$,$ do 
not correlate formally to their meaning in terms of our modern concepts.
For example$,$ in the above displayed equations he would often write
$p = \dot x,\ q = \dot y$ and $r = \dot z.$ This is not equal to the dot 
notation employed in elementary calculus$,$ where $\dot x = \st {\hyper q}.$ 
\bigskip
\leftline{ 9.4. \underbar{ Fluxions and Higher Order Infinitesimals.}}\par
\medskip 
I have often wondered while reading physical arguments that purport to 
derive a particular expression relating measurable quantities what rules$,$ if 
any$,$ govern the elimination or {omission of ``higher order'' increments}.  
Is it only experience or a deep seated intuition that leads to the assumption 
that the omitted quantities will not$,$ indeed$,$ essentially influence the outcome or 
are their other unexpressed factors that have entered into such derivations? 
\par
The actual method used by Newton to establish the majority of the physical 
principles discussed in his  
Mathematical Principles of Natural Philosophy [{\tt Newton [1686]}] is the theory of fluxions and 
a specialized  process of omitting higher order increments.
Newton would expand expressions into finite or infinite series intending to 
remove eventually all expressions involving his infinitely small $o$ through a
division process. This required him to omit various expressions
involving products with $o.$ However$,$ whether or not he omitted expressions 
involving the higher order terms (i.e. terms containing the factors $o^n,$
where $n \geq 2,$) or not depended 
upon the final proposition being sought.\par   
Technically$,$ the process was not fixed; but$,$ was often altered in such a way 
that by omitting 
certain$,$ but not necessarily all$,$ higher order terms from his expressions 
the final results could be  obtained through a division process that eliminated 
all of the remaining infinitely small factors. Moreover$,$ if$,$ due to such 
technical uncertainties$,$ one attempted 
derivation failed to verify a result$,$ then many others were tried.  
This method of derivation is$,$ 
indeed$,$ vague and forced. But$,$ there is  a basic philosophy being expressed 
by such forced procedures. This modeling philosophy is closely resembles 
the concept of \underbar{physical completeness} for a mathematical model. 
{\sl One believes that the descriptive physical content of a 
theory is absolutely 
correct. Then it is believed that the standard mathematical structure 
correlates exactly to the physical theory. Thus if parameters are introduced$,$ 
such as infinitely small numbers$,$ and any ``reasonable'' method is used to 
eliminate these parameters$,$ then$,$ since the final expressions involve only 
terms that correlate to the original physical theory$,$ the final 
expressions ``must'' give an accurate model for physical reality.} This is the 
same reason why the original creators of differential equation models 
were convinced that their equations had solutions$,$ since in their 
philosophy the equation obtained \underbar{must} predict physical behavior and 
the physical behavior does seem to occur.\par
Obviously$,$ Newton's methods were very successful. The time parameter that 
Newton introduced into his dynamic geometry need not$,$ of course$,$ be an 
actually parameter. Thus when some of the infinitely small time increments could not 
be removed by a division process the final expressions could retain the 
increments expressed in terms of $o$ with the understanding that the result was approximate - the 
approximation being relative to the  ``smallest'' of $o.$ The Newtonian 
methods of arbitrarily eliminating higher order increments continues today.
Our concern$,$ however$,$ is with rigorous differential equation modeling 
where such arbitrariness may be partially eliminated.\par
\bigskip
\leftline{9.5. \underbar{What is a tangent?}} \par
\medskip
Our modern analogue for the fluxion concept is captured by the vector notion.
Today$,$ core mathematics courses utilize the theory of free vectors at an 
early stage in the calculus curriculum. The interpretation of such vectors 
as represented by directed line segments ``attached'' to points in $\realp n$ 
yields the associated geometry. In all that follows$,$ it is assume that all 
vectors are limited (i.e. members of ${\cal O}^n.$) Referring back to Example 
4.4.1A or the first page of Appendix 6$,$ a continuous $c\colon [0,1] \to \realp 
n$ is called a curve and is assumed also to be generated by $n$ continuous 
coordinate functions $f_i\colon [0,1] \to \realp n,$ where for each $t \in [0,1];\ 
c(t) = (f_1(t),\ldots,f_n(t)).$ In order to visualize the 
geometry$,$ one may assume that $n = 3$ even though all of the results in this 
section hold for any positive $n \in \nat.$\par
\vskip 18pt
From our NSP-world analytical viewpoint$,$ a curve is represented entirely by a 
collection of hyperpolygonal curves. Can we formulate a geometric definition 
for a possible tangent vector to the curve $c$ at some $c(t),$ where 
$t \in (0,1)$ which 
appropriately generalizes the Euclidean concept? Obviously$,$ requiring a 
tangent line to be unique would eliminate immediately the basic plane 
geometry definition that such a line would intersect the curve in at most one 
point. Since$,$ however$,$ the properties of the curve should be  entailed by the 
collection of hyperpolygonal representations $\{{\cal P}_\Gamma\}$$,$ then a 
seemingly natural 
definition for the existence of a tangent line might involve the collection of 
all hyperline segments that have standard $c(t)$ as one end point and constitutes 
one of the hyperline segments (the {\bf elements}) in some member of $\{{\cal P}_\Gamma\}.$ More 
specifically$,$ this means that a hypervector that would represent such a 
hyperline segment would be of the form $\hyper c(t + dx) - c(t),$ where 
$dx$ is an nonzero infinitesimal. Intuitively$,$ the tangent line and these 
hyperline segments should be infinitesimally near to each other. This 
corresponds to the requirement that the measure of the plane angle between a 
tangent vector and each of these hyperline segments be infinitesimal or 
infinitesimal close to $\pi.$ The next 
definition models this easily grasped geometric picture since the cosine 
of the measure of such an angle is infinitely close to $\pm 1.$ \par
\bigskip
\hrule
\smallskip
{\bf Definition 9.5.1. ({Tangent Vector to a Curve}).} Let continuous 
$c\colon [0,1] \to \realp n,$ and $t \in (0,1).$ A standard unit vector
$\vec T$ is a {\bf TANGENT} at $c(t)$ if for each  $dx \in \dmonad 
0,$ 
$$\vec T \bullet {{\hyper c(t + dx) - c(t)}\over {\Vert\hyper c(t + dx) - 
c(t)\Vert}} \approx \pm 1.$$ \par
\smallskip
\hrule
\vskip 18pt
Definition 9.5.1 is meaningless if there exists some $dx \in \dmonad 0$ such 
that \break $\Vert\hyper c(t + dx) - c(t)\Vert = 0.$ Now writing Theorem 9.1.2 in the 
form $(\hyper c(t + dx) -c(t))/\vert dx \vert= (\hyper T_t)(dx/\vert dx \vert) 
+ \lambda (dx)$ it is obvious that if $c$ is differentiable at $t$ and 
$c^\prime(t) \not= 0,$ then $\hyper c(t + dx) -c(t) \not= 0$ for any $dx \in 
\dmonad 0.$ Definition 9.5.1 points directly to derivative of $c$ as being a 
primary candidate for tangent consideration$,$ as is well-known. In some 
textbooks on this subject a tangent vector is simply defined as the 
derivative. This approach does not show$,$ conclusively$,$ how the tangent vector 
corresponds to either the geometric notions or to the elemental method. In this 
section applications to geometry are being stressed. For this reason$,$ the 
following theorems are established in their entirety within this section 
rather than hiding them in an appendix. The dot product$,$ $\bullet,$ is the Euclidean 
inner product.\par
\vskip 18pt
{\bf Theorem 9.5.1.} {\sl Let continuous $c\colon [0,1] \to \realp n$ be 
differentiable at $t \in (0,1)$ and $c^\prime(t) \not= \vec 0.$ Then 
$c^\prime(t)/\Vert c^\prime(t) \Vert$ is a unit tangent vector at $c(t).$}\par
Proof. From the above 
discussion we know that $\hyper c(t + dx) -c(t) \not= 0$ for any $dx \in 
\dmonad 0.$  Referring 
back to the basic elemental derivation process of example 4.4.1A in appendix
6 or as used in application 8.2.1 it follows that for each $dx \in \dmonad 0,
\ \Vert \hyper c(t + dx) - c(t) \Vert = \Vert c^\prime(t) \Vert \, \vert dx \vert + \vert 
dx\vert \delta,$ where 
$\delta \in \monad 0.$  Thus let $dx \in \dmonad 0.$ Then
$${c^\prime(t)\over {\Vert c^\prime(t) \Vert}} \bullet {{\hyper c(t + dx) - c(t)}\over {\Vert\hyper c(t + dx) - 
c(t)\Vert}} = {c^\prime(t)\over {\Vert c^\prime(t) \Vert}} \bullet {{{\hyper c(t + dx) - 
c(t)}\over \vert dx\vert}\over{\Vert c^\prime(t) \Vert + \delta}} = \eqno (1)$$
$${c^\prime(t)\over {\Vert c^\prime(t) \Vert}} \bullet {{(\pm c^\prime(t)) + 
\lambda}\over {\Vert c^\prime(t) \Vert + \delta}} = w,\eqno (2)$$
where $\lambda \in \monad {\vec 0}.$ Distributing obtains
$$w = {{\pm c^\prime(t) \bullet c^\prime (t) + c^\prime(t) \bullet \lambda}\over
{\Vert c^\prime(t) \Vert (\Vert c^\prime(t) \Vert + \delta)}} \approx \eqno 
(3)$$
$$\pm {{c^\prime(t) \bullet c^\prime(t)}\over {\Vert c^\prime(t) \Vert 
(\Vert c^\prime(t) \Vert + \delta)}}\approx \pm{{c^\prime(t) \bullet 
c^\prime(t)}\over {\Vert c^\prime(t) \Vert 
\Vert c^\prime(t) \Vert }} =\pm 1. \eqno (4)$$\par
\vskip 18pt
{\bf Theorem 9.5.2.} {\sl Let continuous $c\colon [0,1] \to \realp n$ be 
differentiable at $t \in (0,1)$ and $c^\prime(t) \not= \vec 0.$ If $\vec T$ is a 
unit tangent vector at $c(t),$ then $\vec T = \pm(c^\prime(t)/\Vert c^\prime(t) 
\Vert).$}\par
Proof. Let $\vec T = (a_1,\ldots,a_n)$ and $dx \in \dmonad 0.$ Then 
$$1 \approx\left|\vec T \bullet \left({{\hyper c(t + dx) - c(t)}\over{\Vert \hyper c(t 
+ dx)- c(t)\Vert}}\right)\right|\approx\left|\vec T \bullet \left({{c^\prime(t)}
\over{\Vert c^\prime(t)\Vert}}\right)\right| \eqno (5)$$
Thus
$${\left|\sum_{i=1}^na_if_i^\prime(t)\right|\over 
{\sqrt{\sum_{i=1}^nf_i^\prime(t)^2}}}=1\Rightarrow\eqno(6)$$ 
$$\left|\sum_{i=1}^na_if_i^\prime(t)\right|=\sqrt{\sum_{i=1}^nf_i^\prime(t)^2}
=\left(\sqrt{\sum_{i=1}^na_i^2}\right)
\left(\sqrt{\sum_{i=1}^nf_i^\prime(t)^2}\right).\eqno(7)$$
Since neither vector is a zero vector then the Schwarz Inequality states that 
there exists some nonzero $d$ such that $\vec T = d(c^\prime(t)/\Vert 
c^\prime(t) \Vert.$ Thus $\Vert \vec T\Vert = \vert d\vert\bigm \Vert c^\prime(t)/\Vert 
c^\prime(t) \Vert\bigm \Vert \Rightarrow \vert d \vert = 1$ and the proof is 
complete.\par
\vskip 18pt\vfil\eject
\leftline{9.6. \underbar{What is an (k-surface) Osculating Plane?}} \par
\medskip
One of the early concepts used by proponents of differential geometry is the 
notion of ``consecutive'' points in a curve or surface. Modern textbooks$,$ such 
as {\tt Struik [1961]}$,$ emphasize this intuitive approach$,$ even to the point that 
some of their derivations may be considered less than fully justified. The 
definitions for those significant geometric objects associated with the 
tangents are often given without adequate discussion as to the  
geometric content of the definition  and why various constraints are necessary.
Modern infinitesimal analysis eliminates this vagueness.\par
Before extending our results to the $k$-surface case we look at the 
3-dimensional prototype. Let $t_0,t_1,t_2 \in (0,1)$ and assume that the vectors
$\vec c_1 = c(t_1) - c(t_0),\ \vec c_2 = c(t_2) - c(t_0)$ are linearly 
independent. The unique plane containing the two line segment 
representations for $\vec c_1,\ \vec c_2$ attached to $c(t_0)$ is 
$\Pi(c(t_0),\vec c_1,\vec c_2) = \{\vec x \bigm \vert \vec x\in \realp 3 \land
(\vec x = c(t_0) + \lambda\vec c_1 + \nu\vec c_2 )\land (\lambda,\ \nu \in 
\real)\}.$ Elementary linear algebra may also be used to find 
$\Pi(c(t_0),\vec c_1,\vec c_2).$ Consider the matrix                          
$$ A=A(\vec x,c(t_0),\vec c_1,\vec c_2)=\left(\matrix{\vec x -c(t_0)\cr 
                   \vec c_1\cr
                   \vec c_2\cr}\right).\eqno (1)$$
Let nonzero $\vec b$ be a member of the 1-dimensional null space of 
$A(c(t_0),c(t_0),\vec c_1,\vec c_2).$
Then $\Pi(c(t_0),\vec c_1,\vec c_2)= \Pi(\vec b,c(t_0)) = 
\{\vec x \bigm \vert \vec x \in \realp 3 \land (\vec x - c(t_0))\bullet \vec b 
= 0 \}= \{ \vec x\bigm \vert \vec x \in \realp 3 \land \det(A) = 0\}.$  
Once $\vec b$ is selected then the fact that $c(t_0),c(t_1),c(t_2)$ 
are members of this plane can be characterized by the equation
$$f(t) = (c(t)\bullet \vec b) - (c(t_0)\bullet \vec b). \eqno (2)$$ \par                                     
Now let us assume that the second derivation $c^{(2)}$ is continuous on 
$(0,1)$ and that $t_0 <t_1 < t_2.$ Since$,$ in this case$,$ $f(t_0) = f(t_1) = 
f(t_2) = 0$ and $[t_0,t_2] \subset (0,1),$ then  
Rolle's theorem tells us that  
there exists $v_1,\ v_2$ such that
$t_0 \leq v_1 \leq t_1 \leq v_2\leq t_2$ and such that 
$$f^{(1)}(v_1) = c^{(1)}(v_1)\bullet \vec b = 0$$
$$f^{(1)}(v_2) = c^{(1)}(v_2)\bullet \vec b = 0.\eqno (3)$$
Moreover$,$ there exists some $v_3$ such that $v_1 \leq v_3 \leq v_2$ and 
$$f^{(2)}(v_3) = c^{(2)}(v_3)\bullet \vec b = 0. \eqno (4)$$\par 
Embedding the above remarks into the NSP-world then it follows that for 
three $t_0,t_1,t_2 \in \Hyper (0,1)$ and $t_0\approx t_1 \approx t_2\approx 
t^\prime \in (0,1)$ (ordered 
as previously) such that $\hyper c(t_1) - \hyper c(t_0)= \hyper c_1,\ 
\hyper c(t_2) - 
\hyper c(t_0)=\hyper c_2$ are linearly independent in $\hyperrealp 3$ 
and for any nonzero $\vec b \in {\cal O}^3$ in the 1-dimensional null space 
for $\hyper A(\hyper c(t_0),\hyper c(t_0),\hyper c_1,\hyper c_2)$ 
there exist
$v_1,\ v_3 $ such that $v_1 \approx v_3 \approx t_0 \approx t_1\approx t_3$ 
and 
$$(\hyper c(t_0)\bullet \vec b) - (\hyper c(t_0)\bullet \vec b) =\hyper 
c^{(1)}(v_1)\bullet 
\vec b = \hyper c^{(2)}(v_3)\bullet \vec b = 0.\eqno (5)$$\par
But$,$ we also know that $\{v_1,v_3,t_0,t_1,t_2\} \subset \monad {t^\prime}.$ 
Continuity of the derivatives yields that 
$$c^{(1)}(t^\prime)\bullet \vec b \approx 0,\ c^{(2)}(t^\prime)\bullet \vec b 
\approx 0.\eqno (6)$$
and thus 
$$c^{(1)}(t^\prime)\bullet \st {\vec b} = 0, \ c^{(2)}(t^\prime)\bullet \st 
{\vec b} = 0.\eqno (7)$$
\par
Now what would happen if $c^{(1)}(t^\prime),c^{(2)}(t^\prime)$ are linearly 
independent? Let $ B=A(\vec x,c(t^\prime),c^{(1)}(t^\prime),c^{(2)}
(t^\prime))$ and consider the plane $\Pi(\st {\vec b},c(t^\prime)).$ Obviously$,$ 
if $\vec b$ is considered a limited normal to the object
$\Hyper \Pi(\vec b,\hyper c(t_0)),$ then we have normals to the two planes 
that are infinitely close. More significantly$,$ however$,$ is the fact that these 
planes do satisfy an ``infinitely close '' type property. Extending the above 
discussion to $n$-dimensional spaces we have \par
\vskip 18pt
{\bf Theorem 9.6.1.}\par 
{\sl (i) Let $c\colon [0,1] \to \realp n,\ n > 1$ and assume that for some
$k$ such that $1 \leq k \leq n-1$ 
the $k$-derivative$,$ $c^{(k)},$ is continuous on $(0,1).$\par
(ii) Let $\{t_0,\ldots,t_k\} \subset \Hyper (0,1)$ have the properties that
$\{\hyper c_i \bigm \vert \hyper c_i = \hyper c(t_i) - \hyper c(t_0) \land 
1\leq i \leq k\}$ is linearly independent in $\hyperrealp n$ and $t_i \approx 
t_j\approx t^\prime \in (0,1),\break 0 \leq i,j \leq k.$\par
(iii) Let $\{c^{(1)}(t^\prime),\ldots,c^{(k)}(t^\prime)\}$ be linearly 
independent.\par 
\vskip 18pt
Then for every
limited $\vec a \in \Pi_1 = \{\vec x \bigm \vert \vec x \in \hyperrealp n 
\land (\vec x = \hyper c(t_0) + \lambda_1\hyper c_1 + \cdots + 
\lambda_k\hyper c_k)
\land (\lambda_1,\ldots,\lambda_k \in \hyperreal) \}$                         
there exists $\vec c \in \Pi = 
\{\vec x \bigm \vert \vec x \in \realp n \land (\vec x = c(t^\prime)
+ \lambda_1c^{(1)}(t^\prime) + \cdots + \lambda_kc^{(k)}(t^\prime))
\land (\lambda_1,\ldots,\lambda_k \in \real) \}$   
such that $\vec a \approx \vec c.$}\par
Proof. As in the illustrated case$,$ consider the matrix
$$A(\vec x,c(t_0),\vec c_1,\ldots,\vec c_k) =\left(\matrix{\vec x- c(t_0)\cr
                    \vec c_1\cr
                     \cdot\cr
                      \cdot\cr
                     \cdot\cr
                     \vec c_k\cr}\right)\eqno (9)$$
Let $W$ be the $k$-dimensional subset of $\realp n$ spanned by $\{\vec 
c_1,\ldots,\vec c_k\}$ and $W^\bot$ the $(n-k)$-dimensional orthogonal complement 
with respect to the Euclidean inner product. 
Let $\{\vec b_1,\ldots,\vec b_{n-k}\}$ be an orthonormal basis for $W^\bot.$ Now let 
$B$ be the $n \times (n-k)$ matrix formed by considering the vectors $\vec b_i, \ 
1\leq i\leq n-k$ as column vectors. Then consider the (affine) {$k$-plane}
$\Pi(B,c(t_0)) = \{ \vec x\bigm \vert \vec x \in \realp n \land\break (\vec x - 
c(t_0))\cdot B = 0\} = \{\vec x \bigm \vert \vec x \in \realp n \land \exists
z(z \in W \land \vec x = z + c(t_0))\}= \{\vec x \bigm \vert \vec x 
\in \realp n 
\land (\vec x = c(t_0) + \lambda_1 \vec c_1 + \cdots + \lambda_k\vec c_k)
\land (\lambda_1,\ldots,\lambda_k \in \real )\},$ where $t_0 \in (0,1).$ Let
 $$f_i(t) = (c(t) - c(t_0))\bullet \vec b_i,\ 1\leq 1\leq n-k, \eqno (10)$$
where $t \in (0,1).$ 
Applying the same argument as in the above illustration we obtain the 
hyperplane $\Pi_1 = \Hyper \Pi(B,\hyper c(t_0)),$ where $B$ is an $n \times 
(n-k)$ matrix having limited column vectors. Notice that the standard part of $B$ is the 
standard part of each column vector and since the column vectors form an 
*-orthonormal set in $\hyperrealp n$ then the set of standard parts forms an 
orthonormal$,$ hence independent$,$ set of vectors in $\realp n.$ Now 
equations (7) becomes
$$c^{(1)}(t^\prime) \cdot \st B = 0,\ldots, c^{(k)}(t^\prime) \cdot \st B = 
0.\eqno (11)$$
Hence$,$ $\Pi = \{\vec x\bigm \vert \vec x \in \realp n \land (\vec x -
c(t^\prime))\cdot \st B = 0 \}.$ Consequently$,$ if limited $\vec a \in \Pi_1,$ 
then $(\vec a-\hyper c(t_0)) \cdot B = 0 \Rightarrow
(\st {\vec a}- \st {\hyper c(t_0)})\cdot \st B = (\st {\vec a} - 
c(t^\prime))\cdot \st B= 0 \Rightarrow \vec c = \st {\vec a} \in
\Pi$ and that $\vec a \approx \vec c.$ \par
\vskip 18pt
Any set of $(k+1)$ points that satisfies part (ii) of Theorem 9.6.1 is called a 
set of {\bf {$(k+1)$-consecutive} points in the curve $c.$} If all of the 
hypotheses of Theorem 9.6.1. are satisfied$,$ then $\Pi$ is called the 
{\bf { osculating $k$-plane}} and this object is infinitely 
close$,$ 
in the above sense$,$ to any of the other objects $\Pi_1.$ Observe that if  
$k=1,$ then $\Pi$ is but that tangent line to $c$ at $c(t^\prime).$ \par
\bigskip
\leftline{9.7. \underbar{Curvature.}} \par
\medskip
Continue to assume that for the curve $c\colon [0,1] \to \realp n$ the 
$k$-derivative$,$ $c^{(k)},$ is continuous on $(0,1)$ and for our 
discussion let $k=2\ (n>2),$ even though the following results 
also appear to hold for $k=3.$ As is customary for technical purposes$,$ 
consider $c$ re-expressed in terms of the arc length parameter$,$ $s.$ 
Therefore$,$ assume that $c\colon [0,L] \to \realp n$ has continuous
$k$-derivatives on $(0,L),$ where $L$ is the 
length of the arc. For 
$s^\prime \in (0,L),$ let $\{dc(s^\prime)/ds,d^2c(s^\prime)/ds^2\}$ be linearly 
independent. \par
Since $dc(s^\prime)/ds \not= 0$ then Theorem 9.5.1 implies that 
$dc(s^\prime)/ds = T$ is a unit tangent vector at $c(s^\prime).$ It is an easy 
exercise to show that $dT(s^\prime)/ds = d^2c(s^\prime)/ds^2$ is orthogonal to
$T$ and is a linear combination of $\{c^{(1)}(t^\prime),c^{(2)}(t^\prime)\}$ 
for the corresponding $t^\prime \in (0,1).$ Thus$,$ $dT(s^\prime)/ds$ is a 
member of the osculating $k$-plane at $c(t^\prime)$ and for this useful reason
$dT(s^\prime)/ds$ is selected as a (principle) normal vector to the curve in 
that it is orthogonal to $T.$ Normalizing$,$ we obtain $\vec k = \kappa\, \vec N,$ 
where $\vec N$ is a {\bf {unit normal vector to the curve}.} What might 
$\kappa$ represent$,$ where we have that $\kappa>0$?\par
As in section 9.6$,$ but in terms of the arc length parameter$,$ let 
$s_0,s_1,s_2 \in (0,L)$ and assume that the vectors
$\vec c_1 = c(s_1) - c(s_0),\ \vec c_2 = c(s_2) - c(s_0)$ are linearly 
independent. The unique $k$-plane containing the two line segment 
representations for $\vec c_1,\ \vec c_2$ attached to $c(s_0)$ is 
$\Pi(c(s_0),\vec c_1,\vec c_2) = \{\vec x \bigm \vert \vec x\in \realp n \land
(\vec x = c(s_0) + \lambda\vec c_1 + \nu\vec c_2 )\land (\lambda,\ \nu \in 
\real)\}.$ Consider a $k$-sphere of radius $r$ with center $\vec c \in \Pi(c(s_0),
\vec c_1,\vec c_2).$ A general equation for such a $k$-sphere may be written as
$$f(s) =(\vec x - \vec c)\bullet (\vec x - \vec c) - r^2 = 0, \eqno (1)$$
where we have assumed$,$ as with equation (2) in section 9.6$,$ that (1) 
is a function in $s$ 
and that it is satisfied by the points $s_0,s_1,s_2.$ Following the exact same 
process that led from equation (2) through equation (5) of section 9.6 yields,
for $s_0,s_1,s_2 \in \Hyper (0,L),\ s_0\approx s_1\approx s_2 \approx s^\prime
\in (0,L),\ \vec a \in \Pi_1$ and $r_1 \in \hyperreal,$ 
the hyperreal numbers $v_1\approx v_3\approx s^\prime$ such that
$$(\hyper c(s_0) - \vec a)\bullet(\hyper c(s_0) - \vec a) - r_1^2 = 0$$
$$(\hyper c(s_0) - \vec a)\bullet \Hyper (dc/ds)(v_1) =0\eqno (2) $$
$$(\hyper c(s_0) - \vec a)\bullet \Hyper (d^2c/ds^2)(v_3) +  
\Hyper (dc/ds)(v_3)\bullet \Hyper (dc/ds)(v_3)=0.$$\par
The third equation in (2) implies that $\vec a\in {\cal O}^n$. 
The first equation shows that $r_1 \in {\cal O}$ if and only if 
$\vec a\in {\cal O}^n.$ The standard part operator obtains
$$(c(s^\prime) - \st {\vec a})\bullet(c(s^\prime) - \st {\vec a}) - (\st 
{r_1})^2 =0 $$ 
$$(c(s^\prime) - \st {\vec a})\bullet  (dc/ds)(s^\prime) =0$$
$$(c(s^\prime) - \st {\vec a})\bullet  (d^2c/ds^2)(s^\prime) +  
(dc/ds)(s^\prime)\bullet (dc/ds)(s^\prime)=\eqno (3)$$
$$(c(s^\prime) - \st {\vec a})\bullet  (d^2c/ds^2)(s^\prime) + 1=0.$$
Now let $\vec c =\st {\vec a}$ and $r = \st {r_1}$. By Theorem 9.6.1$,$ $\vec 
c$ is a member of the osculating $k$-plane and thus $c(s^\prime) - \vec c= 
\lambda(dc(s^\prime)/ds) + \nu (d^2c(s^\prime)/ds^2).$ Employing equations (3) we 
find that $\lambda = 0$ and $-1 = \nu (d^2c(s^\prime)/ds^2)\bullet 
d^2c(s^\prime)/ds^2)=\nu \kappa^2 .$ Consequently$,$ 
$$ c(s^\prime)- \vec c= (1/\kappa)\vec N, \eqno (4)$$
and the radius of the $k$-sphere or the $k$-circle formed by the intersection of the 
osculating $k$-plane and this $k$-sphere is $1/\kappa.$ As is well-known 
this $k$-circle is called the {\bf {osculating $k$-circle}} 
 and gives a measure to the 
curvature of $c$ at $c(s^\prime).$ \par
Summarizing these concepts from the NSP-world view point$,$ we have  
hypercircles in the hyperplane $\Pi_1$ of Theorem 9.6.1 determined by the three
consecutive points $\hyper c(s_0),\hyper c(s_1),\hyper c(s_2)$ that are 
infinitely close to $c(s^\prime)$ with the center $\vec a \in \Pi_1$ and 
radius $r_1 \in \hyperreal.$ The continuity of the various derivatives 
involved implies that the center of any such 
hypercircle is limited and infinitely close to the center of a unique circle 
in $\Pi$ with 
coordinates $c = c(s^\prime) + (1/\kappa)\vec N$ and having radius $1/\kappa.$ 
Notice that if $\kappa < 0$ the same result holds except that the radius of 
this unique circle is $\vert 1/\kappa \vert.$ 
The 2-dimensional case is obtained by considering the 3-dimensional case but 
restricting the curve to a coordinate plane. Originally$,$ the concept 
of the curvature$,$ as we know it today$,$ was not the approach used by Newton. His 
derivation was for the radius of curvature in the plane. He obtained the 
correct expression but his proof was filled with serious logical 
contradictions. Intuitively$,$ differential geometers describe the osculating 
circle as the 
unique circle containing $c(s^\prime)$ in the osculating plane and 
``containing three consecutive points'' from the curve. The nonstandard 
approach has eliminated these vague concepts with their often forced 
derivations and replaces them with comprehensible and rigorous ideas. 
It would be significant to undergraduate education if the entire subject 
matter of elementary differential geometry$,$ with its wide area of 
physical application$,$ was rewritten and presented in this more easily 
visualized NSP-world approach. Once the basic notions are mastered$,$ then 
the intuition developed is an indispensable asset when they are extended 
to differentiable manifolds.   
\vfil\eject
\centerline{Chapter 10.}
\medskip    
\centerline{\bf THE DIFFERENTIAL}              
\centerline{\bf AND PHYSICAL MODELING}  
\bigskip 
\leftline{10.1. \underbar{Basic Properties.}}\par
There is one principle that is utilized continually in order to pass from the 
monadic environment to the standard world without implication of the standard 
part operator. The principle is called the {\bf {Cauchy Principle}} and is 
the critical step required to establish the major results in this section.\par
\vskip 18pt
{\bf Theorem 10.1.1. (Cauchy Principle)} {\sl Let $P(x)$ be a bounded formal first-order 
expression in one fee variable and employing internal constants (i.e. a 
bounded formal internal property as discussed in Appendix 4$,$ within the proof 
of Theorem 4.2.2.) If $P(\eps)$ holds for each  
$\eps \in \nmonad {\vec 0},$ then there 
exists an open ball $B$ about $\vec 0$ such that $P(\vec x)$ holds for each 
$\vec x \in \hyper B.$}\par
\vskip 18pt
In the last section$,$ it was assumed that we had established the concept of the 
differential formally. Obviously$,$ this was not done and the situation needs to  
be rectified prior to our brief applications to physical problems. First$,$ a 
slight notational change is beneficial. Let $\vec v \in \realp n.$ Then to 
prevent confusion in this chapter$,$ we denote $\monad {\vec v}$ by $\nmonad 
{\vec v}.$ I point out that what follows can be easily extended to normed 
linear spaces over the real or complex numbers.\par
\vskip 18pt
\hrule
\smallskip
{\bf Definition 10.1.1 ({Local Map}).} Let $A \subset \hyperrealp n$ and 
internal $f\colon A \to \hyperrealp m.$ Then $f$ is a {\bf LOCAL MAP} if
$f[\nmonad {\vec 0\,}] \subset \mmonad {\vec 0\,}.$ This means that $f$ 
{\bf {preserves infinitesimals}}.\par
\smallskip
\hrule
\vskip 18pt
\hrule
\smallskip
{\bf Definition 10.1.2. ({Equivalent  Maps}).} Let $f$ and $g$ be local maps. 
If for each nonzero $\vec x \in \nmonad {\vec 0\,}$
$${f(\vec x) - g(\vec x)\over{\Vert \vec x \Vert}} \in \mmonad {\vec 0\,},$$ 
then we write $f \sim g$ and say that the two maps are {\bf EQUIVALENT.}\par
\smallskip
\hrule
\vskip 18pt
{\bf Theorem 10.1.2.} {\sl Two internal maps $f$ and $g$ are equivalent if and 
only if
there is a local map $\alpha$ such that for each nonzero $\vec x \in \nmonad {\vec 0\,}$
$$f(\vec x) = g(\vec x) + (\alpha (\vec x))\Vert \vec x \Vert.$$}\par
\vskip 18pt
Notice that since $\nmonad {\vec 0}$ and $\mmonad {\vec 0}$ are external sets
and the domain$,$ $D_f,$ and range$,$ $f[D_f],$ of an internal map$,$ $f$$,$ 
are internal sets then $\nmonad {\vec 0}$ and $\mmonad {\vec 0}$ are$,$ 
respectively$,$ proper subsets of $D_f$ and $f[D_f].$ 
Probably the simplest type of local map that one can associate with a 
linear space would be one that preserves linearity.\par
\vskip 18pt
\hrule
\smallskip
{\bf Definition 10.1.3. ({Locally Linear})} A local map $f$ is 
said to be {\bf {LOCALLY LINEAR}} (i) if for each $\vec x,\ \vec y \in 
\nmonad {\vec 0\,},$ then $f(\vec x + \vec y) = f(\vec x) + f(\vec y)$ and 
(ii) for each $r \in {\cal O}$ and each $\vec x \in \nmonad {\vec 0\,}$ it 
follows that $f(r\vec x) = rf(\vec x).$\par
\smallskip
\hrule
\vskip 18pt
The interesting thing about locally linear maps is that even though they are 
restricted to linearity on infinitesimals they are uniquely associated with an 
internal linear map defined on the entire space $\hyperrealp n.$\par
\vskip 18pt
{\bf Theorem 10.1.3.} {\sl Suppose that $f$ is any locally linear map. Then 
there exists a unique internal linear map $T_f\colon \hyperrealp n \to 
\hyperrealp m$ such that $\Vert T \Vert \in {\cal O}$ and there exists an open 
$E \subset \realp n,$ where $\vec 0 \in E,$ such that 
$f(\vec x) = T_f(\vec 
x)$ for each $\vec x \in \hyper E.$}\par
\vskip 18pt
For the present$,$ we are not particularly concerned with the domain nor range 
of a locally linear map only that it satisfies the infinitesimal requirements. 
What the next few theorems indicate may be slightly surprising. \par
\vskip 18pt
{\bf Theorem 10.1.4. } {\sl Let $f,\ g$ be locally linear maps. Then 
$f\sim g$ if and only if $T_f(\vec x) \approx T_g(\vec x)$ for each $\vec x \in 
\hyperrealp n$ such that $\Vert \vec x \Vert = 1.$}\par
\vskip 18pt
Previously$,$ the differential was defined for $f\colon [a,b] \to \real$ at a 
point $c \in (a,b)$ where $f^\prime$ exists to be $f^\prime(c)\, dx$ 
for each $dx \in \monad 0.$ This can be viewed as the restriction of the 
nonstandard extension of the standard function $D\colon [a,b] \to \real$ 
defined by $D(x) = f^\prime(c)\, x.$ \par
\vskip 18pt
\hrule
\smallskip
{\bf Definition 10.1.4 ({Differential}).} A local linear map $f$ is a 
{\bf DIFFERENTIAL} if there exists a (standard) bounded linear transformation 
$F\colon \realp n \to \realp m$ such that $\hyper F(\vec x) = f(\vec x)$ 
for each $\vec x \in \nmonad {\vec 0\,}.$\par
\smallskip
\hrule
\vskip 18pt
{\bf Theorem 10.1.5.} {\sl If $f$ is a differential$,$ then there exists a 
unique bounded linear transformation $F\colon \realp n \to \realp m$ and an open 
set $E\subset \realp n$ such that $\vec 0 \in E$ and $\hyper F(\vec x) = f(\vec 
x)$ for all $\vec x \in \hyper E$ and$,$ in particular$,$ $F(\vec x) = f(\vec x)$ 
for each $\vec x \in E.$}\par
\vskip 18pt
What Theorem 10.1.5 does is to lift the internal differential $f$ up to the 
status of a standard linear transformation on some standard nonempty set $E.$ 
This may be a slight surprise since it is only assumed that this is the case 
for the infinitesimals. For applications$,$ the behavior of a local map with 
respect to the infinitesimals is of paramount importance. Further$,$ 
differentials are essentially unique with respect to equivalence.\par
\vskip 18pt
{\bf Theorem 10.1.6.} {\sl If $f$ and $g$ are differentials  and $f \sim g,$ 
then  
there exists some open $E \in \realp n$ such that $\vec 0 \in E$ and 
$f(\vec x) = g(\vec x)$ for each $\vec x \in \hyper E.$}\par
{\bf Corollary 10.1.6.1.} {\sl If $f$ and $g$ are differentials$,$ $f \sim g,$
and $F,\ G$ are$,$ respectively$,$ the corresponding (standard) linear transformation 
representations for $f,\ g,$ then $F=G.$}\par
\vskip 18pt
We now establish the basic relation between differentials and the 
derivative.\par
\vskip 18pt\vfil\eject
\hrule
\smallskip
{\bf Definition 10.1.5. ({Local Increment}).} Let nonempty open 
$G \subset \realp n,$ and $\vec c \in G.$ For any $f\colon G \to \realp m$ the 
{\bf LOCAL INCREMENT} is the map $\Delta \hyper f_{\vec c}$ defined on 
$\nmonad {\vec 0\,}$ by letting $$\Delta \hyper f_{\vec c}(\vec b\,)= 
\hyper f(\vec c + \vec b \,) - 
f(\vec c\,),$$ where $\vec b \in \nmonad {\vec 0\,}.$\par
\smallskip
\hrule
\vskip 18pt
{\bf Theorem 10.1.7.} {\sl Let nonempty open 
$G \subset \realp n,$ and $\vec c \in G.$ A function $f\colon G \to \realp m$ 
is continuous at $c \in G$ if and only if the local increment map 
$\Delta \hyper f_{\vec c}$ is a local map.}\par
\vskip 18pt
Theorems 10.1.5 and 10.1.6 show that there is a standard open set $E$ about 
$\vec 0$ such that a differential is unique on $E$ up to equivalence and is 
exactly equal to some linear transformation restricted to 
$E.$ However$,$ the 
behavior of a differential was only characterized with respect to 
infinitesimals. These results coupled with Theorem 10.1.7 suggest the follow 
definition.\par
\vskip 18pt
\hrule
\smallskip
{\bf Definition 10.1.6. ({Differentiable}).} Let nonempty open 
$G \subset \realp n,$ and $\vec c \in G.$ A function $f\colon G \to \realp m$  
is {\bf DIFFERENTIABLE AT c} if there exists a differential$,$ denoted by 
$d_{\vec c}f,$ that is  equivalent to $\Delta \hyper f_{\vec c}.$\par
\smallskip
\hrule
\vskip 18pt                                                                              
Obviously$,$ from the above results if $d_{\vec c}f$ exists then it is 
representable by a linear transformation $L$ on some standard open 
neighborhood of $\vec 0.$ To complete our basic results$,$ 
all that is needed is to determine the appearance of an $m\times n$-matrix 
representation $A$ for such a transformation. However$,$ this is immediate from 
the fact that if $f,$ of Definition 10.1.5$,$ is differentiable at $c$$,$ then for 
each nonzero $\vec b \in \nmonad {\vec 0\,}$
$$ \Delta \hyper f_{\vec c}(\vec b\,) = \hyper A(\vec b\,) + 
(\alpha(\vec b\,))\Vert \vec b \Vert.\eqno (1)$$
                                                            
\noindent Writing this in coordinate function and column vector form yields
$$(\hyper f_1(\vec c + \vec b\,) - \hyper f_1(\vec c\,),\ldots,
\hyper f_m(\vec c +\vec b\,)-f_m(\vec c\,))^T =(\hyper A){\vec b\,}^T+ (\alpha 
(\vec b\,))\Vert \vec b \Vert,\eqno (2)$$
where nonzero $\vec b \in \nmonad {\vec 0}.$
Now considering any
$$\vec b = \overbrace{(0,\ldots,\eps,\ldots,0)}^j, \ \eps \in \dmonad 
0,\eqno (3)$$
then $\hyper f_{i}( \vec c+  \vec b\, ) - f_{i}(\vec c\, )= a_{i}\eps + (\alpha ( 
\vec b\, ))\vert \eps \vert$ implies that $a_{i} \approx (\hyper f_{i}(\vec 
c + \vec b\,) - 
f_{i}(\vec c\, ))/\eps.$ 
Thus Definition 9.1.1 and the fact that each $a_{i}$ is 
standard yields $a_{i} = \partial f_{i}(\vec c\, )/\partial x_j.$
The elementary propositions about differentiable functions follow readily 
from Definition 10.1.6 and equation (1).\par
With respect to applications equation (1) or 
$$\Delta \hyper f_{\vec c}(\vec b\,) \approx \hyper A(\vec b\,) \eqno (4)$$ 
is the most useful. 
Further$,$ $\Delta \hyper f_{\vec c}(\vec b\,) $ and $\hyper A(\vec b\,)$ are 
related by first-order ideals in the sense that 
$$ \Vert\Delta \hyper f_{\vec c}(\vec b\,) - \hyper A(\vec b\,)\Vert = \Vert
\alpha(\vec b\,)\Vert \,\Vert\vec b \Vert.\eqno (5)$$ 
implies that 
$$ \Vert\Delta \hyper f_{\vec c}(\vec b\,) - \hyper A(\vec b\,)\Vert \in 
o(\Vert \alpha (\vec b\,)\Vert) = o(\eps),\ \eps \in \monad 0,\eqno (6)$$
from Theorem 8.3.5.\par
As far as the notion of the  {n-dimensional derivative} is concerned$,$ simply 
consider equation (1) written for each nonzero $\vec b \in \nmonad {\vec 0}$
as
$${{\Delta \hyper f_{\vec c}(\vec b\,)}\over {\Vert \vec b \Vert}} = 
\hyper A\left({{\vec b}\over {\Vert \vec b \Vert}}\right) + \alpha(\vec 
b\,).\eqno 
(7)$$
For the directional derivative$,$ assuming that $f$ is differentiable at $\vec 
c \in G,$ start with a fixed standard unit vector $\hat u,$ let $\eps \in 
\dmonad 0$ and write (1) as
$$ \Delta \hyper f_{\vec c}(\eps{\hat u}) = \hyper A(\eps{\hat u}) + 
(\alpha(\eps{\hat u}))\Vert \eps{\hat u}\Vert.\eqno (8)$$ 
This leads to                                                              
$$ \Delta \hyper f_{\vec c}(\eps{\hat u}) = \hyper A(\eps{\hat u}) + 
(\alpha(\eps{\hat u}))\vert \eps\vert;\eqno (9)$$ 
which may be written in the more familiar form
$$ {{\Delta \hyper f_{\vec c}(\eps{\hat u})}\over {\eps}}=
\hyper A({\hat u}) \pm  \alpha(\eps{\hat u})\Rightarrow\eqno (10)$$ 
$$\St {\left({{\Delta \hyper f_{\vec c}(\eps{\hat u})}\over {\eps}}\right)}=
A({\hat u}).\eqno (11)$$
Note that if (8) holds for all $\hat u \in \hyperreal^n$ and $\eps \in \monad {0},$ then (1) holds.\par
 There is one other useful concept relative to applications and vectors.
 Two nonzero $\vec v,\  \vec w \in \hyperrealp n$ are said to be 
{\bf {almost parallel}} if $\vec u = \vec v/\Vert \vec v \Vert \approx \pm
\vec w/\Vert \vec w \Vert.$ In an appropriate infinitesimal microscope$,$ two 
almost parallel vectors appear to be parallel and$,$ for applications$,$ yield the 
same effects as do parallel vectors. The hypotheses of the next theorem are 
often realized in applied problems.\parm
{\bf Theorem 10.1.8.} {\sl Let noninfinitesimal $\vec a \in {\cal O}^n$ and 
suppose that $\vec a - \vec b \in \nmonad {\vec 0}.$ Then $\vec a$ is almost 
parallel to $\vec b.$}\par
\bigskip
\leftline{10.2. \underbar{Some General Observations}}\par 
\medskip
In this basic manual$,$ it is only possible to give a few cursory illustrations 
of how to apply the 
derivative to physical problems since there are but a few general procedures 
that can be followed. Moreover$,$ most such applications require specific  
knowledge relative to the refined behavior patterns one associates with the 
development of a natural system. It is more appropriate to concentrate upon 
the modeling of such phenomena within the confines of future more specific 
manuals rather then in this introductory one.\par
It is instructive to return to the thoughts of Newton. When Newton modeled his 
concept of ``instantaneous velocity'' he allowed the scalar velocity $v$ to 
behavior over ``infinitely small'' time periods as if it was constant.
Thus$,$ if $x(t)$ represents the scalar distance and the derivative exists at 
$t$$,$ then for each $\eps \in 
\dmonad 0$ it would follow from Newton's viewpoint that $x(t +\eps) = 
x(t) + v\, \eps.$ Or that $v = (x(t+ \eps) - x(t))/\eps$ - the Galilean 
average velocity. Therefore$,$ $\st v = \
st {(x(t+ \eps) - x(t))/\eps} = 
x^\prime(t) =v.$ But$,$ 
it also holds that if $v \approx (x(t+ \eps) - x(t))/\eps,$ then $v(t) = 
x^\prime(t).$ The same analysis applies to the partial derivative$,$ the 
directional derivative and expression (2) of section 10.1 since $\hyper A = 
A,$ the standard part operator distributes over matrix multiplication$,$ and
each such $\vec b/\Vert \vec b \Vert$ is limited. Also$,$ observe that if 
the velocity is considered to be a vector quantity$,$ $\vec v,$ that is nonzero 
and continuous at $t,$ then $\vec v(t + \eps)$ is almost parallel to $\vec 
v(t).$ {\bf Almost parallel vectors have the same behavior in the NSP-world as 
parallel vectors have in the N-world.} 
For our basic illustrations$,$ our functions are restricted to maps from $\realp 
n$ into $\real.$ \par
(i) One of the primary  principles employed when determining simplistic 
behavior within the monadic environment is to consider$,$ in the natural world$,$ 
how constant rates of change (i.e. linearly varying quantities) or other 
constant quantities are related$,$ and 
then to pass such relations  over to the 
infinitesimal world. When this is properly done one may discover that the end 
result is simply the replacing of $=$ with $\approx,$ parallel vector 
properties replaced by almost parallel$,$ etc. \par
(ii) The {\bf general} physical laws for such simplistic behavior \underbar{may
be} passed over directly to the NSP-world by *-transfer. In certain  
cases$,$ *-transfer will not be sufficient since the basic NSP-world 
law may need to be written in terms of external notation.  \par
Step (ii) may seem slightly vague; but$,$ this is necessary since we cannot 
experiment within the NSP-world monadic environment and determine what the 
actual relation might be. This knowledge can only be indirectly obtained after 
the differential equation is derived$,$ solved and used as a predictor of natural 
system behavior.\par
(iii) Due to the operational restrictions of infinitesimals and limited numbers to 
ring properties$,$ it is often necessary to represent a transferred natural 
law in a manner different from the customary form. \par
In the next section$,$ an attempt is made to use these three general 
observations to derive the partial differential equation for a vibrating 
membrane$,$ with the obvious extension$,$ and the equation for 3-dimensional heat 
transfer. The reader can determine whether or not these derivations are 
more rigorous and more convincing than those that appear in our present day 
texts.\par 
\bigskip
\leftline{10.3. \underbar{Vibrating Membrane}}\par 
\medskip 
Elementary nonrelativistic dynamical problems are closely associated with Newton's second law of 
motion. Little concern is given to whether or not this law is expressed by 
such an equation as $\vec F = m\, \vec a$ or by $\vec F/m = \vec a.$ However$,$ 
when this law is passed over to the monadic environment it may be necessary$,$ due to 
the ring nature of the infinitesimals and limited numbers$,$ to be more selective 
as to the specific expression utilized. [This is discussed more fully in Section 
10.5.]\par
\vskip 18pt
\hrule
\smallskip
{\bf Definition 10.3.1. ({Monadic Second Law of Motion})} Let internal $\vec v \colon 
A \to \hyperrealp m,$ where $A \subset {\hyperreal}^{n+1},$ ($\vec v(x_1,\ldots,x_n,t)
= \vec v(t)),$ denote an internal function representing an 
internal velocity vector. Consider another internal $\vec F \colon 
A \to \hyperrealp m,$ ($\vec F(x_1,
\ldots,x_n,t)= \vec F(t)$)$,$ $t^\prime \in \hyperreal,\ dt \in \dmonad 0,\ 
0\not = m \in \hyperreal,$ and $t^\prime + dt \in A.$ Then
$${{\vec F(t^\prime)}\over{m}} \approx {{\vec v(t^\prime +  dt)- \vec 
v(t^\prime)}\over {dt}}.\eqno (1)$$\par
\smallskip
\hrule
\vskip 18pt 
The usual constraints placed upon a vibrating membrane such as (i) it is 
perfectly elastic$,$ (ii) it is  attached along its entire boundary to a plane$,$ 
(iii) has a ``very small'' deflection (compared to its size) and only in a direction 
perpendicular to the plane$,$ and (iv) it is vibrating about an equilibrium 
position$,$ can be successfully idealized 
in the following manner.\par
At a moment of time$,$ let $A,\, B,\, C,\, D$ be four points in the interior 
of the membrane such that the configuration $A,\, B,\, C,\, D$ forms a 
rectangle. Select the Cartesian coordinate system such that $\overline{AD}$ is 
parallel to the $x$-axis$,$ $\overline{AB}$ is parallel to the $y$-axis and the 
rectangle lies in the $xy$-plane. The line segments $\overline{AB},\ 
\overline{BC},\ \overline{DC},\ \overline{AD}$ are contained in secants to 
membrane surface curves. The tensions  that yield membrane motion are 
measured in terms of forces per unit length and separated into four nonzero constant 
vectors$,$ $\vec F_1,\ \vec F_2,\ \vec F_3,\ \vec F_4.$ 
These tensions are assumed to be exterior 
in nature$,$ $\vec F_1$ is applied to $\overline{AD}$ in the $-z$ direction, $\vec 
F_2$ 
is applied to $\overline{BC}$ in the $z$ direction,  $\vec F_3$ is applied to 
 $\overline{AB}$ in the $-z$ direction, and $\vec F_4$ is applied to $\overline{DC}$ in the $z$ direction. Since perfect nontearing  
vibration occurs then two 
forces have been selected to counteract the effect of the remaining two.\par   
\vskip 18pt
\hrule
\smallskip
\centerline{\bf {Deformable Body Rule}}\par
One of the basic concepts used to model infinitesimally the behavior of a
deformable body is that in the monadic environment its physical elements  
exhibit the same dynamic behavior as that of a rigid body.\par
\smallskip
\hrule
\vskip 18pt
The scalar force $F$ that will accelerate the parallelogram in the 
$\hat k$ direction$,$ when it is above the equilibrium position$,$ is 
$$F = -\Vert \vec F_1 \Vert \vert \overline{AD} \vert + 
\Vert \vec F_2 \Vert \vert \overline{AD}\vert - \Vert \vec F_3 \Vert \vert \overline{AB} \vert + 
\Vert \vec F_4 \Vert \vert \overline{AB}\vert.\eqno (2) $$
\par
\vskip 18pt
%==========TABLE STRUTS

%24 POINT DEEP ROWS
\newbox\medstrutbox
\setbox\medstrutbox=\hbox{\vrule height14.5pt depth9.5pt width0pt}
\def\medstrut{\relax\ifmmode\copy\medstrutbox\else\unhcopy\medstrutbox\fi}

\bigskip

\newdimen\leftmargin
\leftmargin=0.0truein
%for textbook use 1.1truein
\newdimen\widesize
\widesize=3.0truein
\advance\widesize by \leftmargin       
%\moveleft\leftmargin
\hfil\vbox{\tabskip=0pt\offinterlineskip
\def\tablerule{\noalign{\hrule}}

\halign to \widesize{\medstrut\vrule#\tabskip=0pt plus2truein

&\hfil\quad#\quad\hfil&\vrule#

\tabskip=0pt\cr\tablerule

%&\multispan3 \hfil {\tablerm TABLE TITLE} \hfil&\cr\tablerule
%\biggerstrut&ENTRY A&&ENTRY B&\cr\tablerule

&$\Downarrow$ IMPORTANT $\Downarrow$&\cr\tablerule
}} \par
\bigskip
$\Rightarrow$ The following derivations of some significant partial 
differential equation models for natural system behavior are stated in a
protracted mode. Within these derivations$,$ I have presented certain  general 
modeling concepts$,$ gleaned from indirect evidence associated with simplistic 
NSP-world behavior. The physical behavior of the NSP-world relative to each 
particular N-world event must be individually investigated and is usually 
characterized as the intuitive *-transfer of the simplest N-world 
behavior. In classical applications to continuum models$,$ this simplistic 
behavior is often conceived of as intuitively discontinuous and piecemeal 
in character. When this is the case$,$ the entire complex N-world effect is but 
a filtered composition of an infinite replication of what occurs within a 
single monadic environment. It is intuitively disjoint in character since the
composition can be characterized as a ``jumping'' from one monad to another  
rather than some type of  continuous joining of the effects. On the other 
hand$,$ as recently shown [{\tt Herrmann [1989]}]$,$ if we assume that certain N-world 
behavior is fractal or even finitely discontinuous in nature$,$ 
then the NSP-world behavior that would yield such effects may be viewed as
ultrasmooth and ultracontinuous in character. A representation for this 
ultrasmooth NSP-world 
behavior would necessary need to satisfy a differential equation expressed in 
terms of the hyperreals. However$,$ by considering a 
nonstandard model of a nonstandard model it might be possible to show that  
such NSP-world differential equation models are once again  a composition of 
simplistic NSP-world behavior but on a much deeper level.$\Leftarrow$ \par    
\vskip 18pt
{\bf Application 10.3.1.} {\it A differential equation model for a vibrating 
membrane.} \par
All functions are assumed to be defined on an open neighborhood$,$ $G,$ of 
$(a,b,c,t_0)$ where$,$ it appears to be necessary$,$ to assume that $G \cap \real^3$ contains the closed membrane 
 and$,$ at the least$,$ the functions are continuous where defined. Let $\vec T= \vec T(x,y,z)$ represent the internal restoring 
tension per unit length.
Let $z= z(x,y,t)$ represent the configuration of the membrane$,$ as a 
surface$,$ and assume that 
${\partial^2z}/{\partial t^2},\ \partial(\Vert \vec T \Vert(\partial 
z/\partial x))/\partial{x},\
\partial(\Vert \vec T \Vert(\partial z/\partial y))/\partial{y}$ exist  
and that 
$z_x,\ z_y$ are continuous at $(a,b,t_0).$ Let $\rho= 
\rho(x,y,)$ represent the mass per unit area. Finally$,$ 
let $h= h(x,y,z,t)$ represent the per unit area lumped 
load function that yields the $\hat k$ direction  external forces that may 
be applied to the membrane. For any such $(x,y,z,t) = (a,b,c,t_0),$ a 
partial differential equation model for the 
vibrating membrane configuration $z$ is
$$(\rho(x,y))\left({{\partial ^2z}\over{\partial t^2}}\right)= 
{{\partial(\Vert \vec T \Vert(\partial z/\partial x))}\over{\partial{x}}}
+ {{\partial(\Vert \vec T \Vert(\partial z/\partial y))}\over{\partial{y}}}
+ h(x,y,z,t).\eqno (3)$$ \par
{\sl Derivation.} We now infinitesimalize the behavior of an infinitesimal 
portion of the membrane assuming that it will have the same effect in the 
N-world as an infinitesimal vibrating parallelogram. The key to the basic 
modeling technique lies in modeling intuitive statements (iii) and (iv).  
{\bf (A) Expressions that involve the physical notion that something is ``small'' 
compared with something else tend to be the basis for differential equation 
modeling.} First$,$ since (iv) states that this ``small'' deflection (i.e. the 
value of $z$ at $(a,b,t_0)$)  takes place at equilibrium  then this 
should  entail that the tension that produces such a deflection is also
``small'' in the $\hat k$ direction. Now select the plane of attachment for the membrane as the 
$xy$-plane and let the plane represent the equilibrium position. 
By physical observation$,$ as the deflection is made ``smaller'' and 
``smaller'' then the $ABCD$ rectangle discussed above would have its 
normal more nearly parallel to $\hat k.$ The concept of ``small'' deflection 
is now related to the time and the tension. Again by physical observation$,$ it 
can be assumed that to obtain a ``small'' deflection near the $xy$-plane all 
one needs to consider is that the tension itself is ``small'' in the $\hat k$ 
direction and is applied 
for but a ``small'' period of time $\Delta t.$ \par
The above concept of ``smallness'' is now embedded into the monadic world.
The infinitesimal parallelogram representation for the membrane surface is 
assumed to move parallel to the $z$-axis. Thus let $A = (a,b,c), \ 
B= (a,b +dy,c),\ C= (a+dx,b+dy,c),\ D=(a+dx,b,c),$ where positive 
$dx,\ dy \in \monad 0.$ Considering only the $\hat k$ component of the 
tension $\vec T,$ we write 
this component  as $\Vert \vec T \Vert \cos \gamma,$ where $\gamma$ is a 
function continuous on $G_1.$ By continuity 
$\Vert \Hyper {\vec T}(x,y,z,)\Vert\hyper 
\cos {\hyper \gamma (x,y,z,)}\approx \Vert \Hyper {\vec T}(a,b,c,)\Vert\hyper 
\cos {\hyper \gamma (a,b,c,)}$ for each $(x,y,z,)  \in \nmonad {(a,b,c)},\ 
n = 3.$\par
 For the $\hat k$ component$,$ what infinitesimal property will characterize 
the standard modeling 
notion that the deflection must be ``small'' compared to the size of the 
membrane? {\bf (B) Certainly$,$ the smallest possible N-world deflection would occur if 
there were  NO apparent standard deflection at all. But$,$ to analyze the problem  
nontrivially it would need to be assumed that deflection did take place in 
the monadic world.}\par
\vskip 18pt
\hrule
\smallskip
\centerline{{Monadic Restrictions of Physical Processes}}
The idea that a physical process may be restricted to a monadic environment in 
such a way that it does not appear to occur in the standard world but does 
occur in the monadic environment is yet 
another significant infinitesimal modeling technique.\par
\smallskip
\hrule
\vskip 18pt
In order to apply the monadic restriction technique to this problem and without actually evaluating the functions$,$ simply let 
$$\st {\Vert \Hyper {\vec T}(a,b,c)\Vert\hyper 
\cos {\hyper \gamma (a,b,c)}} = 0, \ \vec T(a,b,c) \not= \vec 0.\eqno 
(4)$$\par
{\bf (C) Notice the important modeling technique that from the standard world 
point of view the 
infinitesimal rectangle $ABCD$ appears as if it is the single point $(a,b,c).$}
From N-world physical observation$,$ and indeed the definition of the 
equilibrium position$,$ a single point will vibrate 
only if it is displaced either above or below the 
equilibrium position. By *-transferring this observation to the infinitesimal 
rectangle $ABCD,$  it follows that in order to have a ``small'' 
deflection take place ``near'' the equilibrium position 
$ABCD$ should be placed either above or 
below the 
hyper-$xy$-plane and all members of $ABCD$ should be infinitely close 
to the hyper-$xy$-plane. Let $c = z(a,b,t_0).$\par
The hyperline segments $\overline{AB}$ and $\overline{AD}$ are geometric 
elements for membrane (surface) curves. This fact and our definition (9.5.1) for 
tangents to space curves leads to the conclusion that the hyperslope of 
$\overline{AD} \approx z_x(a,b,t_0)$ and hyperslope of $\overline{AB} 
\approx z_y(a,b,t_0).$ Since every point in $ABCD$ is infinitely close to the 
hyper-$xy$-plane$,$ then slope $\overline{AB}$ and $\overline{AD}$ are both 
infinitely close to 0. Observe that since by continuity there 
exists some open ball about $(a,b,c)$ such that ${\vec T} \not= {\vec 0}$ for any 
member of this ball then $\Hyper {\vec T}(x,y,z,) \not= \vec 0$ for any
$(x,y,z,) \in \nmonad {(a,b,c)}.$ Expression (4) now implies that
$$\hyper \cos {\hyper \gamma (a,b,c)} \approx \hyper \cos {\hyper \gamma 
(a+dx,b,c)}\approx \hyper \cos {\hyper \gamma (a,b+dy,c)} \approx 0.\eqno 
(5)$$
Further$,$ using the continuity of $z_x,\ z_y$ we obtain  
$$z_x(a,b,t_0) \approx z_x(a+dx,b,t_0) \approx z_y(a,b,t_0) \approx 
z_y(a,b+dy,t_0) \approx 0.\eqno (6)$$\par
The problem we now face is how to select the proper combination of entities 
from (5) and (6). The following expressions (7)$,$ (8)$,$ (9) and (10) are obtained 
by means of the observation that$,$ {\bf (D) mathematically$,$ the coordinate spaces are 
independent one from the other and$,$ physically$,$ tensions parallel to the 
$z$-axis along one boundary of the hyperrectangle do not effect 
substantially the adjacent boundaries.} Since $\Hyper {\vec T}$ is a limited vector 
for each $(x,y,z,)\in \nmonad {(a,b,c)}$ then 
$$\Vert \Hyper {\vec T}(a,b,c)\Vert \hyper \cos {\hyper \gamma (a,b,c)}\approx 
\Vert \Hyper {\vec T}(a,b,c)\Vert z_x(a,b,t_0),\eqno (7)$$ 
$$\Vert \Hyper {\vec T}(a,b,c)\Vert \hyper \cos {\hyper \gamma (a,b,c)}\approx 
\Vert \Hyper {\vec T}(a,b,c)\Vert z_y(a,b,t_0),\eqno (8)$$ 
$$\Vert \Hyper {\vec T}(a+dx,b,c)\Vert \hyper \cos {\hyper \gamma 
(a+dx,b,c)}\approx $$
$$\Vert \Hyper {\vec T}(a+dx,b,c)\Vert\hyper z_x(a + dx,b,t_0),\eqno (9)$$
$$\Vert \Hyper {\vec T}(a,b+dy,c) \Vert \hyper \cos {\hyper \gamma 
(a,b+dy,c)}\approx $$
$$\Vert \Hyper {\vec T}(a,b+dy,c) \Vert\hyper z_y(a,b+dy,t_0).\eqno (10)$$\par 
Expression (2) must now be infinitesimalized. The major method used for such 
infinitesimalizing is call {\bf (E) {the (Differential Equation) Method of 
Maximum and Minimum}. In this case$,$ it leads to the conclusion that forces may be considered 
as attached to the $A,\ B,\ D$ vertices of this physical (geometric) element.}
Before proceeding$,$ might there be a purely  physical NSP-world
argument for selecting the values for $T$ at the points 
$A,\ B,\ D$ - an argument 
that is substantially different from the usual one that such a selection is 
done since it ``works''? Recently {\tt Simhony [1987]} has advanced a theory that 
claims that all of the material universe is composed of combinations of the 
configuration he has termed as the {electron-positron lattice}. This is a simple 
{lattice structure} with electrons and positrons located at the lattice nodes.
For such a theory$,$ these objects would be the building blocks of the 
material universe and all forces could be considered as attached to the 
lattice nodes. Our infinitesimal representation $ABCD$ is uniquely determined 
by the vertices $A,\ B,$ and $D.$ Even though Simhony's investigations of the 
possible lattice nature of the material universe 
have interesting conclusions that could by *-transfer lead to an effective 
lattice structure within  the NSP-world$,$ such a structure can also evolve 
from mathematical considerations.\par
The function $\Vert {\vec T}(x,y,z)\Vert  
\cos {\gamma (x,y,z)}$ is continuous on every standard rectangle determined 
by vertices $A = (a,b,c),\ B_1 = (a,b+b_1,c),\ C_1 = (a+a_1,b+b_1,c),\ D_1= 
(a+a_1,b,c),$ where $0 < a_1 \leq r_1,\ 0 < b_1 \leq r_2$ for some $r_1,\, r_2 \in 
\real.$ Thus$,$ since the line segment $\overline{AB_1}$ is compact$,$  
then there exists points $(a,b_m,c),\ (a,b_M,c) \in \overline{AB_1}$ 
such that $$\Vert {\vec T}(a,b_m,c)\Vert  
\cos {\gamma (a,b_m,c)}\leq $$
$$\Vert {\vec T}(a,y,c)\Vert  
\cos {\gamma (a,y,c)}\leq \Vert {\vec T}(a,b_M,c)\Vert  \cos 
{\gamma (a,b_M,c)}\eqno (11)$$
 for each $(a,y,c) \in \overline{AB_1}.$ 
Hence$,$ $$\Vert {\vec T}(a,b_m,c)\Vert  
\cos {\gamma (a,b_m,c)}\,\vert\overline{AB_1}\vert\leq $$
$$\Vert {\vec T}(a,y,c)\Vert  
\cos {\gamma (a,y,c)}\,\vert\overline{AB_1}\vert\leq $$ $$\Vert {\vec T}(a,b_M,c)\Vert  
\cos {\gamma (a,b_M,c)}\,\vert\overline{AB_1}\vert,\eqno 
(12)$$
for each $(a,y,c) \in \overline{AB_1}.$                                                                                     
Whether it be obtained by observation$,$ consideration of the integral or by axiomation$,$ method 
(E) states$,$ for this application$,$ that there exists a net scalar tension  
$F^\prime$ such 
that 
$$\Vert {\vec T}(a,b_m,c)\Vert  
\cos {\gamma (a,b_m,c)}\,\vert\overline{AB_1}\vert\leq $$ 
$$F^\prime
\leq \Vert {\vec T}(a,b_M,c)\Vert  
\cos {\gamma (a,b_M,c)}\,\vert\overline{AB_1}\vert\eqno 
(13)$$                                                                                      
and $F^\prime$ may be considered the force that accelerates the rigid body 
$\overline{AB_1}.$ 
Assuming continuity on connected $\overline{AB_1}$ then there exists some $(a,\overline b,c)$ such that
$$\Vert {\vec T}(a,\overline b,c)\Vert  
\cos {\gamma (a,\overline b,c)}\,\vert\overline{AB_1}\vert 
= F^\prime,\ (a,\overline b,c)\in \overline {AB_1}\eqno (14)$$
and in like 
manner for the other three tensions. Following (2) we have that the total 
force of motion may be considered the sum of these four scalar $\pm \hat k$ 
directed tensions. The results in this paragraph are now *-transferred to the 
NSP-world$,$ in  which case (14) becomes
$$\Vert \Hyper {\vec T}(a,\overline b,c)\Vert  
\hyper \cos {\hyper \gamma (a,\overline b,c)}\,dy 
= F^\prime,\eqno (15)$$
where $b \leq \overline b \leq b + dy.$\par
Continuity at $(a,b,c)$ implies that $\Vert \Hyper {\vec T}( a,\overline b,c)\Vert  
\hyper \cos {\hyper \gamma (a,\overline b,c)} \approx \Vert \Hyper {\vec 
T}(x,y,c)\Vert  
\hyper \cos {\hyper \gamma (x,y,c)},$ where $(x,y,c) \in 
\{(a,b,c),(a+dx,b,c),a,b+dy,c),(a+dx,b+dy,c)\}= A.$ Now selection of the  
appropriate members of 
$A$ comes from repeating the above argument and *-transfer of the behavior. Applying (15) and (2) this implies that the effective scalar 
tension$,$ $\hyper F_T,$ that moves the hyperrectangle $ABCD$
has the property that  \par
\centerline{$\hyper F_T \approx\Vert \Hyper {\vec T}(a+dx,b,c)\Vert\hyper \cos {\hyper 
{\gamma (a + dx,b,c)}}\,dy
-$}
$$ \Vert \Hyper {\vec T}(a,b,c)\Vert \hyper \cos {\hyper {\gamma (a,b,c)}}\,dy
 +$$ $$\Vert \Hyper {\vec T}(a,b+dy,c) \Vert\hyper \cos {\hyper {\gamma 
(a,b+dy,c)}}\,dx - $$ $$
\Vert \Hyper {\vec T}(a,b,c)\Vert \hyper \cos {\hyper {\gamma (a,b,c)}}\, 
dx.\eqno (16)$$
Applying (7)$,$ (8)$,$ (9) and (10) yields
$$\hyper F_T \approx \Vert \Hyper {\vec T}(a+dx,b,c)\Vert\hyper z_x(a + 
dx,b,t_0)\,dy -$$
$$ \Vert \Hyper {\vec T}(a,b,c)\Vert \hyper z_x(a,b,t_0)\,dy
 +$$
$$\Vert \Hyper {\vec T}(a,b+dy,c) \Vert\hyper z_y(a,b+dy,t_0)\,dx - 
\Vert \Hyper {\vec T}(a,b,c)\Vert \hyper z_y(a,b,t_0)\, dx.
 \eqno (17)$$
Applying (E) again and by *-transfer$,$ we note that there exists points 
$(\overline a,\overline b),\ (c,d,f,g)$ such that the mass$,$ $m,$ of the hyperrectangle 
is $\hyper \rho (\overline a,\overline b)\,dx\,dy$ and the entire lumped load $L$ can be represented
by $\Hyper h(c,d,f,g)\,dx\,dy,$ where $(\overline a,\overline b)\approx (a,b)$ and $(c,d,f,g) 
\approx (a,b,c,t_0).$ By continuity$,$ $\hyper \rho (\overline 
a,\overline b)\,dx\,dy \approx
\rho (a,b)\,dx\,dy;\ \Hyper h(c,d,f,g)\,dx\,dy =\Hyper h\,dx\,dy\approx h(a,b,c,t_0)\, dx\, 
dy.$ Writing expression (1) in a scalar form for motion parallel to the 
$z$-axis obtains $\hyper F(t_0)/m \approx (\hyper
v(t_0 + dt) - v(t_0))/dt$ and the total scalar force producing a change in 
the instantaneous velocity is $\hyper F_T + \Hyper h(c,d,f,g).$ Hence,
$${{\hyper v(t_0 + dt) - v(t_0)}\over{dt}}\approx$$
$${{\Vert \Hyper {\vec T}(a+dx,b,c)\Vert\hyper z_x(a + dx,b,t_0)\,dy
- \Vert \Hyper {\vec T}(a,b,c)\Vert z_x(a,b,t_0)\,dy}\over{\hyper \rho (\overline 
a,\overline b)\,dx\,dy}} +$$
$${{\Vert \Hyper {\vec T}(a,b+dy,c) \Vert\hyper z_y(a,b+dy,t_0)\,dx}\over{\hyper \rho 
(\overline a,\overline b)\,dx\,dy}} - $$ $$
{{\Vert \Hyper {\vec T}(a,b,c)\Vert z_y(a,b,t_0)\, dx +
\Hyper h\,dx\,dy}\over{\hyper \rho 
(\overline a,\overline b)\,dx\,dy}}=$$
$${{\Vert \Hyper {\vec T}(a+dx,b,c)\Vert\hyper z_x(a + dx,b,t_0)
- \Vert \Hyper {\vec T}(a,b,c)\Vert z_x(a,b,t_0)}\over{\hyper \rho (\overline 
a,\overline b)\,dx}} +$$
$${{\Vert \Hyper {\vec T}(a,b+dy,c) \Vert\hyper z_y(a,b+dy,t_0) - 
\Vert \Hyper {\vec T}(a,b,c)\Vert z_y(a,b,t_0)}\over{\hyper \rho 
(\overline a,\overline b)\,dy}}+$$ 
$$ {{\Hyper h}\over{\hyper \rho (\overline a,\overline b)}}\approx
{{\hyper z_t(t_0 + dt) - z_t(t_0)}\over{dt}}.\eqno (18)$$ \par
The derivation is completed by applying the standard part operator to 
expression (18) and employing the hypothesized requirement that the various 
partials exist.\par
\vskip 18pt
The following are some useful observations relative to the above derivation.
If one assumes that the tension $\vec T$ and the density $\rho$ are constant 
along with $h\equiv 0,$ then the usual elementary 2-dimensional wave equation 
follows from (3). More to the point$,$ however$,$ is the$,$ not necessarily obvious,
fact that the derivation was not obtained by attempting to {comprehend second 
order rates of change.} Berkeley was correct in his criticism that such things 
are difficult to mentally perceive.\break
 {\bf (F) The general overall approach in the 
derivation the vibrating membrane equation was to consider first-order 
rates of change of other rates of change that had already been investigated 
and found to be representable by infinitely close (locally linear) functions 
expressed in partial derivative form.} This general overall approach is almost 
always the best and least confusing. \par
\vfil\eject
\leftline{10.4. \underbar{Internal Heat  Transfer}}\par 
Consider a vertical (rectangular) slab of {material that uniformly 
conducts heat} where ``heat flow'' is perpendicular to the left and right- hand surfaces and ``flows'' to the right. We quote from {\tt Sears} and {\tt Zemansky [1952]}: {\it
``The figure represents a slab of material of cross section $A$ and thickness
$L.$ Let the whole of the left face of the slab be kept at a temperature
$U_2,$ and the whole of the right face at a lower temperature $U_1.$ The 
direction of heat current is then from left to right through the slab.\par
After the faces of the slab have been kept at the temperatures $U_1$ and 
$U_2$ for a sufficient length of time$,$  the temperature at points within the 
slab is found to decrease uniformly {\rm [linearly]} with distance from the hot to the cold 
face. At each point$,$ however$,$ the temperature remains constant with time. The 
slab is said to be in a ``linear steady state.''\par
It is found by experiment that the rate of flow of heat through the slab in 
the {steady state} is proportional to the area $A,$ proportional to the 
temperature difference $(U_2 - U_1),$ and inversely proportional to the 
thickness $L.$''}\par
Let $Q$ be the usual measure of the {quantity of heat}$,$ $0 < \Delta t\in \real$ 
the time$,$  
$\not= 0< K\in \real$ the constant of proportionality (the {\it thermal 
conductivity}) and $0 < A \in \real$ the 
area.  
Then for this very special and simplistic 
case$,$ it follows that
$$Q = {{KA(\Delta t)(U_2- U_1)}\over{L}}. \eqno (1)$$   
Notice that since this is a rectangular solid and the distance $L$ is measured 
perpendicular to the faces then the convention of heat flow is taken normal to 
the faces.\par
Before the modern theory of {heat conduction}$,$ heat was thought of as a 
invisible weightless fluid called {\it {caloric}} which was produced when a 
substance burned and which could be transmitted by conduction from one body to 
another. The analogy of ``{flowing heat}'' has persisted and due to the 
above steady state law this direction is taken normal to this section. 
For the moment$,$ assume that the temperature $U(a,y,z,t_0)$ is only dependent 
upon time and thickness$,$ $a,$ of the slab. Now alter the position of the 
coordinate system relative to the rectangular slab and 
 let $a + L,\  L > 0,$ 
represent the thickness of the of the slab and $\Delta t = t - t_0.$ 
Our first natural law is now expressed  as
$$Q = {{KA(\Delta t)(U(a,y,z,t_0) -  U(a + L,y,z,t_0))}\over{L}}. \eqno 
(2)$$\par
Let the specific heat$,$ $\sigma,$ and the density$,$ $\rho,$ for a given 
``small'' rectangular solid$,$ $R,$ of volume$,$ $V,$ be constants. At time $t_0$ 
let the temperature of the entire solid be a constant $U(t_0).$ After a time 
change $\Delta t$ and if there is but a ``small'' positive temperature change 
$U(x,y,z,t_0 + \Delta t)$ which appears also to be constant throughout $R,$ then a 
second 
physical law appears to emerge. This law states that the amount of heat$,$ $H,$ 
necessary to achieve such a temperature change is  
 $$H = \rho\,\sigma\,V\,(U(x,y,z,t_0 + \Delta t) - U(x,y,z,t_0)),\eqno (3)$$
where$,$ for this simplistic case$,$ $U$ is not dependent upon position.\par
We need two more requirements prior to our derivation of the heat equation.
First$,$ recall that Count Rumford suggested that heat was really energy in another 
form and Sir James Joule experimentally verified this natural law. In what 
follows$,$ we use the concepts of the conservation and additivity of heat energy. 
The second requirement is Theorem 10.4.1. For nonempty open $G \subset \realp 
n,$ let $f\colon G \to \real.$ Denote any $\vec v \in G$ by 
$\vec v = (x_1,\ldots,x_n).$ Let $y = x_i,\ 1\leq i \leq n.$ 
Recall that a set $C_y$ 
is said to be 
{\bf {convex in the direction $y$}} if for any two points $\vec p=(p_1,\ldots,p_n), \ 
\vec q = (p_1,\ldots,p_i + h_i,\ldots,p_n) \in C_y$ and any $k \in \real$ 
such that $\vec w=(p_1,\ldots,p_i + k,\ldots,p_n)$ is a member of the line segment 
with end point $\vec p,\ \vec q,$ then $\vec w \in C_y.$
\vskip 18pt 
{\bf Theorem 10.4.1.} {\sl Let $f\colon G \to \real,$ where 
nonempty open $G\subset \realp n$ and standard $\vec v \in G.$ 
Suppose that $y = x_i,\ 1\leq i\leq n,\ f_{y}$ is defined on $G$ 
and continuous at $\vec v.$ 
Let $C_{y} \subset \Hyper G,$  where internal $C_{y}$ is 
*-convex in the direction $y.$ If for $h \in \monad 0,$ such that 
$\vec p = (p_1,\ldots,
p_i,\ldots,p_n),\vec q = (p_1,\ldots, p_i + h,\ldots,p_n)\in C_{y}$ and 
$\vec v \approx \vec p,$ then 
there exists $\eps \in \monad 0$ such that
$$\hyper f({\vec q}\,) - \hyper f({\vec p}\,) = f_{y}({\vec v}\,)\, 
h + h\, \eps.$$} \par
{\bf Application 10.4.1} {\it A differential equation model for 3-dimensional 
heat transfer.} \par
All functions are assumed to be defined on a open neighborhood$,$ G$,$ of $(a,b,c,t_0)$ where$,$ it appears necessary$,$ to assume that $G \cap \real^3$ contains the closed solid and$,$ at the least$,$ the functions are continuous where defined. 
Let $U(x,y,z,t)$ represent the {temperature}$,$ never zero $\rho (x,y,z,)$ 
the density 
of the material$,$ never zero $\sigma (x,y,z,)$ the {specific heat} and a 
constant {thermal conductivity }$K.$
 It is assumed that $\sigma, \rho$ do not vary with 
temperature. This constraint implies that our differential equation model is 
only valid for small temperature changes. Even though no sources or sink are 
assumed for this application, they can be easily adjoined to the derivation in a 
manner analogous to the lumped load of application 10.3.1. Let $U_x,\, U_y, \,U_z,$ be continuous on $G$ and $U_t$ be 
defined on $G$ and continuous at $(a,b,c,t_0).$ Let
$U_{xx},\, U_{yy},\, U_{zz},$ exist at $(a,b,c,t_0).$ 
For any such $(x,y,z,t) = (a,b,c,t_0)$ a differential 
equation 
model for internal heat transfer is
$$\rho(a,b,c)\,\sigma(a,b,c){{\partial U}\over{\partial t}}= 
K\left({{\partial^2U}\over{\partial{x^2}}}
+ {{\partial^2U}\over{\partial{y^2}}} +
{{\partial^2U}\over{\partial{z^2}}}\right).$$ \par              
{\sl Derivation.}
Equation (2) assumes that we are dealing with a uniform temperature change in 
that we have a steady state condition. Still retaining the idea that 
the temperature is constant on a section$,$ (2) must be generalized to the case
where the temperature is not uniformly decreasing with the thickness of the slab.
We make a strong appeal to the flow analogy and$,$ thus$,$ consider the temperature 
evaluation to take place along a flow line parallel to $\hat i.$ With this in 
mind$,$ *-transfer yields      
$$\hyper Q ={{K\hyper A(\Delta t)(U(a,b,c,t_0) 
 -  \Hyper U(a + \eps,b,c,t_0))}\over{\eps}}, \eqno (4)$$
where $0 <L = \eps \in \monad 0.$ Can we consider the thermal conductivity as  
a nonconstant function? It appears that if $K$ were not constant$,$ then the 
differential equation for heat transfer would not be derivable. Fortunately$,$ 
in most physical cases$,$ the constancy of $K$ can be verified. 
 Assume that $\partial U/\partial x$ 
exits at $(a,b,c,t_0).$ Then (4) leads to the nonsteady state expression
$$Q=\st {\hyper Q} = \St {\left({{K\hyper A(\Delta t)(U(a,b,c,t_0) -  
\Hyper U(a + 
\eps,b,c,t_0))}\over{\eps}}\right)} =$$ $$- KA(\Delta t)\left({{\partial 
U(a,b,c,t_0)}\over {\partial 
x}}\right).\eqno (5)$$\par  
Equation (5) is conceived of as representing the quantity of heat after the 
time period $\Delta t$ that ``flows'' through a rectangular shaped section cut 
from a conducting solid; but$,$ we are still assuming that the temperature is 
constant on each section.
[Remark: {\it It is claimed by some authors that equation (5) expressed in words with
this flow convention has been experimentally obtained. They then use (5)$,$ along 
with integral concepts and the Divergence Theorem to formulate the heat equation. I 
shall only take (2) as fundamental.}]\par
Once again applying the directed flow analogy$,$ the external view of the net 
amount of heat that 
would accumulate in the interior of the slab from this $\hat i$ direction flow
would be 
$$H_{\hat i}= -KA(\Delta t)\left({{\partial U(a,b,c,t_0)}\over {\partial 
x}}\right) -$$ $$\left(-KA(\Delta t)\left({{\partial U(a+L,b,c,t_0)}\over 
{\partial x}}\right)\right).\eqno (6)$$
Now letting $V$ be the volume of the slab notice that
$$H_{\hat i}/(V(\Delta t)) = K(U_x(a+L,b,c,t_0)-U_x(a,b,c,t_0))/L.\eqno 
(7)$$
Unfortunately$,$ (7) still refers to a temperature that is constant on the 
parallel faces of the slab. To generalize to the nonconstant case$,$ 
first infinitesimalize the slab itself. By *-transfer$,$ (7) holds for an 
infinitesimal rectangular solid $R$ having determining vertices $A = (a,b,c),\ B = 
(a + dx,b,c),\ C =(a,b+dy,c),\ D = (a,b,c + dz),\ 0<dx,\, dy,\, dz \in \monad 
0$ and for an infinitesimal time period $0< dt \in \monad 0.$ Hence$,$ (7) 
becomes $$\hyper H_{\hat i}/(dx\,dy\,dz)(dt)) = K(\Hyper U_x(a+dx,b,c,t_0)-
U_x(a,b,c,t_0))/dx.\eqno (8)$$\par                                                   
Next we correct for the specialized constant temperature case by assuming that 
heat flow takes place only along parallel flow lines. Concentrating upon (7)$,$ 
where the temperature is not assumed to be constant over the sections$,$ we 
apply the maximum and minimum  method (E) from the last section and *-transfer 
the result to the faces of $R.$ Thus$,$ there exists some $(b_1,c_1) \in \monad {(b,c)}$ 
such that the total heat per unit volume-time accumulated by the $\hat i$ 
direction flow is
 $$\hyper H_{\hat i}/(dx\,dy\,dz)(dt)) = K(\Hyper 
U_x(a+dx,b_1,c_1,t_0)-\Hyper U_x(a,b_1,c_1,t_0))/dx.\eqno (9)$$\par
  Applying 
Theorem 10.4.1 yields
 $$\hyper H_{\hat i}/(dx\,dy\,dz)(dt)) = 
KU_{xx}(a,b,c,t_0) + \eps_1,\, (\eps_1 \in \monad 0),\eqno (10)$$
and repeating the argument for the 
$\hat j$ and $\hat k$ directions (and applying the *-transfer of the 
additivity of energy etc.) yields that the total heat accumulated per volume-
time for $R,$ as viewed externally$,$ is 
$$\hyper H/(dx\,dy\,dz)(dt) \approx K(U_{xx}(a,b,c,t_0) + 
U_{yy}(a,b,c,t_0) + U_{zz}(a,b,c,t_0)). \eqno (11)$$
Now consider the 
simplistic natural law represented by equation (3). For a more general case$,$ 
where the temperature$,$ density$,$ and even the specific heat vary by position 
(but not temperature) and are continuous on an open neighborhood $G,$ one can 
and apply the maximum and minimum method (E). Viewing $\sigma,\ \rho,\ U$ as a 
functions of position only$,$ where $t_0$ and $\Delta t$ are assumed fixed. Then 
(dependent upon $t_0$ and $\Delta t$) there exists some $\vec p \in {\rm the \ solid}$ such 
that $$H= V\rho ({\vec p}\,) \sigma ({\vec p}\,) (U(\vec p,t_0 + \Delta t)-
U(\vec p,t_0)). \eqno (12)$$                                 Equation (12) 
can be written for nonzero $\Delta t$ in the form $$H/(V\rho ({\vec p}\,) 
\sigma ({\vec p}\,)\Delta t) =  (U(\vec p,t_0 + \Delta t)-U(\vec 
p,t_0))/\Delta t. \eqno (13)$$ By *-transfer$,$ there exists some $\vec q \in 
R$ such that                                $$\hyper H/(dx\,dy\,dz\hyper \rho 
({\vec q}\,) \hyper \sigma ({\vec q}\,)dt) =  (\Hyper U(\vec q,t_0 + 
   dt)-\Hyper U(\vec q,t_0))/dt. \eqno (14)$$ Applying Theorem 10.4.1 and 
the fact that $\hyper \rho ({\vec q}\,)$ and $\hyper \sigma ({\vec q}\,)$ are limited$,$                                 
we obtain                                       
$$\hyper H/(dx\,dy\,dz)(dt) \approx 
\hyper \rho ({\vec q}\,) \hyper \sigma ({\vec q}\,)\,U_t(a,b,c,t_0). \eqno (15)$$\par
The final step is clear. From expressions (11) and (15) we have
$$\hyper \rho ({\vec q}\,) \hyper \sigma ({\vec q}\,)\,U_t(a,b,c,t_0) \approx
K(U_{xx}(a,b,c,t_0) + $$ $$U_{yy}(a,b,c,t_0) + 
U_{zz}(a,b,c,t_0))\eqno (16)$$ and application of the standard part operator coupled 
with the continuity of $\rho$ and $\sigma$ at $(a,b,c)$ implies 
the result that
$$\rho(a,b,c)\,\sigma(a,b,c){{\partial U}\over{\partial t}}= 
K\left({{\partial^2U}\over{\partial{x^2}}}
+ {{\partial^2U}\over{\partial{y^2}}} +
{{\partial^2U}\over{\partial{z^2}}}\right).\eqno (17)$$\par 
\leftline{10.5. \underbar{Concluding Remarks}}\par 
\medskip
For some derivations$,$ an {alternative to the Monadic Second Law of Motion} may 
be necessary due to the ring nature of the infinitesimals and the limited 
numbers. Let us compare the Newtonian view discussed in section 10.1$,$ 
Definition 10.3.1 and the fundamental expression (1) of section 10.4.
If one lets the scalar velocity $v(t)$ linearly chang with respect to a 
``small'' time change $\Delta t$
(i.e. uniformly changing in the old terminology)$,$ then observation does indeed 
confirm that the constant scalar force that produces such a change is 
$F(\Delta t) = m\, ((v(t_0 +\Delta t)- v(t))/\Delta t_0).$ If this last statement is 
assumed to hold for all nonnegative time changes less than or equal to $\Delta t$
then$,$ as was done in order to obtain equation  (4) of section 10.4 $,$ 
*-transfer yields
$$\hyper F(\eps) = \hyper m \, ((\hyper v(t_0 +\eps)- v(t))/\eps),\,(\eps \in 
\monad 0).\eqno (1)$$
From (1) and application of the standard part operator Newton's Second Law of 
Motion follows. If one removes the standard part operator$,$ then all that can 
be stated is that 
$$\hyper F(\eps) \approx \hyper m \, ((\hyper v(t_0 +\eps)- v(t))/\eps),\,(\eps \in 
\monad 0).\eqno (2)$$
Equation (1) is the alternate NSP-world form of Newton's Second Law of Motion.
It may be necessary to utilize this alternate form for those cases were 
division by infinitesimals and other similar algebraic processes are used 
- processes that do not preserve $\approx.$ Further$,$ the use of (E) and 
Theorem 10.4.1 are useful in order to retain expressions that do not include 
the external notion of $\approx.$ For the derivation of the heat 
equation it was necessary to retain equational expressions until steps (15) 
and (16).\par
Is it necessary to make the derivation of integral or differential equation 
models more rigorous or have I wasted you time? Some mathematicians claim that 
such an exercise is without merit since all one needs is a good guess; and if the 
magically obtained equation is solvable and reasonably predicts natural system 
behavior$,$ then this is all that is required. Not with standing such 
pronouncements$,$ other very productive mathematicians believe that there are 
a few obvious advantages to a more rigorous approach.\par
There is an advantage in education. In the physical sciences$,$ beginning 
students observe simplistic behavior and follow the well-established path of 
attempting to comprehend complex behavior as a composition of the simplistic. 
Extending these intuitively grasped concepts to the NSP-world leads more 
directly to integral or differential equation models and gives them a stronger 
incentive to seek solution methods and investigate the mathematical 
structures. Motivation still remains one of our major educational devices.\par
Scientific theories are not immutable. New experimental scenarios$,$ new 
measuring devices and$,$ indeed$,$ new insights into the foundations of natural 
system behavior have$,$ historically$,$ led to new theories that appear to predict
more accurately system development. Many new insights and theories abound$,$ 
today$,$ within certain scientific journals that specialize in such speculation. 
\par
In the early part of this century$,$ two of the great achievements of the human 
mind were the Special and General Theories of Relativity. Yet, as has been shown [{\tt Herrmann [1995]}]$,$ the major experimentally verified conclusions of 
the Special Theory are obtainable by simple infinitesimal modeling$,$ without 
the use of such concepts as frames of reference and the like. I believe 
that this fact 
will tend to eliminate some of the controversy that still surrounds the 
Special Theory and eventually lead to a better understanding of the basic
nature of electromagnetic radiation. I have previously mentioned how the 
recent concept of fractal behavior has been successfully modeled within the 
NSP-world environment.\par
Indirectly$,$ we are learning the simplistic NSP-world 
behavior that may be assumed to be the cause for certain natural 
world effects that are measured by standard instruments. 
As is well-known$,$ predicting accurately the behavior of many macroscopic 
natural systems$,$ even with our present day theories$,$ is often notably 
unsuccessful.  
It may be that the causes for such behavior are objectively real 
in character within a actual substratum world. And$,$ hence$,$ the more proficient 
we become in analyzing and applying infinitesimal techniques the more likely 
we are  to development new mathematical models that improve upon such 
predictions.\par
\vskip 18pt
\centerline{\bf NOTES}
\noindent [1]  In section 10.1 the concept of two vectors being almost 
parallel was introduced. The following theorem is sometimes useful.\par
\vskip 18pt
{\bf Theorem N.1.} {\sl Two unit vectors$,$ $\vec v,\ \vec u \in \hyperrealp n$ 
have the property that $\vec v \approx \pm \vec u$ if and only if $\vec v \bullet 
\vec u \approx \pm 1.$}\par
\vskip 18pt
\noindent [2] Relative to the statement on page 67 as to the relation between 
the work done along a hyperline segment and the curve if we add to the 
hypotheses of Application 8.2.1 the requirement that for each $t \in [a,b],\ 
c^\prime (t) \not= \vec 0$ and that $F$ is uniformly continuous on $E,$ then
using the notation in that section it follows that there exists some $\eps_j 
\in \monad 0$ such that $\hyper F(\ell_j(t_j^\prime))\bullet \hyper {\vec v_j} 
= \hyper F(\hyper c(t_j^\prime))\bullet \hyper  {\vec v_j} + \eps_j \Vert \hyper 
{\vec v_j} \Vert.$ From this we also have that $\sum_{j=1}^\Gamma 
\hyper F(\ell_j(t_j^\prime))\bullet \hyper {\vec v_j} \approx \sum_{j=1}^\Gamma 
\hyper F(\hyper c(t_j^\prime))\bullet \hyper  {\vec v_j}.$ The proof is at the 
end of Appendix 10.\vfil\eject
\centerline{Appendix For Chapter 2.}
\vskip 18pt
\indent {\bf Throughout these proofs lower case Greek letters will always 
denote infinitesimals.} [Remark: The complete proofs for many of these
fundamental propositions appear for the first time in these appendixes.]
\vskip 18pt
\indent {\bf Theorem 2.1.1.} {\sl Assumption (II) holds for $\hyperreal$ if and
only if $\hyperreal$ is not Archimedean.}\par
\indent Proof. Assume that (II) holds. Since $\hyperreal$ is a field then 
${\eps}^{-1} \in \hyperreal.$ Let $0\not= n\in \nat$. Then $0\not=\eps < 1/n.$
Thus $n < {\eps}^{-1}$ implies that $\hyperreal$ is not Archimedean.\par 
\indent Conversely$,$ assume that $\hyperreal$ is not Archimedean.
Then there exists 
some $\Gamma\in \hyperreal$ such that for each $n\in \nat,\ \vert \Gamma \vert \geq 
n.$ Since $\vert \Gamma\vert \in \hyperreal,$ then it follows
that for each $ n\in \nat,\ \vert \Gamma \vert \geq n+1 > n\Rightarrow 
\vert {\Gamma}^{-1}\vert < 1/n,\ \forall n \in \nat,\ n\not=0$ and the result 
follows.\par                                  
\vskip 18pt
\indent {\bf Theorem 2.2.1.} {\sl The set of limited numbers$,$ $\cal O,$ is a 
subring of $\hyperreal,\ \real \subset \cal O$ and $\cal O$ is not a 
field.}\par
\indent Proof. Let $x,y\in \cal O$. Then there exists some $r_1,r_2 \in \real$ 
such that $\vert x \vert < r_1,\ \vert y \vert < r_2.$ Thus $
\vert x\pm y\vert \leq \vert x\vert + \vert y \vert < r_1 + r_2 \in \real 
\Rightarrow x \pm y \in {\cal O},\ \vert xy \vert = \vert x\vert \vert y\vert 
<r_1r_2 \in \real \Rightarrow xy \in \cal O$.  If $ 0\not= \eps,$ then $\vert 
\eps \vert < r^{-1},\ \forall r \in {\real}^+ \Rightarrow {\eps}^{-1} \notin 
\cal O.$ The result follows. \vskip 18pt \indent {\bf Theorem 2.2.2.} {\sl The 
infinitesimals$,$ $\monad 0,$ form a subring of $\cal O$ and $\Gamma $ is 
infinite if and only if there is some nonzero $\eps \in \monad 0$ such that 
$\Gamma = 1/{\eps}.$}\par Proof. Let $x,y \in \monad 0.$ Then $\forall r\in 
{\real}^+,\ \vert x \vert < r/2,\ \vert y\vert < r/2,\ \Rightarrow \vert x \pm 
y\vert < r.$ Now we also know that $\vert x \vert < \sqrt r,\ \vert y\vert 
<\sqrt r.$ Hence$,$ $\vert xy \vert = \vert x\vert \vert y \vert < r.$ Since $r$ 
is positive and arbitrary then the first result follows. The second follows as 
in the proof of Theorem 2.2.1.\par \vskip 18pt \indent {\bf Theorem 2.2.3.} 
{\sl The set of infinitesimals$,$ $\monad  0,$ is an ideal of $\cal O.$}\par 
Proof. Let $b \in \cal O$. Then $\exists t\in {\real}^+$ such that $\vert b 
\vert < t.$ Also $\forall r \in {\real}^+,\ \forall x\in \monad 0,\ \vert 
x\vert < r/t.$ Hence $\vert xb \vert < r$ and the result follows.\par {\bf 
Theorem 2.2.4.} {\sl The binary relation $\approx$ is an equivalent relation 
on $\hyperreal .$}\par Proof. (i) Since $0 \in \cal O$ then $x-x=0\in \monad 0 
\Rightarrow x\approx x.$\par (ii) Let $x\approx y.$ Then $ x-y \in \monad 0 
\Rightarrow -(x-y) = y-x \in \monad 0 \Rightarrow y \approx x.$\par (iii) Let 
$x\approx t,\ y\approx z.$ Then $x-y = \eps,\ y-z =\delta.$ Then $\eps + 
\delta \in \monad 0 \Rightarrow x-z \in \monad 0 \Rightarrow x\approx z.$\par 
\vskip 18pt {\bf Theorem 2.2.5.} {\sl For each $x,y\in \real,$\hfil\break 
\indent\indent (i) $\monad x \cap \monad y = \emptyset$ if and only if $ x 
\not= y,$\hfil\break \indent\indent (ii) ${\cal O} = \bigcup \{\monad r\vert 
r\in \real\}.$}\par Proof. (i) follows since each monad is an equivalence 
class for the equivalence relation $\approx.$\par (ii) As previously noted 
$\real \subset \cal O.$ Note that $\monad r \cap  \monad {r+1} =\emptyset.$ 
Let $y\in \monad r,\ z \in \monad {r +1 }.$ Then $y = r + \eps,\ z = 
r+1+\delta.$ Assume that $z\leq y.$ Then $r+1+\delta \leq r + \eps \Rightarrow 
1\leq \eps - \delta \in \monad 0$: a contradiction. Thus $y < z.$ But $r+1\in 
\monad {r+1} \Rightarrow y< r+1.$ In like manner$,$ $r-1 <y.$ Thus $y\in {\cal 
O} \Rightarrow \monad r \subset {\cal O} \Rightarrow \bigcup\{\monad r\vert 
r\in \real \} \subset \cal O.$ Now let arbitrary $ x \in \cal O.$ and define 
$S = \{y\vert y\in \real \land y < x\}.$ The set $S$ is bounded above and 
nonempty since $\exists r\in {\real}^+$ such that $\vert x\vert < r,$ and $-r 
\in S.$ By completeness $\sup S \in \real.$ For each $p \in{\real}^+,\ \exists 
y\in S$ such that $\sup S - p<y<x.$ Also $x\leq \sup S + p,$ for otherwise we 
have 
 that $\sup S + p\in S,$ which contradicts the sup definition. Thus 
$\vert x- \sup S \vert \leq p.$ Hence $x - \sup S \approx 0 \Rightarrow
x\in \monad {\sup S} \Rightarrow {\cal O} \subset \bigcup\{\monad r\vert r\in
\real \}$ and the result follows.\par
{\bf Corollary 2.2.4.1.} {\sl If $x,y \in \real,\ x < y,\ z\in \monad x,\ 
w\in\monad y,$ then $z < w.$}\par
Proof. Since $x<y$ then $\exists r\in{\real}^+$ such that $y = x + r.$ Now  
consider the obvious portion of the above argument applied to the monads
$\monad x,\ \monad {x+r}.$\par
\medskip
{\bf Corollary 2.2.4.2.} {\sl If $x,y\in \monad r,\ z\in\hyperreal,\ x<z<y,$ 
then $z\in \monad r.$}\par
Proof. Since $\vert y-z \vert < \vert y-x \vert < r,\ \forall r\in 
{\real}^+$ then $y-z \in \monad 0.$ Thus $r \approx y,\ y \approx z, 
\Rightarrow r \approx z \Rightarrow z \in \monad r.$\par
\vskip 18pt
{\bf Theorem 2.2.6.} {\sl Each monad and the set of limited numbers are bounded 
above {\rm [}resp. below{\rm ]}$,$ but do not possess a least upper bound
 {\rm [}resp. greatest lower bound{\rm ]}.}\par
Proof. This is established for the ``bounded above'' case only. Let $\monad r$ 
be any arbitrary monad. From Corollary 2.2.4.1$,$ the set $\monad r$ is 
bounded above by the real number $r+1.$ Assume that $\sup {\monad r} = s \in 
\hyperreal.$ First$,$ assume that $s \in \real.$ In this case$,$ it is necessary 
that $s \not= r$ for we know $\exists w >0$ and if $ 
s= r,$ then $s + w \in \monad r,\ s < s + w,$ contradicts the bounding 
character of $s.$ Thus $r < s \Rightarrow \forall x \in \monad r,\ \forall
y \in \monad s;\ x < y.$ But $-w + s \in \monad s \Rightarrow \forall x \in 
\monad r; x < -w +s < s.$ Thus $s$ would not be the $\sup {\monad r}.$ 
We must$,$ therefore$,$ pass to the second possible case that $s \in \hyperreal -
\real.$ Obviously$,$ $r < s < r+1$ and thus $s \in \cal O.$ Consequently$,$ 
$\exists r_1 \in \real$ such that $s \in \monad {r_1}.$ Let $r_1 \not=r.$ Then 
$r < r_1$ and $\forall x\in \monad r,\ \forall y \in \monad {r_1},\ x < y.$
But again using $-w +s$ we would have that $s\not=\sup {\monad r}.$
The final case$,$ requires $s \in \monad r.$ But if this were so$,$ then
$s < s+w \in \monad r$ yields that $s$ would not be an upper bound. Thus no 
such $s$ exists and the proof is complete for monads.\par
Let $s = \sup {\cal O}.$ Then $s\notin \cal O$ since if not$,$ then $s+1 \in \cal O$.
On the 
other hand$,$ if $s$ is infinite$,$ then $s > 0 \Rightarrow s-1 $ is infinite 
and that $s-1 >0.$ Thus no such $s$ exists in $\hyperreal$ and the proof in 
complete.\par 
\vskip 18pt
\indent {\bf Theorem 2.3.1.} {\sl For any $n\in\nat,\ n\geq 1,$  the set of limited 
vectors in $\hyperrealp n$ is equal to $${\cal O}^n = \overbrace{{\cal 
O}\times\cdots\times{\cal O}}^{n\rm\;factors}.$$}\par
Proof. See Theorem 2.3.3\par
\vskip 18pt
{\bf Theorem 2.3.2.} {\sl For any $n\in\nat,\ n\geq 1,$ then monad of the 
vector $\vec v=(x_1,\ldots ,x_n)\in \realp n$ is equal to $$\monad {\vec v} = 
\overbrace{\monad {x_1}\times\cdots\times\monad {x_n}}^{n\rm\;factors}.$$}\par 
Proof. See Theorem 2.3.3\par
\vskip 18pt  
{\bf Theorem 2.3.3.} {\sl For any $n\in\nat,\ n\geq 1,$ and for each ${\vec 
v},{\vec w}\in 
\realp n,$\hfil\break
\indent\indent (i) $\monad {\vec v} \cap \monad {\vec w} = \emptyset$ 
f and only if $ {\vec v} 
\not= {\vec w},$\hfil\break
\indent\indent (ii) ${\cal O}^n= \bigcup \{\monad {\vec v}\vert \vec v
\in \realp n\}.$}\par 
Proof. Theorems 2.3.1 and 2.3.2 follow from the *-transform of the basic
Euclidean norm fact that 
if ${\vec v} = (x_1,\ldots,x_n) \in \realp n,$ then $\forall i,\ 1 \leq i 
 \leq n,\ \vert x_i \vert \leq \Vert \vec v \Vert \leq \sum_{i=1}^n\vert x_i 
 \vert$ and
the definitions. Theorem 
2.3.3 can be established directly or it follows easily from the basic
set-theoretic properties for the cross product.\par
\vskip 18pt
{\bf Theorem 2.4.1.} {\sl Let $x,y\in \cal O$. Then\hfil\break
\indent\indent (i) $x \approx y$ if and only if $\st x = \st y,$\hfil\break
\indent\indent (ii) $x \approx \st x,$\hfil\break
\indent\indent (iii) if $x\in \real,$ then $\st x = x,$\hfil\break
\indent\indent (iv) if $x \leq y,$ then $\st x \leq \st y,$\hfil\break
\indent\indent (v) if $\st x \leq \st y,$ then either $x \leq y$ or $x-y \in
\monad 0$ with unknown order.}\par
Proof. 
Observe that for the limited hyperreal number $x,$ that $x = \st x + \eps.$\par
(i) $x,y \in {\cal O} \Rightarrow x = \st x +\eps; \ x \approx y \Rightarrow
x = y + \delta \Rightarrow y+\delta = \st x + \eps \Rightarrow y = \st x + (\eps 
- \delta) \Rightarrow y \in \monad {\st x} \Rightarrow \st x = \st y.$ \par
If $\st x = \st y,$ then $ x = \st x + \eps,\ y = \st y + \delta \Rightarrow
x- \eps = y -\delta \Rightarrow x \approx y.$\par
(ii) $ x = \st x + \eps \Leftrightarrow x \approx \st x.$\par
(iii) For $r \in \real,\ r \in \monad r \Rightarrow r = \st r.$ \par
(iv) Let $ x \leq y.$ Then $x,y \in \monad 0 \Rightarrow \st x = \st y.$
If $x \in \monad r,\ y \in \monad {r_1},\ r\not= r_1,$ then $r > r_1 
\Rightarrow x> y.$ Thus contradiction yields that $\st x = r < \st y = 
r_1.$\par
(v) Assume that $\st x = \st y.$ Then $x - y \in \monad 0.$ If $\st x < \st 
y,$ then $x \in \monad {\st x},\ y \in \monad {\st y} \Rightarrow x < y.$\par
\vskip 18pt
{\bf Theorem 2.4.2.} {\sl Let $x,y\in \cal O$. Then \hfil\break
\indent\indent (i) $\st {x\pm y} = \st x \pm \st y,$ \hfil\break
\indent\indent (ii) $\st {xy} = (\st x)(\st y).$}\par
Proof. Notice that in general$,$ if $ r\in \real,$ then $\st {r+\eps} = r.$
Let $r = \st x,\ s = \st y.$ \par 
(i) $\st {x\pm y} = \st {(r + \eps) \pm  (s + \delta)} = \st {(r \pm s) + 
(\eps \pm \delta)} = r \pm s.$\par
(ii) $\st {xy} = \st {rs + (r\delta + s\eps + \delta\eps)} = rs.$\par
\medskip
{\bf Corollary 2.4.2.1} {\sl Let $x,y\in \cal O$. Then \hfil\break
\indent\indent (i) if $\st y \not= 0,$ then $\st {x/y} = {\st x}/{\st 
y}.$\hfil\break
\indent\indent (ii) if $y = \root n \of x,$ then $\st y = {(\st x)}^{(1/n)}$$,$ 
where it is always the case that if $x\in \cal O$$,$ then $\root n \of x \in 
\cal O$.}\par
Proof. If $\st y \not = 0,$ then $y \notin \monad 0 \Rightarrow y^{-1} \in {\cal 
O} \Rightarrow (x/y) \in \cal O.$ Then $\st x = \st {(x/y)y} =
(\st {x/y})(\st y).$\par
(ii) Obviously$,$ $x\geq 0.$ Also$,$ if $x \geq 0,\ x \in \cal O,$ then $x < r \in 
{\real}^+ \Rightarrow \root n \of x < \root n \of r \Rightarrow \root n \of x 
\in  \cal O.$ Then $ y = \root n \of x \Rightarrow y^n = x \Rightarrow (\st 
y)^n = \st x$ and from this the proof is completed.\par
\vskip 18pt
{\bf Theorem 2.4.3.} {\sl For $\vec v\in {\cal O}^n$$,$ let $\st {\vec v} =
(\st {x_1},\ldots,\st {x_n}).$ The set ${\cal O}^n$ forms a vector space with 
respect to the ring $\cal O$ (i.e. a module) and as such the map {\rm st} 
distributes over the vector space algebra as well as the component defined dot 
and for $n = 3$ the cross product operators.}\par
Proof. This follows immediately from the above standard part operator 
properties.\par
\vskip 18pt
{\bf Theorem 2.4.3.}  {\sl The set $\monad 0$ is a maximal ideal in $\cal O$ and the 
quotient ring ${\cal O}/\monad 0$ is isomorphic to $\real.$}\par
Proof. Let $\Im $ be an ideal such that $\monad 0 \subset \Im \subset \cal O$ and 
assume that $\monad 0 \not= \Im.$ Let $x \in \Im - \monad 0.$  We know that 
$x^{-1} \in \cal O.$ Since $\Im$ is a ideal in $\cal O$ then $x(x^{-1}) =1 \in 
\Im .$ Thus implies that $\Im = \cal O.$ The standard part operator is 
obviously the isomorphism and this completes the proof.\par 
\vfil
\eject
\centerline{Appendix For Chapter 3.}
\vskip 18pt
\indent {\bf Throughout these proofs lower case Greek letters will always 
denote infinitesimals$,$ while upper case will denote infinite numbers.} \par
Even though it is possible to restrict our model to objects such as $\Re$ it 
is useful for future considerations to actually construct 
a general superstructure as originally envisioned by {\tt Robinson and Zakon 
[1969]}. Let $X_0 = \real$ and define by induction $X_{n+1} = {\cal P}(\bigcup 
\{X_j \vert 0\leq j\leq n \}.$  Recalling that it is assumed that $\real$ is 
a set of {atoms} (i.e. they are not sets). Then ${\cal H} =
\bigcup \{X_n \vert n \in \nat \}$ is the 
{\bf {universe}} or {\bf {carrier}} for a structure ${\cal M} = ({\cal 
H},\eps,=).$ In this particular case$,$ this universe it termed a {\bf 
{superstructure} on $\real.$} In the literature there are slightly different
constructions of such superstructures {\tt [Davis 1977].} An object $b \in 
X_0 $ is an individual and otherwise if $b \in X_p$ then $b$ is an entity. 
Thus $b$ is an individual or entity if and only if  $b \in \cal H.$ It will be shown
that $\exists 
p\in \nat$ such that $\Re \in 
X_p.$ A superstructure has may different properties. The following is a brief
list of some of the more important of these$,$ where the proofs are
straightforward and left to the reader. [Note: Later I may change notation and use the usual mathematical practice of considering symbols such as $x,y,z$ in two contexts. One, the a mathematical variable, and the other a formal first-order language variable.]
\medskip
{\bf Proposition 1.} {\sl Let $\cal H$ be a superstructure on $\real.$\par
\indent\indent (i) $X_p \subset X_n$ for $n,p \in \nat,\ 1\leq p \leq n.$ $X_0 \cap X_n = \emptyset, \ n \geq 1.$\par
\indent\indent (ii) $X_0 \cup X_n = X_0\cup \cdots \cup X_n, \ X_{n+1} = \power {X_0 \cup X_n},\ n\geq 0.$\par
\indent\indent (iii) $X_p \in X_{n+1},\ \forall n \in \nat,\ n \geq 1,\ 0\leq 
p\leq n.$\par
\indent\indent (iv) If for $n \geq 1,\ a\in b\in X_n,$ then $a\in X_0 \cup 
X_{n-1}.$\par
\indent\indent (v) If $(a_1,\ldots, a_n) \in b \in X_p\ (p\geq 1),$ then 
$a_1,\ldots, a_n \in X_0 \cup X_{p-1}.$}\par
\indent\indent (vi) If $b \in X_n$ and $a \subset b$$,$ then $a\in X_n.$\par 
\medskip
By use of concepts in abstract model theory$,$ the nonstandard model
$\hyper {\cal M}= (\hyper {\cal H},\eps,=),\ \hyper {\cal H}=\bigcup 
\{\hyper X_n \vert n \in \nat \} $ is constructed through application
of the notion of internal 
individuals or entities and$,$ in particular$,$ (iv) of Proposition 1. [For a 
very clear and concise discussion of this please see page 120 of {\tt Robinson 
and Zakon [1969].}] With respect to our notation$,$ this construction
also leads 
to the following significant proposition that gives a general relation between 
internal objects and elements they contain.\par
{\bf Proposition 2.} {\sl Within our set theory }\par
\indent\indent {\sl (i) If $a \in Q \in \hyper X_{n+1},$
 then $a \in \hyper X_0 \cup \hyper X_n.$ (Every member of an internal entity is internal. Also note that this is the set-theoretic ``or'' statement which is a set even if it contains atoms.)}\par
\indent\indent {\sl (ii) The set $\Hyper {\cal H}$ is closed under finitely many applications of the basic set operation $\cup.$}\par
\indent\indent {\sl (iii) Any finite collection of members from $\Hyper {\cal H}$ is a member of $\Hyper {\cal H}$}\par
\indent\indent {\sl (iv) The set $\Hyper {\cal H}$ is closed under n-tuple formation.}\par
\indent\indent {\sl (v) If $(a_1,\ldots,a_p) \in Q \in \hyper X_{n+1},$ then 
$a_i \in \hyper X_0 \cup \hyper X_n,\ 1\leq i \leq p.$}\par

\indent\indent {\sl (vi) The set $\Hyper {\cal H}$ is closed under the basic set-theoretic operators $\cap, \ -,\ \times,$ among others, for internal sets.\par}
Proof. (i) Sometimes this result is established during construction. However, it follows easily by *-transfer of Proposition 1 (iv). \par
(ii) Note that for any $n\geq 0$ the set $X_n \cup X_0 \in X_{n+1}$. Hence, by Theorem 3.1.3 (i),  $\Hyper (X_n \cup X_0)$ is internal. Now we can use the independently established part (vi) d of Theorem 3.1.3 and obtain that $\hyper X_0 \cup \hyper X_n = \Hyper (X_n \cup X_0)$. Let $A$ and $B$ represent internal entities. For any $n \geq 0,$ the sentence 
$$\forall x \forall y((x \in X_{n+1})\land (y \in X_{n+1}) \to \exists z((z\in X_{n+1}) \land \forall w(w \in (X_n \cup X_0) \to$$ $$( w \in z \iff (w\in x)\lor (w\in y))))$$
holds in $\cal M$; hence, in $\Hyper {\cal M}.$ Now since $A$ and $B$ are internal, there is an $n$ such that$A,\ B \in \hyper X_{n+1} \cup \hyper {X_0} = \hyper (X_{n+1} \cup X_0).$ Thus, there is a $C \in \hyper {X_{n+1}} $ that contains the same elements of $\hyper {X_n} \cup \hyper {X_0}$ as those in ``$A \cup B$.'' But by  (i), all the elements in $A,\ B,\ C$ are members of $\hyper {X_n} \cup \hyper {X_0}.$ Thus $C = A \cup B \in \hyper {\cal H}.$ This proof is completed by induction.\par
(iii) Theorem 3.1.3 (vi)b, shows that this holds for an empty set. The notion of what constitutes the ``number'' of members of a nonempty finite set is an intuitive  metamathematical notion the corresponds this idea to the ``number'' of $\lor$ symbols. [Note: This is not the only way this can be established. Indeed, we could use finite set generating operators.] For each $n \geq 0,$ the sentence 
$$\forall x(x \in (X_{n} \cup X_0) \to \exists y((y \in X_{n+1}\land $$ $$\forall z((z \in (X_n \cup X_0)\land z \in y) \iff z = x))))$$ 
holds in $\cal M$; hence, it holds in $\Hyper {\cal M}.$ Now given any internal $a \in \Hyper {\cal H},$ then $a \in \hyper X_0 \cup \hyper X_n$ for some $n \geq 0.$ Thus, for any $a \in \Hyper {\cal H}$ the set that contains $a$ and only $a$ exists and is an internal. Now to establish this for the notion of a two element with two distinct members, we extend the above sentence. Thus for any $n \geq 0$ the sentence 
$$\forall x\forall y(x \in (X_{n} \cup X_0)\land y \in (X_{n} \cup X_0)\land x \not= y \to \exists y((y \in X_{n+1}\land $$ $$\forall z((z \in (X_n \cup X_0)\land z \in y) \iff (z = x)\lor (z = y))))$$
holds in $\cal H;$ hence, it holds in $\Hyper {\cal H}$. But, we know that for any two internal $a,\ b$ there exists some $n$ such that $a \in \Hyper X_n \cup \Hyper X_0.$ Then every set that contains two and only two distinct internal entities is internal by *-transfer. To apply this to any set of finitely many internal entities we use induction or application of (ii). \par
(iv) Within our set-theoretic definitions for  $X_0$ or $\Hyper {X_0}$ order pair formation is given for specific $a,\ b.$ We use the abbreviation $(a,b)$ to denote the finite set $\{\{a\},\{a,b\}\}.$ This is then extended by induction to the n-tuple $(a_1,\ldots,a_{n-1},a_n) = ((a_1,\ldots,a_{n-1}),a_n).$ Thus any n-tuple is just a finite set, of finite sets, of finites sets,$\ldots$ constructed in the finite manner from the finite set $\{a_1,a_2,\ldots a_n\}.$ For example, consider internal $a,b,\ \not= b$. Then they are both members of some $\hyper {X_n} \cup \Hyper {X_0}.$ Then $\{a\}$ and $\{a,b\}$ are internal and both members of some $\Hyper {X_p}.$  Then the set $\{\{a\},\{a,b\}\} \in \Hyper X_{p+1}.$ Apply this idea to an n-tuple and you have that the n-tuple is internal. \par 
(v) This follows from the definition of the n-tuple, repeated applications of (i) and Proposition 1 (i). \par
(vi) These are established in the same manner as (ii). The $\cap$ is established as an example after Theorem 3.4.1.
\par\medskip

{\bf Theorem 3.1.1.} {\sl It is permissible to assume that\par
\indent\indent (i) if $A\subset \real,$ then $A\subset 
\hyperreal,$\hfil\break
\indent\indent (ii) if $A\subset \realp n,$ then $A\subset \hyperrealp n,$ 
\hfil\break
\indent\indent (iii) if $A\subset (\realp n)\times(\realp m),$ then 
$A \subset (\hyperrealp n)\times(\hyperrealp m).$}\par
\vfil\eject Proof.\par 
(i) This follows from the construction of the nonstandard 
model $\hyper {\cal M} = (\hyper {\cal H},\eps,=)$  by means of the
superstructure technique where each member of the carrier
$\real$ is treated as a constant sequence {\tt [Robinson and Zakon [1969],
Stroyan and Luxemburg [1976]$,$ Davis [1977]$,$ etc.]} Basically this follows
from the fact that each member of $\real$ is an atom or urelement within
our basic set theory.\par
(ii) Let $(a_1,\ldots,a_n)\in A\subset \realp n.$ By construction, $\realp n \in {\cal H}.$ From our basic definition 
of an n-tuple as a finite set of finite sets etc.$,$ then this result is 
immediately established from independently proved Theorem 3.1.3 parts (i) and 
(iii) by finitely distributing the * throughout these finite sets and using 
(i) of this theorem.  If $A = \emptyset,$ then result follows from Theorem 3.1.3 
part (vi).\par
(iii) This is proved in the same manner as part (ii).\par
\vskip 18pt
{\bf Theorem 3.1.2.} {\sl If $A \subset \real$ {\rm [}resp. $\realp n,
\ (\realp n)\times(\realp m)${\rm ]}$,$ then there exists $\hyper A\subset \hyperreal$
 {\rm [}resp. $\hyperrealp n,\ (\hyperrealp n)\times(\hyperrealp m)${\rm ]}
such that \par
\indent\indent (i) $A \subset \hyper A$ and \hfil\break
\indent\indent (ii) $A = \hyper A$ if and only if $A$ is finite.}\par
Proof. From the independently established Theorem 3.4.1 it follows by *-transfer 
that $\hyper A \subset \hyperreal.$\par                                      
(i) This follows from Theorem 3.1.3 part (i) and$,$ as discussed in Theorem 
3.1.1$,$ from the construction of $\hyper {\cal M}.$ \par
(ii) The model $\hyper {\cal M}$ is assumed$,$ at least$,$ to be an enlargement.
Thus consider
the binary relation $R = \{(x,y)\vert (x\in A)\land (y\in A)\land (x\not= 
y)\}.$ Now let $\{(x_1,y_1),\ldots,(x_n,y_n)\} \subset R$ and suppose that $A$ 
is infinite. Then there exists some $b \in A$ such that $x_i \not= b,\ 1\leq 
i\leq n.$ From the concurrency of $R$$,$ this implies that there exists some
internal $c \in \hyper A$ such that $\forall a \in {^\sigma A}= A,\ a\not=c.$ [The ${}^\sigma$ operator is that defined in {\tt Stroyan and Luxemburg, [1976]}. In general, it is the set of all constant sequence $\cal U$-equivalence classes contained in $\Hyper A.$]
The result now follows from Theorem 3.1.3.\par
Conversely$,$ if $A = \emptyset$ or $A$ is finite$,$ then result follows from 
Theorem 3.1.3.\par\medskip [Note: Since $X_{n+1} = \power {X_0 \cup X_n},$ then $(X_0 \cup X_n) \in X_{n+1}$ implies $(X_0 \cup X_n) \in {\cal H}.$ \par\bigskip

{\bf Theorem 3.1.3.} {\sl Unless otherwise
stated all constants represent individuals or entities (i.e. members of 
$C({\cal H})$ where we do not notationally differentiate between the name of an 
object and the object itself.)\par
\indent\indent (i) $a\in A$ if and only if $\hyper a\in \hyper A;\ A\not= B$ 
if and only if $\hyper A \not= \hyper B.$\par
\indent\indent (ii) $A\subset B$ if and only if $\hyper A \subset \hyper 
B.$\par
\indent\indent (iii) $\hyper {\{A_1,\ldots,A_k\}} = \{\hyper {A_1},\ldots,
\hyper {A_k} \}.$\par
\indent\indent (iv) $\hyper {\eskip (A_1,\ldots,A_k)} = (\hyper {A_1},\ldots,
\hyper {A_k}).$\par
\indent\indent (v) $ (A_1,\ldots,A_k) \in A$ if and only 
$(\hyper {A_1},\ldots,\hyper {A_k}) \in \hyper A.$ \par
\indent\indent (vi) Let $A,B$ be sets [resp subsets of $\cal U$].
 Then (d) $\hyper {\eskip (A \cup B)} = \hyper A
\cup \hyper B$$,$ (c) $\hyper {\eskip (A \cap B)} = \hyper A \cup \hyper B$$,$ (a)
$\hyper {\eskip (A - 
B)} = \hyper A - \hyper B$$,$ (e) $\hyper {\eskip (A\times B)} = \hyper A \times \hyper 
B$$,$  (b) $\hyper {\eskip \emptyset} = \emptyset.$\par
\indent\indent (vii) If $A\in \real,$ then $\hyper A = A.$}\par
Proof. (i) The sentences 
$$a\in A,\ \neg (A = B)$$
hold in ${\cal M}$ if and only if the *-transfers 
$$\hyper a\in\hyper A,\ \neg (\hyper A = \hyper B)$$  hold in
$\hyper {\cal M}.$\par
(vi) [Clearly$,$ by construction ${\cal U} \in  {\cal H}.$  
Consequently$,$ $A,B$ are entities by Proposition 1.] 
(a) We know that there exists a $C\in C({\cal H})$ such that $C = A - B.$
The sentence $$\forall x(x \in C 
\leftrightarrow (x\in A)\land(\neg (x\in B)))$$  holds in $\cal M$ if and only if 
$$\forall x(x \in \Hyper C
\leftrightarrow (x\in \hyper A)\land(\neg (x\in \hyper B)))$$ holds in
$\hyper {\cal M}.$ Then the set $\Hyper C = \hyper A - \hyper B.$ Now (b)
is obtained by letting $\emptyset = \real - \real.$ Then from (a)$,$ $\Hyper 
{\emptyset} = \hyperreal - \hyperreal = \emptyset.$ (c) We also know that 
there is a $C \in \cal H$ such that $C= (A\cap B) = A - (A - B).$ Once again by 
(a); (c) follows since $\Hyper (A\cap B) = \hyper A - (\hyper A - \hyper B) =
  \hyper 
A \cap \hyper B.$ Then for (d) there is a $A\cup B=C\in \cal H$. Now simply consider the sentence $\forall x((x \in C) \iff (x\in A) \lor (x \in B)).$ This result follows by *-transfer. (Another proof for (c) also follows from *-transfer.)  
(ii) This follows from parts (a) and (b) of (vi) since $A\subset B$ iff
$A-B=\emptyset.$ (vi) (e) The Cartesian product statement is a little more
  difficult 
to establish. In this elementary approach$,$ it is useful to characterize the 
notion of the singleton and doubleton sets.  The set $S(A)$ of all singleton subsets of $A$ and the set $D(A)$ of all 
doubleton subsets of $A$ are members of $X_{n+2}$
(i.e. $D(A) \subset X_{n+1}).$ The following sentences characterize the
singleton and 
doubleton subsets of $A \subset X_n \cup X_0.$
 $${\Phi_1} = \forall x(x \in D(A) \to 
\exists y \exists w(w \in A \land y \in A \land y \in x \land$$ $$ w \in x 
\land \forall z(z \in X_n \cup X_0 \land z \in x \to z = y \lor z = w))).$$ 
[Note: $P \land Q  \to S \equiv P \to (Q \to S).$] $${\Phi_2} = \forall x(x \in 
S(A) \to  \exists y(y \in A \land y \in x \land $$ $$\forall z(z \in X_n \cup 
X_0 \land z \in x \to z = y))).$$ The *-transfer of these sentences becomes 
$$\Hyper {\Phi_1} = \forall x(x \in \hyper D(A) \to \exists y \exists w(w \in 
\hyper A \land y \in \hyper A \land y \in x \land$$ $$ w \in x \land \forall 
z(z \in \hyper X_n \cup \hyper X_0 \land z \in x \to z = y \lor z = w))).$$ 
$$\Hyper 
{\Phi_2} = \forall x(x \in \hyper S(A) \to  \exists y(y \in \hyper A \land y 
\in x \land $$ $$\forall z(z \in \hyper X_n \cup \hyper X_0 \land z \in x \to 
z = y))).$$ 
[Note: $\hyper D(A) = \hyper (D(A)),\ \hyper S(A) = \hyper (S(A)).$] In what follows, the constants used to name various internal 
objects are in our extended internal language. What these *-transforms 
indicate is that the set $\hyper S(A)\subset S(\hyper A),\ \hyper D(A)\subset 
D(\hyper A).$ Indeed$,$ the same would hold true for any finitely numbered 
subset of $A$. Observe first that $\hyper D(A),\ \hyper S(A) \subset \hyper 
X_{n+1}$ and $\hyper A \subset \hyper X_n \cup \hyper X_0.$ For the
 doubleton case$,$ let $Q \in 
\hyper D(A).$ Thus each such $Q$ is internal since it is a member of $\hyper 
{X_{n+1}}.$ There then exists internal $c,d \in \hyper A$ and for any $f \in 
\hyper X_n \cup \hyper X_0$ and $ f \in Q$ it follows that $f = c$ or $f = d.$ Now since
 $\hyper A 
\subset \hyper X_n \cup \hyper X_0$ then this last statement also includes any $f \in \hyper 
A.$ Proposition 2 (i) states  that \underbar{every} $z\in Q$ is internal and a 
member of $\hyper {X_n}\cup \hyper {X_0}.$ Thus only $c,d \in Q \Rightarrow Q 
\in D(\hyper A).$ Hence $\hyper D(A) \subset D(\hyper A).$ \par
   Now consider 
the sentences $$\Phi_3=\forall x \forall y(x \in A \land y \in A \to \exists 
z(z \in D(A) \land 
 x \in z \land y \in z \land$$ $$ \forall w(w \in X_n \cup X_0 \land w \in z
 \to w = x \lor w = y)))$$
 $$\Phi_4 = \forall x(x \in A \to \exists z(z \in S(A) \land x \in z \land$$
 $$\forall w(w \in X_n \cup X_0 \land w \in z \to w = x))).$$  Once again 
consider the *-transfer 
$$\Hyper {\Phi_3}=\forall x \forall y(x \in \hyper A \land y \in \hyper A 
\to \exists z(z \in \hyper D(A) 
\land
 x \in z \land y \in z \land$$ $$ \forall w(w \in \hyper X_n \cup \hyper X_0
\land w \in z \to w = x \lor w = y)))$$
 $$\Hyper {\Phi_4} = \forall x(x \in \hyper A \to \exists z(z \in \hyper S(A)
\land x \in z \land$$
 $$\forall w(w \in \hyper X_n \cup \hyper X_0 \land w \in z \to w = x))).$$
Now let $W,U$ be objects within our set theory. Let 
$W \in D(\hyper A).$ Then $\exists a,b \in \hyper A$ such that for any $U\in 
W,$ in particular $U \in \hyper X_n \cup \hyper X_0,\ U = a$ or $U = b.$ The 
*-transfer states that there exists some $Q \in \hyper (D(A))$ such that
$a,b \in Q.$ Thus $W \subset Q.$ Postulate 2 states that if any $z \in Q,$ 
then $z \in \hyper X_n \cup \hyper X_0$ and *-transfer yields that $z = a$ or 
$z = b\Rightarrow W = Q \Rightarrow D(\hyper A) \subset \hyper D(A) 
\Rightarrow \hyper D(A) = D(\hyper A).$ In like manner$,$ $\hyper S(A) = 
S(\hyper A).$ (This is not the only way to establish these equalities. I have chosen to use the operators $S$ and $D$ since these and similar operators could prove useful in later investigations.) \par
\baselineskip 12pt
 In order to establish that $\Hyper (A\times B) = \hyper A \times \hyper B$ 
similar techniques are employed. First$,$ noting that $C= A \times B \in X_n$ for 
some $n > 1$ we have that the sentence
$$\Phi_5 = \forall z(z \in C \to \exists x \exists y \exists w \exists u(
x \in S(A) \land y \in D(A\cup B) \land$$ $$ w \in A \land u \in B \land w \in x \land
w \in y \land u \in y))$$ 
holds in $\cal M$ if and only if (iff)
$$\Hyper {\Phi_5} = \forall z(z \in \Hyper C \to \exists x \exists y \exists w 
\exists u(x \in S(\hyper A) \land y \in D(\hyper A\cup \hyper B) \land $$ 
$$w \in 
\hyper A \land u \in \hyper B \land w \in x \land
w \in y \land u \in y))$$ 
holds in $\hyper {\cal M}.$ This implies that $\Hyper C \subset \hyper A \times
\hyper B.$ For the converse$,$ it is clear from the result that $S(\hyper F) = 
\hyper S(F),\ D(\hyper F) = \hyper D(F)$ and *-transfer that the sentence 
$$\Hyper {\Phi_6} =\forall x\forall y \forall w \forall z(x \in S(\hyper A) 
\land y \in D(\hyper A \cup \hyper B) \land w \in \hyper A \land$$ 
$$ z \in
\hyper B \land w \in x \land w \in y \land z \in y \to \exists v(v \in 
\Hyper C
\land x \in v \land y \in v \land $$ $$\forall u(u \in \hyper {X_0} \cup \hyper 
{X_{n-1}} \land u \in \Hyper C \to u = v))).$$
holds in $\Hyper {\cal M}.$ 
Let $W \in \hyper A \times \hyper B.$ Then through application of Proposition  
2 and arguing as above there is a unique $U \in \Hyper {(A \times B)}$ such 
that $U = W.$ From this the result follows.\par
{\it Proofs such as part (e) of this theorem using the above method are 
considered tedious and often do not appear in the literature. In  modern 
nonstandard analysis they are considered trivial examples of the Leibniz 
Principle and are often left to the reader to obtain. I disagree with such 
sentiments. It is easy for an author to say these concepts ``are expressible 
by appropriate first-ordered statements and the result easily follows.'' It 
is$,$ however$,$ not so easy for the untrained to actually do so for  certain
subtle requirements such as applications of Proposition 2 could easily allude 
the neophyte.}\par 
\baselineskip 14pt 
(iii) This is proved by induction with respect to the set $N_1 = \{n\vert n
\in \nat \land n \geq 1\}.$ Let $n = 1,\ A_1 \in \cal H.$ Then there exists
some 
$p \in \nat,\ (p\geq 1)$ such that $C = \{A_1\} \in X_p.$ Thus $\Hyper C 
\subset \hyper X_0 \cup \hyper X_{p-1}.$ The sentence $$\forall x(x \in 
X_0 \cup X_{p-1} \to (x \in C \leftrightarrow x = A_1))$$               
holds in $\cal M$ if and only if $$\forall x(x \in 
\hyper X_0 \cup \hyper X_{p-1} \to (x \in \Hyper C \leftrightarrow x = \hyper 
A_1))$$ 
holds in $\hyper {\cal M}.$     
Thus $\hyper A \in \Hyper C$ and $d \in \Hyper C \Rightarrow d \in \hyper X_0
 \cup \hyper X_{p-1} \Rightarrow d = \hyper A_1$ yields that $\Hyper C = 
\{\hyper A_1 \}.$ Assume that result holds for some $n \in N_1$ and consider 
the set $\{A_1,\ldots,A_n,A_{n+1}\}.$ Note that 
$\{A_1,\ldots,A_n,A_{n+1}\}\in X_p$ for some $p \in \nat$ since it is a finite 
set. Now simply note from the previous case and part (vi) that
$\{A_1,\ldots,A_n,A_{n+1}\} = \{A_1,\ldots,A_n\} 
\cup \{A_{n+1}\} \Rightarrow \Hyper {\{A_1,\ldots,A_n\}} \cup \Hyper 
{\{A_{n+1}\}} = \{\hyper A_1,\ldots,\hyper A_n\} \cup \{\hyper A_{n+1} \}= \Hyper \{A_1,\ldots , A_n\}.$ 
\par
(iv) This follows by induction$,$ the definition of n-tuples and part (iii).\par
(v) This follows by the *-transfer of sentences composed only of constants 
such as in part (i) and from (iv).\par
(viii) See proof of part (i) of Theorem 3.1.1. This completes this proof.\par
\vskip 18pt
Part (vi) (e) of Theorem 3.1.3 has some very useful consequences that allows 
for a great simplification of our first-order language. First$,$ by induction$,$ 
one has immediately that $\Hyper (A_1 \times \cdots \times A_n) = (\hyper A_1 
\times \cdots \times \hyper A_n).$ Also since any object in
$\Hyper (A_1 \times \cdots \times A_n)$ is an actual set-theoretic n-tuple 
from the basic definition then for any $a_i \in \hyper A_i,\ 1\leq 
i\leq n$ {\bf the symbol $(a_1,\ldots,a_n)$ may be used as the correct
set-theoretic abbreviation for the unique set termed an n-tuple.
This simplification is utilized continually throughout the remainder of 
these manuals.}\par
[Note: There are numerous equivalent forms that can be used to represent a statement formally in the required ``bound'' form so that each quantified variable is forced to vary only over members of $\cal H$ or $\cal \Hyper {\cal H}.$ In this regard, notice that $\cal H$ [resp. $\Hyper {\cal H}$] is closed under the basic set-theoretic operations for standard [resp. internal] objects. Thus, for example, if $A \in \cal H$ and $B \in \cal H$, then a statement such as $\forall x ((x \in A)\land (x \in B) \cdots$ satisfies the bound requirement since this is equivalent to $\forall x (x \in A\cup B \cdots)$ and $A \cup B \in \cal H.$] 
 \par\medskip
{\bf Theorem 3.1.4.} {\sl Let $\emptyset \not= R\subset A_1 \times \cdots 
\times A_n,\ (n>1),$\par
\indent\indent (i) $\hyper {P_i(R)} = P_i(\hyper R).$\par
\indent\indent (ii) If $R$ is a binary relation (i.e. n = 2)$,$ then $\hyper 
{R^{-1}} = (\hyper R)^{-1}.$}\par 
Proof. (i) This is established by induction on the number$,$ $m,$ of Cartesian 
products. First note that$,$ in general$,$ $R \subset A_1 \times \cdots \times A_n
\Rightarrow \hyper R \subset \hyper A_1 \times \cdots \times \hyper A_n$ and
 $P_1(R) \subset A_1 \Rightarrow \Hyper {(P_1(R))} \subset \hyper A_1.$  Let
$m = 1.$ Consider the *-transfer of the following sentence, where the symbol $(x,y)$ is an abbreviation for an obvious but more complex first-order expression.
$$\forall x(x \in A_1 \to (x \in P_1(R) \leftrightarrow \exists y(y \in A_2 
\land (x,y) \in R)).$$
Then $q \in \Hyper {(P_1(R))} \Rightarrow \exists s \in \hyper A_2$ such that
$(q,s) \in \hyper R \Rightarrow \Hyper {(P_1(R))} \subset P_1(\hyper R).$
On the other hand$,$ if $q \in P_1(\hyper R),$ then $q \in \hyper A_1$ and 
$\exists s \in \hyper A_2$ such that $(q,s) \in \hyper R \Rightarrow q \in 
\hyper {(P_1(R))}.$ Thus $P_1(\hyper R) = \hyper P(R).$ In like manner for the 
projection $P_2.$\par
Assume that result holds for an arbitrary $m \in N_1,$ (i.e. for $P_i,\ 
 1 \leq i \leq m+1$). Consider any $R\subset A_1 \times \cdots \times 
A_{m+2}.$ Now for a little general set theory. From the definition  
$(A_1 \times \cdots \times A_{m+2}) = (A_1 \times \cdots \times A_{m+1}) 
\times A_{m+2} = C \times A_{m+2}$ and $R \subset C \times A_{m+2}.$ Let 
$R_C = P_C(R)$. Then for each $i,\ 1\leq i \leq m+1,\ P_i(R_C) = P_i(R).$ These 
set-theoretical facts hold for any n-ary relation including $\hyper R.$ Thus$,$ 
from the induction hypothesis$,$ $\forall i,\ 1\leq i\leq m+1,\ \Hyper 
{(P_i(R))} = P_i(P_C(\hyper R)) = P_i(\hyper R).$ Then $\Hyper {(P_{m+2}(R))}
 = P_{m+2}(\hyper R)$ follows from the first case where $m = 1.$ \par
(ii) This follows from the *-transfer of the sentence
$$\forall x\forall y(x \in A_1 \land y \in A_2 \to ((x,y) \in  R 
\leftrightarrow (y,x) \in R^{-1})).$$\par\vskip 16pt
{\bf Proposition 3.} {\sl If internal $a\in \hyper {\cal H}$ and $b \in X_0,$ 
then $a \notin b.$} \par
Proof. This is often proved directly during the construction of $\hyper {\cal 
M}$ but it also follows by *-transfer of each sentence$,$ where $n\in \nat,$
 $$\forall x\forall y (x \in X_0 \land y \in X_n \to \neg (y \in x)).$$\par
\bigskip
{\bf Theorem 3.1.5.} {\sl Let $A$ be an entity. Then\par
\indent\indent (i) $\hyper {\{(x,x)\vert x\in A\}} = \{(x,x)\vert x \in \hyper 
A \}.$\par
\indent\indent (ii) $\hyper {\{(x,y)\vert x\in y\in A\}}= \{(z,w)\vert z\in 
w\in \hyper A\}.$}\par
Proof. (i) Let $I_A$ be the identity relation defined on $A.$  Then by *-
transfer of the sentence $$\forall z(z \in I_A \to \exists x(x \in A \land 
(x,x) = z))$$ it follows that $\Hyper (I_A) \subset I(\hyper A).$ On the other 
hand$,$ *-transfer of the sentence $$\forall x(x \in A \to (x,x) \in I_A)$$
yields that $I(\hyper A) \subset \Hyper (I_A)$ and thus $I(\hyper A)= \Hyper 
(I_A).$\par
(ii) Let $E= \{(x,y)\vert x\in y\in A\},$ where $A \in X_p.$ If $p\leq 1,$ then 
$E = \emptyset$ by Proposition 3 and result is obtained. Hence let $p\geq 2,\ 
a\in b\in A \in X_p \Rightarrow a\in b \in X_{p-1} \Rightarrow a\in X_0 \cup 
X_{p-2}.$ Now let $K =
\{(z,w)\vert z\in w\in \hyper A \}.$ The sentence 
$$\forall z(z \in \hyper E \to \exists x\exists y(x \in \hyper X_0 \cup
 \hyper X_{p-2} \land  y \in \hyper X_{p-1} \land $$ $$x \in y \land 
 y \in \hyper A  \land z = (x,y)))$$ holds in $\hyper {\cal M}.$
 Thus $\hyper E \subset K.$
On the other hand$,$ consider the *-transform $$\forall x\forall y (x \in 
 \hyper X_0 \cup \hyper X_{p-2} \land y \in \hyper X_{p-1}
\land x \in y \in \hyper A \to (x,y) \in \hyper E).$$ Now if $(c,d) \in K,$
then 
$c\in d\in \hyper A \Rightarrow  c\in  d \in \hyper X_{p-1} \Rightarrow
 a \in \hyper X_0 \cup \hyper X_{p-2}$ by 
Proposition 2. This implies that $K \subset \hyper E$ and completes the 
proof.\par
\vskip 18pt
In the next demonstration we meet for the first time what I call a partial *-
transfer. What this means is that for some arbitrary object within a set a 
sentence is discussed and *-transformed. After this  discussion the result
is obtained through application of the metalogic axiom of generalization.\par   
\vskip 18pt
{\bf Theorem 3.2.1.} {\sl Definition 3.2.1 (i.e. $S\to \vec r$) is equivalent to the
limit definition for convergence of a sequence 
(i.e. $\lim_{n\to \infty} S(n) = \vec r$).} \par
Proof. For the sufficiency$,$ let arbitrary $s \in {\real}^+$ and assume that
 $\lim_{n\to \infty} S(n)  = \vec r.$ Then $\exists N \in  \nat$
such that $\forall n > N,\ \Vert S(n) -  \vec r \Vert < s.$
Hence the sentence 
$$\forall x(x \in  \hypernat \land x > N \to \Vert \hyper S(x) -
\vec r\Vert <s)$$
holds in $\hyper {\cal M}.$ In particular let $\Gamma \in {\nat}_\infty.$ 
Then $\Vert \hyper S(\Gamma) - \vec r \Vert < s.$
But$,$ $\Gamma > N,\ \forall N \in 
\nat \Rightarrow \Vert \hyper S(\Gamma) - \vec r \Vert < s,\ \forall 
 s\in {\real}^+.$
Hence $\hyper S(\Gamma) \in \monad {\vec r}.$
 Since $\Gamma$ is also an arbitrary 
member of ${\nat}_\infty$ the sufficiency is obtained.\par
For the necessity$,$ assume that $\forall \Gamma \in {\nat}_\infty,\ 
 \hyper S(\Gamma) \in \monad {\vec r}.$ Let arbitrary
$s \in {\real}^+.$ Since
${\nat}_\infty \not= \emptyset$ then the statement 
$$\exists z(z \in \hypernat \land \forall x(x > z \land x \in \hypernat \to
 \Vert \hyper S(x)- \vec r \Vert < s))$$ holds in
 $\hyper {\cal M}.$ Thus by 
removing the * notation 
$$\exists z(z \in \nat \land \forall x(x > z \land x \in \nat \to
 \Vert S(x) - \vec r \Vert < s))$$ holds in ${\cal M}.$ 
But$,$ this is the limit 
definition for convergence and the proof is complete.
\vskip 18pt 
{\bf Theorem 3.2.2.} {\sl Let the sequence $S\colon \nat \to \realp n,
\ (n\geq 1).$ Then $S\to \vec r \in \realp n$ if and only if for each
$\Gamma \in \nat_\infty,\ \st {\hyper  S_\Gamma} = \vec r.$}\par
Proof. This follows from the definition of the standard part operator.
\vskip 18pt
{\bf Theorem 3.4.1 (The Leibniz Principle).} {\sl A sentence $\Phi$ holds true
for members of $\cal H$ if and only if the sentence $\Hyper {\Phi}$ holds true for 
members of $\hyper {\eskip {\cal H}}.$}\par                                       
Proof. The model theory portion of this proof 
  is not elementary and requires additional effort. 
However$,$ proofs can be found in {\tt Stroyan and Luxemburg [l976]$,$ Davis 
[1977].} \par From this all of our previous theorem involving *-transfer follow.
The proofs that the internal objects are closed under the basic
 set-theoretic operation is very tedious. There are obtained by basic
*-transfer and  
application of the previous results. The method is the same in  all  cases.
Here is one example that shows that if you have a statement like 
$$\Phi =\forall x (x \in A \to \exists y(y \in B \land x \subset y)),$$
then it holds in $\cal M$ if and only if the *-transfer holds in $\hyper {\cal M},$ where 
we do not * the relation symbol $\subset y.$ However, we must always translate this as stating in the *-transfer mode ``internal x a subset of internal y.'' This is so since this is an abbreviation for the statement $x \in (\hyper {{\cal P})}(y),$ where the power-set formation is considered as an operator.
Consider the binary relation $E = \{(x,y)\vert x \in A \land y \in B \land 
x \subset 
y\}.$ Then $\exists p \in \nat$ such that $E \subset X_p$ and this implies 
also that $E \subset X_p \Leftrightarrow \hyper E \subset \hyper X_p.$ By 
taking the *-transfer of 
$$ \forall z(z\in X_p \to (z \in E \leftrightarrow \exists x \exists y(
x \in A \land y \in B \land $$ $$\forall w (w \in X_0 \cup X_{p-2} \land 
 w \in x \to w \in y))),$$
applying Proposition 2 and the usual argument one obtains 
$\hyper E = \{(x,y)\vert x \in \hyper A \land y \in \hyper B \land x 
 \subset y\}.$ From this the *-transform may be written as
$$ \Hyper \Phi =\forall x (x \in \hyper A \to \exists y(y \in \hyper B \land x 
 \subset y)).$$
Thus as long as the *-transfer is applied to sentences where the quantifiers 
are bounded (or assumed bounded) by standard entities then objects are closed 
under the basic \underbar{binary} set-theoretic operators. This even includes
the more
general $\bigcup X$ and $\bigcap X$ where $X$ is an internal entity. However$,$ 
it does not include all set-theoretic operators. In particular$,$ if $A$ is an 
infinite entity$,$ then ${\cal P}(\hyper A)$ is not internal although $\hyper {\power {A}} = \Hyper {(\power {A})}= (\hyper {\cal P})(\hyper {A})$. As seen by the proof of Theorem 3.4.3 every element of $\Hyper {(\power {A})}$ is an entity and *-transfer of the characterizing statement it is also a subset of $\hyper A$. Thus in those
set-theoretic expressions where the power set operator is used$,$ it must
carry a * notation under *-transfer. \par  
\vskip 18pt
{\bf Theorem 3.4.2.} {\sl The set of infinite natural numbers ${\nat}_\infty =
\hypernat - \nat.$}\par
Proof.  Previously proved in Example 3.4.3.\par
\vskip 18pt
{\bf Theorem 3.4.3. (Leibniz' Principle Restated).} {\sl  A sentence $\Phi$
holds true
for members of $\Re = {\cal U} \cup {\cal P}({\cal U})$ if and only if the
sentence $\Hyper {\Phi}$ holds true for the members of $\hyper {\eskip \cal U}$
or internal subsets of $\hyper {\eskip \cal U}.$}\par    
Proof. Observe that $\Hyper \Re = \Hyper {\cal U} \cup \hyper {\cal P}({\cal 
 U}).$ The set ${\cal U} \in X_p \Rightarrow {\cal U} \subset X_0 \cup X_{p-1}
\Rightarrow {\cal P}({\cal U}) \subset {\cal P}(X_0 \cup X_{p-1}) = X_p 
\Rightarrow {\cal P}({\cal U}) \in X_{p+1}.$  Hence, $\Hyper {( \power {\cal U})} \in \hyper {X_{p+1}}.$ Since, by *-transfer, no internal set is a member of $\hyper X_0$, then by Proposition 2 (i) each member of $\hyper {\cal P}({\cal 
 U})$ is an internal subset of ${\hyper {\cal U}}$
and the proof is complete.
\vskip 18pt                     
{\bf Theorem 3.5.1.} {\sl Let $S\colon \nat \to \real.$ If $S$ is an 
increasing {\rm [}resp. decreasing{\rm ]} sequence and there exists some
$\Gamma \in {\nat}_\infty$ such that $\hyper S_\Gamma \in \cal O$$,$ then $S
\to \st {\hyper S_\Gamma.$}}\par
Proof. Assume that $S\colon \nat \to \real$ is increasing. Then by *-transfer
$\hyper S\colon \hypernat \to \hyperreal$ is increasing. Hence for $\Gamma
 \in {\nat}_\infty$ and for $m \in \nat$ it follows since $\Gamma > m$ that
$S(m) \leq \hyper S(\Gamma).$ Consequently$,$ since $\hyper S(\Gamma) \in \cal 
 O$ the $\st {S(m)} = S(m) \leq \st {\hyper S(\Gamma)} = r \in \real.$ Thus the 
sentence 
$$ \forall x(x \in \nat \to S(x) \leq r)$$
holds in $\cal M$; hence in $\Hyper {\cal M}.$  Thus  $\forall \Omega \in 
 {\nat}_\infty,\ \hyper S(\Omega) \leq r \Rightarrow \forall \Omega \in 
 {\nat}_\infty,\ \hyper S(\Omega) \in \cal O.$ Hence $\forall \Omega  \in 
 {\nat}_\infty,\ \st {\hyper S(\Omega)} \in \real.$ So$,$ let arbitrary 
$\Omega \in 
{\nat}_\infty,\ \Omega \not= \Gamma,\ \st {\hyper S(\Omega)} = r_0.$
 Then the same sentence as 
above with $r_0$ replacing $r$ yields that $\hyper S(\Gamma) \leq r_0 
\Rightarrow \st {\hyper S(\Gamma)} = r \leq r_0.$ Consequently$,$ $\st {\Hyper 
 S(\Omega)} = r.$ Since $\Omega$ is arbitrary then $\hyper S(\Omega)
 \in \monad r,\ \forall \Omega \in {\nat}_\infty \Rightarrow S \to r.$\
\vfil\eject
\centerline{Appendix For Chapter 4.}
\vskip 18pt
{\bf Theorem 4.1.1.} {\sl Definition 4.1.1 is equivalent to the limit 
definition for a sequence $S \to +\infty.$ }\par
Proof. First$,$ the method of partial transfer is used. Assume that 
$\lim_{n \to \infty} S_n = +\infty.$ Let arbitrary
$B \in {\real}^+.$ Then from the classical definition there exists some $M_B \in 
\nat$ such that the sentence $$\forall x(x \in \nat \land x \geq M_B \to S_x 
\geq B))$$ holds in $\cal M.$ By *-transfer this particular statement implies 
since $\forall \Gamma \in {\nat}_\infty,\ \Gamma >M_B$ that $\hyper S_\Gamma \geq 
B.$ No matter what $B \in {\real}^+$ is chosen the same conclusion for each
 $M_B$ Thus in all cases $\forall \Gamma \in {\nat}_\infty,\ \hyper S_\Gamma 
\geq B.$ But $B$ is any arbitrary positive real number. Thus 
$\forall \Gamma \in {\nat}_\infty,\ \hyper S_\Gamma \in 
{\real}_\infty^+.$\par
For the converse$,$ assume that $\forall \Gamma \in {\nat}_\infty,\ 
 \hyper S_\Gamma \in {\real}_\infty^+.$ Let $B \in {\real}^+$ and $\Omega \in
 {\nat}_\infty$ which exists since ${\nat}_\infty \not= \emptyset.$ Obviously$,$ 
if $\Lambda \geq \Omega,$ then  $\Lambda \in {\nat}_\infty \Rightarrow
 \hyper S_\Lambda > B.$ Hence$,$ the sentence $$\exists x(x \in \nat \land
 \forall y(y \in \nat \land y \geq x \to S_y > B))$$ holds in $\Hyper 
{\cal M}$ and 
thus holds in $\cal M.$ But this is the classical definition and the result 
follows.\par
\vskip 18pt
{\bf Theorem 4.1.2.} {\sl Let $r \in {\hyperreal}^+,\ S\colon \nat \to \real$
 and $S\to +\infty.$ Then for each
$\Omega \in {\nat}_\infty$ there exists some $\Delta \in {\nat}_\infty$ such 
that $\hyper S(\Omega) + r \leq \hyper S(\Delta).$} \par
Proof. From Theorem 4.1.1 simply consider the *-transfer of the complete 
definition for $\lim_{n\to \infty} S(n) = +\infty.$  Take $\hyper S(\Omega) + 
r$ as the bound. Then $\exists M \in \hypernat$ such that $\forall \Gamma \in
 \hypernat,\ \Gamma \geq M \Rightarrow \hyper S(\Gamma) \geq \hyper S(\Omega) + 
r.$ Obviously$,$ if $M\notin {\nat}_\infty,$ then any member of ${\nat}_\infty$
will suffice and this completes the proof. 
\vskip 18pt
{\bf Theorem 4.2.1. (The Extended Standard Definition Principle)}\par 
 {\sl \indent\indent (i) A set $A \in \Hyper {\cal H}$ is an extended standard 
set (i.e. there exists some $B \in {\cal H}$ such that $A = \hyper B$) if 
and only if there exists some standard set $D$ (i.e. $D \in {\cal H}$) and a standard bound formula 
$\Phi (x)$ in one free variable where each constant in $\Phi (x)$ is a member
of $C({\cal  H})$ and
$$ A = \{x\vert x \in \hyper D \land \Hyper {\Phi}(x) \}.$$ \par
\indent\indent (ii) A set $A$ is an extended standard n-ary relation ($n > 1$)
if and only if there exist n standard sets $D_1,\ldots,D_n$ and a standard bound 
formula $\Phi (x_1,\ldots,x_n)$ in n free variables where each constant in 
$\Phi (x_1,\ldots,x_n)$ is a member of $C({\cal H})$ and 
$$A = \{(x_1,\ldots,x_n)\vert
 x_1 \in \hyper D_1 \land \cdots \land x_n \in \hyper D_n \land \Hyper 
 {\Phi (x_1,\ldots,x_n)}\}.$$}\par
Proof. (i) Let $\Phi (x)$ be any bound formula in which all constants are 
members of $C({\cal H}).$ Let $B = \{x\vert x \in D \land \Phi (x) \},$ 
where $D$ and all 
constants in $\Phi (x)$ represent members of $\cal H.$ In our general set 
theory the set  $A$ exists. Let $a_1,\ldots,a_m$ be the list of all constants 
in $\Phi (x).$ Then there is a finite set $\{X_i\}$ such that each $a_i$ is a
member 
of some $X_i,$ and some $X_p,\ p\geq 1$ such that $D \in X_p$ and since $B \subset D$ then 
from Proposition  1 (vi) $B \in X_p.$ Observe that this also implies that the 
formula $\Phi (x)$ is defined within the superstructure $\cal H$ (i.e. the 
vocabulary of $\Phi (x)$ is interpretable in $\cal H$). Even though once it is 
understood that $\Phi (x)$ has meaning in $\cal H$ the following expanded form 
is not considered necessary$,$ it is presented in preparation for the proof of 
Theorem 4.2.2. The sentence $$a_1 \in X_1\land \cdots \land a_m \in X_m
\land D \in X_p \land B \in X_p \land \forall x(x \in B \leftrightarrow
(x \in D) \land \Phi (x))$$
holds in $\cal M$ if and only if  the *-transform holds in $\Hyper {\cal M}.$ 
Consequently$,$ $A = \hyper B = \{x\vert x \in \hyper D \land \Hyper 
{\Phi (x)}\}.$\par
Conversely$,$ let $A = \hyper B.$ Simply consider the formula
$\Phi (x) = x \in B.$ Then $\hyper B \in \hyper X_p$ and $\hyper B$ be a set.
If $x \in \hyper B,$ then $x \in \hyper X_0 \cup \hyper
 X_{p-1}.$ Further$,$ $X_0 \cup X_{p-1} \in X_p \Rightarrow \hyper X_0 \cup 
\hyper  X_{p-1} \in \hyper X_p.$ Thus $A =\hyper B = \{x\vert x \in \hyper X_0 
 \cup \hyper  X_{p-1} \land \Hyper \Phi (x) \}$ from the first part. \par
(ii) This proof is but a simple modification of part (i). 
Consider the formula $\Phi (x_1,\ldots,x_n),\ n>1$ in $n$ free variables and assume as 
before that the non-quantifier bounding constants represent members of 
$\cal H.$ Adjoin the fact that each $D_i \in X_{p(i)}$ and that there exists 
some $X_k$ such that $B \in X_k,\ k \geq 1$ to the other facts.
Then the sentence $$a_1 \in X_1\land \cdots \land a_m \in X_m \land
 D_1 \in X_{p(1)} \land \cdots \land D_n \in X_{p(n)} \land B \in X_k \land$$
$$\forall x_1 \cdots \forall x_n (x_1 \in 
D_1 \land \cdots \land x_n \in D_n \to \forall z(z \in B \leftrightarrow $$ $$z = (x_1,\ldots,x_n) \land \Phi 
(x_1,\ldots,x_n)))$$ holds in $\cal M$ if and only if  the *-transform holds in $\Hyper
 {\cal M}.$ Consequently$,$ $A = \hyper B = \{(x_1,\ldots,x_n)\vert
 x_1 \in \hyper D_1 \land \cdots \land x_n \in \hyper D_n \land \Hyper 
 {\Phi (x_1,\ldots,x_n)}\}.$\par
The converse follows in the same manner as in part (i) and the proof is 
complete.\par \medskip
 {\bf Corollary 4.2.11.1} {\sl For the $\Phi$ in Theorem 4.2.1, if $B= \{x\vert x \in D \land {\Phi}(x) \}$ [resp. 
$ \{(x_1,\ldots,x_n)\vert
 x_1 \in D_1 \land \cdots \land x_n \in D_n \land  
 {\Phi (x_1,\ldots,x_n)}\}$], then $\hyper B = \{x\vert x \in \hyper D \land \Hyper {\Phi (x)}\}$ [resp. $\{(x_1,\ldots,x_n)\vert
 x_1 \in \hyper D_1 \land \cdots \land x_n \in \hyper D_n \land \Hyper 
 {\Phi (x_1,\ldots,x_n)}\}.$} \par
\vskip 18pt
\baselineskip=12pt
{\bf Theorem 4.2.2. (The Internal Definition Principle)}\par 
 {\sl \indent\indent (i) A set $A$ in our set theory is an internal 
set if and 
only if there exists some internal set $D$ and an internal bound formula $\Phi (x)$ in 
one free variable where each constant in $\Phi (x)$ is a member of 
$C(\Hyper {\cal H})$ and
$$ A = \{x\vert x \in D \land \Phi (x) \}.$$ \par
\indent\indent (ii) A set $A$ is an internal n-ary relation ($n > 1$)
if and only if there exist n internal sets $D_1,\ldots,D_n$ and a internal bound 
formula $\Phi (x_1,\ldots,x_n)$ in $n$ free variables where each constant in 
$\Phi (x_1,\ldots,x_n)$ is a member of $C(\Hyper {\cal H})$ and
 $$A = \{(x_1,\ldots,x_n)\vert
 x_1 \in D_1 \land \cdots \land x_n \in D_n \land \Phi (x_1,\ldots,x_n) \}.$$}\par
Proof. (i) Let $a_1,\ldots,a_m$ be the
list of all constants 
in $\Phi (x).$  Then each $a_i$ is a member of some $\hyper X_i$ and there exists
$X_p$ such that $D \in \hyper X_p.$ Now replace every distinct $a_i$ at each
occurrence 
with a distinct variable $y_i$ yielding the n+1 placed formula 
$\Phi (x,y_1,\ldots,y_m).$ The following sentence
$$\forall y_1 \cdots \forall y_m \forall z(y_1 \in X_1 \land \cdots \land
 y_m \in X_m \land z \in X_p \to \exists u(u \in X_p \land$$
$$\forall x(x \in X_0 \cup X_{p-1} \to (x \in u \leftrightarrow x \in z \land
 \Phi (x,y_1,\ldots,y_m)))))$$ 
holds in $\cal M$ if and only if  it holds in $\Hyper {\cal M}.$ Consequently$,$ by 
considering the specific internal objects denoted by $a_i$ and by $D.$ Then 
the corresponding internal $u$ that exists in $\hyper X_p$ is precisely the 
set  $\{x\vert x \in D \land \Phi (x) \}.$ \par
The converse follows from Proposition 2 in the same manner as in the proof of 
part (i) of Theorem 4.2.1.\par
(ii) For this part follow the same procedure as in the second part of the 
proof of Theorem 4.2.1. Modify the appropriate subformula of the formal 
sentence that appears in the proof of this theorem part (i) in the same manner 
as the sentence was modified in the proof of Theorem 4.2.1 part (ii). With 
this the proof is complete.\par
\baselineskip=14pt
\vskip 18pt
{\bf Theorem 4.3.1.} {\sl Let $n > 1.$ Then for each internal $R \subset
 \hyper A_1 \times \cdots \times \hyper A_n$ and for each 
$i,\ 1\leq i \leq n$ it follows that $\hyper P_i(R) = P_i(R),$ where 
$P_i(R)$ is an internal subset of $\hyper A_i.$}\par
Proof. Observe that for each $i,\ 1\leq i\leq n$ the {\it i}th projection may be 
considered a map $P_i\colon {\cal P}(A_1 \times \cdots \times A_n) \to 
{\cal P}(A_i).$ Thus $\hyper P_i$ is a map from the internal subsets of 
$\hyper A_1 \times \cdots \times \hyper A_n$ to the internal subsets of 
$\hyper A_i.$ Hence$,$ if internal $Q \subset \hyper A_1 \times \cdots \times 
 \hyper A_n,$ then $(\hyper P_i)(Q) \subset \hyper A_i$ is internal.
   Further$,$ in general, there exists some $X_p$ such that if $B \subset
A_1 \times \cdots \times A_n,$ then $B \in X_p.$ The remainder of this proof 
is by induction on the number of Cartesian products$,$ $m,$ as was done in the proof 
of Theorem 3.1.4. Let $m=1.$ The following sentence
$$\forall z(z \in \power {A_1 \times A_2} \to \forall x(x \in A_1
 \to $$ $$ (x \in P_1(z) \leftrightarrow \exists y(y \in A_2 \land (x,y) \in 
z))))$$ holds in $\cal M$ if and only if  it holds in $\Hyper {\cal M}.$ Suppose that
internal $Q \subset \hyper A_1 \times \hyper A_2.$ Let $q \in (\hyper P_1)(Q)$ 
Then $q \in \hyper A_1 \Rightarrow \exists s \in \hyper A_2$ such that
$(q,s) \in Q \Rightarrow (\hyper P_1)(Q) \subset P_1(Q),$ where $P_1$ is a set-theoretic projection for n-ary relations $Q \in \hyper {\cal H}.$ On the other hand$,$ 
let $q \in P_1(Q).$ Then $q \in \hyper A_1$ and $\exists s \in \hyper A_2$ 
such that $(q,s) \in Q.$ This yields$,$ from the *-transform$,$ that $q \in 
(\hyper P_1)(Q)$ and$,$ hence$,$ $(\hyper P_1)(Q) = P_1(Q),$   The remainder of this 
proof follows in exactly the same manner$,$ appropriately modified$,$ as does the 
last part of the proof of Theorem 3.1.4.\par
\vskip 18pt
{\bf Theorem 4.3.2.} {\sl Let $F(B)$ be the set of all finite subsets of $B$. 
Then nonempty internal $A \subset \hyper B$ is hyper finite if and only if $A \in (\hyper 
F)(B)$, where if there can be any confusion, we continue to denote the *-transform of such operators as $(\hyper F).$}\par
Proof. As mentioned, the notion used for finite sets is a first-order statement in our set-theory with ``finitely'' many $\lor$ symbols and constants. This is not the mapping notion used for the definition of hyper finite. However, it is a simple exercise to show by induction that since $\nat \subset X_0$ that for any nonempty finite set $A \in {\cal H},$ there exists in $\cal H$ a bijection $f \colon [0,n] \to A.$ Then $f$ corresponds within our set-theory to a set that can be characterizes as having the same ``number'' of $\lor$ symbols and constants when the range is completely described. Consequently, whenever it is necessary, we characterize finite sets by such mappings. Next$,$ it is necessary to establish by induction that the mapping concept used for 
the hyperfinite set definition is equivalent to the *-transform of the 
standard set-theoretic nonempty finite set mapping correspondence. We show that 
$\forall n \in \nat$ if $g \in B^{[0,n]}$ and $g[[0,n]] = A,$ then $\exists
 k_n\in \nat$ and an injection $f\in B^{[0,k_n]}$ such that $f[[0,k_n]] = A.$
For $n = 0$ then $g$ itself is such an injection. Suppose the hypothesis holds
for an arbitrary $n\in \nat.$ Let $g \in B^{[0,n+1]}$ and $g[[0,n+1]] = 
A^\prime$ 
Our set theory states that there exists a map $g_1 = g\vert [0,n]$ and 
$g_1[[0,n]] \subset A^\prime.$ By induction $\exists k_n \in \nat$ and an injection
$f_1 \in B^{[0,k_n]}$ such that $f_1[[0,k_n]] = g_1[[0,n]].$ If $g_1[[0,n]] = 
A^\prime,$ then $f_1$ suffices.  If $g_1[[0,n]] \not= A^\prime,$ then $g(n+1) \notin 
g_1[[0,n]] = f_1[[0,k_n]]$ and $A^\prime = g_1[[0,n]] \cup \{g(n+1)\}.$ In this case$,$ 
simply define $f = \{(k_n +1,g(n+1))\} \cup f_1.$ Then $f$ is an injection 
defined on $[0,k_n +1]$ onto $A^\prime$. The hypothesis holds by induction. 
Consequently$,$ if $\exists n \in \nat$ such that $g \in B^{[0,n]},$ then 
$g[[0,n]] = A$ is finite in the standard set-theoretic mapping sense. The converse is 
obvious and this implies that our slightly modified mapping definition is 
equivalent to the *-transform of the standard mapping ``definition'' for nonempty finite 
sets.\par
For a given entity $B$ there exists some $X_p$ such that $F(B) \in X_p$ 
Let ${\cal A} = {\cal P}({\nat})$ and ${\cal B} = \{B\}.$ Then as discussed 
just prior to example 4.3.1 there exists some $q \in \nat$ such that for each 
$A \in {\cal A}$ and $B \in {\cal B},\ B^A \in X_q$ and$,$ hence$,$ since the 
members of $B^A$ are not individuals each $f \in B^A$ is a member of $X_{q-
1}.$ From the above proof it follows that the follow sentence
$$\forall x(x \in X_p \to (x \in F(B) \leftrightarrow \exists y\exists z
(y \in \nat \land z \in X_{q-1} \land $$ $$z \in B^{[0,y]} \land z[[0,y]] = x)))$$
holds in $\cal M$ if and only if  it holds in $\Hyper {\cal M}.$ The result follows by *-
transfer.\par
\vskip 24 pt
{\bf Theorem 4.3.3} {\sl Any nonempty finite set of internal individuals or entities is 
internal and hyperfinite.}\par
Proof. This follows from Proposition 2 (iii) and Theorem 4.3.2.\par

\vskip 18pt
{\bf Theorem 4.3.4.} {\sl Let $A \in \Re$ and assume that $A$ is infinite.
Then there exists a hyperfinite set $F$ such that $F \not = A,\ F \not= 
\hyper A$  and  $A \subset F \subset \hyper A.$}\par     
Proof. Let infinite set $A \in X_p.$ Consider the standard binary 
relation $R = \{(x,y)\vert x \in y \land y \in F(A)\
 \land x \in A\}.$ Observe that the domain of $R$ is the set $A.$ Consider
$\{(x_1,y_1),\ldots,(x_n,y_n)\} \subset R.$ Then the set $y_1 \cup \cdots \cup 
y_n \in F(A)$ and $\{(x_1,y),\ldots,(x_n,y)\} \subset A.$ Thus $R$ is concurrent 
on $A$. Thus there exists some $F \in \hyper F(A)$ such that $\forall x \in 
^{\sigma}A,\ x \in F.$ Since $A \in \Re$ then $^{\sigma}A = A \subset F$ and 
by *-transfer $F \subset \hyper A.$ Since $A$ is not internal and $F$ is 
internal then $F \not= A.$ Further$,$ if $F = \hyper A,$ then the sentence
$$\exists x(x \in F(A) \land x = A)$$ holds in $\cal M$ since it holds in 
$\Hyper {\cal M}.$ This would contradict the fact that $A$ is not finite and 
the proof is complete.\par
\vskip 18pt
{\bf Theorem 4.3.5.} {\sl An infinite series$,$ $\sum_{i=0}^\infty a_i,$ converges 
to $r \in \real$ if and only if for each $\Gamma \in {\nat}_\infty$ it follows 
that $\sum_{i=0}^\Gamma a_i \in \monad r.$}\par
Proof. This is an immediate consequence of Theorem 3.2.1.\par
\vskip 18pt
{\bf Theorem 4.4.1.} {\sl Definition 4.4.1 for continuity and uniform 
continuity is equivalent to the $\delta - \eps$ definition.}\par 
Proof. (i) (Continuity). Assume Definition 4.4.1. Thus for nonzero $n,m \in 
\nat,$ nonempty $A \subset \realp n,$ the mapping $f\colon A \to \realp m$ 
has the property at $p\in A$ that $\hyper f[\monad p \cap \hyper A] \subset 
\monad {f(p)}.$ Let $ r\in {\real}^+.$ If $p$ is an isolated point$,$ then 
$\monad p \cap \hyper A = \{p\}$ and $\hyper f[\monad p \cap \hyper A] = \{p\}
\subset \monad {f(p)}.$ Suppose that $p$ is not isolated. Then $\exists 0\not= 
\eps \in \monad 0$ and $p + \eps = q \in \monad p \cap \hyper A.$ Moreover$,$ 
$0 < \Vert q-p \Vert < \Vert \eps \Vert$ and if $s \in \hyper A$ and 
$0 < \Vert s - p \Vert <\Vert \eps \Vert,$ then $s \in \monad p \Rightarrow
\hyper f(s) \in \monad {f(p)}.$ Combining these two cases it follows by *-
transfer that the sentence
$$\exists x(x \in {\real}^+ \land \forall y(y \in A \land 0\leq \Vert p -y 
\Vert < x \to \Vert f(y) - f(p) \Vert <r))$$
holds in $\cal M.$ Hence$,$ $f$ is continuous at $p \in A.$\par
For the converse$,$ let $r \in {\real}^+.$ Then we know that there exists some 
$r_0 \in {\real}^+$ such that 
$$ \forall x(x \in A \land 0\leq \Vert x - p \Vert <r_0 \to \Vert f(x) - f(p) 
\Vert < r)$$
holds in $\cal M.$ Hence the *-transfer holds in $\Hyper {\cal M}.$
Now let $q \in \monad p \cap \hyper A.$ Then $\Vert q - p \Vert <r_0$ for any
$r_0 \in {\real}^+.$ Thus for any $r \in {\real}^+,\ \Vert \hyper f(q) -f(p) 
\Vert < r.$ Since $r$ is an arbitrary positive real number then it follows 
that for any $q \in \monad p \cap \hyper A,\ \hyper f(q) - f(p) \in \monad 0$ 
and $\Rightarrow \hyper f[\monad p \cap \hyper A] \subset \monad {f(p)}.$\par
(ii) (Uniform continuity). Assume Definition 4.4.1. Thus for the map $
f\colon A \to \realp m$ and if $p,q \in \hyper A,\ p -q \in  \monad 0,$ then 
$\hyper f(p) -\hyper f(q) \in \monad 0.$ Suppose that arbitrary $r \in 
{\real}^+.$ We know that there exists a nonzero $\eps \in \monad 0$ and for any
$p,q \in \hyper A$ such that $\Vert p - q \Vert < \eps,\  \hyper f(p) -\hyper 
 f(q) \in \monad 0 \Rightarrow \Vert \hyper f(p) - \hyper f(q) \Vert <r.$ 
Hence$,$ the sentence
$$ \exists x(x \in {\real}^+ \land \forall y \forall z(y \in A \land z \in A 
\land \Vert y - z \Vert <x \to \Vert f(y) - f(z) \Vert <r))$$ holds in 
$\cal M.$ Thus $f$ is uniformly continuous on $A.$\par
For the converse$,$ assume uniform continuity and let arbitrary $r \in 
{\real}^+.$ Then we know that there exists some $r_0 \in {\real}^+$ such that 
for each $x$ and each $y$ such that $x,y \in A$ and $\Vert x - y \Vert < r_0$ 
then $\Vert f(x) - f(y) \Vert <r.$ By *-transfer we have that for any
$p,q \in \hyper A$ and $ p - q \in \monad 0$ it follows that $\Vert \hyper f(p) 
- \hyper f(q) \Vert < r.$ Since $r$ is an arbitrary  nonzero positive real 
number then $\hyper f(p) - \hyper f(q) \in \monad 0.$ This completes the proof. 
\vskip 18pt
{\bf Theorem 4.4.2.} {\sl Definition 4.4.2 for compactness is equivalent to 
the standard definition utilizing open covers.}\par
Proof. Even though there are proofs of this that$,$ in this case$,$ do not 
require the enlargement property (see {\tt Herrmann [l978]}) the enlargement 
property leads to a much shorter demonstration. \par
Suppose $A$ is compact and that there is a $q \in \hyper A$ such that 
$q \notin \bigcup \{\monad r\vert r \in A\}.$ We use the characterization 
that $\monad r = \bigcap \{\Hyper G\vert r\in  G \in  {\cal T} \},$ 
where $\cal T$ is the set of all open 
subsets of $\realp n.$ Thus for each $r \in A$ there exists some $G_r \in 
{\cal T}$ such that $q \notin \Hyper G_r.$ Since $\{G_r\vert r \in A\}$ is an 
open cover for $A$ then there exists a finite $\{G_1,\ldots,G_n\} \subset
\{G_r\vert r \in A\}$ such that $A \subset G_1 \cup \cdots \cup G_n.$ Hence,
$\hyper A \subset \Hyper G_1 \cup \cdots \Hyper G_n \Rightarrow \exists i,\ 
1\leq i\leq n$ such that $q \in \Hyper G_i \in \{G_r\vert r \in A\}.$ But this 
contradiction implies that standard compactness yields Definition 4.4.2.\par
Conversely$,$ assume that $A$ is not compact and let $\cal G$ be an open cover 
of $A$ that does not contain a finite subcover. As usual all of these objects 
are members of $\cal H.$ Consider the 
binary relation $$R = \{(x,y)\vert x\in {\cal G} \land y \in A \land y \notin 
x\}.$$ Let $\{(x_1,y_1),\ldots,(x_n,y_n)\} \subset R.$ Then there exists some 
$q \in A$ such that $q \notin x_i$ for any $i,\ 1\leq i\leq n.$ Hence$,$ $R$ is 
concurrent and the domain of $R = {\cal G}.$ Thus there exists some $b \in 
\hyper A$ such that $(x,b) \in \hyper R$ for each $x \in {^\sigma{\cal G}}$. 
But $x \in {^\sigma{\cal G}} \leftrightarrow x = \Hyper G,\ G \in {\cal G}.$ 
Consequently$,$ since $\cal G$ contains at least one open neighborhood for each 
$r \in A$ then $b \notin \bigcup \{\monad r\vert r \in A\}$ and this completes 
the proof. \parm
\centerline{\bf NOTATION}\parm
Notice that the notation has been generalized slightly. So that there was no possibility of confusion when not carefully read, hyper-function notation such as $\hyper f$ has been presented in two forms. The $(\hyper f) (A)$, to identify the function $\hyper f$ and $\hyper f(A)$, meaning $\Hyper (f(A))$ when $A$ is a standard object. The notational form $(\hyper f) (A)$ is unnecessary since which meaning applies depends upon the argument. In this example, when $A$ is standard, one needs only write $\hyper f(A)$. It is obvious that this notation has only the one possible meaning $\Hyper (f(A))$ for in this form the hyper-function is also being displayed since $\Hyper (f(A)) = \hyper f(\hyper A).$  On the other hand, if $A$ is simply  specified as an internal object, then the notation $\hyper f(A)$ can only mean $(\hyper f)(A).$ The notations are equivalent if $A = \hyper B$. In some cases under our identification process, the notations are also equivalent where we let $\hyper A = A,$ where $A$ is identified with the standard object via the constants sequences of atoms (urelements). \par  
\vfil\eject\centerline{Appendix For Chapter 5.}
\vskip 18pt
{\bf In all that follows the set of natural numbers$,$ $\nat,$ is assumed to 
contain zero and all functions from a  domain that is a subset of $\realp n$ 
into $\real$ are bounded.}\par
For an n-dimensional space $n \geq 1$$,$ the closed set $R = [a_1,b_1] \times 
\cdots \times [a_n,b_n],\  a_i < b_i,\ 1\leq i \leq n$ is called 
{\bf a {rectangle}.} Of course$,$ if $n = 1,$ then a ``rectangle'' is but a 
closed interval. As usual$,$ consider for each $[a_i,b_i]$ a {\bf {partition}} 
 $P_i$ as a finite set of members of $[a_i,b_i]$ such that $a_i,b_i \in P_i$ 
and where $P_i$ is considered as ordered. This is often explicitly written as 
$P_i = \{x_{i0},\ldots, x_{ik}\},\  a_i = x_{i0} < x_{i1} <\cdots<  x_{ik} = b_i.$ 
This determines the closed one-dimensional subintervals 
 $[x_{i(p-1)}, x_{i{p}}],\  1\leq p \leq k.$ In brief$,$ this 
process obtains a partition $P = P_1 \times \cdots \times P_n$ of $R$ and 
a finite collection of closed n-dimensional subrectangles  $S$ obtained by 
considering $([x_{10}, x_{11}] \cup \cdots \cup [x_{1{k-1}},x_{1k}]) \times 
\cdots \times ([x_{n0}, x_{n1}] \cup \cdots \cup [x_{n{m-1}},x_{nm}]).$ Each 
$S$ has a measure$,$ $v(S)\in \real$$,$ assigned to it which  
is intuitively the product of the lengths of the sides. For simplicity of 
notion the definition of the measure $v(S)$ is left intuitively 
understood.\par
Probably the simplest partition to consider would be the one termed a 
{\bf {simple}} partition. These are formed by selecting $n$ nonzero natural 
numbers $m_1,\ldots m_n$ and dividing each interval $[a_i,b_i]$ into an equal 
length partition by adding to each successive partition point the number
$(b_i -a_i)/(m_i).$ This concept is extended to the nonstandard world by 
selecting $n$ infinite natural numbers $\Gamma_1,\ldots\Gamma_n$ and 
generating for each interval $[a_i,b_i]$ an internal hyperfinite partition,
$P_i,$ each subinterval of which has positive infinitesimal length 
$(b_i -a_i)/(\Gamma_i)= dx_i.$ Then the partition $P= P_1 \times \cdots 
\times P_n$ is a {\bf {simple fine 
partition}} of $R.$ Such a partition yields an internal set of hyperrectangles 
$R_q$ such that 
$\hyper v(R_i) = dx_1 \cdots dx_n=dX  \in \monad 0^+.$ 
You could be much more general and consider the {\bf {fine}} partitions which 
are internal collections of hyperfinitely many members of $\Hyper [a_i,b_i]$
such that the length of any subinterval is an infinitesimal.  
Let $P$ be a partition of the rectangle $R$ and assume that $P$ determines 
the set of subrectangles $\{R_q\vert 1\leq q \leq m\}.$  An {\bf {intermediate 
partition}}$,$ $Q$$,$ is any finite sequence of vectors $\{\vec v_q \},$ where 
$\vec v_q \in R_q$ for each $q$ such that $1 \leq q \leq m.$
Let $\cal P$ be the set of all simple partitions of $R$ and ${\cal P}_S$ any 
nonempty subset of $\cal P.$ Then there exists a 
mapping $\cal Q$ on the set ${\cal P}_S$ such that for each $P\in {\cal P}_S$
the image ${\cal Q}(P)$ is the set of all intermediate partitions for $P.$  
From a notational view point subrectangles are denoted by various symbol 
strings. In particular$,$ $S,\ R_q,\ R_S.$\par
For \underbar{any} partition $P$ of $R$ let ${\cal S}(P)$ denote the set of 
all subrectangles generated by $P.$ If $S \in {\cal S}(P),$ then let
$\Vert S\Vert$ denote the length of the diagonal of the subrectangle $S.$
As usual$,$ define the {\bf mesh($P$)}$= \Vert P\Vert =\max \{\Vert S\Vert \vert S \in {\cal 
S}(P)\}.$ Now let $L(f,P) = \sum_{S\in {\cal S(P)}} m_Sv(S)$ be an {\bf lower 
sum} and $U(f,P) = 
\sum_{S\in {\cal S(P)}} M_Sv(S)$ be an {\bf upper sum}$,$ where $m_S 
= \inf \{f(\vec x)\vert \vec x \in S\}$ and $M_S = \sup \{f(\vec x)
\vert \vec x \in S\}$ As is well-known$,$ for any set of partitions of $R$ the 
set of lower and upper sums is a bounded set. Further$,$ let $L(f) =\sup \{L(f,P)\vert P\ {\rm is\ 
 any\ partition\ of\ }R\}$  and  $U(f) =\inf \{U(f,P)\vert P\ {\rm is\ 
 any\ partition\ of\ }R\}.$  The function $f$ is {\bf {Darboux}} integrable if
$L(f) = U(f)$ and the value of this integral is $D\int_R f(\vec x)\, dX = L(f) 
= U(f).$ As is well-known $f$ is Darboux integrable if and only if $f$ is Riemann 
integrable in the sense of Riemann Sums and has the same value. Now the 
concept of the mesh extends to the *-mesh and for any internal partition $Q$ of 
$\hyper R$ the $\Hyper \Vert Q \Vert$ exists by *-transfer of the finitary 
statement dealing with maximum values that exist within a finite set of real 
numbers.  Also$,$ $Q$ is a fine partition if and only if $\Hyper \Vert Q \Vert \in \monad 
0.$ \par
Notice that if ${\cal P}_S$ is any nonempty subset of the set of all 
partitions of $R,$ then the lower and upper sums we may consider maps 
$L(f,\bullet)\colon {\cal P}_S \to \real$ and $U(f,\bullet)\colon 
{\cal P}_S \to \real.$ These maps extend to the maps $\hyper L(f,\bullet)
\colon \hyper {\cal P}_S \to \hyperreal$ and $\Hyper U(f,\bullet)\colon 
\hyper {\cal P}_S \to \hyperreal.$    
\vskip 18pt     
\hrule
\smallskip
{\bf Definition A5.1.} Let bounded $f\colon R \to \real$ and ${\cal P}_I$ be any nonempty 
set of hyperfinite partitions of $\hyper R$ such that $\Hyper \Vert P \Vert
\in \monad 0$ for each $P \in {\cal P}_I.$ Then $f$ is said to be 
{\bf {Z-integrable}} for ${\cal P}_I$ if there exists some $P \in {\cal P}_I$ such that 
$$\st {\hyper L(f,P)} = \st {\Hyper U(f,P)} = Z\int_Rf(\vec x)\, \Hyper \Vert 
P \Vert.$$
\smallskip
\hrule
\vskip 18pt
{\bf Proposition 5.1.} {\sl Suppose that bounded $f\colon R \to \real.$ 
Let $Q =  Q_1\times \cdots \times Q_n$ be a hyperfinite partition  of $\hyper 
R,$ where each $Q_i$ is a hyperfinite partition of $[a_i,b_i],$ for 
$1\leq i\leq n.$ If $\Hyper {\rm mesh}(Q_i) \in \monad {0},$ then 
$$0\leq L(f) - \hyper L(f,Q) \in \monad {0},$$
$$0\leq \Hyper U(f,Q) - U(f)\in \monad {0}.$$}\par
Proof. Consider $P_i= \{a_i = x_{i0},\ldots,x_{i(m+1)}=b_i\}$ 
a partition of $[a_i,b_i].$ Then since it is a finite 
set of real numbers it follows that $\hyper P_i = P_i.$ Now let $Q=Q_1 \times 
\cdots \times Q_n$ be a hyperfinite partition of $\hyper R,$ where 
each $Q_i$ is hyperfinite and $\Hyper {\rm mesh}(Q_i) \in \monad {0}.$ 
[Note that $\Hyper {\rm mesh}(Q_i) \in \monad 0$ for $1\leq i\leq n$ if and only if 
$\Hyper {\rm mesh}(Q) \in \monad 0.$] Let $Q_i= \{a_i = y_{i0},\ldots,y_{i(\nu +1)}=b_i\}.$
For each $j = 1,\ldots,m,$ we know there are finitely many members of $Q_i$
such that $y_{ik(j)} \leq x_{ij} < y_{i(k(j) + 1)}.$ This holds for each $i = 
1,\ldots,n.$ Thus there is a \underbar{finite} set ${\cal R}$ of infinitesimal 
subrectangles $R_q$ such that $\vec x \in P$ if and only if there exists some
$R_q \in {\cal R}$ such that $\vec x \in R_q.$ Further$,$ if $R_q \in {\cal 
R},$ then $\hyper v(R_q) \in \monad 0.$  Since ${\cal R}$ is also a 
hyperfinite set then by *-transfer we can rearrange the hyperfinite sum
$\sum _{S \in {\cal S}(Q)} m_Sv(S)$ as follows: 
$$\hyper L(f,Q) = \sum_{S \in {\cal R}}m_S\hyper v(S) + 
\sum_{S \in ({\cal S}(Q)-{\cal R})}\!\!\!\!\!\!\!\!m_S\hyper v(S).$$                                                
This all implies since a finite sum of infinitesimals is infinitesimal that 
$$\hyper L(f,Q) \approx \sum_{S \in({\cal S}(Q)-{\cal R})}\!\!\!\!\!\!\!\!m_S\hyper v(S).$$                                                
Moreover$,$ with respect to the common refinement $K= P\cup Q$ the partition points in $P$ only finitely partition the infinitesimal subrectangles $S \in {\cal R}.$ 
Hence$,$ in a similar manner the hyperfinite sum $\hyper L(f,K)$ can be 
rearranged so that$$\hyper L(f,K) \approx \sum_{S \in({\cal S}(Q)-
{\cal R})}\!\!\!\!\!\!\!\!m_S\hyper v(S).$$ Thus $\hyper L(f,K) \approx \hyper L(f,Q).$\par
 Now $K$ is a hyperfinite partition of $R$ that 
*-refines $P.$ By *-transfer of the elementary properties of lower sums and 
partitions it follows that $\hyper L(f,K) \geq L(f,P).$
Let $r \in \real^+.$ Then there is a standard partition $P_1$ such that
$L(f,P_1) \geq L(f) -r.$ From this we have that
$$L(f) \geq \hyper L(f,Q) \approx \hyper L(f,K) \geq L(f,P_1) \geq L(f) -r.$$
Since $r$ is arbitrary then taking the standard part operator we have that
$\st {\hyper L(f,Q)} = L(f).$ Thus $0 \leq L(f) - \hyper L(f,Q) \in 
\monad  0.$ In similar manner the second conclusion follows and this completes 
the proof.\par
\vfil
\eject
{\bf Proposition 5.2.}\par {\sl (i) A bounded function $f\colon R \to \real$ is 
Darboux integrable
iff it is Z-integrable for some ${\cal P}_I$ if and only if it is Z-integrable for all 
${\cal P}_I.$\par
(ii) If bounded $f\colon R \to \real$ is Darboux integrable or Z-integrable for 
${\cal P}_I$$,$ 
then  $$D\int_Rf(\vec x) dX = Z\int_Rf(\vec x)\, \Hyper \Vert P\Vert.$$}\par
Proof. These results are indeed immediate from Proposition 5.1 and application 
of the standard part operator.\par
\vskip 18pt
{\bf Proposition 5.3.} {\sl Let bounded $f\colon R\to \real$ and ${\cal P}_I$ 
be any nonempty 
set of hyperfinite partitions of $\hyper R$ such that $\Hyper \Vert P \Vert
\in \monad 0$ for each $P \in {\cal P}_I.$ If $f$ is Z-integrable and
$P,Q \in {\cal P}_I,$ then $$\st {\hyper L(f,P)} = \st {\Hyper L(f,Q)}= Z\int_Rf(\vec x)\, \Hyper \Vert P\Vert =
Z\int_Rf(\vec x)\, \Hyper \Vert Q\Vert.$$}\par
Proof. Since $f$ is Darboux integrable then once again application of 
Proposition 5.1 and the standard part operator completes the proof.\par
\vskip 18pt
{\bf Theorem 5.1.1.} {\sl A bounded function $f\colon R \to \real$ is 
integrable if and only if  it is integrable in the sense of Darboux and Riemann Sums.}\par
Proof. Let $\cal P$ be the set of all simple partitions of $R.$ Then there exists a 
mapping $\cal Q$ on the set $\cal P$ such that for each $P\in {\cal P}$ the 
image ${\cal Q}(P)$ is the set of all intermediate partitions for $P.$ For 
each $P \in {\cal P}$ and each $Q \in {\cal Q}(P)$ let 
$$S(f,P,Q) = \sum_{S \in {\cal S}(P), \vec x \in S \cap 
Q}\!\!\!\!\!\!\!\!\!\!f(\vec x)v(S).$$
First$,$ assume that $f$ is integrable as defined in by Definition 5.1.1. Then 
we know that there is some $r \in \real,$ some $P^\prime \in \hyper {\cal P}$ 
such that for each $Q^\prime \in \hyper {\cal Q}(P^\prime),\ \hyper 
S(f,p^\prime,Q^\prime) \in \monad r.$ Hence for arbitrary $A \in \real^+$
$$r -A < \hyper S(f,P^\prime,Q^\prime) < r + A.$$
By *-transfer$,$ there exists a standard simple
partition $P\in {\cal P}$ such that for each $Q\in {\cal Q}(P)$ it follows 
that 
$$r -A < S(f,P,Q) < r + A.$$
From this it follows that $r -A \leq L(f,P) \leq U(f,P) \leq r +A.$ Thus
$r -A \leq L(f) \leq U(f) \leq r + A,\ \Rightarrow \vert U(f) -r \vert \leq A$
and
$\vert L(f) - r \vert \leq A.$ Since $A$ is arbitrary then this implies that
$U(f) = L(f) = r.$ Consequently$,$ $f$ is Darboux integrable and as is 
well-known 
$f$ is Riemann integrable in the sense of Riemann Sums. \par
 Conversely$,$ assuming that $f$ is Darboux (or Riemann) integrable then $f$ is 
Z-integrable for$,$ ${\cal P}_I,$ the set of all fine partitions of 
$\hyper R.$ Thus $P^\prime \in {\cal P}_I$ and $\Hyper \Vert P^\prime \Vert = dX,$ 
then
$$\st {\hyper L(f,P^\prime)} = \st {\Hyper U(f,P^\prime)} = Z\int_Rf(\vec 
 x)dX = D\int_Rf(\vec x)dX.$$
But from *-transfer$,$ for each $Q^\prime\in \hyper {\cal Q}(P^\prime),$
$$\hyper L(f,P^\prime) \leq \hyper S(f,P^\prime,Q^\prime)\leq \Hyper 
U(f,P^\prime).$$
 The result follows by taking the standard part operator.\par
\vskip 18pt
{\bf Theorem 5.1.2.} {\sl If bounded $f\colon R \to \real$ is integrable$,$ then
there exists a unique $r \in \real$ such that  
for every fine partition $P = \{\vec x_0,\ldots \vec x_\Omega \},\ \Omega \in 
{\nat}_\infty$ and every internal intermediate partition $Q = \{\vec v_q\},\ 
1 \leq q \leq \Gamma \in {\nat}_\infty$ it follows that
$$\sum_{k=1}^\Gamma \hyper f(\vec v_q)\hyper v(R_q) \in \monad r.$$}\par
Proof. Since $f$ is Darboux integrable then it is Z-integrable. Let
${\cal P}_I$ be the set of all fine partitions of $R.$ Then every member of 
${\cal P}_I$  is internal. Result follows by repeating the proof of the 
converse of Theorem 5.1.1 for ${\cal P}_I.$ \par
\vskip 18pt
Obviously$,$ the restriction of fine partitions to the simple fine partitions is 
not necessary since the value of the integral defined in Chapter 5$,$ as well as 
the equivalent Darboux or Riemann integral$,$ is independent of the choice of the 
fine partition. Clearly$,$ the simplicity of the simple fine partitions is often 
useful in applications and the customary integral notation of $dX$ implies 
that each infinitesimal subrectangle has *-measure = $dX.$ However$,$ for much 
that follows $\cal P$ is often assumed only to be an infinite set of 
partitions $P = P_1 \times \cdots P_n$ and in this case$,$ letting 
${\cal S}(P)$ the set of all subrectangles of $R \subset \realp n$ generated by 
$P \in {\cal P},$ we define the set ${\cal C}_{PSR}^\prime = \{S\vert 
\exists P(P\in {\cal P} \land S \in {\cal S}(P)\}.$ Extend the definition of 
simply additive to ${\cal P}.$\par
{\bf Proposition 5.4.} {\sl Let $f \colon R \to \real.$ Let $\cal P$ 
be any infinite set of partitions (including the trivial one $\{R\}$) of 
$R,\ B\colon {\cal C}_{PSR}^\prime \to R,$ and $B$ simply additive on each $P\in 
{\cal P}.$ Suppose that there exists a fine 
partition $Q \in \Hyper {\cal P}$ (i.e. $\Hyper \Vert Q \Vert \in \monad 0)$
 such that for every $S \in \hyper {\cal 
S}(Q)$ there exists 
some $\vec p \in S$ such that
$$\hyper B(S)/(\hyper v(S)) \approx \hyper f(\vec p),$$
then for any $r \in {\real}^+$
$$-rv(R) + \hyper L(f,Q) < B(R) < \Hyper U(f,Q) + rv(R).$$}\par
Proof. Let $r \in {\real}^+$ and assume that for the fine partition
$Q \in \Hyper {\cal P}$ $$\hyper (B(R)) = \hyper B(\hyper R) \geq \Hyper U(f + r,Q).$$
First$,$ we make the following standard observations. Let $P \in \cal P$ 
be any standard partition of $R.$ Assume that $B(R) \geq U(f+r,P).$
The simple additivity of $B$ yields that there must exist some $S \in {\cal S}(P)$ 
such that 
$$B(S) \geq M_ Sv(S),\ M_S = \sup \{f(x) + r\vert x \in S\}.$$
By *-transfer$,$ there exists some $T \in \hyper {\cal S}(Q)$ such that
$$\hyper B(T) \geq M_T\hyper v(T),\ M_T =
 \hyper \sup \{\hyper f(\vec x) + r\vert x \in T\}.$$ 
Hence for each $\vec p \in T$ $$\hyper B(T) \geq (\hyper f(\vec p) +r)\hyper v(T).$$
Therefore since $\hyper v(T) > 0$ this implies that for each 
$\vec p \in T,$ $$ \hyper B(T)/(\hyper v(T)) -\hyper f(p) \geq r.$$
This$,$ however$,$ contradicts the hypotheses. Thus $\hyper (B(R)) = \hyper 
B(\hyper R) = B(R)< \Hyper U(f + r,Q) =
\Hyper U(f,Q) + rv(R).$\par
In like manner$,$ it follows that $-rv(R)  + \hyper L(f,Q) < B(R).$
\vskip 18pt
{\bf Proposition 5.5.} {\sl Let $f \colon R \to \real.$ Let $\cal P$ 
be any infinite set of partitions (including the trivial one $\{R\}$) of 
$R,\ B\colon {\cal C}_{PSR}^\prime\to R,$ and $B$ simply additive on each $P\in 
{\cal P}.$
Suppose that there exists a fine 
partition $Q \in \hyper {\cal P}$ such that for every $S \in \hyper {\cal 
S}(Q)$ there exists 
some $\vec p \in S$ such that
$$\hyper B(S)/(\hyper v(S)) \approx \hyper f(\vec p),$$
then $f$ is Z-integrable for $\{Q\}$ and $$B(R) = Z\int_Rf(\vec x) \Hyper \Vert Q 
\Vert.$$}\par
Proof. From Proposition 5.4 we have that for any $r \in {\real}^+,\ 
-rv(R) + \hyper L(f,Q) < B(R) < \Hyper U(f,Q) + rv(R).$ Taking the standard 
part operator and using Proposition 5.1 it follows that 
$$-rv(R) + L(f) \leq B(R) \leq U(f) + rv(R).$$
But $r$ an arbitrary positive real number  $\Rightarrow$ 
$$L(f) = U(f) = Z\int_R f(\vec x) \Hyper \Vert Q \Vert.$$ \par
\vskip 18pt
{\bf Theorem 5.2.1. (An Infinite Sum Theorem.)} {\sl  Let bounded 
$f\colon R \to \real$ 
and simply additive $B\colon {\cal C}_{PSR} \to \real.$ If there exists a 
simple fine partition $\{R_q \vert 1 \leq q \leq \Gamma \}$ and for each $R_q$ 
there exists some $\vec p \in R_q$ such that 
$$ \hyper B(R_q)/dX \approx \hyper f(\vec p),\eqno(*)$$
then $f$ is integrable and 
$$B(R) = \int_R f(\vec x)\, dX.$$}\par
Proof. Simply apply Proposition 5.5.\par
\vskip 18pt
\hrule
\smallskip
{\bf Definition A5.2. ({Jordan-Supernearness}).} Let $J$ be a nonempty  
Jordan-measurable subset of $R \subset \realp n.$ A map 
$B\colon {\cal C}_{J} \to \real,$ where ${\cal C}_{J}$ is a nonempty set of  
Jordan-measurable subsets of $J,$ is 
{\bf {JORDAN-SUPERNEAR}} to bounded $f\colon J \to \real$ if for every
$K \in \Hyper {\cal C}_{J}$ such that $\hyper v(K\cap R_S) \in {\monad 
0}^+,$ where $R_S$ is an infinitesimal subrectangle of $\hyper R,$ 
and every $\vec p \in K\cap R_S$ it follows that
$$ \hyper B(K \cap R_S)/\hyper v(K\cap R_S) \approx \hyper f(\vec p),\eqno (**)$$
where $v(A)$ is the Jordan content for any Jordan-measurable $A \subset J$ 
and $\hyper B$ is defined on $K \cap R_S.$\par
\smallskip
\hrule
\vskip 18pt
{\bf Proposition 5.6.} {\sl Let $J$ be a Jordan-measurable 
subset of $\realp n.$ If bounded $f\colon J 
\to \real$ is uniformly 
continuous on $J$ and ${\cal C}_J$ is any nonempty set of connected 
Jordan-measurable 
subsets of $J,$ then there exists a map $ B\colon {\cal C}_J \to \real $ that is 
Jordan-supernear to $f.$}\par
Proof. Consider any rectangle $R$ such that $J \subset R$ and as usual define for each $A \in {\cal C}_J,$ 
$$B(A) = \int_A f(\vec x)\, dX = \int_{R}\hat f(\vec x)\, dX.$$ 
Let $A\in {\cal C}_J$ and $R_0\subset  R$ be any subrectangle such that $A \cap R_0 
\not= \emptyset.$ Then we know that there exists some 
$\vec s \in A \cap R_0$ such that $\int_{A\cap R(0)} f(\vec x)\, dX = f(\vec s) v(A\cap 
R_0).$ Suppose that $Q\in \Hyper {\cal C}_J,\ \hyper v(Q \cap R_S) \in 
{\monad 0}^+$ and $B$ is defined on $ Q \cap R_S,$ where $R_S$ is an 
infinitesimal subrectangle. Then $Q \cap R_S \not= \emptyset$ and *-transfer 
implies that there exists some  $\vec q \in Q\cap R_S$ such that 
$\hyper B(Q\cap R_S) = \hyper f(\vec q)\hyper v(Q\cap R_S).$ Consequently,
$$(\hyper B(Q\cap R_S)/\hyper v(Q\cap R_S)\approx \hyper f(\vec p).$$
Since $Q\cap R_S \subset R_S,$ then $\vec x,\vec y \in Q \cap R_S \Rightarrow
\vec x\approx \vec y,$ and  uniform continuity of 
$f\ \Rightarrow \hyper f(\vec x) \approx \hyper f(\vec y).$ Thus for 
any $\vec p \in Q\cap R_S,$
$$\hyper B(Q\cap R_0)/\hyper v(Q \cap R_S) \approx \hyper 
f(\vec q) \approx \hyper f(\vec p).$$
This completes the proof.\par
In order for Jordan-supernearness to yield continuity on $J,$ it appears 
necessary to select both of the sets 
$J$ and ${\cal C}_J$ more carefully.\par
{\bf Proposition 5.7.} {\sl Let the rectangle $R \subset \realp n,\ f\colon R \to 
\real$ and ${\cal C}_{SR} =\{R_q\vert R_q\ {\rm is\ a\ subrectangle\ of}\ R\}.$ 
Suppose that $B\colon {\cal C}_{SR} \to \real$ is Jordan-supernear to $f,$ then 
$f$ is uniformly continuous on $R.$}\par
Proof. Let $\vec p \not= \vec q,\ \vec p,\vec q \in \hyper R,\ \vec p \approx 
\vec q.$ Assuming that
$\vec p= (x_1,\ldots,x_n),\ \vec q = (y_1,\ldots, y_n)$ then $dx = \max \{
\vert x_i - y_i \vert \vert 1 \leq i\leq n\} >0.$ Now letting $R_S =
\{(z_1,\ldots,z_n)
\vert \forall i(1\leq i \leq n \to (x_1 -dx \leq z_i\leq x_n +dx)\land (z_i 
\in \hyperreal) \}$ then for this infinitesimal subrectangle we have that
$R_S \in \Hyper {\cal C}_{SR},\ \vec p,\vec q \in R_S$ and $\hyper v(R_S) \in 
{\monad 0}^+.$ Jordan-supernearness implies that for each $\vec x \in R_S$
$$\hyper B(R_S)/\hyper v(R_S) \approx \hyper f(\vec x).$$   
Since $f$ is bounded on $R$ and $\vec p,\vec q \in R_S$ then 
$$\hyper B(R_S)/\hyper v(R_S) \approx \hyper f(\vec p)\approx \hyper f(\vec 
q).$$ 
Consequently $f$ is uniformly continuous on $R.$\par
{\bf Theorem 5.2.2.} {\sl A bounded function $f\colon R \to \real$ is 
continuous if and only if there exists a map $B\colon {\cal C}_{SR} \to \real$ 
that is supernear to $f.$} \par
Proof. Obvious from above propositions. \par
{\bf Corollary 5.2.2.} {\sl Suppose that $f\colon R \to \real$ is continuous.
For each $R_S \in {\cal C}_{SR}$ define $B(R_S) = \int_{R_S} f(\vec x)\, dX.$
Then $B$ is supernear to $f.$} \par
Proof. This is established in the proof of Proposition 5.6\par
\vskip 18pt

{\bf Theorem 5.2.3.} {\sl Let bounded $f\colon R \to \real.$ If $B\colon {\cal 
C}_{SR} \to \real$ is supernear to $f$ and simply additive on each simple 
partition of each $R_S,$ then $f$ is continuous on $R$ and 
$$B(R_S) = \int_{R_S} f(\vec x)\, dX$$ for each $R_S \in {\cal C}_{SR}.$}\par
Proof. Proposition 5.7 implies continuity of $f$ implies that the mapping $B^\prime \colon {\cal 
C}_{SR} \to \real$ defined by $B^\prime (R_S) = \int_{R_S}f(\vec x)\, dX$ is 
supernear to $f.$ Now $f$ is bounded on $R_S.$ Let ${\cal P}_I$ be the set 
of all simple partitions of $R_S.$ Then $B$ is defined on all subrectangles 
generated by members of ${\cal P}_I$ and simply additive on ${\cal S}(P)$ for 
each $P \in {\cal P}_I.$ Hence there exists a fine partition $Q \in \Hyper
{\cal P}_I$ and $B$ restricted to ${\cal C}_{P(I)SR(S)}$ is also (Jordan) supernear 
to $f$. Of course$,$ $f$ is bounded on each $R_S.$ 
Thus for each $R_S,\ B$ satisfies the hypotheses of Proposition 5.5. 
Therefore,
$$B(R_S) = Z\int_{R_S} f(\vec x) dX = \int_{R_S} f(\vec x) dX = 
B^\prime(R_S).$$\par
{\bf Corollary 5.2.3.1} {\sl Let bounded $f\colon R \to \real.$ If $B\colon 
{\cal C}_{SR} \to \real$ is supernear to $f$ and additive on ${\cal C}_{SR},$ 
then $f$ is continuous on $R$ and 
$$B(R_S) = \int_{R_S} f(\vec x)\, dX$$ for each $R_S \in {\cal C}_{SR}.$}\par
{\bf Corollary 5.2.3.2} {\sl Let bounded $f\colon R \to \real.$ There exists one 
and only one map $B\colon {\cal C}_{SR} \to \real$ that is supernear to $f$ 
and either simply additive on each simple partition of each $R_S$ 
or additive on ${\cal C}_{SR}$}.\par
\vfil\eject\centerline{Appendix For Chapter 6.}
\vskip 18pt
Recall that a curve is a continuous map 
$c\colon [0,1] \to \realp n.$ This is equivalent to considering $c$ as 
determined by n continuous coordinate
functions 
$x_i = f_i(t),\ 1 \leq i \leq n$ each defined on $[0,1] \subset \real.$ 
Of course$,$ the geometric curve $C$ 
determined by these functions is usually considered as the set 
$\{(x_1,\ldots,x_n)\vert t \in [0,1]\}.$ The *-transform of these 
defining functions leads to the functions $x_i = \hyper f_i(t),\ 1 \leq i 
\leq n,$ each defined on $\Hyper [0,1] \subset \hyperreal$ and they generate the 
``hypercurve" $\Hyper C \subset \hyperrealp n.$\par
\vskip 18pt                              
{\bf {Example 4.4.1.A.}}  Let $\Gamma \in {\nat}_\infty.$ Then $F =\{t_i\vert t_i 
= i/\Gamma \land 0 \leq i \leq \Gamma \}$ is an internal and hyperfinite 
subset of $\Hyper [0,1].$ By *-transfer $F$ behaves like an ordered partition 
of the interval $[0,1]$ as defined in the standard sense. Such a set is 
termed a {\bf {fine partition}} (i.e. hyperfinitely many members of $[0,1]$ 
generating subintervals that are infinitesimal in length). The internal set 
$F$ generates the 
internal set of ``points'' $P = \{(\hyper f_1(t_i),\ldots,\hyper 
f_n(t_i))\vert t_i \in F\}$ that are members of the hypercurve $\Hyper C.$
Now for each $i = 0,\ldots,\Gamma - 1,$ and each $j,\ 0\leq j \leq n$ 
let $\hyper f_j(t_{i+1}) - \hyper 
f_j(t_i) = d(j,i).$ (If $c$ is continuous$,$ then each $d(j,i) \in \monad 0.$) 
 For each $i \in \hypernat$ such that $0 \leq i \leq 
\Gamma -1,$ the internal set $\ell_i = \{(x_1,\ldots,x_n)\vert
\forall j \in \hypernat,\ 0\leq j \leq n,\  x_j =
\hyper f_j(t_i) + t(d(j,i)) \land t \in 
\Hyper [0,1] \}$ is a hyperline segment connecting the two points 
$(\hyper f_1(t_i),\ldots,\hyper f_n(t_i)),\ (\hyper f_1(t_{i+1}),
\ldots,\hyper f_n(t_{i+1}))$ on the curve $\Hyper C.$ From this one obtains 
the internal hyperpolygonal curve ${\cal P}_\Gamma = \bigcup \{\ell_i \vert 
0\leq i \leq \Gamma - 1\}.$ As to the length of ${\cal P}_\Gamma $ simply 
extend the concept of length in the classical sense by defining for each 
$i = 0,\ldots,\Gamma -1$ the vector $\vec v_i = (d(1,i),\ldots,d(n,i)) \in 
\hyperrealp n.$ Then let the hyperfinite sum $\sum_{i=0}^{\Gamma - 
1} \Vert \vec v_i \Vert = \vert {\cal P}_\Gamma \vert \in \hyperreal.$ Even though$,$ 
in general$,$ you would have a different hyperpolygon with a different hyperreal 
length for $\forall \Gamma \in {\nat}_\infty$ we show$,$ using the above 
terminology and notation that:\par

{\sl If 
$c$ is continuously differentiable$,$ then for all $\Gamma \in {\nat}_\infty, \
\vert {\cal P}_\Gamma \vert \in \monad r$ and the real number 
$$r= \int_a^b \sqrt{\sum_{j=1}^n  f_j^\prime(t)^2}\, dt.$$}\par
Proof. Consider an arbitrary $\Gamma \in {\nat}_\infty$ and the internal 
partition $F$ of $[0,1].$ This is a simple fine partition with $t_{i+1} - t_i = dt 
= 1/\Gamma.$ ({\it Actually the following proof holds for any fine 
partition.})  
By *-transfer of the standard mean value theorem for the 
derivative$,$ it follows that for each $i = 0,\ldots, \Gamma -1,$ and each $j 
= 1,\ldots, n$ there exists some $t_{ji}^\prime \in (t_i, t_{i+1})$ such that 
$$d(j,i) = \hyper f_j^\prime(t_{ji}^\prime)\, dt.$$
Thus $\Vert \vec v_i \Vert = \sqrt{\sum_{j=1}^n  d(j,i)^2} = 
\sqrt{\sum_{j=1}^n (\hyper f_j^\prime(t_{ji}^\prime)dt)^2} = 
\sqrt{\sum_{j=1}^n (\hyper 
f_j^\prime(t_{ji}^\prime ))^2}\, dt.$ Since each $f_j^\prime$ is uniformly continuous then 
$\hyper f_j^\prime (t_i) = \hyper f_j^\prime(t_{ji}^\prime) + \delta_{ji},\ 
 \delta_{ji} \in 
\monad 0.$ Thus $\Vert \vec v_i \Vert = \sqrt{\sum_{j=1}^n (\hyper 
f_j^\prime (t_i)+ \delta _{ji})^2}\, dt = \sqrt{\sum_{j=1}^n (\hyper 
f_j^\prime (t_i))^2+ \eps_i}\, dt,\ \eps_i \in \monad 0.$ However,
$\sum_{j=1}^n (\hyper f_j^\prime (t_i))^2+ \eps_i \approx
\sum_{j=1}^n (\hyper f_j^\prime (t_i))^2 \Rightarrow 
\sqrt {\sum_{j=1}^n (\hyper f_j^\prime (t_i))^2+ \eps_i } \approx
\sqrt {\sum_{j=1}^n (\hyper f_j^\prime (t_i))^2 }.$ Hence$,$
$\Vert\vec v_i \Vert =\sqrt {\sum_{j=1}^n (\hyper f_j^\prime 
(t_i))^2}\, dt
+ \delta_i\,dt,\ \delta_i \in \monad 0.$ By *-transfer of the finite case$,$ 
$\max \{\vert \delta_i\vert \vert 0\leq i \leq \Gamma -1 \} = \delta \in \monad 0.$ Also 
note that the function  $\sqrt{\sum_{j=1}^n f_j^\prime(t)^2}$ is 
continuous on $[a,b].$  Putting the above together we have that
$$\sum_{i=0}^{\Gamma - 1} \Vert \vec v_i \Vert = \vert 
{\cal P}_\Gamma \vert =\sum_{i=0}^{\Gamma -1} {\sqrt {\sum_{j=1}^n 
(\hyper f_j^\prime (t_i))^2 }\, dt} +\sum_{i=0}^{\Gamma -1} \delta 
_i\, dt.$$
But,
$$ \vert \sum_{i=0}^{\Gamma -1} \delta _i\, dt \vert \leq \sum_{i=0}^{\Gamma -1}
\vert \delta \vert\, dt =  \delta (b-a) \in \monad 0.$$
Therefore$,$ 
$$\sum_{i=0}^{\Gamma - 1} \Vert \vec v_i \Vert = \sum_{i=0}^{\Gamma -1} {\sqrt {\sum_{j=1}^n 
(\hyper f_j^\prime (t_i))^2 }\, dt} + \lambda,\ \lambda \in \monad 
0.$$ Since $\sqrt{\sum_{j=1}^n f_j^\prime(t)^2}$ is bounded and 
integrable then theorem 5.1.2 yields that 
$${\rm st}\Bigl(\sum_{i=0}^{\Gamma -1} {\sqrt {\sum_{j=1}^n 
(\hyper f_j^\prime (t_i))^2 }\, dt}\Bigr) = \int_a^b 
\sqrt{\sum_{j=1}^n f_j^\prime(t)^2}\, dt,$$
and the proof is complete.\par
\vskip 18pt
What the above definition shows is that in this case 
our definition $\vert {\cal 
P}_\Gamma \vert$ is independent of the particular $\Gamma \in {\nat}_\infty$ 
chosen  (as mentioned within  the proof it is actually independent 
of any fine partition chosen) and also coincides with the classical one. 
Indeed$,$ we have the following added proposition that shows that the classical 
concept of the 
rectifiable curve and the existence of $r\in \real$ such that 
$\vert {\cal P}_P \vert\in \monad r$ for every fine partition $Q$ of $\Hyper [0,1]$
are equivalent concepts. \par
\vskip 18pt
{\bf Proposition 6.1.} {\sl Consider continuous $c\colon [0,1] \to \realp n.$ Then 
$c$  is rectifiable if and only if  for every fine partition $Q$ of 
$\Hyper [0,1]$
$$\st {\Hyper \vert {\cal P}_Q \vert} = L.$$}\par
Proof. (Sketch) Let $P =\{0 = u_0 < u_1 < \cdots u_q = 1\}$ be a partition  
and let ${\cal P}_P$ denote the polygonal curve generated by the functions
$f_1,\ldots,f_n$ and $\vert {\cal P}_P \vert$ the standard length for the 
polygonal curve. Consider any fine partition $Q$ of $\Hyper [0,1].$ Let 
$I =[t_i,t_{i+1}]$ be a subinterval generated by $Q.$ Then for any $t \in I,\ 
 \vert \hyper f_j(t_i) -\hyper f_j(t) \vert,\ \vert \hyper f_j(t) - \hyper 
f_j(t_{i+1})\vert \in \monad 0$ for each $j$ since each $f_j$ is uniformly 
continuous on $[0,1].$ Note that by *-transfer $\Hyper \vert {\cal P}_Q \vert$ is 
but the length of the interval hyperpolygonal curve ${\cal P}_Q$ as defined in 
Example 4.4.1.A. Now there are but finitely many infinitesimal subintervals 
generated by $Q$ that contain the partition points from $P.$ Considering the 
standard properties of the Euclidean norm$,$ the internal common refinement 
$P\cup Q$ and the fact that the finite sum of infinitesimals is infinitesimal 
this leads to 
$$\vert {\cal P}_P \vert \leq \Hyper \vert {\cal P}_{P \cup Q} \vert \leq 
\Hyper \vert {\cal P}_Q \vert + 
\eps,\ \eps \in \monad 0.$$
Since $c$ is rectifiable then given any $r \in {\real}^+$ there exists a 
partition $P(r)$ of $[0,1]$ such that $0\leq L - \vert {\cal P}_{P(r)} \vert< r.$
 By 
*-transfer we also have that $\Hyper \vert {\cal P}_Q \vert \leq L.$ Consequently$,$ 
$$0\leq L - \Hyper \vert {\cal P}_Q \vert < r+\eps.$$
Thus $\st {\Hyper \vert {\cal P}_Q \vert } = L.$ \par
Conversely$,$ let $P$ be any standard and $Q$ any fine partition of $\Hyper [0,1]$ 
Then internal partition $P \cup Q$ is a fine 
partition. Thus $$\vert {\cal P}_P \vert \leq \Hyper \vert {\cal P}_{P\cup Q}
\vert \in \monad L.$$ Hence$,$
$$\st {\vert {\cal P}_P \vert} = \vert {\cal P}_P \vert \leq \st {\Hyper \vert 
{\cal P}_{P\cup Q}}\vert = L$$ and $L$ is an upper bound for the standard 
$\vert {\cal P}_P \vert$ and$,$ thus$,$ an upper bound for all internal
partitions of 
$\Hyper [0,1].$ Moreover$,$ if $s \in 
{\real}^+,$ then $ 0\leq L- \Hyper \vert {\cal P}_{P\cup Q} \vert  < s.$   Thus by *-
transfer$,$ for every $s \in {\real}^+$ there exists some partition $K$ of 
$[0,1]$ such that $ 0\leq L- \vert {\cal P}_K \vert  < s.$ 
Thus $c$ is rectifiable. \par
\vskip 18pt
In Proposition 6.1$,$ the number $L$ is the length of the curve. We now proceed 
to establish the other Chapter 6 theorems.\par
\vskip 18pt
{\bf Theorem 6.2.1.} {\sl Let $A$ be a compact subset of $\realp n$ and the 
infinitesimal subrectangle $R \subset \hyper A.$ Then there exists some $\vec 
p \in A$ such that $R \subset \monad {\vec p}.$}\par
\vskip 18pt
Proof. Since $R \not= \emptyset$ let $\vec x,\vec y \in R.$
 Then $\vec x \approx \vec y.$ Since $A$ is compact then  
there exists some $\vec p \in \realp n$ such that $y \in \monad {\vec p}.$
Hence for each $\vec x \in R,\ \vec x \approx \vec p \Rightarrow R \subset 
\monad {\vec p}.$ \par
\vskip 18pt
\centerline{IR3}
(1) We wish to measure a quantity $M$ for a compact Jordan-measurable set $J 
\subset R \subset \realp n,$ where $M$ is defined on and$,$ at least$,$ additive 
over members of the set ${\{\cal C}_{SR},R-J,J\}.$ Further$,$ if subrectangle 
$S \subset R - J,\ M(S) = 0,$ and $M(R-J) = 0.$ 
Let $v(J)$ denote the Jordan content.\par
(2) There is a generating function $f(\vec x)$ that is related to the 
functional 
$M$ in the following manner:\par
\indent\indent (i) The functions $f$ is continuous on $J.$\par
\indent\indent (ii) Let $P$ be some simple fine partition$,$ 
$S \in \hyper {\cal S}(P)$ and $K = \hyper J \cap S \not= \emptyset.$ 
Then there exist $\vec x_m \in K$ and $\vec x_M 
\in K$ such that $\hyper {f_m}=\hyper {f(\vec x_m)} = 
\Hyper {\inf} \{\hyper f(\vec 
x)\vert \vec x \in K\} = \inf \{\hyper f(\vec x)\vert \vec x \in K\}$ and
$\hyper {f_M} =\hyper {f(\vec x_M)} = \Hyper {\sup} \{\hyper f(\vec 
x)\vert \vec x \in K\} = \sup \{\hyper f(\vec x)\vert \vec x \in K\}.$ \par 
\indent\indent (iii) $(\hyper {f_m)}\, \hyper v(S) \leq \hyper M(S) 
\leq (\hyper f_M)\, \hyper v(S).$ [Note: this is the case where $\hyper L$ 
is the identity map.]\par
\vskip 18pt
{\bf Theorem 6.2.2.} {\sl If IR3  holds$,$ then $$M(J) = \int_Jf(\vec x)\, 
dX.$$}\par    
Proof. Continuity of $f$ on $J$ yields$,$ since the  set of discontinuities of 
$\hat f$ has  
Lebesgue measure zero$,$ that $\hat f$ is integrable on $R.$ 
Thus $\int_J f(\vec x)\, dX = \int_R \hat f(\vec x)\, dX.$  Let $P$ be a 
a simple fine partition and any $S \in {\cal S}(P)$ such that $K = \hyper J \cap S 
\not= \emptyset.$ Then we have from (iii) that 
$$\hyper f_m \hyper v(S) = \hyper f_m\ dX \leq \hyper M(S)\leq \hyper f_M 
\hyper v(S) = \hyper f_M\ dX.$$\par
Assume that $K \neq \emptyset.$
Since $f =  \hat f$ is uniformly continuous on $J$ and $\vec x_m, \vec x_M \in 
K \Rightarrow \vec x_m \approx \vec x_M$ then $\hyper {\hat 
f(\vec x_m)} = \hyper f(\vec x_m) \approx \hyper f(\vec x_M) = \hyper {\hat 
f(\vec x_M)}\Rightarrow$ 
$$\hyper M(S)/dX \approx \hyper{\hat f(\vec x_m)}.\eqno(*)$$ \par

For the case that $\hyper J \cap S = \emptyset$$,$ it follows that $S \subset R-J;$
 which implies that $\hyper M(S) = 0.$ Hence in 
this case $\hyper M(S) = \hyper {\hat f(\vec x)},$ for any $\vec x \in 
S\Rightarrow$ expression (*). By application of the Infinite Sum Theorem
$$M(R) = \int_R \hat f(\vec x)dX = \int_J f(\vec x)\,dX=M(R-J) + M(J) 
= M(J).$$ \par 
\vskip 18pt
{\bf Theorem 6.2.3.} {\sl Let compact Jordan-measurable $J \subset R \subset 
\realp n.$ If continuous $f\colon J \to \real,$ then for any partition $P$ of 
$R$ and any $S \in {\cal S}(P),$ where $K = J \cap S \not= \emptyset$ 
there exist $\vec x_m \in K$ and $\vec x_M 
\in K$ such that $f_m=f(\vec x_m) = \inf \{f(\vec 
x)\vert \vec x \in K\}$ and
$f_M = f(\vec x_M) = \sup \{f(\vec 
x)\vert \vec x \in K\}.$} \par
Proof. This follows immediately since each such $K$ is compact and $f$ is 
continuous on it.\par         
\vskip 18pt
{\bf Theorem 6.2.4. ({Self-evident Max. and Min.})} {\sl Let the rectangle $R\subset \realp 
n$ and suppose that compact Jordan-measurable $J \subset R.$ Let $M$ be 
defined as in (1) of IR3$,$ continuous 
$f\colon J\to \real,\ {\cal P}$ an acceptable set of partitions and any 
$P \in {\cal P}.$ If for any $S \in {\cal S}(P)$ such 
that $J \cap S \not= \emptyset$ it follows that 
$(f_m)v(S) \leq M(S) \leq (f_M)v(S),$ 
then the infinitesimalizing process IR3 holds.}\par
\vskip 18pt
Proof. Everything stated in the hypothesis can be written in our first-order 
set-theoretic language. Noting such things as the *-minimum [resp. *-max.]
value of $\hyper f$ on a *-compact $\hyper J \cap S,$ where $S$ is an 
infinitesimal subrectangle$,$ is the same as the minimum [resp. max.]
value of $\hyper f$ on $\hyper J \cap S$ by *-transfer and the fact that 
some simple fine partition exists the result follows from Theorem 6.2.2.\par
\vskip 18pt
\centerline{IR4--{Method of Constants}}
(1) In what follows$,$ let for any $A \subset \realp n\ $ ``int'' denote the 
interior of $A.$ Let ${\cal A}=\{{\cal C}_{SR},\{{\rm int}(J\cap S)
\not= \emptyset\bigm\vert S \in {\cal C}_{SR}\}\}.$
We wish to measure a quantity $M$ for a Jordan-measurable set $J 
\subset R \subset \realp n,$ where $M$ is$,$ at least$,$ defined on and   
additive over the members of the set $\{{\cal A}, R-J,J\}$ and
for a subrectangle $S,\ S \subset R -{\rm int}(J)$ it follows that
$\ M(S) = 0$ and $M(R-J) = 0.$ Let $v(J)$ denote the Jordan content.\par
(2) There is a generating function $f(\vec x)$ that is related to the 
functional 
$M$ in the following manner:\par
\indent\indent (i) The function $f$ is bounded on $J.$\par
\indent\indent (ii) Let $P$ be any arbitrary simple fine partition$,$ 
$S \in \hyper {\cal S}(P)$ an arbitrary infinitesimal subrectangle and 
$\emptyset \not=  K = \Hyper {\rm int}(\hyper J \cap S).$ \par
\indent\indent (iii) There exists some $\vec x \in K$ such that 
$\hyper M(S)= \hyper f(\vec x)\hyper v(S)$ or  $\hyper M(S)/\hyper v(S) 
\approx \hyper f(\vec x).$  
\par
\vskip 18pt
{\bf Theorem 6.4.1.} {\sl If IR4 holds$,$ then $$M(J) = \int_Jf(\vec x)\, dX.$$}
Proof. First$,$ $f$ is bounded if and only if $\hat f,$ as it is defined on 
$R,$ is bounded. We confine our attention to the function $\hat f\colon 
R\to \realp n.$ From (iii) of IR4$,$ letting $\emptyset \not= 
 K = \Hyper {\rm int}(\hyper J \cap S)$ for $S \in \hyper {\cal S}(P)$ there is 
some $\vec x \in K,$ hence in $S,$  such that 
$$\hyper M(S) = \hyper f(\vec x)\, v(S) = \hyper {\hat f}(\vec x)\, dX
\ {\rm or}$$
$$\hyper M(S)/\hyper v(S) \approx \hyper f(\vec x) = \hyper {\hat f}(\vec 
x).$$
For the case that $K =\emptyset$ then $S \subset \hyper R - \Hyper {\rm 
int}(\hyper J) \Rightarrow \hyper M(S) = 0 =\hyper {\hat f}(\vec x)\, dX$ 
for any $\vec x \in S.$ 
Thus in all cases,
$\hyper M(S)/dX \approx \hyper {\hat f}(\vec x),\ \vec x \in S.$ Thus
from the Infinite Sum Theorem
$$M(R) = \int_R \hat f(\vec x)\, dX = \int_J f(\vec x)\, dX=M(R-J) + M(J) 
= M(J)$$ 
and the proof is complete.\par
\vskip 18pt
Prior to the next proof I point out one aspect of the Jordan-measurable 
subsets that will tacitly appear throughout many of these proofs. 
Let $J$ be a Jordan-measurable. Then it is part of the definition that $J$ is 
a bounded subset of $\realp n.$ As is well-known this implies that the 
boundary of $J,\ \partial J$ is Jordan-measurable and that $v(\partial J) = 
0.$ Further$,$ $\partial ({\rm int}(J)) \subset \partial J\Rightarrow$ that 
${\rm int}(J)$ is Jordan-measurable. Since ${\rm int}(\partial J) \cap {\rm 
int}({\rm int}(J)) = \emptyset$ and $J \subset \partial J \cup {\rm int}(J) 
\Rightarrow v(J) \leq v({\rm int}(J)) \Rightarrow v({\rm int}(J)) = v(J).$                      
\vskip 18pt
{\bf Theorem 6.4.2. ({Self-evident Method of Constants })} 
{\sl Let the rectangle $R\subset \realp 
n$ and suppose that Jordan-measurable $J \subset R.$ Let $M$ be defined as in 
(1) of IR4$,$ continuous 
 $f\colon R\to \real,\ {\cal P}$ an acceptable set of partitions of $R$ and any
$P \in {\cal P}.$ If for any $S \in {\cal S}(P)$ such that $\emptyset \not= K= 
{\rm int}(J \cap S)$ there exists some $\vec x \in K$ and some $\vec y \in S$ 
such that (i) $M(K)= f(\vec x)\, v(K)$ and (ii) $M(S)= f(\vec y)\, v(S),$  
then the infinitesimalizing process IR4 holds for $f$ restricted to $J.$}\par
Proof. Let $k = f\vert J$ and $\ Q \in \hyper {\cal P}$ be a simple fine 
partition. Assume that $\emptyset \not= K =\Hyper {\rm int}(\hyper J \cap S)$ 
for $S \in Q.$ Note that $k$ is bounded. By *-transfer of the hypotheses there 
exists some $\vec x \in K$ (thus $\vec x \in S$ also) and some $\vec y \in S$ 
such that 
(i) $\hyper M(K) = \hyper k(\vec x)\hyper v(K) = \hyper f(\vec x)\hyper v(K)$ 
and (ii) $\hyper M(S) =  \hyper f(\vec y)\hyper v(S).$ Thus since $\hyper v(K) 
\not= 0$ then
$$\hyper M(K)/ \hyper v(K) = \hyper f(\vec x)$$
$$\hyper M(S)/\hyper v(S) = \hyper f(\vec y).$$
But$,$ both $\vec x,\ \vec y \in S\ \Rightarrow \vec x \approx \vec y.$ Uniform 
continuity of $f$ yields that $\hyper f(\vec x) \approx \hyper f(\vec y)
\Rightarrow \hyper M(S)/\hyper v(S) \approx \hyper f(\vec x) =\hyper k(\vec x)
= \hyper {\hat k(\vec x)},$ where $\vec x \in S$ and thus 
$$M(R) = \int_R \hat k(x)\, dX = \int_J f(\vec x)\, dX.$$\par
\vskip 18pt
{\bf Theorem 6.4.3. (Extended Self-evident Method of Constants)}
{\sl Let the rectangle $R\subset \realp n$ and suppose that 
Jordan-measurable $J \subset R.$ Let $M$ be defined as in (1) of IR4$,$ continuous 
 $f\colon R\to \real,$ continuous $g\colon R\to \real,\ {\cal P}$ an acceptable set of 
partitions of $R$ and any
$P \in {\cal P}.$ If for any $S \in {\cal S}(P)$ and $\emptyset \not= 
K= {\rm int}(J \cap S)$ there exists some 
$\vec x_1,\ \vec x_2 \in K$ and some $\vec y_1,\ \vec y_2 \in S$ such that 
(i) $M(K)= f(\vec x_1)\, g(\vec x_2)\, v(K)$ and (ii) $M(S)= 
f(\vec y_1)\, g(\vec y_2)\, v(S),$ then the infinitesimalizing process IR4 
holds for $fg$ restricted to $J.$}\par
Proof. Let $k = f\vert J$ and $h = g\vert J.$ By *-transfer there exists a simple fine partition $Q \in \hyper {\cal 
P}$ with all the indicated properties. Let $S \ \hyper {\cal 
S}(Q),\ \emptyset \not= K = \Hyper {\rm int}(\hyper J \cap S).$ Then there 
exist $\vec x_1,\ \vec x_2 \in K, (\vec x_1,\vec x_2 \in S)$ and 
$\vec y_1,\ \vec y_2 \in S$ such that
$$\hyper M(K)/\hyper v(K) = \hyper k(\vec x_1)\hyper h(\vec x_2)=\hyper f(\vec 
x_1)\hyper g(\vec x_2)$$
$$\hyper M(S)/\hyper v(S) = \hyper f(\vec y_1)\hyper g(\vec y_2).$$
Uniform continuity of $f$ and $g$ imply that $\hyper f(\vec x_1) \approx
\hyper f(\vec y_1)$ and $\hyper g(\vec y_2) \approx \hyper g(\vec x_1).$
Using the fact that $f$ and $g$ are also bounded we have that
$$\hyper f(\vec y_1) \hyper g(\vec y_2) \approx \hyper f(\vec x_1)
\hyper g(\vec x_1).$$
Hence$,$ $$\hyper M(S)/\hyper v(S) \approx \hyper {\hat k(\vec x_1)} 
\hyper {\hat h(\vec x_1)},\ \vec x_1 \in S $$
and the result follows.\par 
\vskip 18pt
[{\bf Remark.} It is obvious that in the above two theorems (i) is not really 
needed if we know that the $y$'s are members of $K.$ The reason that (i) is 
included in the hypothesis is that it is easier to argue for the location of 
the $x$'s when applications are considered.] \par
For a given rectangle $R,$ recall that $F({\cal C}_{SR})$ is the set of all 
finite sets of subrectangles of $R.$  Extend the definition of the maps 
$L(f,\bullet),\ U(f,\bullet)$ for any nonempty collection ${\cal F} 
\in F({\cal C}_{SR}),$ where $\emptyset \notin {\cal F}.$ Let $P$ be any 
standard partition of $R.$ For a given 
$J \subset R\subset \realp n,$ where to avoid trivialities we always assume 
that ${\rm int}(J) \not= \emptyset$$,$  
define the following 
set theoretic operators each defined on the set of all partitions.
The boundary subrectangles $\partial (P) = \{S\bigm \vert S \in  
{\cal S}(P) \land S \cap J 
\not= \emptyset \land S \cap (R - J) \not= \emptyset\};$  
the exterior subrectangles ${\rm ext}(P) = \{S\bigm \vert S \in
{\cal S}(P) \land S \subset (R - J) \},$  and the inner 
subrectangles ${\rm inn}(P) = {\cal S}(P) - (\partial (P) 
\cup {\rm ext}(P)).$ Notice that $S \in {\rm inn}(P)$ if and only if  
$S\subset J.$ Let ${\cal P}$ be a set of partitions of 
$R$ and $I_{\cal P} = \{\cup A\bigm \vert \exists P(P \in {\cal P} \land 
A = {\rm inn}(P))\}\bigcup \{S\bigm \vert \exists P(P \in {\cal P} \land S 
\in {\rm inn}(P))\}.$
 When considering the nonstandard extension of these 
operators$,$ we use  the notation $\partial, \ 
{\rm inn}$ and ${\rm ext}$ rather than the notation $\hyper 
 \partial, \, \hyper {\rm ext}$ and $\Hyper {\rm inn}.$  
Since  any fine partition $Q$ is hyperfinite and the sets $\partial (Q),\,
{\rm inn}(Q)$ and ${\rm ext}(Q)$ are internal subsets of $Q$ then they are 
hyperfinite and by *-transfer of the standard finite case it follows that 
for bounded $f\colon J \to \real$
$$\hyper L(\hat f,Q)= \hyper L(\hat f,\partial (Q)) + \hyper L(\hat f,{\rm ext}(Q)) +
\hyper L(\hat f,{\rm inn}(Q)).$$ 
$$\Hyper U(\hat f,Q)= \Hyper U(\hat f,\partial (Q)) + \Hyper U(\hat f,{\rm ext}(Q)) +
\Hyper U(\hat f,{\rm inn}(Q)).$$\par
Applying a method similar to that used to establish Proposition 6.1,
it follows from the definition of Jordan-content that 
$J$ has Jordan-content (i.e. is Jordan-measurable) if and only if  for every fine 
partition $Q$ of $\Hyper R,$ 
$$\sum_{S \in \partial (Q)}\!\!\!\!\!\! \hyper v(S) \approx 0.$$\par
{\bf Proposition 6.2.} {\sl Let Jordan-measurable $J \subset R \subset \realp 
n,$ bounded $f\colon J \to \real,$ and $Q$ be any fine partition of $\Hyper 
R.$ Then
$$ \hyper L(\hat f,\partial (Q))\approx 0,\ \Hyper U(\hat f,\partial (Q))
\approx 0,$$
$$\hyper L(f,{\rm inn}(Q))\approx \hyper L(\hat f,Q), \Hyper U(f,{\rm 
inn}(Q))\approx \Hyper U(\hat f,Q).$$}\par
Proof. There exist some $m,\, M \in \real$ such that
for each $\vec x \in R,\  m\leq \hat f(\vec x\,) \leq M.$ Thus
$$m\!\!\!\!\!\sum_{S \in \partial (Q)}\!\!\!\!\!\! \hyper v(S) \leq 
\hyper L(\hat f,\partial (Q)) \leq \Hyper U(\hat f,\partial (Q) \leq 
M\!\!\!\!\!\sum_{S \in 
\partial (Q)}\!\!\!\!\!\! \hyper v(S).$$
From the above observation relative to the Jordan-content$,$ it follows that
$ \hyper L(\hat f,\partial (Q))\approx \Hyper U(\hat f,\partial (Q))\approx 0$ 
and the second part of the conclusion follows from the fact that 
$\Hyper {\hat f(\vec x\,)} = 0$ for each $\vec x \in S$ and for each $S \in {\rm ext}(Q).$\par
\vskip 18pt  
{\bf Proposition 6.3.} {\sl Suppose that Jordan-measurable $J \subset R \subset 
\realp n$ and bounded $f\colon J \to \real.$ 
Let $Q =  Q_1\times \cdots \times Q_n$ be an internal partition  of $\hyper 
R,$ where each $Q_i$ is an internal partition of $[a_i,b_i],$ for 
$1\leq i\leq n.$ If $\Hyper {\rm mesh}(Q_i) \in \monad {0}$ (i.e. $Q$ is a 
fine partition)$,$ then 
$$ L(\hat f) \approx \hyper L(f,{\rm inn}(Q)),\  U(\hat f)\approx \Hyper 
U(f, {\rm inn}(Q)).$$}\par
Proof. From Proposition 5.1 of Appendix 5$,$ $L(\hat f) \approx 
\hyper L(\hat f,Q)$ and $U(\hat f) \approx \Hyper U(\hat f,Q)$ since 
$\hat f\colon R \to \real$ is bounded. The result follows by application of 
Proposition 6.2.\par
\vskip 18pt
{\bf Proposition 6.4.} {\sl Suppose that $J \subset R \subset 
\realp n$ and bounded $f\colon J \to \real.$ Let $\cal P$ be any infinite set 
of partitions$,$ $B\colon I_{\cal P} \to \real,$ 
and for each $P \in {\cal P},\ B$ is simply additive on each ${\rm inn}(P)).$ 
Suppose that 
there exists a fine partition $Q \in \hyper {\cal P}$ such that for every 
$S \in {\rm inn}(Q)$ there exists some $\vec p \in S$ 
such that                                                                             
$$\hyper B(S)/(\hyper v(S)) \approx \hyper f(\vec p\,),$$
then for any $r \in {\real}^+$ and $A = {\rm inn}(Q)$
$$ -r\hyper v(\cup A) + \hyper L(f,A) <
\hyper B(\cup A) < \Hyper U(f,A) +  r\hyper v(\cup A).$$}\par
Proof. First$,$ note that by *-transfer that $\hyper v$ is defined for 
$\cup  A,\ \hyper v(\cup A) = \sum_{S \in A} \hyper v(S)$  
and that $\hyper B$ is defined on $\cup  A.$ Now repeat the proof of 
Proposition 5.4 with ${\rm inn}(P)$ substituted for $P,\ \cup\, {\rm inn}(P)$ or
$\cup\, {\rm inn}(Q)$ substituted for $R$ and $\cup\, {\rm inn}(Q)$ for $Q.$ \par
\vskip 18pt
{\bf Proposition 6.5. {(An Infinite Sum Theorem)}} {\sl Suppose that Jordan-measurable $J \subset R \subset 
\realp n$ and bounded $f\colon J \to \real.$ Let $\cal P$ be any infinite set 
of partitions$,$ $B\colon I_{\cal P} \to \real,$ 
and for each $P \in {\cal P},\ B$ is simply additive on ${\rm inn}(P).$ 
Suppose that 
there exists a fine partition $Q \in \hyper {\cal P}$ such that for every 
$S \in {\rm inn}(Q)=A$ there exists some $\vec p \in S$ 
such that                                                                             
$$\hyper B(S)/(\hyper v(S)) \approx \hyper f(\vec p\,),$$
then $ f$ is integrable and 
$$\st {\hyper B(\cup  A)} = \int_J f(\vec x\,)\, dX.$$}\par
Proof. From Proposition 6.4$,$ we have that for any $r \in {\real}^+,\ 
 -r\hyper v(\cup A) + \hyper L(f,A) <
\hyper B(\cup A) < \Hyper U(f,A) +  r\hyper v(\cup A).$  By *-transfer of 
the standard case$,$ then $\hyper v(\cup  A) \leq \hyper v(\hyper J) = v(J) \leq 
\hyper v((\cup  A)\cup (\partial (Q))) = \hyper v(\cup  A) + \hyper v(\cup 
(\partial (Q))$ implies that 
$\st {\hyper v(\cup  A)} \leq v(J) \leq \st {\hyper v(\cup  A)}$ or that 
$\st {\hyper v(\cup  A)} = v(J).$ Consequently$,$ 
$$ -r\,\st {\hyper v(\cup A)} + \st {\hyper L(f,A)} \leq
\st {\hyper B(\cup A)} \leq$$ $$\st {\Hyper U(f,A)} +  r\,\st {\hyper v(\cup A)}$$
yields by Proposition 6.3 that
$$ -r\,v(J) + L(\hat f) \leq
\st {\hyper B(\cup A)} \leq U(\hat f) +  r\,v(J).$$
Since $r$ is an arbitrary member of $\real^+$ then 
$$L(\hat f) = U(\hat f) = \int_R \hat f(\vec x\,)\, dX = \int_J f(\vec x\,)\, 
dX= \st {\hyper B(\cup A)}.$$\par                                           
\vskip 18pt
Obviously$,$ the above Infinite Sum Theorem almost yields the result being 
sought. But$,$ just as obviously the functional $B$ would need some additional 
property akin to a simple property displayed by the basic Jordan-measure $v.$ 
Assume  
that a functional such as $B$ is defined on $J$ and say that $B$ is 
{{\bf Jordan-like}} if for the fine partition 
$Q$ of Proposition 6.5 $\hyper B(\cup A) \approx B(J).$ It is not difficult to 
show if for each $r \in\real^+$ there exists a $\delta \in \real^+$ such that 
for each $P\in {\cal P}$ with $\Vert P\Vert < \delta$ it follows that
$\vert B({\rm inn}(P)) - B(J) \vert < r,$ then $B$ is Jordan-like for each 
fine partition $Q \in {\cal P}.$ (Of course$,$ it is assumed that for each 
$\delta$ mentioned in this last statement that such a partition $P$ exists 
such that $\Vert P \Vert < \delta.$) If $B$ satisfies this property for a set 
of partitions ${\cal P}$ then we say that $B$ {{\bf has an ordinary Jordan-like 
property.}}\par
\vskip 18pt
\centerline{IR5} \par
(1) We wish to measure a quantity $M$ for a compact Jordan-measurable set $J 
\subset R \subset \realp n,$ where for simplicity assume that $M$ is defined 
on and$,$ at least$,$ additive 
over members of the set $\{I_{\cal P},J\},$ and ${\cal P}$ is the set of all 
simple partitions of $R.$ As usual$,$ $v(J)$ denotes the Jordan-content.\par
(2) There is a generating function $f(\vec x\,)$ that is related to the 
functional 
$M$ in the following manner:\par
\indent\indent (i) The function $f$ is continuous on $J.$\par
\indent\indent (ii) For a some simple fine partition$,$ $Q,$ and each  
$S \in {\rm inn}(Q)$ there exist $\vec x_m \in S$ and $\vec x_M 
\in S$ such that $\hyper {f_m}=\hyper {f(\vec x_m)} = 
\Hyper {\inf} \{\hyper f(\vec 
x)\vert \vec x \in S\} = \inf \{\hyper f(\vec x\,)\vert \vec x \in S\}$ and
$\hyper {f_M} =\hyper {f(\vec x_M)} = \Hyper {\sup} \{\hyper f(\vec 
x)\vert \vec x \in S\} = \sup \{\hyper f(\vec x\,)\vert \vec x \in S\}.$ \par 
\indent\indent (iii) $(\hyper {f_m)}\, \hyper v(S) \leq \hyper M(S) 
\leq (\hyper f_M)\, \hyper v(S).$\par 
\indent\indent (iv) For $Q,\ M$ is Jordan-like.\par
\vskip 18pt
{\bf Proposition 6.6.} {\sl Suppose there exists a fine partition such that  IR5
holds for $Q.$ Then $$M(J) = \int_Jf(\vec x\,)\, 
dX.$$}\par                                                                     
Proof. Consider the fine partition $Q$ from (iii) of IR5. Then for 
each $S \in {\rm inn}(Q)$
$$ (\hyper f_m)\, \hyper v(S) \leq \hyper M(S) \leq (\hyper f_M)\, \hyper 
v(S).$$
Now $f$ is uniformly continuous on $J$ and $\vec x_m,\ \vec x_M \in S$ imply 
that $\vec x_m\approx \vec x_M$ and that $\hyper f_m \approx \hyper f_M.$ 
Therefore$,$ 
$$\hyper M(S)/\hyper v(S) \approx \hyper f(\vec x_m).$$
Consequently$,$ from (iv)
$$\int_R \hat f(\vec x\,)\, dX = \int_J f(\vec x\,)\, 
dX= \st {\hyper B(\cup A)}= B(J).$$\par        
\vskip 18pt
We also have a self-evident theorem that goes with Proposition 6.6.\par
\vskip 18pt           
{\bf Proposition 6.7. ({Self-evident Max. and Min.})} 
{\sl Let the rectangle $R\subset \realp 
n$ and suppose that compact Jordan-measurable $J \subset R.$ Let $M$ be 
defined as in (1) of IR5$,$ continuous 
$f\colon J\to \real,\ \cal P$ is the set of all simple 
partition on $R$ and any 
$P \in {\cal P}.$ If for any $S \in {\rm inn}(P)$  
it follows that 
$(f_m)v(S) \leq M(S) \leq (f_M)v(S),$ and $M$ has the ordinary Jordan-like 
property,
then the infinitesimalizing process IR3 holds.}\par
Proof. For each fine partition $Q \in {\cal P}$ the functional $M$ is Jordan-
like. Since there exists a fine partition $Q \in \hyper {{\cal P}}$ then 
the *-transfer of Theorem 6.2.3 and the hypotheses yields 
that (ii)$,$ (iii)$,$ and (iv) hold for $Q.$  The conclusion follows from 
Proposition 6.6.\par
\vskip 18pt
The Method of Constants can be greatly simplified for functionals $M$ that are 
Jordan-like.\par
\vskip 18pt        
\centerline{IR6--{Method of Constants}}
(1) We wish to measure a quantity $M$ for a Jordan-measurable set $J 
\subset R \subset \realp n,$ where for simplicity assume that $M$ is defined 
on and$,$ at least$,$ additive 
over members of the set $\{I_{\cal P},J\},$ and ${\cal P}$ is the set of all 
simple partitions of $R.$ As usual$,$ $v(J)$ denotes the Jordan-content.\par
(2) There is a generating function $f(\vec x)$ that is related to the 
functional 
$M$ in the following manner:\par
\indent\indent (i) The function $f$ is bounded on $J.$\par
\indent\indent (ii) For some simple fine partition$,$ $Q,$ and each  
$S \in {\rm inn}(Q)$ there exists some $\vec x \in S$ such that 
$\hyper M(S)= \hyper f(\vec x)\hyper v(S)$ or  $\hyper M(S)/\hyper v(S) 
\approx \hyper f(\vec x).$\par
\indent\indent (iii) The functional $M$ if Jordan-like for $Q.$ \par
\vskip 18pt
{\bf Proposition 6.8.} {\sl Suppose there exists a fine partition such that IR6
holds for $Q.$ Then $$M(J) = \int_Jf(\vec x\,)\, 
dX.$$}\par                                                                     
Proof. This is immediate from Proposition 6.5 and that $M$ is Jordan-like for
$Q.$\par
\vskip 18pt
{\bf Proposition 6.9. ({Self-evident Method of Constants})} 
{\sl Let the rectangle $R\subset \realp 
n$ and suppose that Jordan-measurable $J \subset R.$ Let $M$ be 
defined as in (1) of IR6$,$ bounded  
$f\colon J\to \real,$ where $\cal P$ 
is the set of all simple partitions on $R$ and any 
$P \in {\cal P}.$ If for any $S \in {\rm inn}(P)$  
there exists some $\vec x \in S$  such that 
$M(S) = f(\vec x)\, v(S),$ and $M$ has the ordinary Jordan-like 
property$,$ then the infinitesimalizing process IR3 holds.}\par
Proof. Again this is immediate from *-transfer and Proposition 6.5.\par
The advantages of having a Jordan-like functional are obvious when one 
compares our last proposition in the section with Theorem 6.4.3. \par
{\bf Proposition 6.10. (Extended Self-evident Method of Constants).} 
{\sl Let the rectangle $R\subset \realp 
n$ and suppose that Jordan-measurable $J \subset R.$ Let $M$ be 
defined as in (1) of IR6$,$ continuous  
$f,\, g\colon J\to \real,$ where $\cal P$ 
is the set of all simple partitions on $R$ and any 
$P \in {\cal P}.$ If for any $S \in {\rm inn}(P)$  
there exists some $\vec x_1 \in S$  and some $\vec x_2\in S$ such that 
$M(S) = f(\vec x_1)g(\vec x_2)\, v(S),$ and $M$ has the ordinary Jordan-like 
property$,$ then the infinitesimalizing process IR6 holds.}\par
Proof. The proof follows from Proposition 6.5 and the observation that 
for a fine partition $Q$ and $S \in {\rm inn}(Q)$ if $\vec x_1,\, \vec x_2
\in S,$ then $\vec x_1 \approx \vec x_2$ and $\hyper f(\vec x_1) \hyper
g(\vec x_2) \approx \hyper f(\vec x_1) \hyper g(\vec x_1).$\par               
\vfil\eject
\centerline{Appendix For Chapter 7.}\par
\bigskip
{\bf Theorem 7.2.2.} {\sl Let Jordan-measurable $\ J \subset R 
\subset \realp n,$ bounded $f\colon J \to \real,$ and $Q$ be any fine 
partition of $\hyper R.$ Let hyperfinite $\hyper {\cal S}(Q) = \{S_i 
\bigm \vert 0\leq i \leq \Gamma -1\}.$  Assume that there exists a 
hyperfinite sequence $U_i,$ where $U_i \in S_i$ for each $S_i \in {\cal S}(Q).$ 
Then 
$$\st {\!\!\!\!\!\!\sum_{S(i) \in {\rm inn}(Q)}\!\!\!\!\!\!\hyper f(U_i) \hyper v(S_i)}= 
\st {\!\!\!\!\!\sum_{S(i) \in {\cal S}(Q)} \!\!\!\!\!\!\hyper 
{\hat f}(U_i) \hyper v(S_i)}.$$}\par
Proof. First$,$ recall the following definitions.
The boundary subrectangles $\partial (Q) = \{S\bigm \vert S \in \hyper 
{\cal S}(Q) \land S \cap \hyper J 
\not= \emptyset \land S \cap (\hyperrealp n - \hyper J) \not= \emptyset\};$  
the exterior subrectangles ${\rm ext}(Q) = \{S\bigm \vert S \in\hyper 
{\cal S}(Q) \land S \subset (\hyperrealp n - \hyper J) \},$  and the inner 
subrectangles ${\rm inn}(Q) = {\cal S}(Q) - (\partial (Q) 
\cup {\rm ext}(Q)).$ Notice that we have used the notation $\partial, \ 
{\rm inn}$ and ${\rm ext}$ rather than the notation $\hyper 
 \partial, \, \hyper {\rm ext}$ and $\Hyper {\rm inn}$ even though these sets 
are generated by obvious extended standard set-theoretic operators. 
It follows directly from the definition of Jordan-content that 
$J$ has Jordan-content if and only if  for every fine partition $Q$
$$\sum_{S(i) \in \partial (Q)}\!\!\!\!\!\! \hyper v(S_i) \approx 0.$$
Since there exists $m,\ M \in \real$ such that for all
$\vec x \in R,\ m \leq \hat f(\vec x )\leq M$ then noting that 
$\partial (Q)$ is internal and$,$ hence$,$ a hyperfinite set we have from *-
transfer of the finite case that  
$$m\!\!\!\!\!\sum_{S(i) \in \partial (Q)}\!\!\!\!\!\! \hyper v(S_i) \leq 
\sum_{S(i) \in \partial (Q)}\!\!\!\!\!\! 
\hyper {\hat f}(U_i)\hyper v(S_i) \leq M\!\!\!\!\!\sum_{S(i) \in 
\partial (Q)}\!\!\!\!\!\! \hyper v(S_i).$$
Thus $\sum_{S(i) \in \partial (Q)} 
\hyper {\hat f}(U_i)\hyper v(S_i) \approx 0.$\par
The sets ${\rm inn}(Q)$ and ${\rm ext}(Q)$ are also hyperfinite sets and from 
the properties of finite addition and the fact that $\hyper {\hat f(U_i)} = 0$ for 
each $U_i \in S_i \in {\rm ext}(Q)$ it follows that
$$\sum_{S(i) \in {\cal S}(Q)} \!\!\!\!\!\!\hyper {\hat f}(U_i) \hyper v(S_i) =
\!\!\!\!\!\!\sum_{S(i) \in {\rm inn}(Q)} \!\!\!\!\!\!\hyper f(U_i) 
\hyper v(S_i)\, +
\!\!\!\!\!\!\sum_{S(i) \in \partial (Q)} \!\!\!\!\!\!\hyper {\hat f}(U_i) 
\hyper v(S_i)\, +$$
$$\sum_{S(i) \in {\rm ext}(Q)} \!\!\!\!\!\!\hyper {\hat f}(U_i) \hyper v(S_i) 
\approx \!\!\!\!\!\!\sum_{S(i) \in {\rm inn}(Q)} \!\!\!\!\!\!\hyper f(U_i) \hyper 
v(S_i).$$
Since 
$$m\!\!\!\!\!\sum_{S(i) \in {\cal S}(Q)}\!\!\!\!\!\! \hyper v(S_i) \leq 
\sum_{S(i) \in {\cal S} (Q)}\!\!\!\!\!\! 
\hyper {\hat f}(U_i)\hyper v(S_i) \leq M\!\!\!\!\!\sum_{S(i) \in 
{\cal S}(Q)}\!\!\!\!\!\! \hyper v(S_i)$$
and $$\sum_{S(i) \in {\cal S}(Q)}\!\!\!\!\!\! \hyper v(S_i) =v(R)$$
then $$\sum_{S(i) \in {\cal S} (Q)}\!\!\!\!\!\! 
\hyper {\hat f}(U_i)\hyper v(S_i)$$
is a limited number and the result follows. 
\vfil\eject
\centerline{Appendix For Chapter 8.}\par
\vskip 18pt
{\bf Theorem 8.3.1.} {\sl For each $\eps \in \monad 0$ the set
$o(\eps) = \{\eps\, h \bigm \vert h \in \monad 0 \}$ is an ideal in $\monad 
0.$}\par
Proof. First$,$ we show that $o(\eps)$ is a subring of the ring $\monad 0$. Since 
$o(\eps) \subset \monad 0,$ then consider arbitrary $\eps \, h_1 \in o(\eps)$ 
and $\eps \, h_2 \in o(\eps).$ Observe that $(\eps \, h_1)-(\eps \, h_2)=
\eps\, (h_1 - h_2) = \eps\,  h_3,$ where $h_3 \in \monad 0,$ and $(\eps \, 
h_1)(\eps \, h_2)= \eps \,(\eps \, h_1 h_2),$ where $(\eps \, h_1 h_2)\in 
\monad 0).$ Thus $o(\eps)$ is a subring of $\monad 0.$ Let $\delta \in \monad 0.$
Then $\delta\, (\eps \, h_1) = \eps \,(\delta \, h_1)$ implies that $o(\eps)$ 
is an ideal in $\monad 0.$  \par

\vskip 18pt                                               
{\bf Theorem 8.3.2.} {\sl Let $\eps \in {\monad 0}^+.$ Suppose that
$w \in \hyperreal$ and $0\leq w\leq \eps\, h \in o(\eps).$ Then 
$w \in o(\eps).$}\par
Proof. Obviously$,$ if $w =0$ or $w = \eps\, h,$ then result holds. Hence$,$ 
suppose that $0 < w < \eps\,h.$ Now $0 < w/\eps < h$ implies that
$w/\eps \in \monad 0$ implies that $\eps\, (w/\eps) = w \in o(\eps).$ \par

\vskip 18pt
{\bf Theorem 8.3.3.} {\sl Let $\eps,\, \delta \in {\monad 0}^+.$ If 
$0 \leq \delta \leq \eps,$ then $o(\delta) \subset o(\eps)$ and $o(\delta)$ is a 
ideal in $o(\eps).$}\par
Proof. We know that $o(\delta)$ is an ideal in $\monad 0.$ Thus let $\delta\,h 
\in o(\delta)$ and $\eps\,h_1 \in o(\eps).$ Then $\delta\,(\eps\, hh_1) \in 
o(\delta)$ implies that $o(\delta)$ is an ideal in $o(\eps).$ To show that 
$o(\delta) \subset o(\eps)$ we need only consider the nonnegative elements.
Let $0\leq \delta \, h \in o(\delta).$ Then $0\leq \delta\, h < \eps\, h$ 
implies by Theorem 8.3.3 that $\delta\, h \in o(\eps).$\par

\vskip 18pt
Are there order ideals such that $o(\delta) \subset o(\eps)$ and $o(\delta) \not= 
o(\eps)$?\par
(1) Let $0 < \delta \leq \eps,\ \delta,\eps \in {\monad 0}^+.$ Then 
$o(\delta\eps) \subset o(\eps)$ and $o(\delta\eps) \not= o(\eps).$\par
Proof. $\delta\, \eps \in o(\eps),$ but $\delta\, \eps \notin o(\delta\, 
\eps)$ for  $\delta\, \eps\, h = \delta\, \eps\Rightarrow h = 1.$ Indeed$,$ if 
$r \in {\cal O}- \monad 0,$ then $\delta\,\eps\,r \in o(\eps),$ but
$\delta\, \eps\,r \notin o(\delta\,\eps).$

\vskip 18pt
(2) Let $n \in {\nat}^+ = \nat - \{0\},\ \eps \in {\monad 0}^+.$ Then 
$o(\eps^n) \subset o(\eps^{n-1}) \subset \cdots \subset o(\eps)$ and
$o(\eps^i) \not= o(\eps^j);\ 1 \leq i,j \leq n;\ i \not= j.$\par
Proof. From (1).

\vskip 18pt
Let $\{\eps_1,\ldots,\eps_k\} \subset \monad 0.$ Define
$o(\eps_1,\ldots,\eps_k) = \{\eps_1\, h_1 + \cdots + \eps_k\, h_k \bigm \vert
h_i \in \monad 0 \land 1\leq i\leq k\}.$\par

\vskip 18pt
{\bf Theorem 8.3.4.} {\sl Let $\eps =\max{\{\vert \eps_1\vert,\ldots,\vert 
\eps_k \vert\}}.$ Then $o(\eps_1,\ldots,\eps_k) = o(\eps).$}\par
Proof. Assume that $x = \sum_{i=1}^n\eps_i\,h_i \in o(\eps_1,\ldots,\eps_n).$
Let $\delta \in \monad 0$ and $ \delta \geq \max\{\vert h_1\vert,\ldots, \vert 
h_n \vert \}.$ Then 
$$0\leq \vert \sum_{i=1}^n\eps_i\,h_i \vert \leq \sum_{i=1}^n\vert \eps_i 
\vert\, \vert h_i\vert \leq n\eps\, \delta = \eps\, (n\delta) \in o(\eps).$$
Thus $o(\eps_1,\ldots,\eps_n) \subset o(\eps).$ \par
On the other hand$,$ let $\eps\, h \in o(\eps).$ Since $\eps = \eps_j$ for some 
$j$ such $1\leq j\leq n$ then $\eps\, h \in o(\eps_1,\ldots,\eps_n)\Rightarrow
o(\eps) \subset o(\eps_1,\ldots,\eps_n).$\par

\vskip 18pt
{\bf Theorem 8.3.5.} {\sl Let $\eps =\max{\{\vert \eps_1\vert,\ldots,\vert 
\eps_k \vert\}}.$ Then
 $$o(\eps_1,\ldots,\eps_k) = o(\sqrt {\eps_1^2 + \cdots
+\eps_k^2}) = o(\eps).$$}\par
Proof.  By *-transfer $\sqrt{\sum_{i=1}^n \eps_i^2} \leq \sum_{i=1}^n \vert
\eps_i\vert.$ Now for $h \in \monad 0$ it follows that
$\vert h\vert \sqrt{\sum_{i=1}^n \eps_i^2}\leq n\vert h\vert \,\eps;$ 
in which case by Theorem 8.3.3
$\vert h\vert\sqrt{ \sum_{i=1}^n \eps_i^2} \in o(\eps)\Rightarrow 
h\sqrt{ \sum_{i=1}^n \eps_i^2} \in o(\eps).$ Hence$,$ 
$o(\sqrt{\sum_{i=1}^n \eps_i^2}) \subset o(\eps).$\par
Conversely$,$ let $\eps\, h \in o(\eps).$ Then $\eps\, \vert h\vert \in o(\eps).$
But $\vert \eps \, h \vert = \eps\, \vert h\vert = \vert \eps_j \vert \vert 
h\vert \leq \vert h\vert \sqrt{\sum_{i=1}^n \eps_i^2}.$ Once again Theorem 8.3.3
implies that $\eps\, h \in o(\sqrt{\sum _{i=1}^n\eps_i^2}).$\par
 
\vskip 18pt
{\bf Theorem 8.4.1.} {\sl Let $1 \leq n \in \nat.$ Suppose that $f^{(n-
1)}\colon [a, nb] \to \real$ and that $f^{(n)}\colon (a,nb) \to \real,$ where
$f^{(k)}$ denotes the {\it k}th derivative of $f.$ Then there exists some $t \in
(a,nb)$ such that $\Delta^n f(a,b) = f^{(n)}(t)\, b^n.$}\par
Proof. By induction. We know the result holds for $n = 1,$ therefore$,$ assume 
that result holds for $n-1$ and let $g(x) = f(x + b) - f(x),\ x \in [a, a+(n-
1)b].$ Now $g^{(n-1)}(x) = f^{(n-1)}(x + b) - f^{(n-1)}(x) \in \real$ for 
each $x \in [a, a+(n-1)b].$ Thus$,$ by induction$,$ there exists some
$t_0 \in (a,a + (n-1)b)$ such that $\delta^{(n-1)} g(a,b) = g^{(n-
1)}(t_0)b^{n-1}= (f^{(n-1)}(t_0 + b) - f^{(n-1)}(t_0))\, b^{n-1}.$
Applying the mean value theorem yields that there exists some $t \in
(t_0,t_0 + b)$ such that $g^{(n-1)}(t_0) = f^{(n-1)}(t_0 + b) - f^{(n-
1)}(t_0)=f^n(t)\, b.$ Thus$,$ $\Delta^nf(a,b) = \Delta^{(n-1)}g(a,b) =
f^n(t)b^n,$ where $t \in (a,a + nb).$\par
{\bf Corollary 8.4.1.1} {\sl  Let $1 \leq n \in \nat.$ Suppose that $f^{(n-
1)}\colon [a, b] \to \real$ and that $f^{(n)}\colon (a,b) \to \real,$ then for 
each $dx\in {\monad 0}^+$ and $c \in \Hyper [a,b),$ there exists some $t  \in 
(c,c +ndx)$ such that $\Delta^n\hyper f(c,c+dx) = \hyper f^{(n)}(t)\,(dx)^n.$}\par   
Proof. This follows from *-transfer and the fact that $[c,c+n\,dx] \subset 
\Hyper [a,b].$\par

\vskip 18pt
{\bf Theorem 8.4.2.} {\sl Let $1 \leq n \in \nat.$ Suppose that $f^{(n-
1)}\colon [a, b] \to \real$ and that $f^{(n)}\colon (a,b) \to \real.$ 
If $c \in (a,b),$ then for each $dx\in \monad 0,\ dx \geq 0 $ {\rm [}resp. $dx 
<0${\rm ]}
$$f^n(c)\, (dx)^n \approx \Delta^n \hyper f(c,c+dx), {\rm [resp.}f(c+dx,c){\rm 
]} \pmod {o((dx)^n)}.$$}\par
Proof. Let $c \in (a,b)$ and consider any $r \in \real^+$ such that
$[c,c+n\,r] \subset [a,b].$ Define $g\colon [c,c+(n-1)r] \to \real$ by
$g(y) = f(y + r) - f(y).$ Theorem 8.4.1 implies that there exists some
$t_1 \in (c,c+(n-1)r)$ such that $\Delta^{(n-1)}g(c,r) =g^{(n-1)}(t_1)
\, r^{n-1}.$
By *-transfer it follows that for $dx \in {\monad 0}^+$ there exists
some $t \in (c,c + (n-1)dx)$ such that $\Delta^{(n-1)}\hyper g(c,dx) =
\hyper g^{(n-1)}(t)\, dx^{(n-1)}.$ Consequently$,$ by *-transfer$,$ one obtains 
$$\Delta^n\hyper f(c,dx) = {{\hyper f^{(n-1)}(t + dx) - \hyper f^{(n-
1)}(t)}\over{dx}}\, dx^n.$$
Consequently,
$${{\Delta^n\hyper f(c,dx)}\over{dx^n}}= {{\hyper f^{(n-1)}(t + dx) - 
\hyper f^{(n- 1)}(t)}\over{dx}} =$$
\eject
$$\left({{\hyper f^{(n-1)}(t + dx) - f^{(n- 1)}(c)}\over{(t+dx) - 
c}}\right)\left({{(t+dx)-c}\over{dx}}\right)\ + \ \ \ \ \ \ \ \ \ \ \ \ \ \ \ \ \ \ \ \ $$
$$\ \ \ \ \ \ \ \ \ \ \ \ \ \ \ \ \ \ \left({{f^{(n-1)}(c) - \hyper f^{(n-1)}(t)}\over{c-t}}\right)\left({{c-
t}\over{dx}}\right) =$$
$$\left({{\hyper f^{(n-1)}(c + dx_1) - f^{(n- 1)}(c)}\over{dx_1}}\right)
\left({{dx_1}\over{dx}}\right)\ + \ \ \ \ \ \ \ \ \ \ \ \ \ \ \ \ \ \ \ \ $$
$$\ \ \ \ \ \ \ \ \ \ \ \ \ \ \ \ \ \ \ \left({{\hyper f^{(n-1)}(c+dx_2) - f^{(n-
1)}(c)}\over{dx_2}}\right)\left({{-dx_2}\over{dx}}\right),$$ 
where $dx_1 = t + dx -c,\ dx_2 = t - c \in \monad 0.$
We observe that $1 < dx_1/dx < n$ and that $0 < dx_2/dx < n-1.$ Recall 
that $$\St {\left(\left({{\hyper f^{(n-1)}(c + dx_1) - f^{(n- 1)}(c)}\over{dx_1}}
\right) \left({{dx_1}\over{dx}}\right)\right)}= f^n(c)\St 
{\left({{dx_1}\over{dx}}\right)}$$ 
 $$\St {\left(\left({{\hyper f^{(n-1)}(c+ dx_2) -  f^{(n-
1)}(c)}\over{dx_2}}\right)\left({{-dx_2}\over{dx}}\right)\right)}= f^n(c)
\St {\left({{-dx_2}\over{dx}}\right)}$$ 
Therefore,
$$\left({{\hyper f^{(n-1)}(c + dx_1) - f^{(n- 1)}(c)}\over{dx_1}}\right)
\left({{dx_1}\over{dx}}\right)\ + \ \ \ \ \ \ \ \ \ \ \ \ \ \ \ \ \ \ \ \ $$
$$\ \ \ \ \ \ \ \ \ \ \ \ \ \ \left({{\hyper f^{(n-1)}(c+dx_2) - f^{(n-
1)}(c)}\over{dx_2}}\right)\left({{-dx_2}\over{dx}}\right)\approx$$ 
$$f^n(c)\left(\St {\left({{dx_1}\over{dx}}\right)} + \St 
{\left({{-dx_2}\over{dx}}}\right)\right) = $$ $$f^n(c)\left(\St {\left({{dx_1}\over{dx}} 
- {{dx_2}\over{dx}}\right)}\right) =f^n(c).$$ \par
\vskip 18pt
{\bf From page 71.} {\sl A simple proof shows that if $r,\,r_0 \in \Hyper [-1,1],$ then the best we can 
say is that $$\Vert \hyper c(t + r\,\eps) - \vec g(t + r_0\, \eps)\Vert 
\in o(\eps).\eqno (6)$$}\par
Proof.$\Vert \hyper c(t + r\,\eps) - \vec g(t + r_0\, \eps)\Vert 
= \Vert(h_1,\ldots,h_p)r\,\eps \Vert = \lambda \, \eps, \ \ \lambda \in 
\monad 0.$ \par

\vskip 18pt
{\bf From page 72.} {\sl Considering any $(s,t)\, (s_0,t_0) \in \Hyper [-1,1] 
\times \Hyper [-1,1]$ a simple proof yields that                    
$$\Vert \Hyper {\vec r}(u_0 + s\,\delta,v_0 + t\,\eps) - 
\Hyper {\vec k}(u_0 + s_0\,\delta,v_0 + t_0\,\eps)\Vert \in o(\max 
\{\delta,\eps\}) = o(\lambda).\eqno (4)$$} \par
Proof. $\Vert \Hyper {\vec r}(u_0 + s\,\delta,v_0 + t\,\eps) - 
\Hyper {\vec k}(u_0 + s_0\,\delta,v_0 + t_0\,\eps)\Vert = 
\sqrt {s^2\delta^2h_s^2 + t^2\eps^2h_t^2}.$ Hence
$$\vert s\,\delta\,h_s\vert\leq \sqrt {s^2\delta^2h_s^2 + t^2\eps^2h_t^2} \leq
\vert s\,\delta\,h_s\vert +\vert t\,\eps\,h_t\vert,$$                          
$$\vert t\,\eps\,h_t\vert\leq \sqrt {s^2\delta^2h_s^2 + t^2\eps^2h_t^2} \leq
\vert s\,\delta\,h_s\vert +\vert t\,\eps\,h_t\vert,$$                          
implies that
$$\left| {{s\,\delta\,h_s}\over{2}}\right|+\left| {{t\,\eps\,h_t}\over{2}}\right|
\leq\sqrt {s^2\delta^2h_s^2 + t^2\eps^2h_t^2} \leq
\vert s\,\delta\,h_s\vert +\vert t\,\eps\,h_t\vert.$$ 
Consequently$,$ 
$$\lambda_s\,\delta + \lambda_t\,\eps \leq \sqrt {s^2\delta^2h_s^2 + t^2\eps^2h_t^2}
\leq \gamma_s\,\delta + \gamma_t\,\eps,$$
where $\lambda_s,\,\lambda_t,\,\gamma_s,\,\gamma_t \in \monad 0.$
The result follows from Theorem 8.3.4.\par                         

\vskip 18pt
{\bf From page 75.} Prior to discussing the McShane integral$,$ a general 
theorem that establishes that $\lambda$-fine partitions exist for each 
$L$-gauge is a useful addition to this appendix.\par

\vskip 18pt
{\bf Proposition 8.1.} {\sl Let $X$ be a nonempty  connected topological 
space. Let $\cal C$ be a collection of subsets of $X$ with the following two 
properties.\par
(i) If $x \in X,$ then there exists some $C \in {\cal C}$ such that
$x  \in {\rm int}(C).$\par
(ii) If $C_1,\ C_2 \in {\cal C}$ and $C_1 \cap C_2 \not= \emptyset,$ then 
$C_1 \cup C_2 \in {\cal C}.$\par
The following two statements are equivalent. \par
(A) The space $X$ is compact.\par
(B) The set $X$ is an element of every collection of subsets of $X$ that 
satisfy (i) and (ii).}\par
Proof. Assume that $X$ is compact and the set ${\cal C}$ satisfies (i) and 
(ii) above. Then there exist a finite subset $\{C_1,\ldots,C_n\}$ of 
${\cal C}$ such that $X = \bigcup \{{\rm int}(C_i)\bigm\vert
1\leq i \leq n\}= \bigcup\{C_i\bigm\vert 1\leq i\leq n\}$ from (i) and 
compactness. Consider $C_1$ and assume that $n > 1.$ Then there exists some 
$C_k$ such that $C_1 \cap C_{k_1} \not= \emptyset$ from the connectedness of 
$X.$ From (ii) $A_1 = C_1 \cup C_{k_1} \in {\cal C}.$ If $n >2,$ then there 
exists some $C_{k_2},$ where $C_{k_2} \notin \{C_1,C_{k_1}\},$ such that 
$A_1 \cap C_{k_2} \not= \emptyset.$ Again from (ii) $A_2 =A_1 \cup C_{k_2} 
\in {\cal C}.$ Continuing this finite process leads to the conclusion that
$X = \bigcup\{C_i\bigm\vert 1\leq i\leq n\} \in {\cal C}.$\par
Conversely$,$ assume that (B) holds and let $\cal G$ be an open cover of $X.$ 
Recall that $F({\cal G})$ is the set of all finite subsets of ${\cal G}.$
Let ${\cal C} = \{A\bigm \vert A\subset X \land \exists Y(Y \in F({\cal 
G})\land (A \subset \bigcup Y))\}.$ Obviously$,$ since ${\cal G}\subset {\cal C},$ (i) 
holds. Moreover$,$ if $C_1,\ C_2 \in {\cal C}$ and $C_1 \subset \bigcup Y_1$ and 
$C_2 \subset \bigcup Y_2,$ then $Y_1 \cup Y_2 \in F({\cal G})$ and 
$C_1 \cup C_2 \subset \bigcup \{Y_1 \cup Y_2\}$ imply that
$C_1 \cup C_2 \in {\cal G}.$ Thus (ii) holds. Hence$,$ $X \in {\cal C}$ implies 
that $X$ is compact. This completes the proof. \par

\vskip 18pt
I point out that the partitions that are used in the gauge type integrals and 
other generalizations of the Riemann integral often are not considered to 
generate closed subrectangles; but$,$ rather are considered to generate  
``left-closed'' or ``right-closed'' subrectangles. This is what is done by 
McShane in his definition of the McShane integral. Thus referring to Chapter 
5$,$ Section 5.1$,$ replace the set of closed subrectangles obtained by 
considering the expansion of 
$([x_{10}, x_{11}] \cup \cdots \cup [x_{1{k-1}},x_{1k}]) \times 
\cdots \times ([x_{n0}, x_{n1}] \cup \cdots \cup [x_{n{m-1}},x_{nm}])$ by the 
collection ${\cal A} = \{A_i\bigm \vert 1 \leq i\leq m\}$ of n-dimensional 
right-closed subrectangles obtained by considering the expansion of $((x_{10}, x_{11}] 
\cup \cdots \cup (x_{1{k-1}},x_{1k}]) \times 
\cdots \times ((x_{n0}, x_{n1}] \cup \cdots \cup (x_{n{m-1}},x_{nm}]).$  
This gives a collection of right-closed rectangles that subdivides  
$(x_{10},x_{1k}] \times \ldots \times (x_{n0},x_{nm}] = 
R_{rc}$ into nonoverlapping right-closed subrectangles. However$,$ when 
considering an  
intermediate partition $Q,$ then it is often (but not always) required 
that each $\vec x_i \in Q$ also has the property that $ \vec x_i\in 
\overline{A}_i.$ This is the case with the $\delta$-fine partitions discussed
in {\tt Mawhin [1985]}$,$ and {\tt Swartz and Thomson [1988]}. This is \underbar{not} the 
case with the $\lambda$-fine partition as defined by McShane. \par
Let the closed rectangle $R = \overline{R}_{rc}\subset \realp n.$ Recall that 
$B \subset  
\realp n$ is a nonempty {\bf open n-cell} if $B = (a_1,b_1) \times \cdots \times 
(a_n,b_n),\  a_i < b_i,\ 1 \leq i \leq n.$ Let ${\cal B}$ be the set of all 
open n-cells for $\realp n.$ In order to establish results with the least 
amount of notation$,$ call a map  $\lambda \colon \realp n \to {\cal B}$ a 
{\bf $L$-gauge} if for each $\vec x \in \realp n,\ \vec x \in \lambda (\vec x).$ 
(In all cases and without further mention$,$ when one concentrates upon a domain 
such as $R$ then the gauge is 
considered to be the above map restricted to $R.$) Intuitively$,$ think of a 
 $L$-gauge as simply carving out specific open n-cell neighborhoods for each 
member of $\realp n.$ Now there is associated with gauge integration the 
pairs $(Q,P)$ which are termed a {\bf P-partition of $R$} if\par 
 (i) $P= \{A_1,\ldots,A_k\}$ is a set of 
right-closed subrectangles that subdivides $R_{rc}$ into nonoverlaping 
right-closed subrectangles in the usual manner and\par
(ii) the finite sequence $Q = \{\vec 
x_1,\ldots,\vec x_k\}$ has the properties that each $\vec x_i \in R.$ 
More significantly$,$ however$,$ are the {\bf $\lambda$-fine partitions} that are 
P-partitions with the additional property that\par
(iii)  each $A_i \subset \lambda(x_i),$ for $1 \leq i\leq k,$ where $\lambda$ 
is a $L$-gauge.\par
 What needs to be determined$,$ however$,$ is whether 
or not there exists at least one $\lambda$-fine partition for a given 
$L$-gauge. 

\vskip 18pt
{\bf Proposition 8.2.} {\sl Let $\lambda \colon \realp n \to {\cal B}$ be an 
$L$-gauge and $R \subset \realp n.$ Then there exists at least one 
$\lambda$-fine partition of $R.$}\par
Proof. Let ${\cal C}$ be the set of all (closed) subrectangles of compact $R$ 
for which 
there exists an $\lambda$-fine partition. First$,$ we show that $\cal C$ 
satisfies (i) of Proposition 8.1. Let $\vec x \in R.$ Then $\vec x \in 
\lambda (\vec x)$ and $\lambda (\vec x)$ is an open neighborhood of $\vec x.$
Hence$,$ $R \cap \lambda(\vec x)$ is a nonempty $R$-open neighborhood of $\vec 
x.$ Moreover$,$ since there exists a closed subrectangle $I \subset \realp n$ 
such that $\vec x \in {\rm int}_{\srealp n}(I) \subset I \subset 
\lambda(\vec x),$   
it follows that 
$\vec x \in R \cap {\rm int}_{\srealp n}(I)\subset {\rm int}_R(R \cap I)
\subset I \subset \lambda (\vec x).$ 
Now simply consider the $\lambda$-fine partition $\{(\vec x,R_{rc}\cap 
I_{rc})\}$ of $(R \cap I)_{rc}.$ 
Since $R\cap I= \overline{(R \cap I)}_{rc}$ is a (closed) subrectangle of $R$ 
then $R\cap I \in {\cal C}$ 
and $\cal C$ satisfies (i) of Proposition 8.1.\par

Let $C_1,\ C_2 \in {\cal C}$ and $C_1 \cap C_2 \not= \emptyset.$ 
Let $P_{\lambda}= \{A_1,\ldots,A_k\}$ and $P_{\lambda}^\prime
= \{A_1^\prime,\ldots,A_p^\prime\}$ be the two subdividing collections of
right-closed subrectangles and $Q_{\lambda},\ Q_{\lambda}^\prime$ 
the corresponding 
finite sequences for $C_1$ and $C_2,$ respectively. These exist from the 
definition of ${\cal C}.$ Now 
$P_{\lambda}= \{A_1,\ldots,A_k\}$ and $P_{\lambda}^\prime 
= \{A_1^\prime,\ldots,A_p^\prime\}$ are generated by partition $P$ and 
$P^\prime,$ respectively. Considering the common refinement $P^*$ then it 
follows$,$ 
since any new right-closed subrectangle generated by this common refinement 
is a subset of some member of  $P_{\lambda}\cup P_{\lambda}^\prime,$ that a 
finite sequence of members of $Q_{\lambda} \cup Q_{\lambda}^\prime$ may be 
chosen in such a manner that selection yields a possible new pair 
$(P_{\lambda}^*,Q_{\lambda}^*)$ which forms a $\lambda$-fine partition for 
the closed subrectangle
$C_1 \cup C_2.$ Thus $\cal C$ satisfies part (ii) of Proposition 8.1. Since 
$R$ is compact then $R \in {\cal C}$ and the proof is complete.\par

\vskip 18pt
In order to show the existence of a L-micropartition$,$ a special type of 
internal function$,$ $\hat{\lambda},$ termed a {\bf L-microgauge} is 
needed.\par

\vskip 18pt
{\bf Proposition 8.3.} {\sl Let $L(\realp n,{\cal B})$ be the set of all 
$L$-gauges. Then there exists an internal map $\hat{\lambda} \colon
\hyperrealp n \to \hyper {\cal B}$ such that for each $\vec x \in \hyperrealp 
n$ and every $\lambda \in L(\realp n,{\cal B}),$ 
$$\hat{\lambda}(\vec x) \subset \hyper \lambda(\vec x).$$}\par
Proof. Consider the internal binary relation
$$B = \{(x,y)\bigm\vert x,\ y \in \Hyper (L(\realp n,{\cal B}))\land
\forall z(z \in \hyperrealp n \to y(z) \subset x(z))\}.$$
We show that $B$ is concurrent$,$ at least$,$ on the set $^\sigma(L(\realp n,{\cal 
B})).$ Let $\{(x_1,y_1),\ldots,(x_n,y_n)\} \subset B$ and $\{x_1,\ldots,x_n\}
\subset {^\sigma(L(\realp n,{\cal B}))} = \{\hyper \lambda \bigm\vert
\lambda \in L(\realp n,{\cal B})\}.$ By *-transfer of the result that for any 
$\vec x \in \realp n$ the intersection of any nonempty finite subset of  
${\cal B}$ each member of which contains $\vec x$ is itself a member of 
${\cal B}$$,$ we can define     
the internal map 
$f\colon \hyperrealp n \to \hyper {\cal B}$ as follows: For each 
$\vec x \in \hyperrealp n$ let $f(\vec x) 
= x_1(\vec x) \cap \cdots \cap x_n(\vec x).$  This yields  
$f \in \Hyper (L(\realp n,{\cal B}))$ and $\{(x_1,f),\ldots,(x_n,f)\} \subset 
B.$ Consequently$,$ $B$ is concurrent$,$ at least$,$ on $^\sigma(L(\realp n,{\cal 
B})).$ But$,$ we are working in a polysaturated enlargement. This implies that 
there exists some internal $\hat{\lambda} \in \Hyper (L(\realp n,{\cal B}))$ 
such that for each $\vec x \in \hyperrealp n$ and each $\lambda \in 
L(\realp n,{\cal B})$ it follows that $\hat{\lambda}(\vec x) \subset
\hyper \lambda(\vec x)$ and completes the proof.\par

\vskip 18pt
The object $\hat{\lambda}$ that exists by Proposition 8.3 is called 
a {\bf L-microgauge} and has many interesting properties. For example$,$ 
if $\vec p \in R,$ then $\hat{\lambda}(\vec p) \subset \monad {\vec p},$ but
$\hat{\lambda}(\vec p) \not= \monad {\vec p}$ since $\hat{\lambda}(\vec p)$ is 
an internal *-open member of $\hyper {\cal B};$ while $\monad {\vec p}$ is 
an external set which is the union of *-open members of $\hyper {\cal B}.$
[{\tt Herrmann [1978]}] However$,$ our major use of a  L-microgauge
$\hat{\lambda},$ in this appendix$,$ is to obtain a {\bf L-micropartition}
$(Q_{\hat{\lambda}},P_{\hat{\lambda}})$ which now exists by *-transfer of the 
results in Proposition 8.2. Observe that  $(Q_{\hat{\lambda}},
P_{\hat{\lambda}})$ is also a $\hyper \lambda$-fine partition for each 
standard $L$-guage $\lambda$ since each internal $A_i \subset \hat 
\lambda(\vec x_i) \subset \hyper {\lambda}(\vec x_i).$\par
McShane defined his integral in the following manner. Let $\lambda$ be a 
$L$-gauge. Then there exists a map ${\cal L}$ such that ${\cal L}(\lambda)$ 
is the 
nonempty set of all $\lambda$-fine partitions of $R.$ Let $f\colon R \to 
\real.$ Then {\bf $f$ is McShane integrable over $R$} if 
there exists a real number $J$ such that: for each positive real number $r$
there corresponds a $L$-gauge $\lambda$ such that for each  
$(Q_\lambda,P_\lambda) \in {\cal L}(\lambda),$
$$\left| \sum_{i=1}^kf(\vec x_i)\,v(A_i) - J \right| < r,$$
where $v(A_i) = v(\overline{A}_i).$ Since McShane's integral is equivalent to 
the Lebesque integral$,$ our final results in this particular appendix are 
stated in terms of the Lebesque integral.\par

\vskip 18pt
{\bf Proposition 8.4.} {\sl If $f\colon R \to \real$ is Lebesgue integrable 
with $\int_R f(\vec x)\, dX = J,$ then for each L-microgauge $\hat{\lambda}$ 
and for each $(Q_{\hat \lambda},P_{\hat \lambda}) \in 
\hyper {\cal L}(\hat \lambda),$
$$ \sum_{i=1}^\Gamma \hyper f(\vec x_i) \, \hyper v(A_i) \in \monad J.$$}
Proof. Suppose that $f\colon R \to \real$ is Lebesgue integrable 
with $\int_R f(\vec x)\, dX = J.$ Consider arbitrary positive real $r.$
The following sentence
$$\exists x(\forall y(x \in L(\realp n,{\cal B}) \land y \in {\cal L}(x) \to 
\left|\sum_{i=1}^k f(\vec x_i)\, v(A_i) - J \right| < r))$$                
holds in $\Hyper {\cal M}$ by *-transfer. Now let $\hat \lambda$ be a 
L-microgauge and $(Q_{\hat \lambda},P_{\hat \lambda}) \in 
\hyper {\cal L}(\hat \lambda).$
From the above observation $(Q_{\hat \lambda},P_{\hat \lambda})\in 
\hyper {\cal L}(\hyper \lambda)$ for any standard $L$-gauge $\lambda.$ 
Thus$,$ independent of any particular $L$-gauge and hence any particular 
positive real $r$
$$\left|\sum_{i=1}^\Gamma \hyper f(\vec x_i)\, \hyper v(A_i) -J\right| < r,$$
for each $(Q_{\hat \lambda},P_{\hat \lambda}) \in \hyper {\cal L}(\hat 
\lambda).$ Consequently, $\sum_{i=1}^\Gamma \hyper f(\vec x_i)\hyper v(A_i) \in \monad J$
for each $(Q_{\hat \lambda},P_{\hat \lambda}) \in \hyper {\cal L}(\hat 
\lambda)$ and the proof is complete.\par

\vskip 18pt
{\bf Proposition 8.5.} {\sl Suppose that $f\colon R \to \real$ and there 
exists a real number $J,$ and a L-microgauge $\hat \lambda$ such that 
for each $(Q_{\hat \lambda},P_{\hat \lambda}) \in \hyper {\cal L}(\hat 
\lambda),$
$$\sum_{i=1}^\Gamma \hyper f(\vec x_i) \, \hyper v(A_i) \in \monad J.$$
Then $f$ is Lebesgue integrable on $R$  and $\int_R f(\vec x)\, dX 
= J.$}\par
Proof. Let $r$ be any positive real number. Then 
$\left|\sum_{i=1}^\Gamma \hyper f(\vec x_i) \, \hyper v(A_i)- J\right| < r.$
Thus the sentence
$$\exists x(\forall y(x \in \Hyper (L(\realp n,{\cal B})) \land y \in 
\hyper {\cal L}(x) \to 
\left|\sum_{i=1}^k \hyper f(\vec x_i)\, \hyper v(A_i) - J \right| < r))$$                
holds in $\Hyper {\cal M},$ hence$,$ in ${\cal M}$ and the result follows.\par  
\vfil\eject
\centerline{Appendix For Chapter 10.}\par
\vskip 18pt
{\bf Proposition 10.1} {\sl Let standard $A \in {\cal H}.$ Then $^\sigma A$ is 
internal if and only 
if $A$ is finite.}\par                                                          
Proof. If $A$ is empty$,$ then $\hyper A = {^\sigma A} = \emptyset$ and thus
$^\sigma A$ is internal. If $A=\{a_1,\ldots,a_k\},$ a nonempty finite set$,$ 
then $\hyper A = \{\hyper a_1,\ldots,\hyper a_k\} = {^\sigma A}$ is internal. 
Conversely$,$ assume that $A$ is infinite. Then there some $B \subset A$ and a 
bijection $f\colon B \to \nat.$ Then for each $b \in B$ and there exists 
$n \in \nat$ such that $(\hyper b, \hyper n) \in \hyper f.$ Further$,$ for each 
$n \in \nat$ there exists some $b \in B$ such that $(\hyper b,\hyper n) \in 
\hyper f$ from Theorem 3.1.3 part (v). Hence$,$ $\hyper f[{^\sigma B}] = 
{^\sigma N}.$ If $^\sigma A$ is internal$,$ then ${^\sigma A}\cap \hyper B = 
{^\sigma B}$ is internal. But$,$ the image of an internal set under an internal 
map is internal. This contradicts the result that $^\sigma N$ is external.\par
\vskip 18pt
The following is the very important {\bf {Cauchy Principle}} for normed 
linear spaces.\par
{\bf Theorem 10.1.1. (Cauchy Principle)} {\sl Let $P(x)$ be a bounded formal first-order 
expression in one fee variable and employing internal constants (i.e. a 
bounded formal internal property as discussed in Appendix 4$,$ within the proof 
of Theorem 4.2.2.) If $P(\eps)$ holds for each 
$\eps \in \nmonad {\vec 0},$ then there 
exists an open ball $B$ about $\vec 0$ such that $P(\vec x)$ holds for each 
$\vec x \in \hyper B.$}\par
Proof. (In the following proof$,$ as an illustration$,$ we retain the 
$\sigma$ notation rather than assume the identification. This will show how the 
embedding technically allows one to go back-and-forth between the standard 
model and the $\sigma$ model$,$ in many cases. Some modern texts in this subject 
do not utilize our identification process for objects in $\Re.$) 
Let $\Phi(z) = z \in (\hypernat - \{\Hyper 0\}) \land \forall y(y \in 
\hyperrealp n \land \Vert y \Vert < 1/z \to P(y)).$ Then $\Phi(z)$ is a 
bounded internal first-order expression and from the hypothesis $\Phi(z)$ 
holds for each $n \in \nat_\infty.$ Let $A= \{z\bigm \vert z \in \hypernat \land
\neg \Phi(z)\}.$ From the internal definition theorem$,$ $A$ is an internal 
subset of $\hypernat$ and $A \subset {^\sigma \nat}.$ Thus $A = {^\sigma C},$ 
where $C \subset \nat.$ From proposition 10.1$,$ $C$ is finite. Thus $C$ has a 
greatest member$,$ say $k.$ Hence $\Phi(z)$ holds for each $n \in \hypernat$ such 
that $n \geq {\hyper k} + {\Hyper 1}.$ Now simply consider the standard ball
$B$ about $\vec 0$ with radius $1/(k+1).$ Then $\hyper B$ is but the ball in 
$\hyperrealp n$ about $\Hyper {\vec 0}  $ with radius $1/(k + 
1) = 1/(\hyper k + \Hyper 1).$\par
\vskip 18pt
{\bf Theorem 10.1.2.} {\sl Two internal maps $f$ and $g$ are equivalent if and only if
there is a local map $\alpha$ such that for each nonzero $\vec x \in \nmonad {\vec 0\,}$
$$f(\vec x) = g(\vec x) + (\alpha (\vec x))\Vert \vec x \Vert.$$}\par
Proof. Define 
$$\alpha (\vec x) = {{f(\vec x) - g(\vec x)}\over{\Vert \vec x \Vert }},$$ 
for nonzero $\vec x \in \nmonad {\vec 0},$ and $\alpha (\vec 0) = \vec 0.$
Let $f \sim g.$ Then for each nonzero $\vec x \in \nmonad {\vec 0}$ we have 
that $(f(\vec x) - g(\vec x))/\Vert \vec x \Vert \in \mmonad {\vec 
0}\Rightarrow \alpha$ is a local map. \par
Conversely$,$ assume that for nonzero $\vec x \in \nmonad {\vec 0}$ that 
$f(\vec x) = g(\vec x) + (\alpha (\vec x))\Vert \vec x \Vert,$ where 
$\alpha$ is a local map. The result follows from the definition of a local 
map.\par
\vskip 18pt
If the maps $f$ and $g$ of Theorem 10.1.2 are locally linear maps$,$ then 
the equation in the conclusion of Theorem 10.1.2 holds for $\vec x = \vec 0$ 
as well.\par
\vskip 18pt
{\bf Theorem 10.1.3.} {\sl Suppose that $f$ is any locally linear map. Then 
there exists a unique internal linear map $T_f\colon \hyperrealp n \to 
\hyperrealp m$ such that $\Vert T \Vert \in {\cal O}$ and there exists an open 
$E \subset \realp n,$ where $\vec 0 \in E,$ such that $f(\vec x) = T_f(\vec 
x)$ for each $\vec x \in \hyper E.$}\par
Proof. First$,$ we define an internal function by means of the Internal 
Definition Theorem. Let fixed $\Lambda \in \nat_\infty$ and $\eps = 
1/\Lambda.$ 
Notice that for nonzero $\vec x \in \hyperrealp n,\ \eps (\vec x/\Vert \vec 
x\Vert) \in \nmonad {\vec 0}.$ Next$,$ consider
$$T_f = \{(x,y)\bigm \vert x \in \hyperrealp n \land x \not= \vec 0 \land
y \in \hyperrealp m \land $$ $$y = \Lambda \Vert x\Vert f\left(\eps\left({x\over{
\Vert x\Vert}}\right)\right)\} \cup (\vec 0,\vec 0).$$ \par
We know that for any $\vec x \in \hyperrealp n$ there is a nonzero $\alpha
\in \monad 0 $ such that $\alpha \vec x \in \nmonad {\vec 0}.$ Consider any
$\vec x \in \hyperrealp n$ and an $\alpha \in \monad 0$  such that
$\alpha \vec x \in \nmonad {\vec 0}.$ Then
$$T_f(\vec x) = \alpha\Vert \vec x \Vert \left({{f\left(\eps {\vec x \over \Vert 
\vec x \Vert}\right)\over{\alpha\eps}}}\right) = $$
$${{f\left({\alpha\eps\Vert \vec x \Vert \vec x}\over{\Vert \vec 
x\Vert}\right)}\over{\alpha \eps}}= {f(\alpha\eps\vec x)\over {\alpha\eps}}=
{{f(\eps\alpha \vec x)}\over{\alpha \eps}}= {{\eps f(\alpha \vec x)}\over 
{\eps\alpha}} = {{f(\alpha \vec x)}\over \alpha}.\eqno (*)$$
Obviously$,$ if $\vec x \in \nmonad {\vec 0},$ then taking any nonzero
$\alpha \in \monad 0$ we have that 
$$T_f(\vec x) = {{f(\alpha \vec x)}\over{\alpha}} = {{\alpha f(\vec 
x)}\over{\alpha}} = f(\vec x).$$\par
Consider the bounded internal first-order expression
$P(x) = x\in \hyperrealp n \land T_f(x) = f(x).$ Since $P(x)$ holds for all 
$\vec x \in \nmonad {\vec 0}$ then the Cauchy Principle implies that there 
exists an open ball $E$ about $\vec 0$ such that for each $\vec x \in \hyper 
E,\ T_f(\vec x) = f(\vec x).$\par
To establish that $T_f$ is a linear transformation$,$ let $\vec x,\ \vec y
\in \hyperrealp n$ and $\alpha \in \monad 0$ be such that $\alpha \vec x,\ 
\alpha \vec y \in \nmonad {\vec 0}.$ Then using (*) and the locally linear 
property it follows that
$$T_f(\vec x+\vec y) = {1\over \alpha}(f(\alpha \vec x + \alpha y)) = 
{1\over \alpha}(f(\alpha \vec x) + f(\alpha \vec y)) = T_f(\vec x) + T_f(\vec 
y).$$
Now let $r \in \hyperreal.$ Then there is a nonzero $\alpha \in \monad 0$ such 
that $\alpha \vec x,\ r\alpha \vec x \in \nmonad {\vec 0}.$ Consequently$,$ 
$$T_f(r\vec x) = {1\over \alpha}f(r\alpha \vec x) = {r\over\alpha}f(\alpha 
\vec x) = rT_f(\vec x).$$\par
For uniqueness$,$ let $f(\vec x) = T(\vec x) = G(\vec x)$ for each 
$\vec x \in \nmonad {\vec 0},$ where $T,\ G$ are internal linear transformations from 
$\hyperrealp n$ into $\hyperrealp m.$ [Note: We$,$ of course$,$ mean that they 
are linear transformations over the vector space $\hyperrealp n$ into 
$\hyperrealp m$ where the field in question is $\hyperreal.$]
Then by linearity for any $\vec x \in \hyperrealp n$
$$T(\vec x) = \Lambda \Vert \vec x\Vert T\left(\eps\left({\vec x\over{\Vert 
\vec x\Vert}}
\right)\right)= f(\vec x) = \Lambda \Vert \vec x\Vert G\left(\eps\left({\vec x
\over{\Vert \vec x\Vert}}\right)\right)= G(\vec x).$$\par
Lastly$,$ we show that $\Vert T_f \Vert$ is limited. Recall that if $T$ is a 
standard linear transformation$,$ then $\Vert T \Vert = \sup {\{T(\vec x) \bigm 
\vert \vec x \in \realp n \land \Vert \vec x \Vert = 1\}}.$ 
The $\Hyper\sup$ operator extends$,$ in general$,$ to internal subsets of 
$\hyperreal$ and$,$ indeed$,$ is the same operator as a $\sup$ defined on some 
subsets of the ordered field $\hyperreal$ even though it may not exist for 
certain bounded sets such as $\monad 0.$ Assume that $\Vert T_f\Vert$ is not 
limited. Then by *- transfer of the standard case$,$ there is some 
$\vec x \in \hyperrealp n,\  \Vert \vec x\vert = 1,$ such that 
$\Vert T_f(\vec x) \Vert$ is infinite since the set $\{T_f(\vec x) \bigm 
\vert \vec x \in \realp n \land \Vert \vec x \Vert = 1\}$ is internal and 
not *-bounded. Consequently$,$ $\vec x/\Vert T_f(\vec x)\Vert \in \nmonad {\vec 
0}.$ From (*) one obtains
$$T_f(\vec x) = \Vert T_f(\vec x)\Vert f\left({\vec x \over {\Vert T_f(\vec x)
\Vert}}\right),$$                                                                         
which implies that
$$\Vert T_f(\vec x) \Vert= \Vert T_f(\vec x)\Vert \left\Vert f\left({\vec x 
\over {\Vert T_f(\vec x)\Vert }}\right)\right \Vert.$$    
Hence$,$ $\Vert f(\vec x/\Vert \vec x\Vert)\Vert = 1 \notin \monad 0$ implies 
the contradiction that $f$ is not a local map and this completes the proof. \par
                                                                     
\vskip 18pt
{\bf Theorem 10.1.4. } {\sl Let $f,\ g$ be locally linear maps. Then 
$f\sim g$ if and only if $T_f(\vec x) \approx T_g(\vec x)$ for each $\vec x \in 
\hyperrealp n$ such that $\Vert \vec x \Vert = 1.$}\par
Proof. First$,$ assume that $T_f(\vec x) \approx T_g(\vec x)$ for each $\vec x \in 
\hyperrealp n$ such that $\Vert \vec x \Vert = 1.$ Consider any $\vec y \in 
\nmonad {\vec 0}.$ Then from Theorem 10.1.3$,$ it follows that
$${{f(\vec y) - g(\vec y)}\over {\Vert \vec y \Vert}} = {{T_f(\vec y) - 
T_g(\vec y)}\over {\Vert \vec y \Vert}}= T_f\left({{\vec y}\over{\Vert \vec y 
\Vert}}\right) - T_g\left({{\vec y}\over{\Vert \vec y \Vert}}\right)\in 
\mmonad {\vec 0}.$$
Thus $f \sim g.$\par
Conversely$,$ assume that $f \sim g,$  $\Vert \vec x \Vert = 1$and that 
positive $\alpha \in \monad 0.$
Then from theorem 10.1.3$,$ 
$$T_f(\vec x) - T_g(\vec x) = {1\over \alpha}(T_f(\alpha\vec x)-T_g(\alpha\vec 
x))= {{f(\alpha \vec x) - g(\alpha \vec x)}\over{\Vert \alpha \vec x 
\Vert}}\in \mmonad {\vec 0}.$$\par
\vskip 18pt
{\bf Proposition 10.2.} {\sl Let $F,G\colon {\cal V} \to {\cal W}$ be two linear 
transformations and ${\cal V},\  {\cal W}$ be (nontrivial) linear 
spaces over the real or complex fields$,$ where ${\cal V}$ is normed. If there exists an 
open ball $B \subset {\cal V},$ about $\vec 0,$ such that $F(\vec x) = G(\vec x)$ for each $\vec x \in 
V,$ then $F=G.$}\par
Proof. Let nonzero $\vec x \in \cal V.$ and positive $r \in \real$ be the radius 
of the ball $B.$ Then $(r\vec x)/2\Vert \vec x \Vert \in B$ and 
$T((r\vec x)/(2\Vert \vec x \Vert)) =(r/(2\Vert \vec x \Vert))T(\vec x) = 
G((r\vec x)/(2\Vert \vec x \Vert)) =(r/(2\Vert \vec x \Vert))G(\vec x)\Rightarrow
T(\vec x) = G(\vec x).$ Obviously$,$ $\vec 0 = T(\vec 0) = G(\vec 0)$ and the 
proof is complete.     
\par 
\vskip 18pt
{\bf Theorem 10.1.5.} {\sl If $f$ is a differential$,$ then there exists a 
unique bounded linear transformation $F\colon \realp n \to \realp m$ and an open 
set $E\subset \realp n$ such that $\vec 0 \in E$ and $\hyper F(\vec x) = f(\vec 
x)$ for all $\vec x \in \hyper E$ and$,$ in particular$,$ $F(\vec x) = f(\vec x)$ 
for each $\vec x \in E.$}\par
Proof. Since $f$ is a differential then there exists a bounded linear 
transformation $F\colon \realp n \to \realp m$ and a unique internal linear 
transformation $T_f\colon \hyperrealp n \to \hyperrealp m$ such that 
$\hyper F(\vec x) = f(\vec x)= T_f(\vec x)$ for each $\vec x \in \nmonad 
{\vec 0}.$  Consider the bounded internal
first-order expression $P(x) = x \in \hyperrealp n \land \hyper F(x) = f(x)= 
T_f(x).$
Then since $P(x)$ holds for all $\vec x \in \nmonad {\vec 0}$ then the Cauchy 
Principle implies that there exists some open ball $E$ about $\vec 0$ such 
that for each $\vec x \in \hyper E,$ $\hyper F(\vec x) = f(\vec x)= T_f(\vec 
x).$ Hence$,$ in particular$,$ $\hyper F(\vec x) = F(\vec x) = f(\vec x)$ for each 
$\vec x \in E.$ The fact that $F$ is unique comes from the *-transfer of 
Proposition 10.2. 
\par
\vskip 18pt     
{\bf Theorem 10.1.6.} {\sl If $f$ and $g$ are differentials  and $f \sim g,$ 
then there exists some open $E \in \realp n$ such that $\vec 0 \in E$ and 
$f(\vec x) = g(\vec x)$ for each $\vec x \in \hyper E.$}\par
Proof. There are two bounded standard transformations $T,G \colon \realp n\to \realp 
m$ such that $f(\vec x) = \Hyper T(\vec x)$ and $g(\vec x) = \Hyper G(\vec x)$ 
for each $\vec x \in \nmonad {\vec 0}.$ By the uniqueness property of Theorem 
10.1.3 and Theorem 10.1.4 $\Hyper T(\vec y) \approx \Hyper G(\vec y)$ for each 
$\vec y \in \hyperrealp n$ such that $\Vert \vec y \Vert = 1.$ Considering any 
$\vec z \in \realp n$ such that $\Vert \vec z\Vert = 1$ and taking the 
standard part operator this yields $T(\vec z) = G(\vec z).$ Taking any
nonzero $\vec x \in \realp n,$ we have that
$$T(\vec x) = \Vert \vec x \Vert T\left({{\vec x}\over{\Vert \vec x 
\Vert}}\right)=\Vert \vec x \Vert G\left({{\vec x}\over{\Vert \vec x 
\Vert}}\right)= G(\vec x).$$
Thus $f(\vec x) = g(\vec x)$ for each $\vec x \in \nmonad {\vec 0}.$ Once 
again application of the Cauchy Principle yields the result.\par
 
\vskip 18pt
{\bf Theorem 10.1.7.} {\sl Let nonempty open 
$G \subset \realp n,$ and $\vec c \in G.$ A function $f\colon G \to \realp m$ 
is continuous at $c \in G$ if and only if the local increment map 
$\Delta \hyper f_{\vec c}$ is a local map.}\par
Proof. This follows immediately from Definition 4.4.1 and the fact that if 
$G$ is open and $\vec c \in G,$ then $\nmonad {\vec c} \subset \Hyper G.$\par
\vskip 18pt
{\bf Theorem 10.1.8.} {\sl Let noninfinitesimal $\vec a \in {\cal O}^n$ and 
suppose that $\vec a - \vec b \in \nmonad {\vec 0}.$ Then $\vec a$ is almost 
parallel to $\vec b.$}\par
Proof. Since $\vec a \in {\cal O}^n$ and $\vec 0 \not= \vec a$ then there 
exists some nonzero $\vec c \in \realp n$ such that $\vec a \in \nmonad {\vec c}.$
Consequently$,$ $\vec b \in \nmonad {\vec c}.$ Further$,$ $\st {\Vert \vec 
a\Vert},\ \st {\Vert \vec b\Vert} \not= 0$ and $\vert \ \Vert \vec a \Vert  - 
\Vert \vec b \Vert \ \vert \leq \Vert \vec c - \vec b \Vert \Rightarrow
\st {\Vert \vec a \Vert)} = \st {\Vert \vec b \Vert}.$ Hence,
$${{\st {\vec a}}\over {\st {\Vert \vec a \Vert}}} = \St {\left({{\vec a}\over
{\Vert \vec a \Vert}}\right)} = {{\st {\vec b}}\over {\st {\Vert \vec b 
\Vert}}} 
= \St {\left({{\vec b}\over
{\Vert \vec b \Vert}}\right)}$$
and the proof is complete.                              
\vskip 18pt
{\bf Theorem 10.4.1.} {\sl Let $f\colon G \to \real,$ where 
nonempty open $G\subset \realp n$ and standard $\vec v \in G.$ 
Suppose that $y = x_i,\ 1\leq i\leq n,\ f_{y}$ is defined on $G$ 
and continuous at $\vec v.$ 
Let $C_{y} \subset \Hyper G,$  where internal $C_{y}$ is 
*-convex in the direction $y.$ If for $h \in \monad 0,$ such that 
$\vec p = (p_1,\ldots,
p_i,\ldots,p_n),\vec q = (p_1,\ldots, p_i + h,\ldots,p_n)\in C_{y}$ and 
$\vec v \approx \vec p,$ then 
there exists $\eps \in \monad 0$ such that
$$\hyper f({\vec q}\,) - \hyper f({\vec p}\,) = f_{y}({\vec v}\,)\, h + h\, \eps.$$} \par
Proof. Recall that for the hypothesized behavior of $f$ and $f_{y}$
 the mean value theorem for the partial derivative $f_{y}$ 
states that for any set $D_{y} \subset G$ convex in the direction ${y}$ if
$\vec u =  (u_1,\ldots,
u_i,\ldots,u_n),\, \vec w = (u_1,\ldots, u_i + r,\ldots,u_n)\in D_{y},\ r \in
\real $ then 
there exists $\vec s$ in the line segment with end points $\vec u,\ \vec w$
such that $f(\vec w) - f(\vec u) = r\,f_y({\vec s}\,).$ Thus $\vec s \in D_{y}$ 
and $\Vert \vec u - {\vec s}\, 
\Vert \leq \Vert \vec u - \vec w\Vert = \vert \vec r\,\vert.$ By *-transfer
it follows that there exists some $\vec t \in C_{y}$ such that 
$\vec t \in C_{y}$ and $\Vert \vec p - \vec t\, \Vert \leq \vert h \vert$ and
$\hyper f(\vec q) - f(\vec p) = h\,\hyper f_y(\vec t\,).$
But$,$ $\vec t\approx \vec p \approx \vec v$ and the continuity of $f_{y}$ at 
$\vec v$ imply that $\hyper f_{y}({\vec t}\,) \approx f_{y}({\vec v}\,).$ Thus there 
exists some $\eps \in \monad 0$ such that $\hyper f_{y}({\vec t}\,)= f_{y}({\vec 
v}\,) 
+ \eps.$ The result follows by substitution.\par
\vskip 18pt
\centerline{\bf NOTES}
\vskip 18pt
\noindent [1] {\bf Theorem N.1.} {\sl Two unit vectors$,$ $\vec v,\ \vec u \in \hyperrealp n$ 
have the property that $\vec v \approx \pm \vec u$ if and only if $\vec v \bullet 
\vec u \approx \pm 1.$}\par
Proof.  First note that $\vec v,\, \vec u \in {\cal O}^n.$ 
For the necessity$,$ assume that $\vec v \approx \pm \vec u.$ Then $\vec v \pm 
\vec u \in \nmonad {\vec 0\,} \Rightarrow \vec v \bullet \vec v \pm 
\vec u \bullet \vec v = 1 \pm \vec u \bullet \vec v\in \monad 0 \Rightarrow
\pm \vec u \bullet \vec v \approx 1 \Rightarrow  \vec u \bullet \vec v \approx 
\pm 1.$\par
For the converse$,$ assume that $\vec v \bullet \pm\vec u \approx 1 $ but that
$\vec v \not\approx \pm \vec u.$ Thus there exist $\vec 0 \not= \vec r_{\pm} 
\in \realp n$ 
and $\eps_{\pm} \in \nmonad {\vec 0\,}$ such that $\vec v = \pm \vec u +  \vec r_{\pm}
+  \eps_{\pm}.$ Hence$,$ $\vec v \bullet \pm\vec u = \pm \vec u \bullet \pm\vec u
+ \vec r_{\pm} \bullet \pm\vec u + \eps_{\pm}\bullet \pm\vec u = 1  + 
\vec r_{\pm} \bullet \pm\vec u + \eps_{\pm}\bullet \pm\vec u \approx 1  + 
\vec r_{\pm} \bullet \pm\vec u.$\par
 We now show that  $\vec r_{\pm} \bullet\pm \vec u
\not\approx 0.$ We know that $\vec r_{\pm}\bullet\pm \vec u = (1/4)(\Vert
\pm\vec u + \vec r_{\pm}\Vert^2 + \Vert \pm\vec u - \vec r_{\pm}\Vert^2).$ Assume that                                                                           
$\vec r_{\pm} \bullet \pm\vec u \approx 0\Rightarrow
\Vert \pm\vec u + \vec r_{\pm}\Vert^2 \approx \Vert \pm\vec u - \vec r_{\pm}
\Vert^2.$ Now $\vec v = \pm \vec u +  \vec r_{\pm} +  \eps_{\pm}\Rightarrow
\Vert \vec v \Vert^2 = 1 \approx \Vert \pm \vec u +  \vec r_{\pm} 
\Vert^2\Rightarrow 2 \approx \Vert
\pm\vec u + \vec r_{\pm}\Vert^2 + \Vert \pm\vec u - \vec r_{\pm}\Vert^2 =
2\Vert\pm\vec u \Vert^2 + 2\Vert \vec r_{\pm}\Vert^2 = 2 + 2\Vert \vec r_{\pm}\Vert^2 
\Rightarrow 0 \approx \Vert \vec r_{\pm}\Vert;$ a contradiction. Thus 
$\vec r_{\pm} \bullet\pm \vec u \not\approx 0 \Rightarrow 
\vec v \bullet \pm\vec u \not\approx 1.$ The result follows from this 
contradiction.\par 
\vskip 18pt
\noindent [2] Let $c\colon [a,b] \to \realp n$ be a continuous differentiable curve with 
graph $C$ and $c^\prime(t) \not= \vec 0$ for each $t \in [a,b].$ 
Assume that uniformly continuous $F\colon E \to \realp n,$ open $E \supset 
C.$ Then $\hyper F(\ell_j(t^\prime))\bullet \hyper {\vec v_j} = 
\hyper F(\hyper c(t^\prime)) \bullet \hyper {\vec v_j} + \eps_j \Vert \hyper 
{\vec v_j}\Vert, \eps_j \in \monad 0,$ and $\sum_{j=1}^\Gamma
\hyper F(\ell_j(t^\prime))\bullet \hyper {\vec v_j}\approx \sum_{j=1}^\Gamma
\hyper F(\hyper c(t^\prime)) \bullet \hyper {\vec v_j},$ where the symbols 
have the same meaning as in the derivation for Application 8.2.1.\par
Proof. For an interior point of $t\in (a,b)$ we use the concept of uniform 
differentiability and the result [{\tt Stroyan and Luxemburg [1976:94-97]}] that for 
any $t_{j-1},\ t_j \approx t$ that there exists some $\vec {\eps}_j \in
\nmonad {\vec 0}$ such that $\vec v_j =\hyper c(t_j) - \hyper c(t_{j-1}) = 
c^\prime(t)(t_j - t_{j-1}) + \vert t_j - t_{j-1} \vert \vec {\eps}_j.$
Assuming that $\hyper c(t_j) - \hyper c(t_{j-1}) = \vec 0 \Rightarrow
c^\prime(t)(t_j - t_{j-1}) + \vert t_j - t_{j-1} \vert \vec {\eps}_j =
\vec 0 \Rightarrow c^\prime(t)((t_j - t_{j-1})/\vert t_j - t_{j-1} \vert) +  
\vec {\eps}_j = \vec 0 \Rightarrow \Vert c^\prime(t)((t_j - t_{j-1})/
\vert t_j - t_{j-1} \vert) 
\Vert = \Vert c^\prime(t)\Vert = \Vert -\vec {\eps_j}\Vert \Rightarrow
\st {\Vert c^\prime(t)\Vert} = \Vert c^\prime(t)\Vert = \st {\Vert -\vec 
{\eps_j}\Vert} = 0;$ a contradiction. Thus $\vec v_j=\hyper c(t_j) - \hyper c(t_
{j-1}) 
\not = \vec 0.$ Now the same proof that appears in 
{\tt Stroyan and Luxemburg [1976:94-97]} shows that this result also holds if $t = 
a$ or $ t = b.$ Next observe that $\ell_j(t^\prime) \approx \ell_j(t_j) = 
\hyper c(t_j) \approx \hyper c(t^\prime)\approx c(t),$ where $t \in [a,b].$
Thus $\ell_j(t^\prime), \, \hyper c(t^\prime) \in \nmonad {c(t)} \subset 
\Hyper E$ and the uniform continuity of $\hyper F\Rightarrow
\hyper F(\ell_j(t^\prime)) = \hyper F(\hyper c(t^\prime)) + \vec {\lambda}_j.$ 
Hence,
$\hyper F(\ell_j(t^\prime)\bullet(\vec v_j/\Vert \vec v_j \Vert) = 
\hyper F(\hyper c(t^\prime))\bullet(\vec v_j/\Vert \vec v_j \Vert)  +
\vec {\lambda}_j\bullet (\vec v_j/\Vert \vec v_j \Vert)=    
\hyper F(\hyper c(t^\prime))\bullet(\vec v_j/\Vert \vec v_j \Vert) + \eps_j$ 
and the first result follows.\par
For the second result follow the usual method and notice that there 
exists infinitesimal $\delta = \max \{\eps_j\}$ and that 
$\vert \sum _{j=1}^\Gamma\eps_j \Vert \vec v_j \Vert\,\vert \leq \vert\delta\vert 
\sum_{j=1}^\Gamma \Vert \vec v_j \Vert \Rightarrow  \sum _{j=1}^\Gamma\eps_j 
\Vert \vec v_j \Vert \in \monad 0$ since $\st {\sum_{j=1}^\Gamma \Vert \vec 
v_j \Vert }=$ length of the curve $c.$\par
\vfil
\eject
\centerline{REFERENCES}
\noindent {\bf Aerts$,$ D.} [1984]$,$ The missing elements of reality in the description of 
quantum mechanics of the E.P.R. paradox situation$,$ {\sl Helvetica Physica,}
{\bf 57}: 421-428.\hfil\break 
{\bf Apostal$,$ Tom M.} [1957]$,$ {\sl Mathematical Analysis,} Addison-Wesely$,$ 
Reading$,$ MA.\hfil\break
{\bf Barwise$,$ Jon (ed.)} [1977]$,$ {\sl Handbook of Mathematical Logic,} 
North-Holland$,$ Amsterdam.\hfil\break
{\bf Berkeley} [1734]$,$ {\sl The Analyst,} London.\hfil\break
{\bf Cesari} [1956]$,$ {\sl Surface Area}$,$ Princeton University Press$,$ 
Princeton$,$ NJ. \hfil\break
{\bf Cutland N. J.} [1986]$,$ Private communication.\hfil\break
{\bf De Lillo$,$ Nicholas J.} [1982]$,$ {\sl Advanced Calculus with Applications,}
Macmillan$,$ New York$,$ NY. \hfil\break  
{\bf Gauss$,$ Karl F.} [1827]$,$ {\sl General Investigation of Curved Surfaces},
Raven Press$,$ Hewlett$,$ NY. \hfil\break
{\bf Henstock$,$ R.} [1961]$,$ Definitions of the Riemann type of variational 
integral$,$ {\sl Proc. London Math. Soc.,} {\bf (3) 11}: 402-418.\hfil\break
{\bf Herrmann$,$ Robert A.} [1976]$,$ The Q-topology$,$ Whyburn type filters and the 
cluster set map$,$ {\sl Proc. Amer. Math. Soc.,} {\bf 59}: 161-166.\hfil\break 
{\bf Herrmann$,$ Robert A.} [1978]$,$ (2003) {\sl Nonstandard Analysis - A Simplified 
Approach,} \hfil\break
http://www.arxiv.org/abs/math.GM/0310351\hfil\break
{\bf Herrmann$,$ Robert A.} [1980]$,$ A nonstandard approach to pseudotopological 
compactifications$,$ {\sl Z. Math. Logik Grundlagen Math.,} {\bf 26}: 361-384. 
\hfil\break
{\bf Herrmann$,$ Robert A.} [1985]$,$ Supernear functions$,$ {\sl Math. Japanica,}
{\bf 30}: 169-185. \hfil\break
{\bf Herrmann$,$ Robert A.} [1986]$,$ [1987]$,$ (1993) {\sl The Theory of Ultralogics} \hfil\break
http://www.arxiv.org/abs/math.GM/9903081\hfil\break
http://www.arxiv.org/abs/math.GM/9903082\hfil\break
\noindent {\bf Herrmann$,$ Robert A.} [1989]$,$ Fractals and ultrasmooth microeffects$,$ 
{\sl J. Math. Phys.}$,$ {\bf 30}(4)$,$ April 1989: 805-808. \hfil\break
{\bf Herrmann, Robert A.} [1995], {\sl Nonstandard Analysis Applied to Special and General Relativity - The Theory of Infinitesimal Light-Clocks} http://arxiv.org/abs/math/0312189 \hfil\break
{\bf Hurd$,$ A.E. and P.A. Loeb} [1985]$,$ {\sl An Introduction to Nonstandard 
Real Analysis,} Academic Press$,$ Orlando$,$ FL.\hfil\break
{\bf Keisler$,$ H. Jerome} [1986]$,$ {\sl Elementary Calculus - An Infinitesimal 
Approach,} (Second edition)$,$ Prindle$,$ Weber \& Schmidt$,$ Boston. \hfil\break
{\bf Jarnik$,$ J.$,$ J. Kurzweil and S. Schwabik} [1983]$,$ On Mawhin's approach to 
multiple nonabsolutely convergent integrals$,$ {\sl \v Casopis P\v est. Mat.,} 
{\bf 108}: 157-167.\hfil\break
{\bf Jeck$,$ Thomas J.} [1971]$,$ {\sl 
Lectures Notes in St Theory, Lecture Notes in Mathematics \#217,}
Macmillan$,$ Springer-Verlag$,$ New York$,$ NY. \hfil\break 
{\bf Leibniz G. W.} [1701]$,$ {\sl M\'emoire de M.G.G Leibniz touchant son sentiment 
sur le calcul diff\'erentiel$,$ Journal de Tr\'evoux}$,$ Mathematische Schriften$,$ 
ed. C.I. Gerhardt$,$ Vol. 5 (1858).\hfil\break
{\bf Lorentz$,$ H. A.} [1915]$,$ {\sl The Theory of Electrons }$,$ (Dover$,$ New 
York$,$ NY$,$ 1952.)\hfil\break
{\bf Luxemburg$,$ W.A.J.} [1962]$,$ {\sl Non-Standard Analysis - Lectures on A. 
Robinson's Theory of Infinitesimals and Infinitely Large Numbers,} Math. 
Dept.$,$ California Institute of Technology$,$ Pasadena$,$ CA.\hfil\break
{\bf Luxemburg$,$ W.A.J.} [1973] What is nonstandard analysis? in {\sl Papers 
in the Foundations of Mathematics$,$ No. 13 Slaught Memorial Papers,} Amer. 
Math. Monthly {\bf 80}: 38-67.\hfil\break
{\bf Machove$,$ M. and J. Hirschfeld} [1969]$,$ {\sl Lectures on Non-Standard 
Analysis,} Lecture Notes in Mathematics V.94$,$ Springer-Verlag$,$ New York$,$ NY.\hfil\break
{\bf Mawhin$,$ J.} [1986]$,$ Nonstandard analysis and generalized Riemann 
integrals$,$ {\sl \v Casopis P\v est. Mat.,} {\bf 111}:34-47. 
\hfil\break
\noindent {\bf Maxwell James Clark} [1890] {\sl The Scientific Papers of James Clark 
Maxwell}$,$ Cambridge University Press$,$ Cambridge (Dover$,$ New York$,$ 
NY$,$ 1965.)\hfil\break
{\bf McShane$,$ E. J.} [1973]$,$ A unified theory of integration$,$ {\sl Amer. 
Math. Monthly,} {\bf 80}: 349-359.\hfil\break
{\bf Morley$,$ Arthur} [1942]$,$ {\sl Mechanics for Engineers}$,$ Longmans$,$ Green 
and Co. London.\hfil\break
{\bf Newton$,$ Isaac} [Summer 1665]$,$ The calculus becomes an algorithm,
in {\sl The Mathematical Papers of Isaac Newton}$,$ (ed. D.T Whiteside$,$ V. I$,$ 
 Cambridge University Press$,$ New York$,$ NY$,$ 1967): 298-368.\hfil\break
{\bf Newton$,$ Isaac} [Oct. 1665 - May 1666]$,$ The general problem of tangents$,$ 
curvature and limit-motion analysed by the method of fluxions$,$ in {\sl The 
Mathematical Papers of Isaac Newton}$,$ (ed. D.T Whiteside$,$ V. I$,$ 
Cambridge University Press$,$ New York$,$ NY$,$ 1967): 369-399. \hfil\break
{\bf Newton$,$ Isaac} [Oct. 1666]$,$ The October 1666 tract on fluxions,
in {\sl The Mathematical Papers of Isaac Newton}$,$ (ed. D.T Whiteside$,$ V. I$,$ 
Cambridge University Press$,$ New York$,$ NY$,$ 1967): 400-448. \hfil\break
{\bf Newton$,$ Isaac} [Winter 1670 - 1671]$,$ The tract `De methodis serierum et 
fluxionum'$,$ in {\sl The Mathematical Papers of Isaac Newton}$,$ (ed. 
D.T Whiteside$,$ V. III$,$ 
Cambridge University Press$,$ New York$,$ NY$,$ 1969): 32-353. \hfil\break
{\bf Newton$,$ Isaac} [1686]$,$ {\sl Mathematical Principles of Natural 
Philosophy,} (Revised translation by Florian Cajori$,$ University of California 
Press$,$ Berkeley$,$ CA$,$ 1934). \hfil\break
{\bf Robinson$,$ Abraham} [1961]$,$ Non-standard analysis$,$ {\sl Nederl. Akad. 
Wetensch. Proc. Ser. A 64$,$ and Indag. Math.} {\bf 23}: 432-440.\hfil\break
{\bf Robinson$,$ Abraham} [1966]$,$ {\sl Non-standard Analysis}$,$ North-Holland$,$ 
Amsterdam. \hfil\break
\noindent
{\bf Robinson$,$ A. and E. Zakon} [1969]$,$ A set-theoretic characterization of 
enlargements$,$ in {\sl Applications of Model Theory to Algebra$,$ Analysis$,$ and 
Probability} (ed. W.A.J. Luxemburg$,$ Holt$,$ Rinehart  and Winston$,$ New York$,$ 
NY): 109-122.\hfil\break
{\bf Sears$,$ F. W. and M. W. Zemansky} [1952]$,$ {\sl College Physics,} 
Addison-Wesley$,$ Cambridge$,$ MA. \hfil \break
{\bf Simhony$,$ M.} [1987]$,$ {\sl The Electron-Positron Lattice Space,} Physics 
Section 5$,$ The Hebrew University$,$ Jerusalem. \hfil\break 
{\bf Spivak$,$ Michael} [1965]$,$ {\sl Calculus on Manifolds,} W. A. Benjamin$,$ 
New York$,$ NY.\hfil\break
{\bf Stroyan$,$ K.D. and W.A.J. Luxemburg} [1976]$,$ {\sl Introduction to the 
Theory of Infinitesimals,} Academic Press$,$ New York$,$ NY.\hfil\break
{\bf Struik$,$ D. J.} [1961]$,$ {\sl Differential Geometry,} Addison-Wesley,
Reading$,$ MA. \hfil\break
{\bf Suppes$,$ Pactrick} [1960]$,$ {\sl Axiomatic Set Theory,}
D. Von Nostrand$,$ New York$,$ NY. \hfil\break 
{\bf Swartz$,$ Charles and Brian S. Thomson} [1988]$,$ More on the Fundamental 
Theorem of Calculus$,$ {\sl Amer. Math. Monthly,} {\bf 95}: 644-648.\hfil\break
{\bf Synge$,$ J. L  and B. A. Griffith} [1959]$,$ {\sl Principles of Mechanics,}
McGraw-Hill$,$ NY.\hfil\break
\vfil\eject

 \centerline{\hskip -1.50in{\bf Special Symbols}}
\centerline{\hskip -1.50in{(Alphabetically listed by first symbol letter.)}}
\vskip 18pt
\+Symbol\dotfill&Name$,$ if any\dotfill&Page no.\cr
\+${\cal A}$\dotfill&\dotfill&50\cr
\+$BHF$\dotfill&Basic Hyperfinite& \cr
\+&Subsets of $\hypernat$\dotfill&30 \cr
\+${\cal C}(\Hyper {\cal H})$\dotfill&Internal Constants\dotfill&21 \cr
\+${\cal C}({\cal H})$\dotfill&Standard Constants\dotfill&21 \cr
\+${\cal C}_{PSR}$\dotfill&\dotfill&37 \cr
\+${\cal C}_{SR}$\dotfill&\dotfill&38 \cr
\+${\cal C}_{PSR}^\prime$\dotfill&\dotfill&123\cr
\+${\cal C}_J$\dotfill&A Set Of Jordan&&\cr
\+&Measurable subsets\dotfill&125\cr
\+$D\int_Rf(\vec x)\, dX$\dotfill&Darboux Integral\dotfill&125\cr
\+${\rm ext}(Q)$\dotfill&Exterior Subrectangles\dotfill&63, 132\cr
\+$F(B)$\dotfill&Set Of All Finite& \cr
\+&Subsets of $B$\dotfill&31 \cr
\+$\cal H$\dotfill&The Superstructure\dotfill&17 \cr
\+$\Hyper {\cal H}$\dotfill&Hyperstructure\dotfill&21 \cr
\+${\rm int}(A)$\dotfill& Interior Points in $A$\dotfill&50\cr
\+${\rm inn}(Q)$\dotfill&Inner Subrectangles\dotfill&63, 132\cr 
\+$J$\dotfill&Jordan-measurable set\dotfill&42 \cr
\+$L[S]$\dotfill&m-demensional& \cr
\+&Element\dotfill&47 \cr
\+$L(f)$\dotfill&\dotfill&120\cr
\+$L(f,P)$\dotfill&Lower Sum\dotfill&120\cr
\+$L(\realp n,{\cal B})$\dotfill&Set Of $L$-gauges\dotfill&141\cr
\+${\cal L}(\lambda)$\dotfill&$\delta$-fine Partitions\dotfill&142\cr
\+$M_1$\dotfill&Infinitesimals \dotfill&11&\cr
\+$\dmonad {\vec p}$\dotfill&Deleted Monad\dotfill&77\cr
\+$\monad 0$\dotfill&Infinitesimals \dotfill&11&\cr
\+$\monad r$\dotfill&Monad about $r$ \dotfill&12&\cr
\+$\monad {\vec v}=\nmonad {\vec v}$\dotfill&Monad about $\vec v$ 
\dotfill&13$,$ 88&\cr
\+${\monad 0}^+$\dotfill&Nonnegative&\cr
\+&Infinitesimals\dotfill&69\cr            
\+$({\rm mod}\,o(\cdot))$\dotfill&\dotfill&69\cr
\+$m(x_1\delta_1,\ldots,x_p\delta_p)$\dotfill&Magnification Operator\dotfill&72\cr
\+$\hat \delta$\dotfill&Microguage\dotfill&75\cr
\+$\nat$\dotfill&Natural Numbers \dotfill&10&\cr
\+${\nat}_\infty$\dotfill&Infinite Natural No.s\dotfill&19 \cr
\+$o$\dotfill&Infinitesimals \dotfill&11&\cr
\+$\cal O$\dotfill&Limited numbers \dotfill&11\cr
\+${\cal O}^n$\dotfill&Limited Vectors \dotfill&13&\cr
\+$o(\eps)$\dotfill&Order Ideal\dotfill&68\cr
\+$o^n(\eps_1,\ldots,\eps_n)$\dotfill&\dotfill&71\cr
\+${\cal P}_\Gamma $\dotfill&Internal Hyper-& \cr
\+&polygonal Curve\dotfill&33 \cr
\+${\cal P}(W)$\dotfill&Power Set\dotfill&17\cr
\+$P_i(R)$\dotfill&i'th Projection\dotfill&18 \cr
\+$P_i$\dotfill&Hyperfinite& \cr
\+&Partition of $[a_i,b_i]$\dotfill&35 \cr
\+$P = P_1 \times \cdots \times P_n$\dotfill&Simple Fine Partition\dotfill&36\cr
\+$P_\delta$\dotfill&Gauge Partition\dotfill&75\cr
\+$(Q_\delta,P_\delta)$\dotfill&$\delta$-fine partition\dotfill&75\cr
\+$(Q_{\hat \delta},P_{\hat \delta})$\dotfill&Micropartition\dotfill&75\cr
\+$Q$\dotfill&Intermediate Partition\dotfill&120\cr
\+${\cal Q}$\dotfill&Intermediate&&\cr
\+&Partition Map\dotfill&122\cr
\+$\hyperreal$\dotfill&Extended reals \dotfill&10&\cr
\+$\hyperreal$\dotfill&Hyperreals \dotfill&10&\cr
\+$\hyperreal$\dotfill&Star-reals \dotfill&10&\cr
\+$\real ^+$\dotfill&Positive Reals\dotfill&10& \cr
\+$\real$\dotfill&Real Numbers \dotfill&10&\cr
\+$\hyperreal - \cal O$\dotfill&Infinite Hyperreals \dotfill&11&\cr
\+${\hyperreal}_\infty$\dotfill&Infinite Hyperreals\dotfill&11&\cr
\+$\realp n$\dotfill&Euclidean n-space \dotfill&12&\cr
\+$\Re$\dotfill&\dotfill&19 \cr
\+${\real}_\infty^+$\dotfill&Positive Infinite& \cr
\+&Hyperreals\dotfill&26 \cr
\+$R$\dotfill&Rectangle\dotfill&35 \cr
\+$R_q$\dotfill&Subrectangle\dotfill&35 \cr
\+$R_S$\dotfill&Subrectangle\dotfill&37 \cr
\+$\st$\dotfill&Standard Part \dotfill&14&\cr
\+${\cal S}(P)$\dotfill&\dotfill&37 \cr
\+$S$\dotfill&Subrectangle\dotfill&48\cr
\+$S(f,P,Q)$\dotfill&Riemann Sum\dotfill&122\cr
\+$\cal U $\dotfill&Basic Universe\dotfill&17 \cr
\+$U(f)$\dotfill&\dotfill&120\cr
\+$U(f,P)$\dotfill& Upper Sum\dotfill&120\cr
\+$v(R_q)$\dotfill&Volume\dotfill&85 \cr
\+$X_{n+1} $\dotfill&\dotfill&17 \cr
\+$Z\int_R f(\vec x)\, \Hyper \Vert P \Vert$\dotfill&Z-integral\dotfill&121\cr
\bigskip
\centerline{\hskip -1.5in{(Non-alphabetical listing.)}}
\medskip
\+$\approx$\dotfill&Infinitely Close \dotfill&11\cr
\+$\Vert \cdot \Vert$\dotfill&Euclidean Norm \dotfill&11\cr
\+$\partial A$\dotfill& Boundard Points of $A$\dotfill&55\cr\vfil\eject
\+$\partial (Q)$\dotfill&Boundary Subrectangles\dotfill&\cr
\+&\dotfill&62\cr
\+$\Vert P \Vert$\dotfill&\dotfill&120\cr
\+$\Vert S \Vert$\dotfill&\dotfill&120\cr\vfil\eject
\centerline{}
\vskip 1.50in
\centerline{\bf Some Applications of Nonstandard Analysis}
\centerline{\bf to Advanced Undergraduate Mathematics}
\centerline{$\diamondsuit$ Very Elementary Physics$\diamondsuit$}
\vskip .5in
\centerline{\bf Robert A. Herrmann}
\vskip 4.75in
\centerline{A 1991 Instructional Development Project from the}
\centerline{Mathematics Department}
\centerline{United States Naval Academy}
\centerline{572C Holloway Road}
\centerline{Annapolis$,$ Maryland 21402-5002}
\eject
\ \ \par
\vskip 2.75in
\centerline{IMPORTANT NOTICE}\pars
\indent Since the writing of this book was financed entirely by a designated 
grant from the Federal Government that was specifically obtained for this 
sole purpose then a copyright for this specific book cannot be obtained by 
its author. 
Any portion of its contents can be copied and used without seeking permissions 
from the author. However$,$ when such copying or use is made of this material$,$ 
it is 
necessary that the author and the U. S. Naval Academy be indicated as 
the source of the  material being used. 
Further note that certain new results that appear in this book 
will be published under the author's name in scholarly journals. \par
\vfil\eject
\centerline{}
\centerline{\bf CONTENTS}
\bigskip
\indent Chapter 1\par
{\bf Introduction}\par
\line{\indent 1.1 \ \ Brief Comments \leaderfill 165}
\line{\indent 1.2 \ \ Manual Structure \leaderfill 165}
\bigskip
\indent Chapter 2\par
{\bf Mechanics}\par
\line{\indent 2.1 \ \ Instantaneous Velocity\leaderfill 167}
\line{\indent 2.2 \ \ Acceleration\leaderfill 171}
\line{\indent 2.3 \ \ Forces and Newton's Law\leaderfill 172}
\line{\indent 2.4 \ \ Vectors\leaderfill 174}
\line{\indent 2.5 \ \ Energy and Force Fields\leaderfill 175}
\line{\indent 2.6 \ \ General Impulse\leaderfill 176}
\bigskip
\indent Chapter 3\par
{\bf Slightly Less Basic Mechanics}\par
\line{\indent 3.1 \ \ Mass\leaderfill 181} 
\line{\indent 3.2 \ \ Moments and Center of Mass\leaderfill 183}
\line{\indent 3.3 \ \ Point Masses\leaderfill 185}
\line{\indent 3.4 \ \ Standard Rules\leaderfill 187}
\line{\indent \ \ \ \ References \leaderfill 189}
\line{\indent\ \ \ \ Additional Special symbols \leaderfill 190}\vskip 0.25in
\centerline{To be continued by members of the physics community.}
 
\vfil
\eject
\centerline{Chapter 1.}
\medskip
\centerline{\bf INTRODUCTION}\par
\bigskip
\leftline{1.1 \underbar{Brief Comments}}\par 
\medskip
Since the time of Archimedes the major applications of infinitesimal reasoning 
have been in the general discipline of geometry and what we now categorize as 
the subject matter of Physics. All of the applications that appear in 
the first volume in this series -- Infinitesimal Modeling -- are taken from these 
two disciplines. The methods employed within 
this manual are based exclusively upon those that appear in the our major  
reference the {\it Some Applications of Nonstandard Analysis to Undergraduate 
Mathematics -- Infinitesimal Modeling.}   
From time-to-time$,$ a portion of 
certain applications$,$ discussions and conclusions are directly retrieved 
from the Infinitesimal Modeling manual  so that this Elementary Physics manual 
will present$,$ from the view 
point of applications$,$ a continuous and cohesive structure that parallels the 
standard first undergraduate physics course that requires the Calculus as a 
prerequisite. {\bf One important feature of this physics manual is that  
many of the basic rigorous derivations are followed  immediately by additional 
derivations that have been translated into  the  
classical language used in most undergraduate calculus courses. This will 
enable most undergraduate students to more easily comprehend a derivation's logical 
sequence.}    
\par
\bigskip
\leftline{1.2 \underbar{Manual Structure}}\par
\medskip
We will not replicate an actual physics course in this manual but 
rather present \underbar{representative derivations} using the rules 
established within the Infinitesimal Modeling manual for some of the more significant 
integral and differential equation models for the behavior of will-known 
natural systems.\par
 This presentation will only be for mechanics. 
Obviously$,$ we can only make a minute sampling from these very broad 
categories. However$,$ it is hoped that$,$ if care is exercised$,$ the 
examples chosen will lead the instructor to seek more rigorous derivations 
for the more complex and refined aspects of system behavior. The individual 
physics instructor is certainly more intuitively and academically prepared than 
you author for a penetrating and rigorous investigation of the more subtle 
aspects of this subject. \par
 It is the belief of your author that the most expedient approach is to train 
the scientific community in the rudiments of rigorous infinitesimal analysis 
by such devices as the Infinitesimal Modeling manual$,$ the new infinitesimal calculus courses 
that have been introduced throughout the world$,$ and manuals similar to this 
Physical Manual. Once individual 
scientists achieve a working knowledge of the basic principles then those 
who specialize in a given subject area are the appropriate ones to continue a 
more in depth exploration. What is discovered by an in depth rigorous 
infinitesimal approach is that simple fundamental observations lead to simple 
standard expressions. These expressions$,$ after being transferred to the 
nonstandard model$,$ yield a simple view of a new world called the nonstandard 
physical world (NSP-world). These transferred simple processes also lead to 
NSP-world processes that when applied within the NSP-world lead back again to 
standard integral or differential equation models that mirror natural system 
behavior. Usually$,$ one acquires knowledge about the appropriate NSP-world 
processes through observation of simple or idealized natural system 
behavior and then accepts those NSP-world views that lead to verified 
predictions. It is by means of this back-and-forth approach that we gain useful 
knowledge about the NSP-world.\par
The following notation indicates the beginning and ending of each derivation.
For the rigorous analytical derivation using the language of The Basic Manual$,$ the 
beginning and ending are marked by a $\diamondsuit$ symbol. A second or 
third 
derivation is also be denoted by $\diamondsuit;$ but$,$ each is included 
within a subsection marked by $\{${\sl Second Derivation}..........$\}$ or
$\{${\sl Third Derivation}..........$\}$  \par
\vfil\eject
\centerline{Chapter 2.}
\medskip    
\centerline{\bf MECHANICS}              
\bigskip 
\leftline{2.1. \underbar{Instantaneous Velocity.}}\par
\medskip
In 1686$,$ Newton {\tt [Newton [1934]]} gives what he claims is the easily 
comprehended notion of the ``{ultimate velocity},'' or what we now term the
{\bf {instantaneous velocity}}$,$ for an actual real material object. 
{\it But by the same argument it may be alleged 
that a body arriving at a certain place$,$ and there stopping$,$ has no ultimate 
velocity; because the velocity$,$ before the body comes to the place$,$ is not its 
ultimate velocity; when it has arrived$,$ there is none. But the answer is easy; 
for by the ultimate velocity is meant that with which the body is moved$,$ 
neither before it arrives at its last place and the motion ceases$,$ nor after$,$ 
but at the very instant it arrives; that is$,$ the velocity with which the body 
arrives at its last place$,$ and with which the motion ceases.} [{\tt Scholium to 
Lemma XI in Book 1}] For the case of nonzero instantaneous velocity$,$ one might 
gather from this the power of Newton's mental vision and his intuitive 
comprehension of future behavior. Since even though the object may not 
appear to move at the ``instant''(i.e. an instant of time) 
 one observes the hands of clock point at a 
numerical representation for the time$,$ the object did arrive at 
a space location and has the \underbar{capacity} to 
change its position. It is claimed$,$ incorrectly$,$ that this type of change in 
position 
is noted  when a second 
observation is made and the hands of the same clock are assumed to point at a 
different numerical representation for the time. \par 
Newton's modeling of this idea is firmly rooted in his concept of the 
relation between geometry 
(the basic mathematical structure of the 1600's) and its relation to mechanics.
{\it Geometry does not teach us to draw lines$,$ but requires them to be drawn$,$ 
for it requires that the learner should first be taught to describe these 
accurately before he enters geometry$,$ then it shows how by these operations 
problems may be solved. To describe right lines and circles are problems$,$ but 
not geometrical problems. The solution of these problems is required from 
mechanics,....therefore geometry is founded in mechanical practice$,$ and is 
nothing but that part of universal mechanics which accurately proposes and 
demonstrates the art of measure.} [{\tt Newton [1934:xvii]}]
 Newton's claim is that 
our observations  and intuitive comprehension of mechanics comes first in our 
education. These concepts are then abstracted to include the vague notion 
that objects have certain ``capacities or potentials to do things''- 
the {capacity or potential idea}. 
We are told that it is \underbar{after} experimentation$,$ observation and 
reflection that the mathematical structure is evoked and these ``easy'' 
capacity concepts are modeled.\par
The abstract notion of {instantaneous velocity} may have been ``easy'' for Newton 
to grasp$,$ but it was incomprehensible to {Berkeley} and many others who 
believed that such abstractions could not be applied to {actual real material 
objects}. The paramount philosophy of science for Berkeley was a science of the 
material and {directly observed universe}. Any arguments that relied upon 
such abstractions would need to be rejected.\par
How can we communicate such an abstract idea to students who do not necessarily 
possess Newton's obvious mental ability? {\tt Tipler [1982:26]} writes: {\it At first glance$,$ it might seem impossible to define the velocity of a 
particle at a single instant$,$ i.e.$,$  at a specific time. At a time $t_1,$
the particle is at a single point $x_1.$ If it is at a single point$,$ how can 
it be moving? On the other hand$,$ if it is not moving$,$ shouldn't it stay at the 
same point? This is an age-old paradox$,$ which can be resolved when we realize 
that to observe motion and thus define it$,$ we must look at the position of the 
object at more than one time. It is then possible to define the velocity at an 
instant by a limiting process.} Unfortunately$,$ Tipler has reversed Newton's 
original notion$,$ that something exists prior to the modeling of motion and 
this something is the {capacity to move}. This capacity then leads to the need 
to seek various observations from which a numerical value can be ``defined.''
Is it now  possible to \underbar{derive} the well-known derivative expression for 
instantaneous velocity and$,$ at least partially$,$ retain Newton's capacity 
concept by infinitesimally modeling what is indeed easily observed behavior?  
An affirmative answer to this question depends upon your acceptance  
fundamental properties the {Galilean theory of uniform velocity}$,$ infinitesimal 
analysis and its associated interpretations.\par
$\diamondsuit$ Assuming that we are working in the laboratory setting with a  {fixed standard 
for the measure} of (linear) distance and time$,$ let increasing
$d\colon I \to \real,$ where $I = [a,b],\ a\not=b,\ a,\,b \in \real,$ 
represent the {distance an object travels over the time interval} $[a,b].$ 
Of course$,$ a lot has been assumed$,$ even that it makes sense to consider time 
as representable by a continuum such as $[a,b].$ 
[{\it Remark.} Recent work 
[{\tt Herrmann [1989]}] has shown that if time is not a continuum then there are 
internal functions that relate time to a continuum and these functions are 
infinitely close to any discrete (discontinuous) time concept.]  Suppose that 
we extend our ``observations'' of finitely many cases and accept for a very 
simple motion that $d(t_2) -d(t_1)/(t_2 -t_1)= v=\Vert \vec v \Vert $ a constant (the scalar 
velocity)$,$ for any 
$t_1,\, t_2 \in [a,b],\ t_1 < t_2.$ If such a motion persisted$,$ then$,$ 
of course$,$ $v$ can be used to calculate a change in the distance over 
a change in time. Let $t_1 \in (a,b)$ and $\Delta t$ any positive real number 
such that $t_1 + \Delta t \in (a,b].$ Then$,$ assuming a constant scalar 
velocity$,$ 
one has that $d(t_1 + \Delta t) = d(t_1) + v\,\Delta t$ or that
$$d(t_1 + \Delta t) - d(t_1) = v\,\Delta t.\eqno (1)$$ \par
What if the distance expression was  
more complex than the linear type 
expressed by (1)? We are seeking an appropriate definition for $v$ 
that extends this case. First$,$ it follows  
immediately that $v$ would be a function in $t\in [a,b].$ Since a constant 
function is the simplest in the collection of continuous functions$,$ 
then$,$ at this stage of our analysis$,$ we simply require $v$ to be continuous on 
$[a,b].$ However$,$ we know from the Extreme Value Theorem that for any $[t_1,t_2] 
\subset [a,b], t_2 = t_1 + \Delta t,\  \Delta t > 0,$ there exists 
$t_m,\, t_M \in [t_1,t_2]$ such that for each $t \in [t_1,t_2],\ 
v(t_m) \leq v(t) \leq v(t_M).$ Hence$,$ 
$$(t_2-t_1)v(t_m) \leq (t_2-t_1)v(t)\leq (t_2-t_1)v(t_M),\eqno (2)$$
for each $t \in [t_1,t_2].$\par
Now consider the physical processes involved and correspond equation (2) 
to these processes. For the time span $t_2 - t_1,$ it appears reasonable to 
state$,$ using the case where $t_m$ and $t_M$ are constant and our intuitive 
notion of {distance traveled}$,$ that the actual distance 
moved $d(t_2) - d(t_1)$ has the property that
$$ (t_2-t_1)v(t_m) \leq d(t_2) - d(t_1) \leq (t_2-t_1)v(t_M),\eqno (3)$$
If you accept the model for distance expressed by (3)$,$ then from  the 
Intermediate Value Theorem there would necessarily exist
some $t^\prime \in [t_1,t_2]$ such that $d(t_2) - d(t_1)= v(t^\prime)(t_2 - 
t_1).$ What this means is that the distance can be calculated$,$ knowing 
$v(t^\prime),$ as if it were created by a constant scalar velocity. 
Transfer the above results by intuitive *-transfer to the 
infinitesimal NSP-world (i.e. they ``hold'' true for the infinitesimals). 
 Thus$,$ for 
any positive $\eps \in \monad 0$ and any $t_1 \in (a,b),$ it follows that there 
exists some $t^\prime \in [t_1,t_1 + \eps]$ such that 
$$(\hyper d(t_1 +\eps) - d(t_1))/\eps = \hyper v(t^\prime). \eqno (4)$$
The same argument shows that if negative $\eps \in \monad 0,$ then there 
exists some $t^\prime \in [t_1+ \eps,t_1 ]$ such that (4) as well. \par
Unfortunately$,$ we do not know the value of $t^\prime$ in (4). But$,$ once 
again$,$ continuity 
of $v$ at $t_1$ does allow us to write that $\hyper v(t^\prime) \approx 
v(t_1).$ Hence$,$ $$(\hyper d(t_1 +\eps) - d(t_1))/\eps \approx v(t_1). \eqno (5)$$
Obviously$,$ since $\eps$ is an arbitrary nonzero infinitesimal then if there 
exists a distance function$,$ $d,$  that satisfies (3) for all such time 
intervals,
then $d$ must be differentiable at $t_1$ and application of the standard part 
operator implies that $d^\prime (t_1) = v(t_1).\diamondsuit$\par
\medskip
$\{${\sl Second derivation.} In what follows$,$ the above derivation for the instantaneous velocity 
function is reworded into a quasi-classical description using slightly 
modified calculus terminology. The ground rules for this second derivation 
are:\par
(i) As is done in Internal Set Theory$,$ the ``*'' notation is removed from the 
functions since whether they  are 
nonstandard extensions of standard functions is clear from the function's 
argument (i.e. preimage).\par
 (ii) The symbols ``$\approx$'' is translated by the term ``infinitely 
close.'' This relation can be physically characterized by stating that no 
standard machine can measure any difference between the quantity on the left 
and the quantity on the right no matter how small the machine error.
\par
(iii) Except for $\nat_\infty,$  hyperreal numbers are usually limited. Hence$,$ 
simply call such a hyperreal number by the single word term ``number.''\par
(iv) We use the fact that functions defined and continuous on $[a,b]$ 
preserve the infinitely close 
concept for these numbers. That is if $t,\ t_1 \in [a,b]$ and $t \approx t_1,$ then $f(t) \approx 
f(t_1).$ Infinitesimals may be called the ``infinitely or very small.'' These 
numbers 
can be physically characterized as measures that are smaller than any 
standard machine error --- measures that appear to a machine to be zero.\par
(v) Rather than use the standard part operator$,$ where applicable$,$ 
use the simple term ``limit'' in its place$,$ since it has the same 
operative properties. Also use the fact 
that limits of two infinitely close numbers are equal.\par
The modified classical derivation is exactly the same until after equation 
(3). Then it continues as follows:\par
$\diamondsuit$ If you accept the model for distance expressed by (3)$,$ then 
from  the 
Intermediate Value Theorem there would necessarily exist
some $t^\prime \in [t_1,t_2]$ such that $d(t_2) - d(t_1)= v(t^\prime)(t_2 - 
t_1).$ What this means is that the distance can be calculated$,$ knowing 
$v(t^\prime),$ as if it was created by a constant scalar velocity.            
Now (3) and these facts hold for the infinitely small.  Thus$,$ for 
any positive infinitely small $\Delta t$ and any $t_1 \in (a,b),$ 
it follows that there exists some $t^\prime$  such that $t_1 \leq t^\prime 
\leq t_1 + \Delta t$ (i.e.$,$ $t^\prime \in [t_1, t_1 + \Delta t]$ and 
$$(d(t_1 +\Delta t) - d(t_1))/\Delta t = v(t^\prime). \eqno (4)$$
The same argument shows that if $\Delta t$ is a negative infinitely small 
number$,$ 
then there 
is some number $t^\prime $ such that $t_1 + \Delta t \leq t^\prime \leq t_1$ 
and once again (4) holds. \par
Unfortunately$,$ we do not know the value of $t^\prime$ in (4). But$,$ since 
$t^\prime 
\approx t_1$ then continuity 
of $v$ at $t_1$ allows us to write that $v(t^\prime) \approx 
v(t_1).$ Hence$,$ $$(d(t_1 +\Delta t) - d(t_1))/\Delta t \approx v(t_1). \eqno (5)$$
Obviously$,$ since nonzero $\Delta t$ is an arbitrary and  infinitely small$,$ then if there 
exists a distance function$,$ $d,$  that satisfies (3) for all such time 
intervals,
then $d$ must be differentiable at $t_1$ and the limit of the left hand side 
of (5)$,$ [{\it as $\Delta t$ varies could be added$,$ but is not necessary}] 
must equal the limit of the right hand side 
which is the constant $v(t_1).$ This implies that $d^\prime (t_1) = 
v(t_1).\diamondsuit\}$\par
\medskip 
$\{${\sl Third derivation --- entirely classical.}  
As discussed in the Infinitesimal Modeling manual$,$ 
many derivations for differential equation models require an infinitesimalizing 
process for behavior that is only approximated within the standard world for 
"small" quantities. One of the simplest illustrations of this is the 
differential equation model for Newton's Law of Cooling. The necessity for 
this special process comes from the experiential evidence that the observed 
behavior holds only for small measures of the independent variables
and as the measures are reduced such 
behavior is more closely approximated by a standard functional expression. 
This is in direct contrast to an expression such as (3) and the discussion 
that follows where the results appear to hold for all intervals $[t_1,t_2].$
In the case of such concepts as the instantaneous velocity$,$ it is possible to 
present an entirely classical derivation. This classical derivation 
begins with equation (3) and continues as follows: \par
$\diamondsuit$ If you accept the model for distance expressed by (3)$,$ then 
from  the 
Intermediate Value Theorem there would necessarily exist
some $t^\prime \in [t_1,t_2]$ such that $d(t_2) - d(t_1)= v(t^\prime)(t_2 - 
t_1).$ What this means is that the distance can be calculated$,$ knowing 
$v(t^\prime),$ as if it were created by a constant scalar velocity. Thus for
$t_1 \in (a,b)$ and for every positive $\Delta t$ such $t_2 = t_1 + \Delta t 
\in (a,b)$ there exists some $t^\prime $ such that 
$$(d(t_1 + \Delta t) - d(t_1))/\Delta t) = v(t^\prime), \eqno (4)$$
and $t_1 \leq t^\prime \leq t_1 + \Delta t.$\par
Using the Axiom of Choice$,$ we can consider a function $f$ defined on the 
respective $\Delta t$ such that $f(\Delta t) = t^\prime.$ Repeating the 
process for the negative $\Delta t$ such that $[t_1 + \Delta t, t_1] \subset 
(a,b)$ and extending the function $f$ to include these negative $\Delta t$
leads to the conclusion that $\lim_{\Delta t \to 0}f(\Delta t) = t_1.$ 
From the assumed continuity of $v$ it follows that
$$\lim_{\Delta t \to 0}\left({{d(t_1 + \Delta t) - d(t_1)}\over{\Delta 
t}}\right)=\lim_{\Delta t \to 0}v(f(\Delta t))=v(t_1). \eqno (5)$$
This all implies that under the conditions stated $d^\prime (t_1) = v(t_1).$ 
$\diamondsuit \}$\par
How the above derivations improve our comprehension of the concept of 
instantaneous velocity is discussed at the conclusion of this section. 
Returning to the concept of 
the {capacity to move}$,$ Theorem 9.1.2 of the Infinitesimal Modeling manual indicates that for 
each $\eps \in \monad 0$
$$ \Delta \hyper d_{t_1}(\eps) =\hyper d(t_1 +\eps) - d(t_1)= 
v(t_1)\eps + \eps\lambda (\eps)= \Hyper T(\eps) + \eps\lambda (\eps), \eqno 
(6)$$
where $\lambda$ is a local function defined by (6) and linear $T\colon \real 
\to \real.$
Thus for each $\eps \in \monad 0,\  \Delta \hyper d_{t_1}(\eps)$ and 
$v(t_1)\eps$ are not just infinitely close$,$ but they are {infinitely close of the 
first order}. (See the Infinitesimal Modeling manual Section 8.3.)\par
\vskip 18pt
\hrule
\smallskip
{\bf {Definition 2.1. (Infinitely Close of Order One)}.} 
 Two hyperreal valued functions defined on $\monad 0$ are said to be 
{\bf {Infinitely Close of Order One or of the First Order}} if for each 
$\eps \in \monad 0$ there exists some $t\in \monad 0$ such that 
$$f(\eps) = g(\eps) + \eps\, t.$$
In which case$,$ this is denoted by $f \sim_1 g.$ Further$,$ if the  two functions 
are considered to be measuring some physical properties$,$ then we often say 
that the natural world effects of these properties are 
{{\bf indistinguishable (at level one or on the first level)} }.\par
\smallskip
\hrule
\vskip 18pt
It is a simple matter to show that $\sim_1$ is an equivalence relation on the 
set of all hyperreal valued functions defined on $\monad 0.$ The capacity to 
move concept is now represented in the monadic environment by noting that 
$$ \Delta \hyper d_{t_1}\sim_1 \Hyper T= d^\prime(t_1)(\cdot).\eqno (7)$$ 
Or in words$,$ {\sl within the monadic world the distance represented by 
$ \Delta \hyper d_{t_1}$ is indistinguishable from (has the same effect as) 
that produced by a scalar velocity $d^\prime(t_1).$} Of course$,$ the function 
$d^\prime$ is termed the  {{\bf instantaneous velocity}}. As Newton claimed$,$ if 
mechanics leads to geometry$,$ then this is what motivated the geometric 
concepts of the  
rectifiable curve$,$ tangents$,$ curvature and the like that appear in the 
Infinitesimal Modeling 
manual Chapter 7$,$ section 7.2$,$ Chapter 8$,$ sections 8.5 -- 8.6 and Chapter 9$,$ 
sections 9.5 -- 9.7.\par
Within the foundations of any discipline it is often difficult to refine even 
slightly what may have been assumed previously to be an elementary and 
not dissectible assertion. Thus$,$ until now$,$ this has been with the idea of 
instantaneous velocity. However$,$ \par
(i) let the distance function$,$ $d(t),$ and an  unknown continuous scalar 
velocity 
function$,$ $v(t) = \Vert \vec v(t)\Vert,$ be related by expression (3).\par
(ii) Let (3) hold for every time subinterval $[c ,d] =[t_1 + (\Delta t)_1, t_1]\subset 
(a,b)$ and $[c,d] = [t_1, t_1 + (\Delta t)_2],$ where 
$(\Delta t)_1< 0, \ (\Delta t)_2> 0.$\par
(iii) There exists some $t^\prime \in [c,d]$ such that the actual distance 
traveled $d(t_2) - d(t_1) = v(t^\prime)(t_2 - t_1).$ \par
{\bf Then $d$ must be differentiable at $t_1$ and the only standard scalar 
velocity function that satisfies (i)$,$ (ii) and (iii) is the function $v(t_1) =
d^\prime(t_1).$} The instantaneous velocity is not obtained by simply 
postulating a definition but is derived from more elementary observations.\par
\bigskip 
\leftline{ 2.2. \underbar{Acceleration}}
\medskip
For twenty years$,$ {Galileo struggled} with the problem of representing the 
velocity of a falling body in terms of distance {\tt [Gillispie [1960:42]]}.  
After failing in every 
attempt$,$ a new idea began to ferment - an idea that today seems so common 
place since$,$ as illustrated in the previous section$,$ we are taught to think of 
elementary velocity as expressed in a time coordinate. But$,$ it was Galileo's 
revolutionary concept of {time as an independent abstract entity} that led to 
the theory of motion that Newton applied in his dynamical geometry and 
Mathematical Principles of Natural Philosophy [{\tt [Newton [1686]}].  
The remarkable insight exemplified by Galileo should not be underestimated.
His struggle once again indicates the difficulty that scientists often face 
when$,$ through reflection$,$ they alter some well entrenched but erroneous 
elementary precept.\par
 $\diamondsuit$ As seen in the previous section$,$ the only elementary  
standard function the preserves our intuitive understanding about the 
elementary measure of distances traveled is the instantaneous velocity. 
Further$,$ in the monadic environment$,$ the Galilean properties for constant 
or uniform scalar velocities and how they are compounded to yield distances traveled are 
indistinguishable from the actual quantities. But$,$ now that we have accepted 
the instantaneous velocity$,$ $v = \Vert \vec v\Vert,$ as the appropriate 
elementary nonconstant velocity concept$,$ we can certainly apply a section 1
type analysis to $v.$ This requires the strict application of the {Galilean
theory of motion} for a constant acceleration on one hand$,$ and Newton's concept 
of the nonconstant acceleration (produced by a force) that leads to nonuniform
velocity on the other hand.  \par
Thus for time interval $[a,b]$ substitute in expression (3) of section 1,
 $v$ for $d$ and a 
representation for a continuous scalar acceleration $a = \Vert \vec a \Vert$ for $v$
and obtain
$$(t_2-t_1)a(t_m) \leq v(t_2) - v(t_1) \leq (t_2-t_1)a(t_M).\eqno (1)$$
The expression on the left of (1) applies the concept of a constant 
acceleration and measures the minimum possible [linear]
{uniform change in the velocity} and that on the 
right the maximum possible uniform change.
This implies$,$ as in the case of instantaneous velocity and  under parallel 
hypotheses as stated in section 1 (i)$,$ (ii)$,$ that  there exists 
a unique scalar acceleration function $a$ such that $a(t_1) = v^\prime(t_1)= 
d^{(2)}(t_1)$ 
- the {\bf {instantaneous acceleration}}. $\diamondsuit$ As with the case of the 
instantaneous velocity we also have 
$$\Delta \hyper v_{t_1} \sim_1 \Hyper T_a = v^\prime(t_1)(\cdot),\eqno (2)$$ 
which indicates that within the monadic environment that the change in 
velocity in the N-world is indistinguishable on the first level from that 
produced by a constant acceleration. However$,$ application of Theorem
8.4.2 yields that an associated change in distance satisfies the two$,$ 
difficult to visualize$,$ statements 
$$ \Delta^2 \hyper d_{t_1}\sim_1 \Hyper T_a= v^\prime(t_1)(\cdot)=d^{(2)}(\cdot),\ {\rm and}$$ 
$$\Delta^2 \hyper d_{t_1}\sim_2 \Hyper T_a= \Hyper 
d^{(2)}(\cdot),\eqno (3)$$
where $\sim_2$ is defined in the obvious manner by replacing $\eps$ in 
the last term in the equation in definition 2.1 by $\eps^2$ and the second equation in (3) is the best that we 
can state$,$ in generally.\par 
\bigskip
\leftline{ 2.3. \underbar{Forces and Newton's ``Law''}}
The story is told that a student asked {Max Planck} to explain how he 
 perceived nuclear forces? Planck is said to have replied that he 
would perceive them as someone pulling on his coat sleeve. Whether or not this 
story is factual$,$ experience indicates that a change in velocity is better 
comprehended by considering such a change as the effect of a more easily 
sensed ``{force}'' that may be the cause of the change. This brings use the 
Newton's so-called Second Law of Motion.\par
It is not true that Newton formulated his Second Law as it is taught in our 
elementary physics courses where we are told that the scalar force is proportional to the 
instantaneous acceleration. It is also not true that he formulated his Second 
Law in terms of something equivalent to the derivative of the {momentum}. First,
Newton defined the ``quantity of motion'' as follows: {\it The {quantity of 
motion} is the measure of the same$,$ arising from the {\rm [scalar]} velocity and the quantity 
of matter conjointly} [{\tt Newton [1934:1]}] Thus the quantity of motion  is the 
momentum. As indicated by the explanation that 
follows his statement$,$ his Second Law was the observation that {\it The {\rm 
[uniform]} change of {\rm [the quantity of]} motion  is proportional to the 
{\rm [constant]} motive force impressed; and is made in the direction of the
right {\rm [i.e.$,$ straight]} line in which the force is impressed.} [{\tt [Newton
 [1934: 13]}] He proceeds in the Scholium to that section to apply this 
Second Law and his idea that the 
total effect of finitely many constant forces is additive over time to establish 
Galileo's discovery that the {\it ...descent of bodies varies as the square of 
the time.} [{\tt Newton [1934:21]}] In terms of a constant {impulse} notion$,$ Newton's argument 
does not include the mass$,$ but rather leaves the mass as the constant of 
proportionality. In summation notation$,$ the argument has the following form.
The time $[a,b]$ is subdivided into $n$ equal time intervals $\Delta t.$ 
The scalar force during these time intervals is a constant $F.$ Then the force 
times the length of time impressed (i.e. the {\bf {impulse}})is proportional to the uniform change in 
scalar velocity$,$ $\Delta v,$ and all such uniform changes in the scalar velocity are equal in value. 
Hence$,$ the composition of such forces would yield a total effect
$$\sum F(\Delta t) = F \sum \Delta t = F(b-a) \propto \sum \Delta v=
v(b) - v(a)= v,\eqno (1)$$
where $F$ now becomes the {\it {whole force}}$,$ the constant of proportionality is the 
same for each summand and from  the additivity of velocities in a 
straight line motion $v$ is the  {\it {whole velocity}.} Now translating into our 
symbols$,$ Newton writes $F(b-a) \propto v$ as $(b-a) \propto v.$ He then states 
that the {\it spaces} [i.e.$,$ distance traveled] in {\it proportional times are 
as the product of the velocities and times;....} [{\tt Newton [1934:21]}] Thus such a 
distance $d \propto v(b-a) \propto (b-a)^2.$ \par
In his Principles$,$ Newton does not utilize his method of fluxions$,$ even though 
in other communications he does$,$ to develop his theory of motions of 
material bodies. He first presents arguments delineating the ``ultimate ratios'' 
between geometric measures - arguments that employ those intuitive concepts 
acceptable and apparently comprehensible by geometers$,$ expressions such as 
{\it magnitude diminished  \underbar{in infinititum}}. He then correlates time 
to the measure of one leg of a right triangle and velocity to the other leg. 
Considering his geometric notion of the ultimate ratio of the areas of these 
triangles as 
the length of the sides {\it diminish \underbar{in infinititum}}$,$ which he 
previously established was as to the squares of the {\it homologous sides},
he draws the conclusion that {\it the spaces which a body by an finite force 
urging it$,$ whether that force is determined and immutable$,$ or is continually
augmented or continually diminished {\rm [with respect to time]}$,$ are in the 
very beginning of the motion 
to each other as the squares of the times.} [{\tt Newton [1934:34]}] In corollary 
{\it iv} he writes: {\it And therefore the forces are directly as the spaces 
described in the very beginning of the motion$,$ and inversely as the squares of 
the times.} [{\tt Newton [1934:35]}] The expression {\it very beginning} is used to 
describe the ultimate ratio concept and that this is a point force 
associated with an instant of time. Further notice that he has replaced the 
idea of constant forces over a time subinterval with forces that are 
being altered {\it continually}. \par
The logical vagueness of Newton's arguments can be eliminated by modern 
infinitesimal analysis. $\diamondsuit$ Consider the time interval $[a,b]$ and any
$[t_1,t_2] \subset [a,b].$ Assume that the scalar {force $F$ is a continuous 
real valued function} defined on $[a,b].$ We look at Newton's observations 
relative to {constant scalar forces} and there relation to uniform changes in 
the scalar velocity. First$,$ from continuity$,$ there exists $
t_m,\, t_M \in [t_1,t_2]$ such that 
$$F(t_m)(t_2 -t_1) \leq F(t)(t_2-t_1) \leq F(t_M)(t_2-t_1)\eqno (2)$$ 
for each $t \in [t_1,t_2].$ Let $C$ be a constant of proportionality which$,$ 
for this 
particular simplistic problem$,$ is considered absolutely constant in 
character. Let $q$ denote a real valued function that measures$,$ with respect 
to time$,$ the notion of the quantity  of motion. 
Then from the actual stated Second Law
it follows that the minimal possible change in momentum over the $[t_1,t_2]$ 
is $q(t_m)= F(t_m)(t_2 -t_1)$ and the maximal change is $q(t_M)= 
F(t_M)(t_2 -t_1).$ Now the actual velocity function in terms of time 
is$,$ at present$,$ unknown. But$,$ whatever it may be$,$ we consider the actual 
quantity of motion $q$ to be an extension of the constant case and$,$ hence$,$ it 
is proportional$,$ for a given object$,$ to the actual change in velocity or 
to $C(v(t_2) - v(t_1)).$ Assuming that such a change lies somewhere 
between the minimal and maximal changes then applying the continuity of the 
force function one obtains that there exists some $t^\prime \in [t_1,t_2]$ such that
$$F(t^\prime)(t_2-t_1) = C(v(t_2) - v(t_1)).\eqno (3)$$
Letting $t_1 \in (a,b)$ and noticing that the above is assumed to hold for 
{\it any} subinterval $[t_1,t_1 + \Delta t],\ \Delta t>0$ of $(a,b]$ or {\it 
any} subinterval $[t_1 + \Delta t, t_1],\ \Delta t < 0$ of $[a,b)$ then *-transfer yields 
that for nonzero $\eps \in \monad 0$ there exists some $t^\prime \approx t_1$ 
such that
$$\hyper F(t^\prime) = {{C(\hyper v(t_1 + \eps) - 
v(t_1))}\over{\eps}}\approx F(t_1).\eqno (4)$$\par
Taking the standard part operator we arrive at a derivation of our modern 
Newton's Second Law of Motion. That for all of this to occur there must exist 
a velocity function that is differentiable at $t_1$ and the only relation 
between such a point force function and such a velocity function is 
$$F(t_1) = Cv^\prime(t_1) = Ca(t_1),\eqno (5)$$
where $a$ is the {\bf {instantaneous acceleration}.} $\diamondsuit$ As previously$,$ this can be 
further related to the change in momentum for an infinitesimal 
time by the expression
$$C(\Delta \hyper v)\sim_1  F(t_1)(\cdot). \eqno (6)$$\par
Or$,$ as stated in words$,$ that the  {change in momentum over an infinitesimal 
time} is 
indistinguishable on the first level from that produced by a constant force 
applied to infinitesimal time periods. Moreover$,$ the easily grasped 
concept of the impulse was the  starting point in Newton's original 
arguments. Unfortunately$,$ some modern textbooks do not introduce the impulse 
relative to constant forces as a first principle.\par
$\{${\sl Second derivation.}
Start the modification of the above after equation (3).\par
$\diamondsuit$ Letting $t_1 \in (a,b)$ and noticing that the above 
is assumed to hold for 
{\it any} subinterval $[t_1,t_2]$ of $[a,b]$ then all of the above holds for 
infinitely small subintervals. Thus for  nonzero infinitely small $\Delta t$ 
there exists some $t^\prime \approx t_1$ 
such that
$$F(t^\prime) = {{C(v(t_1 + \Delta t) - 
v(t_1))}\over{\Delta t}}\approx F(t_1).\eqno (4)$$\par
Taking the limit of the left hand side [as $\Delta t$ varies] we arrive at a derivation of our modern 
Newton's Second Law of Motion.  
$$C\lim {{(v(t_1 + \Delta t) - 
v(t_1))}\over{\Delta t}}= F(t_1).\eqno (5)$$
Thus for all of this to occur there must exist 
a velocity function that is differentiable at $t_1$ and the only relation 
between such a point force function and such a velocity function is 
$$F(t_1) = Cv^\prime(t_1) = Ca(t_1),\eqno (6)$$
where $a$ is the {\bf {instantaneous acceleration}.} $\diamondsuit \}$\par
The reader may replicate the third derivation along with the 
discussion of properties (i)$,$ (ii) and (iii) as they appear in Section 2.1 for 
the instantaneous acceleration and Newton's Second Law.\par
\vskip 18pt
\hrule
\smallskip
\hrule
\smallskip
Notice that one derivation method -- the maximum and minimum method -- yields 
the instantaneous velocity,
instantaneous acceleration and the modern Second Law of Motion from what may 
be more fundamental observations. \par
\smallskip
\hrule
\smallskip
\hrule
\vskip 18pt
\leftline{ 2.4. \underbar{Vectors}}
\medskip 
For constant forces$,$ Newton's Corollary I to his three laws establishes for 
constant forces the idea that bodies move in $\realp 2$ by the {vector 
composition} of two scalar forces acting simultaneously upon a particle or 
fixed point in a material body.
 [{\tt Newton [1934:14]}] Newtonian mechanics may appear to begin with the idea 
that the position of a particle in an {$n$-dimensional coordinate 
system} 
is dependent upon the {composition of $n$ forces (the cause)}; but$,$ in actually$,$ 
it is the {position (i.e. the effect)} that is the fundamental concept. The 
{position of a particle} is defined in terms of $n$ coordinate functions$,$ each 
expressed in the independent parameter --- time. This leads to the 
{\bf {position vector (i.e. radius vector)}} 
$\vec r = (x_1(t),\ldots,x_n(t)),\ t \in [a,b].$
Applying the previous section to these coordinate functions independently$,$ leads 
to the vector form for the {instantaneous velocity}$,$ {instantaneous 
acceleration}$,$ and {force}.\par
In the Infinitesimal Modeling manual$,$ the geometric concept of the {length of a continuous 
curve} is fully discussed and$,$ where possible$,$ this length is correlated to the 
integral by means of our integral rules. Following Newton's notion of a {dynamic 
geometry}$,$ the same conclusions evolve from the instantaneous velocity. 
Suppose that a point force $\vec F(t),\ t \in [a,b],$ is the cause that 
induces an instantaneous acceleration $\vec a(t)$ upon a particle. Since $\vec F$ 
was assumed continuous on some time interval $[a,b]$ then $\vec v\,^\prime$ is 
continuous on $[a,b].$ Consequently$,$ $\vec v$ and $\Vert \vec v\Vert$ are 
continuous on $[a,b].$ 
Assuming we are in $\realp n$ then$,$ noting that for $t \in [a,b],\ \Vert v(t) 
\Vert = \sqrt {(x_1^\prime(t))^2 + \cdots (x_n^\prime(t))^2} = 
d^\prime(t)/dt,$ we obtain
$$d(b) - d(a) = \int_a^b \sqrt {(x_1^\prime(t))^2 + \cdots (x_n^\prime(t))^2}\, 
dt,\eqno (1)$$
which is the same result obtained in the Infinitesimal Modeling manual by considering 
{hyperpolygonal representations} for the length of a continuously differentiable 
geometric curve. \par
Let the {path of motion} be represented by continuous $c\colon [a,b] \to 
\realp n,$ where $c(t) = (x_1(t),\ldots,x_n(t)).$ 
By Theorems 9.5.1 and 9.5.2 of Chapter 9 and Theorem N.1 at 
the end of Chapter 10 in the Infinitesimal Modeling manual it is shown that if 
$\vec v (t) =\vec c\,^\prime(t) \not= \vec 0,\ t \in (a,b),$ then the unit tangent $\vec T = 
\pm \vec c\,^\prime (t)/\Vert \vec c\,^\prime (t)\Vert$ and that $\vec T$ is 
{almost parallel} to every {hyperpolygonal directed line segment}
$(\hyper c(t + dx) - c(t))/\Vert \hyper c(t + dx) - c(t)\Vert$ for each 
nonzero $dx \in \monad 0.$ This implies that the unit instantaneous velocity 
vector $\vec v(t)/\Vert \vec v(t) \Vert$ is not only one of the 
two possible unit tangent vectors but also has the property of being 
almost parallel to each 
$(\hyper c(t + dx) - c(t))/\Vert \hyper c(t + dx) - c(t)\Vert.$ Since the 
effect of vectors as a model for natural world (i.e. {N-world}) behavior is often 
related to the physical concept of ``{direction}'' then as far as the N-world 
effects are concerned the direction of the standard velocity vector 
attached to the point $c(t)$ is {indistinguishable from} the direction of 
$\pm (\hyper c(t + dx) - c(t)),\ 0\not= dx  \in \monad 0.$ The 
indistinguishableness of such effects is beneficial when energy is to be 
considered.\par
\bigskip
\leftline{2.5. \underbar{Energy and Force Fields}}
\medskip
In the Infinitesimal Modeling manual$,$ the line integral is employed to measure  the {\sl energy 
expended within a forced field while moving along a curve} (Application 8.2.1)
and the customary method of line integral evaluation obtained. We examine here 
the first portion of this derivation. \par
$\diamondsuit$ In experimental physics$,$ the concept of ``{work}'' ({energy expended}) 
is introduce. {\bf All 
one needs to do is to establish its properties for a polygonal curve.}  
Suppose we have continuous force field $F\colon E \to 
\realp n.$ Let ${\cal P}_k\subset \realp n$ be any finite polygonal curve$,$ 
$\ell_j$ one of 
the line segment portions of ${\cal P}_k$ with $\vec v_j$ denoting this line 
segment considered as a directed line segment in the direction of motion 
through the field $F.$ If $F$ is constant on $\ell_j,$ then the work done 
moving along $\ell_j$ is defined as $W(\ell_j) = F \bullet (\vec v_j/\Vert 
\vec v_j \Vert)\Vert \vec v_j \Vert,$ where length of $\ell_j= \Vert \vec v_j 
\Vert.$
 What if the force field is not constant? Consider ${\cal P}_k$ as 
represented by a continuous $\ell\colon [a,b] \to \realp n$ and assume that 
$F$ is defined on $\ell.$ Then for a given 
$\ell_j= \{(x_1(t),\ldots,x_n(t))\bigm \vert t \in [t_{j-1},t_j]\}$ there exists 
some $t_m,\ t_M$ such that $W_m(\ell_j) = F(\ell(t_n))\bullet \vec 
v_j\leq W(\ell_j)=F(\ell(t))\bullet\vec v_j= W_M(\ell_j)=F(\ell(t_M))\bullet 
\vec v_j$ for each $t \in [t_{j-1},t_j].$ 
Let's make the one assumption that the actual amount of work expended moving 
along the line segment $\ell_j$ is $W(\ell_j)$ and that 
$W_m(\ell_j) \leq W(\ell_j) \leq W_M(\ell_j).$
Then from continuity there exists some $h_j^\prime \in [t_{j-1},t_j]$ such that 
$W(\ell_j) = F(\ell(h_j^\prime))\bullet \vec v_j.$ 
The idea of the nonconstant force field over a line segment is embedded into the 
the NSP-world by *-transfer assuming that what has been established above 
holds for all such polygonal curves. Hence$,$ let ${\cal P}_\Omega$ be a 
hyperpolygonal representation for the curve generated by a fine partition$,$ 
$\ell_j$ an hyperline segment in ${\cal P}_\Omega.$ Since $F$ is continuous on 
$E$ then $\hyper F$ is defined on ${\cal P}_\Omega.$ It follows 
that $\Hyper W(\ell_j) = \hyper F(\ell_j (h_j^\prime))\bullet \hyper {\vec 
v_j}.$ 
For polygonal curves$,$ in general$,$ the work done is an additive function. 
Thus for the hyperpolygonal curve ${\cal P}_\Omega$ 
$$\Hyper W({\cal P}_\Omega) = \sum_{j=1}^\Omega \hyper F(\ell_j (h_j^\prime))
\bullet \hyper {\vec v_j}.\eqno (1)$$ \par
The following is established within this derivation and in the 
Infinitesimal Modeling manual. 
Let $c\colon [a,b] \to \realp n$ be a continuous differentiable curve with 
graph $C$ and assume that $c^\prime (t) \not= 0$ for each $t \in [a, b].$ 
Assume that uniformly continuous $F\colon E \to \realp n,$ open $E \supset 
C.$ The work done in moving through the force field on the path $C$ is
$$W(C) = \St {\left(\sum_{j=1}^\Gamma \hyper F(\ell_j(h_j^\prime))\bullet
\hyper {\vec v_j}\right)}= \int_C F\bullet d\vec R,\eqno (2) $$
where the integral is the line integral over $c$ and ${\cal P}_\Gamma$ is 
\underbar{any} {hyperpolygonal representation} for the curve $C.$ Thus$,$ under the 
hypotheses given$,$ the work done in the NSP-world moving along these 
hyperpolygonal curves is indistinguishable$,$ in general$,$ from what is accepted 
as the N-world work in traveling along the curve. $\diamondsuit$\par
$\{${\it Second derivation.} $\diamondsuit$ See the second derivation for the 
general impulse 
in section 2.6 and modify the above accordingly. $\diamondsuit \}$ \par
{\bf It is significant to realize 
that even though the line integral will exist under less constraints and one 
could extend the concept of energy$,$ to say curves that are  rectifiable but not 
smooth$,$ this need not imply that there is a meaningful correlation between 
this extended concept and its simplistic restriction to polygonal curves.} 
On the other hand$,$ the above derivation once again utilizes a {maximum and 
minimum} approach relative to a basic geometric element.\par  
\bigskip
\leftline{2.6 \underbar{General Impulse}}
\medskip
The energy or work done is considered to be the standard part of any one of a 
collection of hyperfinite sums each term of which is modeled by the
amount of energy expended moving along a hyperline segment through a 
constant force field. The {hyperfinite sums} may be  manipulated internally
as if they were finite sums and preserve the intuitive concept of finite 
additivity. Unfortunately$,$ when an object is represented by a simple 
differential equation obtained by such methods as illustrated in section 2.1$,$ 
2.2$,$ and 2.3 then many authors solve such expressions by elementary 
differential equations methods without given full infinitesimal meaning to the process 
involved. The Fundamental Theorem of Integral Calculus is relevant to the 
existence of such solutions; but$,$ it seems$,$ that in order to grasp the NSP-world 
significance of the concept 
being measured that {due consideration} should be given {to the integral 
solution} 
rather than simply expressing the result by means of a rote procedure. As an 
illustration of this consider the extension of the concept of the impulse.\par
Previously the impulse was a natural notion restricted to constant force 
fields. A scalar force $F,$ when applied for a period of time $\Delta t$ 
produces an altered momentum $\Delta p.$ This leads  to the expression
$F\Delta t = \Delta p.$ Then $F \Delta t=I$ is defined as the {impulse}.
The extension to vector notation is obvious $\vec I = \vec F\,\Delta t= \Delta 
\vec p,$ with a meaningful measure being the Euclidean norm. Applying 
the infinitesimalizing process (2)$,$ (3)$,$ (4) and (5) of section 2.3$,$ where 
$C(v(t_2) - v(t_1))$ is replaced by $\vec p(t_2) - \vec p(t_1)$ leads to the 
customary elementary  derivative expression 
$$\vec F = {{d\vec p}\over{dt}}.\eqno (1)$$
Of course$,$ the vector $\vec F= (p_1^\prime,\ldots,p_n^\prime)$ where the $p_i$ 
are the components for the momentum vector $\vec p.$ \par
As far as a generalization for the impulse is consider there are various 
approaches. The Self-evident Max. and Min. Theorem 6.2.4 in the Infinitesimal 
Modeling manual 
applied to the components of $\vec p$ leads to the conclusion$,$ if one has an 
intuitive comprehension of the basic additive 
properties for the impulse concept$,$  that the proper 
expression for the general impulse for a continuous time dependent 
force field $\vec F$ over the time interval $[a,b]$ should be 
$\vec I = \int_a^b \vec F\, dt.$ However$,$ it may be more motivational and 
instructive to consider$,$ once again$,$ hyperpolygonal paths of motion since 
momentum is modeled after the motion concept and the concept of 
hyperfinite summation only. The summation of the measures of elementary physical 
processes still remains a primary modeling procedure that dictates the overall 
physical effect.\par 
$\diamondsuit$ Suppose$,$ as in section 2.5$,$ that uniformly continuous force field 
$F\colon E \to \realp n,$ open $E \subset \realp n,\ F =(F_1,\ldots,F_n).$
For a line segment $\ell^\prime\subset E,$ determine by a linear system of 
equations
over the time interval $[t_1,t_2]$ and for a component $F_i$ equation 
(2) section 2.3 can be re-expressed in terms of an impulse vector
$\vec I = (I_1,\ldots,I_n)$ as $F_i(t_m)(t_2 - t_1) = I_{im} \leq F_i(t)(t_2 - 
t_1) \leq F_i(t_M)(t_2-t_1) = I_{iM}.$ Thus there is some $t_2^\prime \in [t_1,t_2]$ 
such that the actual $i$'th component of the impulse over that line segment is 
$I_i(t_2^\prime) = F_i(t_2^\prime)(t_2 - t_1).$ By *-transfer this holds for 
any hyperline segment. {\bf Over a finite polygonal line $\ell$ it is assumed 
that 
the total $i$'th component of the impulse is the simple sum of the $i$'th 
component over the individual line segments}. Considering the impulse to be a 
function defined over the set of all finite polygonal curves in $E,$ then  
for any hyperpolygonal curve ${\cal P}_\Gamma$
$$\Hyper I_i({\cal P}_\Gamma) = \sum_{j=1}^\Gamma  \hyper F_i(t_j^\prime)(t_j - 
t_{j-1}),\eqno (2)$$                                                          
where $\{t_0,\ldots,t_\Gamma\}$ is the fine partition of $\Hyper [a,b]$ that 
generates ${\cal P}_\Gamma.$\par
Now from uniform continuity$,$ $ \hyper F_i(t_j^\prime) =  \hyper F_i(t_j) + \eps_j,\ \eps_j \in 
\monad 0.$ Since any finite set of real numbers contains a maximum$,$ then any 
hyperfinite set of hyperreal numbers contains a maximum. Thus there exists 
some $\eps = \max \{\eps_1,\ldots,\eps_\Gamma\} \in  \monad 0$ and following 
the {elemental derivation process} from the Infinitesimal Modeling manual$,$ section 8$,$ we have
$$\left| \sum_{j=1}^\Gamma \eps_j(t_j - t_{j-1})\right| \leq \sum_{j=1}^\Gamma 
\vert \eps_j\vert\,(t_j - t_{j-1})\leq \vert \eps \vert \sum_{j=1}^\Gamma 
(t_j - t_{j-1})= \vert \eps \vert\,(b-a)\Rightarrow$$
$$\lambda_i=\sum_{j=1}^\Gamma \eps_j(t_j - t_{j-1})\in \monad 0.
\eqno (3)$$ 
Substitution into equation (2) yields
$$\hyper I_i({\cal P}_\Gamma) = \left(\sum_{j=1}^\Gamma \hyper F_i(t_j)(t_j - 
t_{j-1})\right) + \lambda_i,\eqno (4)$$\par
Now $F_i$ is continuous. Thus $\st {\hyper I_i} = \st {\sum_{j=1}^\Gamma \hyper F_i(t_j)(t_j - 
t_{j-1})}$ exists$,$ by Theorem 5.1.2 of the Infinitesimal Modeling manual$,$ since 
$\{t_1,\ldots,t_\Gamma\}$ is an internal intermediate partition and 
$\st {\hyper I_i} = \int_a^b F_i(t)\, dt.$ Notice that the hyperpolygonal curve 
was only an auxiliary notion since this last result is determine by a
fine partition of $\Hyper [a,b]$ and is$,$ indeed$,$ the same for all such 
fine partitions. Extending the impulse to the vector $\vec I,$ then 
$\vec I =\st {\hyper {\vec I}\,} =  \int_a^b \vec F(t)\, dt.\diamondsuit$\par
$\{${\it Second derivation.} First$,$ it should be obvious that there is no complete classical 
counterpart 
to the nonstandard language used in an elemental derivation. Recall$,$ however$,$ 
how it is possible to us a quasi-classical language in the second derivation 
given in section 2.1. 
(i) As is done in Internal Set Theory$,$ the ``*'' notation is removed from the 
functions since whether they  are 
nonstandard extensions of standard functions is clear from the function's 
argument (i.e. preimage).\par
 (ii) The symbols ``$\approx$'' is translated by the term ``infinitely 
close.'' This relation can be physically characterized by stating that no 
standard machine can measure any difference between the quantity on the left 
and the quantity on the right no matter how small the machine error.
\par
(iii) Except for $\nat_\infty,$  hyperreal numbers are usually limited. Hence$,$ 
simply call such a hyperreal number by the single word term ``number.''\par
(iv) We use the fact that functions defined and continuous on $[a,b]$ 
preserve the infinitely close 
concept for these numbers. That is if $t,\ t_1 \in [a,b]$ and $t \approx t_1,$ then $f(t) \approx 
f(t_1).$ Infinitesimals may be called the ``infinitely or very small.'' These 
numbers 
can be physically characterized as measures that are smaller than any 
standard machine error --- measures that appear to a machine to be zero.\par
(v) Rather than use the standard part operator$,$ where applicable$,$ 
use the simple term ``limit'' in its place$,$ since it has the same 
operative properties. Also use the fact 
that limits of two infinitely close numbers are equal.\par   For the elemental 
derivation method some additions to this quasi-classical language are 
necessary. These additions will necessarily be somewhat 
less precise for the concepts of the first-order property and the 
internal object will not be mentioned.\par
(vi) Let the {region $J$ be any of those studied in elementary calculus}.\par
(vii) Hyperline segments are termed ``{infinitely small line segments},'' which 
means line segments of infinitesimal length.\par
(viii) Hyperpolygonal curves may be called ``{infinitesimal polygonal 
curves}'' 
and  
defined as ``polygonal curves with each 
line segment being an  infinitely small segment.''\par
(ix) Call a fine partition of an interval an ``{infinitesimal partition}''. This 
means that it has a $\Gamma + 1$ number of elements that determine $\Gamma$ 
subintervals each of which is infinitely small in length.\par
(x) Hyperfinite sums are termed as ``special finite sums'' 
and behave as if they are finite sums.\par
(xi) Hyperfinite sets are termed as ``special finite sets'' and they 
also behave as if they are finite.\par
(xii) For the special finite sets$,$ or special finite sums$,$ $\Gamma$  
denotes the ``{number'' of elements or terms}$,$ respectively.\par
(xiii) The product of an infinitesimal and a number [assuming limited] is 
an infinitesimal. \par
(xiv) Any nonnegative number less than or equal to an infinitesimal is an 
infinitesimal.\par
$\diamondsuit$ Suppose$,$ as in section 2.5$,$ that uniformly continuous force field 
$F\colon E \to \realp n,$ open $E \subset \realp n,\ F =(F_1,\ldots,F_n).$
For a line segment $\ell^\prime\subset E,$ determine by a linear system of 
equations
over the time interval $[t_1,t_2]$ and for a component $F_i$ equation 
(2) section 2.3 can be re-expressed in terms of an impulse vector
$\vec I = (I_1,\ldots,I_n)$ as $F_i(t_m)(t_2 - t_1) = I_{im} \leq F_i(t)(t_2 - 
t_1) \leq F_i(t_M)(t_2-t_1) = I_{iM}.$ Thus there is some $t_2^\prime \in [t_1,t_2]$ 
such that the actual $i$'th component of the impulse over that line segment is 
$I_i(t_2^\prime) = F_i(t_2^\prime)(t_2 - t_1).$ This result holds for 
any infinitesimal line segment. {\bf Over a finite polygonal line $\ell$ it is assumed 
that 
the total $i$'th component of the impulse is the simple sum of the $i$'th 
component over the individual line segments}. Considering the impulse to be a 
function defined over the set of all finite polygonal curves in $E,$ then this  
$i$'th component of the impulse is defined on an infinitesimal polygonal curve 
${\cal P}_\Gamma$ with $\Gamma$ many sides and its value is the special 
finite sum
$$ I_i({\cal P}_\Gamma) = \sum_{j=1}^\Gamma F_i(t_j^\prime)(t_j - 
t_{j-1}),\eqno (2)$$                                                          
where $\{t_0,\ldots,t_\Gamma\}$ is an infinitesimal partition of $ [a,b]$ 
with infinitely small subintervals that 
generates ${\cal P}_\Gamma.$\par
Now since the above results hold for these subintervals  and 
$t^\prime \in [t_{j-1},t_j]$ for $1 \leq j \leq \Gamma$ then $t^\prime \approx 
t_j.$ Hence$,$ $ F_i(t_j^\prime) \approx  F_i(t_j)$ and from the 
definition of $\approx$ for each $j$ there exists some infinitely small $\eps_j$ 
such that $ F_i(t_j^\prime) =  F_i(t_j) + \eps_j.$ 
Since any finite set of real numbers contains a maximum$,$ then any 
special finite set of numbers contains a maximum. But$,$ 
$\{\eps_1,\ldots,\eps_\Gamma\}$
is a special finite set. Thus there exists some infinitely small 
$\eps = \max \{\eps_1,\ldots,\eps_\Gamma\}.$ Recall that what holds true for 
finite sums holds true for special finite sums. Thus the triangle 
inequality holds. Consequently $$\left| \sum_{j=1}^\Gamma \eps_j(t_j - t_{j-1})\right| \leq \sum_{j=1}^\Gamma 
\vert \eps_j\vert\,(t_j - t_{j-1})\leq \vert \eps \vert \sum_{j=1}^\Gamma 
(t_j - t_{j-1})= \vert \eps \vert\,(b-a)$$
But$,$ since $\vert \eps \vert\,(b-a)$ is infinitely small then
$$\lambda_i=\sum_{j=1}^\Gamma \eps_j(t_j - t_{j-1})\eqno (3)$$ 
is infinitely small. 
Substitution into equation (2) and again using the fact that a special finite 
sum 
behaves like a finite sum yields
$$ I_i({\cal P}_\Gamma) = \left(\sum_{j=1}^\Gamma  
(F_i(t_j)+\eps_j)(t_j -t_{j-1})\right)=$$ $$\left(\sum_{j=1}^\Gamma  F_i(t_j)(t_j - 
t_{j-1})\right) + \lambda_i,\eqno (4)$$\par
Hence$,$ $$ I_i({\cal P}_\Gamma) \approx 
\left(\sum_{j=1}^\Gamma  F_i(t_j)(t_j - 
t_{j-1})\right).\eqno (5)$$
Now $F_i$ being continuous implies from the definition of the integral 
that $\lim \left(\sum_{j=1}^\Gamma  F_i(t_j)(t_j - 
t_{j-1})\right)$ $= \int_a^b F_i(t)\, dt.$  Taking the limit of expression 
(5) yields that $ I_i=\int_a^b F_i(t)\, dt.$ Extending the 
impulse to the vector $\vec I,$ then 
$\vec I =\lim {{\vec I}\,} =  \int_a^b \vec F(t)\, dt.$ $\diamondsuit \}$\par
There are$,$ at least$,$ two ways used to arrive at the relation  
between the change of momentum $\vec p(b) - \vec p(a)$ and the impulse 
$\vec I.$ The first is the usual one of simply stating that 
$\vec I = \int_a^b \vec F(t)\, dt= \int_a^b (d\vec p/dt)\, dt=\vec p(b) - 
\vec p(a)$ from the Fundamental Theorem of Integral Calculus. There is a 
second method utilizing hyperfinite summation that incorporates the fact that 
$F_i$ may be considered as uniformly continuous on $[a,b].$ Since this notion 
is significant for the 
proper intuitive understanding of the underlying physical processes$,$ we 
give an illustration using the momentum vector $\vec p = (p_1,\ldots,p_n).$ 
\par
$\diamondsuit$ Start with equation (2) and substitute $\hyper p_i^\prime(t_j^\prime)$ for 
$ F_i(t_j^\prime).$ This obtains
$$\hyper I_i({\cal P}_\Gamma) = \sum_{j=1}^\Gamma \hyper p_i^\prime 
(t_j^\prime)(t_j - t_{j-1}),\eqno (5)$$ 
Then by the nonstandard mean value theorem  there exists some $h_j^\prime
\in (t_{j-1},t_j)$ such that 
$$\hyper p_i(t_j) - \hyper p_i(t_{j-1}) = \hyper p_i^\prime(h_j^\prime)(t_j - 
t_{j-1}). \eqno (6)$$
But uniform continuity of $p_i^\prime$ on $[a,b]$ implies that there exists some 
$\delta_j \in \monad 0$ such that
$\hyper p_i^\prime(h_j^\prime) = \hyper p_i^\prime(t_j^\prime) + 
\delta_j^\prime.$                                                          
The elemental derivation process then 
yields 
$$\delta_i +\sum_{j=1}^\Gamma \hyper p_i^\prime (h_j^\prime)(t_j - t_{j-1})= 
\sum_{j=1}^\Gamma p_i(t_j^\prime)(t_j - t_{j-1}).\eqno (7)$$
Substituting (6) into (7) and again using the elemental derivation process yields
$$\delta_i + p_i(b) - p_i(a)  = \sum _{j=1}^\Gamma \hyper 
p_i^\prime(t_j^\prime)(t_j - t_{j-1}) = \hyper I_I({\cal P}_\Gamma),\ 
\delta_i,\lambda_i \in \monad 0. \eqno (8)$$
Consequently$,$ $p_i(b) - p_i(a) = \st {\hyper I_i({\cal P}_\Gamma)}$ implies the 
result sought that $\vec p(b) - \vec p(a) = \st {\hyper {\vec I}({\cal 
P}_\Gamma)}= \vec I.\diamondsuit$\par 
$\{${\it Second derivation.} $\diamondsuit$ Modify the second derivation for the general impulse 
integral.$\}\diamondsuit$\par
Once again$,$ the above illustrates that the total impulse is not 
dependent upon the hyperpolygonal path along which the object moves$,$ but is 
simply the change in momentum. However$,$ by physical intuition$,$ an object has 
changed its momentum by traversing a physical path $C$ and we are using a 
hyperpolygonal representation for such a path of motion. For any two 
such hyperpolygonal curves ${\cal P}_\Gamma,\ {\cal P}_\Lambda,$ it follows 
from the hypotheses that $\hyper {\vec I}({\cal P}_\Gamma) \approx  \hyper 
{\vec I}({\cal P}_\Lambda).$ Now the basic geometric measure for $C$$,$ the 
length$,$ is from the derivation in the Infinitesimal Modeling manual completely determined 
by the *-length of each of these hyperpolygonal representations. Therefore$,$ it 
seems appropriate to consider the standard impulse over the curve $C$ to be 
the unique effect $\st {\hyper {\vec I}({\cal P}_\Gamma)},$ where ${\cal 
P}_\Gamma$ is any hyperpolygonal representation  for $C.$\par
Why does the elemental derivation process work and what is it indicating about 
integral styled quantities? All of the basic quantities in any expression 
prior to considering their hyperfinite sum must not only be infinitely close 
but must be$,$ in the above case$,$ infinitely close of order one. The hyperfinite 
summation of infinitesimals need not be infinitesimal or even limited; but$,$ 
depending on the type of limited objects this type of special infinite 
closeness utilizes$,$ then the elemental derivation process does imply that 
this particular hyperfinite sum of infinitesimals is infinitesimal. 
In order to guarantee that this is the case$,$ strong hypotheses where required 
for the functions involved. Except for 
the possibility of restricting the hyperfinite summation to special sets of 
partitions$,$ such as in the case of the gauge type integral discussed in section 
8.9 of the Infinitesimal Modeling manual$,$ it$,$ at present$,$ has not been possible to obtain 
rigorous derivations for integral expressions without these strong 
hypotheses.\par 
\vfil\eject
\centerline{Chapter 3.}
\medskip    
\centerline{\bf SLIGHTLY LESS BASIC MECHANICS}              
\bigskip 
\leftline{3.1 \underbar{Mass}}\par
\medskip
Since the {Infinite Sum Theorems}$,$ the rules {IR1--IR6} and  the {Self-Evident 
Theorems} that appear in the Infinitesimal Modeling manual$,$ have not been known to the 
physicist previously then almost all of the elementary definitions or derivations that 
involve infinitesimal quantities and the integral have relied upon  {elemental 
methods}. These 
methods refer directly to the ``{particle}''$,$ point charges$,$ point 
masses and other such notions. This process obviously forces the measure 
to have the same properties as the integral and the converse of starting with 
a functional$,$ considering fundamental properties and showing that such a 
functional must be a measure as defined by some integral need not be 
considered. The mathematician might find this converse approach as appealing. 
But$,$ it is more of a global view for a particular scenario and elementary 
instruction in the physical sciences tends to force upon a student an 
atomistic view in the sense that complex observed behavior is 
conceived of as the effects produce by minuscule objects. This is the 
prevalent textbook approach.\par
The {Self-Evident Theorems} that appear in the Infinitesimal Modeling manual require$,$ at least$,$ 
three global assumptions as well as a strong additivity property for the 
measure under consideration.  
The simplifications for regions that have Jordan-content$,$ 
Propositions 6.7 and 6.10 as they 
appear in appendix 6$,$ are very easy to apply in the laboratory setting since 
one need only investigate these properties on a rectangular interior to the 
region. 
However$,$ this does require the assumption that the measure being considered is 
Jordan-like. Theorem 7.2.2 in the Infinitesimal Modeling manual implies that the 
integral over a {Jordan-measurable} region is an {ordinary Jordan-like 
measure}.
Further$,$ all of the customary geometric regions used within elementary 
physics are all compact Jordan-measurable sets and these simplified 
self-evident theorems may be applied. But the technical difficulties of 
establishing these simplified 
self-evident theorems should not be underestimated as amply 
illustrated by the formal proofs that appear in appendix 6 of the 
Infinitesimal Modeling manual. These proof 
methods are well beyond almost 
all undergraduates who require even advance physics courses as part of their 
curriculum.\par 
 In the 
Infinitesimal Modeling manual the Infinite Sum Theorem is the only approach used to obtain 
the mass of an object as determined by a continuous density function and 
this approach does not appear in any present day physics textbook. 
Consequently$,$ it is obvious that$,$ for the present$,$ the {elemental method} with the 
notion of hyperfinite summation is more appropriate. As discussed towards the 
end of section 2.6$,$ unless we restrict integral modeling to rectangular 
regions then$,$ even within the elemental derivation$,$ it is necessary to assume 
that the quantity being measured is Jordan-like. A brief review of this 
concept is useful.\par
Let $R$ be an n-dimensional rectangle and $R \subset \realp n.$  
Suppose that Jordan-measurable  
$J \subset R\subset \realp n,$ where to avoid trivialities we always assume 
that the interior of $J$ is nonempty. Let $\cal P$ be a set of partitions of $R,$ 
where$,$ as usual$,$ if $P \in {\cal P},$ then each $S \in P$ is but an 
n-dimensional subrectangle in $R.$ 
Define on ${\cal P}$ the following map ``inn.''
For each $P \in {\cal P},$ let ${\rm inn}(P) = \{S\bigm\vert S \in P\  
{\rm and}\ S\subset J\}.$ Of course$,$ the set ${\rm inn}(P)$ is called a set of {\bf {inner 
subrectangles.}} Let $I_{\cal P} = \{\cup A\bigm \vert 
\exists P(P \in {\cal P} \land 
A = {\rm inn}(P))\}\bigcup \{S\bigm \vert \exists P(P \in {\cal P} \land S 
\in {\rm inn}(P))\}.$                                                                     
 When considering the nonstandard extension of $\rm inn$ to a {fine partition}
$Q$ we usually do not write this as  $\Hyper {\rm inn}(Q)$ but rather retain 
the original notation ${\rm inn}(Q).$ Since  any fine partition $Q$ is 
hyperfinite and the set ${\rm inn}(Q)$ is an internal subset of $Q$ 
then it is also hyperfinite. For simplicity$,$ as far as the functional 
(i.e. measure) $B\colon \{I_{\cal P},J\} \to \real$ is concerned$,$ 
when the elemental derivation 
method is used it is assumed that $B$ is$,$ at least$,$ {simply additive} on each 
${\rm inn}(P),$ when $P \in {\cal P}.$ This means that if $P \in {\cal P}$ and
$A = {\rm inn}(P),$ then $B(\cup A) = \sum _{S\in A} B(S).$\par
In general$,$ if a generating function for the above functional $B$ is constant
over the entire Jordan-measurable $J,$ then  it is necessary that 
$B$ display certain Jordan-measure properties. In particular$,$ if $v$ denotes 
the Jordan-measure$,$ then for any fine parition $Q$ of $\hyper R$$,$ it follows from   
Theorem 7.2.2 
that $v(J) = \st {\sum_{S \in A} \hyper v(S)} = \st {\hyper v(\cup S)}.$ Thus,
$v(J) \approx \hyper v(\cup A).$ Of course$,$ $\hyper v(S)$ is but the product 
of the infinitesimal lengths of the sides of the $n$-dimensional subrectangle
$S.$ Bearing in mind the basic properties of Jordan-measure$,$ we assume 
that for a fine partition $Q$ that the functional $B$ has the property that
$B(J) \approx \hyper B(\cup A),$ where $A = {\rm inn}(Q).$ This property is 
termed the {\bf {Jordan-like}} property and a standard criterion for the 
functional $B$ to be Jordan-like for all fine partitions is given in the Infinitesimal Modeling manual. Of course$,$ these technical matters need made be 
discussed in most elementary courses. {\sl Indeed$,$ in the second derivations one 
only needs  to state that the region $J$ is one of those studied in the elementary 
calculus.}\par
$\diamondsuit$ Let compact and Jordan-measurable $J \subset R \subset \realp 3$ 
and $S$ be an inner rectangle as defined 
above. Let $\rho \colon J \to \realp 3$ be a continuous {density function}. 
If $\rho (x,y,z)$ is 
constant over a subrectangle $S$ of $J,$ then from the definition of 
$\rho$ the mass $m(S) =
\rho\, v(S),$ where $v(S)$ is the volume of $S.$ Suppose that $\rho$ is not 
constant. Then there exists $(x_m,y_m,z_m)\in S$ and $(x_M,y_M,z_M)\in S$ such 
that $\rho_m =\rho (x_m,y_m,z_m)\leq \rho (x,y,z) \leq \rho (x_M,y_M,z_M)=\rho 
_M$ for each $(x,y,z) \in  S.$ Hence$,$ if the mass of $S= m(S),$ then 
observation and the above property for constant density indicates that$,$ 
at least$,$ from the macroscopic viewpoint 
$$\rho_m\, v(S) \leq m(S) \leq \rho_M\ v(S).\eqno (1)$$\par
The Intermediate Value Theorem for connected subsets of $R$ 
yields that there is some $(x_1,y_1,z_1) \in S$ 
such that $m(S) = \rho (x_1,y_1,z_1)\, v(S).$ If we have more than one such 
inner rectangle$,$ say $\{S_1,\ldots ,S_n\},$ than observation indicates that 
the total mass$,$ $m(A),$  of the configuration $A =S_1 \cup \cdots \cup S_n $ is the finite 
sum 
$$m(A) = \sum_{j=1}^n m(S_j) = \sum_{j=1}^n \rho (x_j,y_j,z_j)\, 
v(S_j).\eqno (2)$$                                             
By *-transfer$,$ these results hold for infinitesimal subrectangles of $J$  and 
hyperfinite summation. Now let $Q$ be any fine partition of $\hyper R.$ Then ${\rm 
inn}(Q)$ is an internal hyperfinite set of infinitesimal subrectangles. 
Hence,
$$\hyper m(\cup \{S_j\vert S_j \in {\rm inn}(Q)\}) = \sum_{j=1}^\Gamma \hyper {\rho} 
(x_j,y_j,z_j)\, \hyper v(S_j).\eqno (3)$$ \par
Even though each point $\vec v_j =(x_j,y_j,z_j) \in S_j,$ this does not imply that the 
intermediate partition $\{\vec v_j\bigm \vert 1 \leq j \leq \Gamma \}$ is 
internal. However$,$ letting $\vec c_j$ denote the corner of $S_j$ nearest or 
equal to $(0,0,0)$ then $\{\vec c_j \bigm \vert 1 \leq j \leq \Gamma \}$ is an 
internal intermediate partition. Uniform continuity of $\rho$ on $J$ yields 
that for each $j$ there exists some $\eps_j \in \monad 0$ such that
$\hyper {\rho} (x_j,y_j,z_j) = \hyper {\rho} (\vec c_j) + \eps_j.$ Hence,
$$\hyper m(\cup \{S_j\vert S_j \in {\rm inn}(Q)\}) = \sum_{j=1}^\Gamma 
(\hyper {\rho} (\vec c_j)+ \eps_j)\, \hyper v(S_j).\eqno (4)$$ 
The elemental derivation process yields
$$\hyper m(\cup \{S_j\vert S_j \in {\rm inn}(Q)\}) \approx \sum_{j=1}^\Gamma 
(\hyper {\rho} (\vec c_j)\, \hyper v(S_j).\eqno (5)$$
Since $\rho$ is integrable over $J$ then taking the standard part and applying 
Theorem 7.2.2 of the Infinitesimal Modeling manual one obtains
$$\st {\hyper m(\cup \{S_j\vert S_j \in {\rm inn}(Q)\})} = \int_J \rho (x,y,z)\, dX.\eqno 
(6)$$\par
As mentioned above$,$ one final step is required. It must be assumed that the concept 
of the mass of an object has the same Jordan-like quality as does the concept 
of the volume. Making this last assumption implies that $m(J) = \int_J \rho 
(x,y,z)\, dX.\diamondsuit$ \par
A physical interpretation of the Jordan-like quality of the mass is also 
possible. For every fine partition $Q,$ the NSP-world effects of the mass of a 
simple internal rectangular configuration $\cup \{S_j\bigm \vert S_j \in {\rm inn}(Q)\}$
is indistinguishable from the effects of the mass of $J$ at the first (order) 
level$,$ if no comparison is to be made with other such measures. When such a 
comparison is made$,$ as 
pointed out in the Infinitesimal Modeling manual$,$ Section 8.3$,$ the effects of
$\hyper m(\cup \{S\bigm \vert S \in {\rm inn}(Q)\})$ and $m(J)$ would be 
infinitely close of order $3$ and indistinguishable at that level.
Of course$,$ all that is being derived for the case $\realp 3$ holds  
for $\hyper m(\cup \{S\bigm \vert S \in {\rm inn}(Q)\})$
and any $\realp n,\ n\geq 1.$\par
\bigskip
\hrule
\smallskip
{\sl In all of the previous derivations$,$ equations of form (1) are the 
essential modeling requirements. The effect of the elemental derivation process 
is to eliminate the additional additivity property for such functionals $m.$ 
However$,$ if such additivity is assumed$,$ then the integral expression follows 
immediately from the self-evident theorems that appear in the Infinitesimal Modeling manual.}
\smallskip
\hrule
\bigskip \par
$\{${\it Second derivation.} Now that we have the essential procedures needed 
to modify a 
rigorous infinitesimally styled derivation and obtained a quasi-classical one,
it is not necessary to present an entire second derivation for 
equation (6). It is useful to conjoin the following to our list of 
alterations. The set $J$ can simply be described as one of the regions 
studied in the elementary calculus and for which an integral expression for 
volume is obtained. For the necessity of expressing ${\rho} 
(x_j,y_j,z_j)$ in equation (3) as ${\rho} (x_j,y_j,z_j) 
= {\rho} (\vec c_j) + \eps_j,$ without invoking the concept of the 
internal intermediate partition$,$ one can simply argue that it is necessary 
to evaluate $\rho$ at a known partition generated point rather than at a less 
explicitly known $(x_j,y_j,z_j)\in S.$ The requirement that $m$ be 
Jordan-like can be described as a common feature that $m$ must share with the 
volume of the region $J.$ The feature being$,$ that for any infinitesimal 
partition the configuration 
composed of all of the infinitely small subrectangles contained within $J$ 
must have its volume infinitely close the volume of $J.$ This same feature 
holds for the mass and other such measures since if $\rho \equiv 1$ on $J,$ 
then $v(J) = m(J).$ I leave to the reader the actual 
construction of the second derivation.$\}$\par
\bigskip
\leftline{3.2 \underbar{Moments and Center of Mass}}
\medskip \par
It is often difficult to decide whether it is more significant to derive  
a specific mathematical model 
completely from fundamental physical observations and the defining properties 
of the mathematical structure or to replace portions of the derivation with 
well-known theorems gleaned from the abstract structure itself. It appears 
that within modern theoretical physics$,$ abstract mathematical results are 
introduced as soon as practical even though they may be couched in a 
quasi-physical language. This procedure follows the routine assumption that the 
rigorous logic displayed when a proposition is proved abstractly is equivalent 
to the logic required for a much longer and more complex derivation that 
utilizes but the fundamental correspondence between the defining properties of 
the structure and the physical terms that describe the physical scenario. The 
following theorem is useful in 
order to illustrate the economy achieved by such an early introduction of well-known mathematical results.  \par
\vskip 18pt
{\bf Theorem 3.2.1. {A Weighted Mean Value Theorem.}} {\sl Let $J \subset 
 E \subset R \subset \realp n,\ J$ be Jordan-measurable$,$ and $E$ be 
compact and connected. Let continuous $f\colon E \to \real,$ 
integrable $g\colon J\to \real$ and $g(\vec v) \geq 0$ for each $\vec v \in J.$
Then there exists some $\vec v_0 \in E$ such that $\int_Jf(\vec v)g(\vec v)\, 
dX = f(\vec v_0)\int_Jg(\vec v)\, dX.$}\par
\vskip 18pt
Even though the elemental derivation process is being stressed throughout this 
physics manual$,$ it is worthwhile to once again mentioned that the Self-Evident 
theorems of the Infinitesimal Modeling manual are always available. The Method of 
Constants$,$ as 
well as the Maximum-Minimum Method$,$ is actually exemplified within the elemental 
process derivations. Immediately following equation (1) of this section is 
the fact that there does exit some $(x_1,y_1,z_1) \in S$ such that 
$m(S) = \rho (x_1,y_1,z_1) v(S).$ This is the explicit requirement for 
application of Proposition 6.9 in Appendix 6 assuming the requisite 
additivity properties for the measure $m.$ If the measure generating 
function is composed of the product of two or more nonconstant functions$,$ then 
the Extended Self-Evident Method of Constants Proposition 6.10 may be 
appropriate.\par
$\diamondsuit$ Notwithstanding our discussion in the above paragraph$,$ let's consider the 
a {moment} generating function $M_{(\cdot)}$ defined as follows: Let nonnegative 
continuous $\rho (x,y,z) \colon J \to \real,$ where $J \subset R \subset \realp 3$ is 
compact and Jordan-measurable. A moment function defined on a 
partition subrectangle $S \subset J$ is $M_{yz}(S) = x\, \rho (x,y,z)\, v(S),$ 
where $x = f(x,y,z)$ is continuous on $\realp 3.$   In the same manner$,$ 
define $M_{xz}(S) =y\, \rho (x,y,z)\, v(S),\         
 M_{xy}(S) = z\, \rho (x,y,z)\, v(S).$ Although it may appear to be 
sufficient to 
consider the point $(x,y,z)$ as an arbitrary member of $S,$ certain special 
selections are necessary. For example,
if $\rho $ is constant$,$ then $(x,y,z)$ might be  selected as  
the center of $S.$  It is well-known that if we consider 
a small 
enough homogeneous 
rectangular solid $S,$ then the  gravitational field of the Earth in relation 
to $S$ is effectively a parallel vector field and the acceleration of gravity 
is a constant $g.$ The center of the rectangle is 
the point of rotational stability within such a field --- the so-called 
{center of gravity}. This is demonstrated in the customary manner by 
considering the expressions $g\, M_{xy}(S),\ g\, M_{yz}(S),\ g\, M_{xz}(S).$
Thus from the Newtonian gravitational point of view$,$ the {moment function} can 
be viewed as a measure of the rotational effect with 
respect to the coordinate planes within such a gravitational field. Notice that 
the necessary observations that lead to this conclusion are relative to the 
actual small size of the objects.\par
 For the general case  
of continuous $\rho$ on $S,$ there exist
 $x_m\rho (x_m, y_m,z_m)$ and 
$x_M\rho (x_M,y_M,z_M),$ where $(x_m,y_m,z_m),$ $(x_M,y_M,z_M) \in S$ such that
$$x_m\rho (x_m, y_m,z_m)v(S)\leq x\, \rho (x,y,z)v(S)\leq x_M\rho 
(x_M,y_M,z_M)v(S).\eqno (1)$$
The usual assumption is now imposed upon our problem. Suppose that the actual moment 
effect $M_{yz}(S)$ lies somewhere between these two extremes. Thus there 
would exist some $(x_1^\prime,y_1^\prime,z_1^\prime) \in S$ such that  
$$M_{yz}(S) = x_1^\prime\, \rho (x_1^\prime,y_1^\prime,z_1^\prime)\, v(S). \eqno 
(2)$$\par
Equation (2) is now extended to a finite system $\{S_j\bigm \vert 1 \leq j 
\leq n\}$ of partition inner subrectangles using the apparent experiential result that the 
total moment effect of the system is the sum of the individual effects. Hence,
$$M_{yz}(\cup S_j) = \sum_{j=1}^n x_j^\prime\, \rho 
(x_j^\prime,y_j^\prime,z_j^\prime)\, v(S_j).\eqno (3)$$\par
By *-transfer$,$ we assume that within the NSP-world the above behavior holds 
for infinitesimal rectangles. Since this is a modeling technique it is not 
necessary to assume that such infinitesimal rectangles exist in some type 
of objective reality. {\sl However$,$ it is possible to describe such behavior 
within a {substratum NSP-world} by considering the simplistic behavior of 
a hyperfinite 
set of NSP-world infinitesimal rectangles as a {superstructure of objects} 
that controls the behavior of a
corresponding system of ``small'' natural world rectangular objects.} As in 
the least section$,$ assuming that $Q$ is a fine partition of $\hyper R$ this leads to the conclusion that 
$$\hyper M_{yz}(\cup\{S_j\vert S_j \in {\rm inn}(Q)\})
= \sum_{j=1}^\Gamma x_j^\prime\, \hyper 
{\rho}(x_j^\prime,y_j^\prime,z_j^\prime)\, 
\hyper v(S_j).\eqno (4)$$ \par
Letting $(x_j,y_j,z_j)$ denote the corner of $S_j$ nearest or 
equal to $(0,0,0)$ then $\{(x_j,y_j,z_j)\bigm \vert 1 \leq j \leq \Gamma \}$ is an 
internal intermediate partition. Uniform continuity of $x\rho$ on $J$ yields 
that for each $j$ there exists some $\eps_j \in \monad 0$ such that
$x_j^\prime\hyper {\rho} (x_j^\prime,y_j^\prime,z_j^\prime) 
=x_j \hyper {\rho} (x_j,y_j,z_j) + \eps_j.$ Applying the elemental 
derivation process one obtains
$$\st {\hyper M_{yz}(\cup \{S_j\vert S_j \in {\rm inn}(Q)\})} = \int_Jx\, \rho 
(x,y,z)\, dX.\eqno (5)$$ 
The additional requirement that the measure of the moment is Jordan-like 
yields $M_{yz}(J) = \int_Jx\, \rho (x,y,z)\, dX.$ \par
Now  $R\subset 
\realp 3$ is compact and connected and $f(x,y,z,) = x$ is continuous on $R.$ 
Thus$,$ by the Weighted Mean Value Theorem there exists some real 
$\overline 
x$ such that
$$M_{yz}(J) = \int_Jx\, \rho (x,y,z)\, dX= \overline{x} \int_J\, \rho (x,y,z)\, dX =
\overline{x}\, m(J). \eqno (6)$$
Repeating the above argument there exist real $\overline{y}$ and $\overline{z}$ such 
$M_{xz}(J) = \int_Jy\, \rho (x,y,z)\, dX= \overline{y} \int_J\, \rho (x,y,z)\, dX =
\overline{y}\, m(J)$ and $M_{xy}(J) = \int_Jz\, \rho (x,y,z)\, dX= \overline{z} \int_J\, \rho (x,y,z)\, dX =
\overline{z}\, m(J).$ Consequently$,$ as far as the moment effects are concerned the 
object $J$ can be consider as represented by the single point $(\overline{x},
\overline{y},\overline{z})$ with the mass number $m(J)$ attached to 
it.$\diamondsuit$\par
$\{${\it Second derivation.} From this point on in this Elementary Physics manual$,$ the 
second quasi-classical derivation will not be given unless it is substantially 
different from our previous examples.$\}$\par
Notice that the above derivation did not start with the concept of the point 
masses and then derive the integral expression for the center of mass. Rather$,$ 
we derived the concept by means of infinitesimal analysis. Additionally$,$ the 
statement that the center of mass is equivalent to the center of gravity 
appears in the above discussion to depend upon the parallel gravitational 
field concept. In the next 
section$,$ we show that the idea of the less substantiated {point masses}$,$ 
{\sl if viewed from the NSP-world,} does lead to the same center of mass 
conclusion.\par
\bigskip
\leftline{3.3 \underbar{Point Masses}}
\medskip \par
In {\tt Tipler [1982]}$,$ the concept of the {point mass} is used to develop the 
{center of gravity} and {center of mass} for such objects$,$ assuming that 
this technique has 
been justified. On page 229 of volume 1$,$ Tipler states: {\it``If the 
center-of-mass coordinates of a continuous body are to be calculated$,$ the sum 
$\sum m_ix_i$ must be replaced by the integral $\int x\, dm,$ where $dm$ 
is an element of mass.''} In the text by {\tt Young$,$ Riley$,$ McConnell$,$ Rogge [1974$,$ p. 281]}$,$ when 
{moments of 
inertia} are discussed$,$ once again the student is instructed that such 
a measure is given by an integral over $dm.$ No further explanation is given 
as to why this particular technique is justified. This vague modeling 
technique can be justified within the NSP-world once equations such as (6) of 
Section 3.2 have been derived.\par
$\diamondsuit$ Assume the hypotheses used to derive (6) of Section 3.2 and  let $Q$ 
once again be a fine partition of $\hyper R.$ Let $S_j\in {\rm inn}(Q)$ Then there 
exists some $(x_j^{\prime\prime},y_j^{\prime\prime},z_j^{\prime\prime})\in 
S_j$ such that $\hyper m(S_j) = \hyper {\rho}(x_j^{\prime\prime},y_j^{\prime\prime},
z_j^{\prime\prime})\, \hyper v(S_j).$ Now $\hyper {\rho}(x_j^{\prime\prime},y_j^{\prime\prime},
z_j^{\prime\prime})\approx \hyper {\rho}(x_j^{\prime},y_j^{\prime},
z_j^{\prime}),$ where $(x_j^{\prime},y_j^{\prime},
z_j^{\prime})$ is as described in the derivation of Section 3.2. 
Hence,
 $x_j^\prime \hyper {\rho}(x_j^{\prime},y_j^{\prime},
z_j^{\prime}) \approx x_j^\prime \hyper {\rho}(x_j^{\prime\prime},y_j^{\prime\prime},
z_j^{\prime\prime}),$ for $x_j^\prime$ is limited. The elemental derivation 
process yields 
$$M_{yz}(J) =\int_Jx\, \rho (x,y,z)\, dX = \St {\left(\sum_{j=1}^\Gamma x_j^\prime
\hyper {\rho} 
(x_j^{\prime\prime},y_j^{\prime\prime},z_j^{\prime\prime})\, \hyper 
v(S_j)\right)}.\eqno (1)$$
Consequently$,$ letting $dm_j = \hyper m(S_j) = \hyper {\rho} 
(x_j^{\prime\prime},y_j^{\prime\prime},z_j^{\prime\prime})\, \hyper v(S_j)$ 
denote that mass of the infinitesimal subrectangle (not a point!) and repeating the argument 
for the other two moments one can describe the moment effects within the 
standard world as follows: The effect is indistinguishable from the effect 
of a hyperfinite sum of mass numbers attached to the points 
$(x_j^\prime,y_j^\prime,z_j^\prime).$ Thus the points 
$(x_j^\prime,y_j^\prime,z_j^\prime)$ within the NSP-world can be viewed as 
point masses. Obviously$,$ the point $(x_j^\prime,y_j^\prime,z_j^\prime)$ is 
not unique since it may be replaced by any $(x,y,z) \in S_j.\diamondsuit$ \par
The reader might be inclined to attempt to *-transfer the above derivation to 
the standard world and arrive at the conclusion that there exists a set 
of point masses that would yield the moment effects expressed by the 
integrals. This 
would be an error$,$ however$,$ since the standard part operator and as well as
$\approx$ are external concepts. The technique of *-transfer$,$ at this stage$,$ 
only allows the hyperfinite sum of point masses to be transferred   
into a statement about finitely many point 
masses the sum of the moments of which would approximate the moment effect
within an given positive $\eps$ for all the partitions of $R$ with mesh less 
than some positive $\delta.$ What this implies is that rather than accepting 
an ad hoc modeling technique that utilizes unrealistic standard world 
point masses to derive the integral expressions for the moment effects$,$ 
it may be more conducive to student comprehension to employ the NSP-world
point masses since their use can be more rigorously justified.
{\sl However$,$ in certain cases once infinitesimal analysis has established 
equation (6) of 
Section 3.2$,$ and the like$,$ then \underbar{standard means} can be applied to 
investigate an effective center of mass for a finite collection of 
objects.} This we do next. \par 
$\diamondsuit$ Suppose that $\{J_1,\ldots,J_p\}$ is a nonempty {\sl finite} 
set of   
of pairwise disjoint$,$ compact$,$  and Jordan-measurable subsets of 
$R \subset \realp 3.$ 
Further$,$ let nonnegative continuous $\rho_i \colon J_i \to \real,$ for 
each $i$ such that $1\leq i\leq p.$ The function $f(x,y,z) = x$ is obviously 
continuous on compact$,$ connected $R.$  The piecewise well-defined function
$h(x,y,z) = \rho_i(x,y,z);\ (x,y,z) \in J_i,\ 1\leq i\leq p$ is continuous on $J = 
J_1 \cup\cdots \cup J_p.$ [Since $\realp 3$ is a normal topological space it 
follows that if $\vec v \in J_i,$ then $\monad {\vec v} \cap J_j = \emptyset,
\ i \not= j.$] Now the set $J$ is compact; hence$,$ closed and bounded. 
Consequently$,$ $h$ is a nonnegative bounded function defined on $J.$ 
The function $f(x,y,z) = x$ is continuous on compact$,$ connected $R$ and 
integrable on $J.$ Consider the moment effect generating function $m_{yz} = 
x\, h(x,y,z).$ Then from the Weighted Mean Value Theorem there exists 
$\overline{x}_j$ such that the moment effect 
$$M_{yz}(J) = \int_Jm_{yz}\, dX = \overline{x}_j\, m(J). \eqno (2)$$
However$,$ 
$$\int_Jm_{yz}\, dX= \sum_{i=1}^p\int_{J_i}x\rho_i(x,y,z)\, dX = 
\sum_{i=1}^p\overline{x}_i\, m(J_i).\eqno (3)$$
Thus 
$$M_{yz}(J) = \sum_{i=1}^p \overline{x}_i\, m(J_i). \eqno (4)$$
Repeating the above derivation yields similar equations as (4) for the other 
two moments. However$,$ I repeat$,$ once again$,$  that this approach is only 
relative to a nonempty 
finite set of disjoint$,$ compact$,$ and Jordan-measurable subsets of 
$R.$ $\diamondsuit$
\bigskip
\leftline{3.4 \underbar{Standard Rules and the Elemental Derivation Process}}
\medskip \par
Although the elemental derivation process is very appealing to the intuition$,$ 
what happens when this process is viewed as a mathematically stated theorem?
If you were to analyze the standard hypotheses needed to model this process$,$ 
then what would be obtained is the rule IR5 as it appears in the 
Infinitesimal Modeling manual. This rule coupled with the hypotheses stated in 
a theorem  such as Proposition 6.7 in the Infinitesimal 
Modeling manual leads to a formal theorem that can be established not by an 
infinitesimal sum theorem but be the elemental derivation method. We develop 
such a theorem next --- a theorem that allows us to eliminate the actual 
elemental derivation process.  
The necessary notation for what follows is defined in this manual. \par 
We will not state what comes next as a formal theorem but state it somewhat 
informally. We point out that it is but a restatement in a slightly expanded 
form of Proposition 6.7 in the Infinitesimal Modeling manual.  
Let $\cal P$ be any  set of partitions of the rectangle $R\subset 
\realp n$ and compact Jordan-measurable $J \subset R.$ Suppose that $\hyper 
{\cal P}$ 
contains a fine partition. Of course$,$ if $\cal P$ is the set of all simple 
partitions of $R,$ then such a fine partition exists. Next let $B$ be a real 
valued 
function(al) defined  on $\{I_{\cal P},J\}$ and$,$ at least$,$ additive on 
$\{I_{\cal P},J\}.$ Let continuous $f\colon J \to \real.$ Suppose that for any 
$P \in \cal P$ and any $S \in {\rm inn}(P),$ it follows that $(f_m)v(S) \leq B(S) \leq 
(f_M)v(S),$ where $(f_m)$ [resp. $(f_M)$] is the minimal [resp. maximal] value 
of $f$ on $S$ and $v(S)$ is the Jordan-measure of $S$ (i.e. its simple 
volume). Then if $B$ has the ordinary Jordan-like property$,$ it follows (from the 
elemental derivation process) that $B(J) = \int_Jf(\vec{x})\, dX.$ \par
What the last paragraph signifies is that in all cases where elementary 
physical measures are concerned one needs only argue for the acceptance of the 
stated hypotheses. Once such hypotheses are accepted as reasonable$,$ then the 
conclusion follows from both the infinite sum theorem or the elemental 
derivation process. \par
The assumption that $B$ has the ordinary Jordan-like property is not difficult 
to accept. A standard criterion appears in the Infinitesimal Modeling 
manual and that property is apparently necessary in order for $B$ to be obtained 
by an integral. 
This comes from the fact that the integral itself when viewed as a 
functional satisfies this property. 
It intuitively signifies that an approximation for the value of $B(J)$ that 
is better that any 
machine error can be obtained by considering  the value $B(C),$ where 
$C\subset J$ is a configuration composed of subrectangles taken from a 
partition with  ``small enough'' mesh. Now  compare this with the 
requirements of IR3 (1) in the Infinitesimal Modeling manual and Theorem 6.2.4, where 
the ordinary Jordan-like property is not assumed. In this case$,$ 
the maximal--minimal assumption is weakened for boundary rectangles. This 
weakening is relative to the value of the functional as extended to the 
boundary rectangles. In an elementary exposition$,$ it may be more reasonable 
to   ``build'' configurations such as $C$ and accept the ordinary 
Jordan-like property for such configurations$,$ then to alter the intuitive 
acceptance of such maximal-minimal statements as (1) on page 22 of this manual. \par 
Thus far$,$ we have needed to include the strong requirements that the 
function $f$ be continuous on $J$ and that $J$ be$,$ at least$,$ compact and 
Jordan-measurable. Can either or both of these requirements  be 
relaxed and an acceptable derivation method developed for integral models? An 
answer to this question will depend upon what one considers as ``acceptable'' 
and the areas of application. We will attempt to answer this question in 
later sections of this manual. \par
\bigskip
To be continued by properly trained members of the physics community.
\vfil\eject
\centerline{REFERENCES}
\bigskip
\noindent {\bf Gillispie$,$ C.C.,} [1960]$,$ {\sl The Edge of Objectivity,} 
Princeton University Press$,$ Princeton$,$ NJ.\pars 
\noindent {\bf Herrmann$,$ R. A.,} [1989]$,$  Fractals and ultrasmooth microeffects$,$ J. 
Math. Physics$,$ {\bf 30} (4): 805-808.\pars
\noindent {\bf Newton$,$ I.}$,$ [1934]$,$ {\sl Mathematical Principles of Natural 
Philosophy}$,$ (Translated by Cajori)$,$ University of Cal. Press$,$ Berkeley. 
\pars
\noindent {\bf Tipler$,$ P. A.,} [1982] {\sl Physics,} Worth Publishers$,$ Inc.$,$ 
New York. \pars
\noindent {\bf Young$,$ Riley$,$ McConnell,} [1974] {\sl Essentials of Mechanics,} 
The Iowa State University Press$,$ Ames$,$ Iowa. \pars
\vfil\eject
\centerline{\hskip -1.50in{\bf Additional Special Symbols}}
\centerline{\hskip -1.50in{(Alphabetically listed by first symbol letter.)}}
\vskip 18pt
\+Symbol\dotfill&Name$,$ if any\dotfill&Page no.\cr
\+$\sim_1$\dotfill&Infinitely Close of&\cr
\+&Order One\dotfill&170\cr
\+$\sim_2$\dotfill&Infinitely Close of&\cr
\+&Order two\dotfill&172\cr
\+$\vec I$\dotfill&Impulse Vector\dotfill&177\cr
\+$m(S)$\dotfill&Mass of $S$\dotfill&182\cr
\+$M_{(\cdot)}$\dotfill& Moment Function\dotfill&184\cr
\+$dm$\dotfill&Mass Of&\cr
\+&Infinitesimal $S=$&\cr
\+&Element Of Mass\dotfill&186\cr
\vfil\eject
\catcode`@=11 % from plain.tex
\newdimen\pagewidth \newdimen\pageheight \newdimen\ruleht
% These values were modified:
   \hsize=6.5in  \vsize=54pc  \maxdepth=2.2pt  \parindent=2pc
   \hoffset=.35in
\pagewidth=\hsize  \pageheight=\vsize  \ruleht=.5pt

%This routine is used by \output; this is different from 
%  the one found in App. E since some are not needed here.
\def\onepageout#1{\shipout\vbox{\offinterlineskip
  \vbox to \pageheight {\makeheadline
          #1 % the content of page
          \makefootline \boxmaxdepth=\maxdepth}}
  \advancepageno}

\output{\onepageout{\unvbox255}}
\newbox\partialpage
\def\begind{\begingroup
  \output={\global\setbox\partialpage=\vbox{\unvbox255\bigskip}}\eject
  \output={\doublecolumnout} \hsize=2.85in \vsize=109pc} %2.85
  % Again the sizes have been changed  !!

\def\doublecolumnout{\splittopskip=\topskip \splitmaxdepth=\maxdepth
  \dimen@=54pc \advance\dimen@ by-\ht\partialpage
  % Need to change the value of \dim@ also...
  \setbox0=\vsplit255 to\dimen@ \setbox2=\vsplit255 to\dimen@
  \onepageout\pagesofar
  \unvbox255 \penalty\outputpenalty}
\def\pagesofar{\unvbox\partialpage
  \wd0=\hsize \wd2=\hsize \hbox to\pagewidth{\box0\hfil\hfil\box2}}
\def\balancecolumns{\setbox0=\vbox{\unvbox255} \dimen@=\ht0
  \advance\dimen@ by\topskip \advance\dimen@ by-\baselineskip
  \divide\dimen@ by2 \splittopskip=\topskip
  {\vbadness=10000 \loop \global\setbox3=\copy0
    \global\setbox1=\vsplit3 to\dimen@
    \ifdim\ht3>\dimen@ \global\advance\dimen@ by1pt \repeat}
  \setbox0=\vbox to\dimen@{\unvbox1}
  \setbox2=\vbox to\dimen@{\dimen2=\dp3 \unvbox3\kern-\dimen2 \vfil}
  \pagesofar}

\begind
                       %but it is slight more spread out.
\def\sub#1{{\leftskip=0.2in \noindent #1 \par}\par}
\noindent Index for pages 1-153.\par \noindent Location is within $\pm$ 1.\par
\noindent Abel  6.\par
\noindent acceptable set of partitions  49.\par
\noindent additive  39.\par
\noindent additive$,$ simply 37.\par
\noindent algebraists$,$ information for  15$,$ 17.\par
\noindent algorithm  78.\par
\noindent almost parallel  91.\par
\noindent Apostal [1957]  35.\par
\noindent applications$,$ simple  23.\par
\noindent Archimedean  10.\par
\noindent Archimedes  6.\par
\noindent {\quad}\par
\noindent {\quad}\par
\noindent {\quad}\par
\noindent Barwise [l977]  21.\par
\noindent basic:\par
\sub {element  47.}
\sub {elementary integral  48.}
\sub {hyperfinite subsets  30.}
\sub {laboratory experimentation  46.}
\sub {region  39.}
\noindent bound formula  29.\par
\noindent boundary subrectangles 62.\par
\noindent boundedness concept  46.\par
\noindent bounding method  23.\par
\sub {redefinition  24.}
\noindent {\quad}\par
\noindent {\quad}\par
\noindent {\quad}\par
\noindent caloric  98.\par
\noindent Cauchy  6.\par
\sub {definition for the integral  35.}
\sub {Principle  88.}
\noindent Cavalieri's  36.\par
\noindent Cesari [1956]  71.\par
\noindent compactness  32.\par
\noindent compactness theorem  21.\par
\noindent comprehend second order rates of change 98.\par
\noindent conduction$,$ heat 98.\par
\noindent conductivity$,$ thermal 99.\par
\noindent consecutive points 84.\par
\noindent constant values$,$ extension fo 50.\par
\noindent {\quad}\par
\noindent {\quad}\par
\noindent constant$,$ use of  16.\par
\noindent constants:\par
\sub {elemental method of 64.}
\sub {extended standard 23.}
\sub {external 23.}
\sub {internal 20.}
\sub {internal 23.}
\sub {method of 50.}
\sub {self-evident method of 51.}
\sub {unstarred 23.}
\sub {use of  28.}
\noindent continuity  32.\par
\noindent continuous distribution of matter  49.\par
\noindent continuous$,$ uniformly 32.\par
\noindent convergence$,$ sequential 19.\par
\noindent convex in the direction $y$  99.\par
\noindent curvature$,$ radius of 87.\par
\noindent curve  31.\par
\sub {tangent to 71.}
\noindent curves$,$ hyperpolygonal 73.\par
\noindent Cutland [1986]  37.\par
\noindent {\quad}\par
\noindent {\quad}\par
\noindent {\quad}\par
\noindent d'Alembert  6.\par
\noindent d'Alembert - Euler  8.\par
\noindent D-world  14.\par
\noindent Darboux  35.\par
\sub {integral  35.}
\noindent de l'Hospital  6$,$ 26.\par
\noindent De Lillo [1982]  42.\par
\noindent definition:\par
\sub {for Euclidean Spaces  11.}
\sub {for infinite  11.}
\sub {for infinitely close  11.}
\sub {of hyperfinite  31.}
\sub {of infinitesimals  10.}
\sub {of limited  11.}
\sub {of monad  12.}
\noindent deformable body rule  93.\par
\noindent deleted monad  77.\par
\noindent $\delta$-fine partition  75.\par
\noindent derivative$,$ n-dimensional 90.\par
\noindent differentiable  90.\par
\noindent differential  89.\par
\noindent differential equation:\par
\sub {method of maximum and minimum  95.}
\sub {models  49.}
\noindent dot notation$,$ Newton's  81.\par
\noindent dynamical (loci) methods  78.\par
\noindent {\quad}\par
\noindent {\quad}\par
\noindent {\quad}\par
\noindent electron-positron lattice  95.\par
\noindent element:\par
\sub {basic 47.}
\sub {hypertrapezoid 60.}
\sub {infinitesimal 64.}
\sub {surface 73.}
\noindent elemental method of constants  64$,$ 76.\par
\noindent elementary:\par
\sub {geometric 46.}
\sub {prototype 47.}
\sub {integral 48. }
\noindent elements  48.\par
\sub {geometric 59$,$ 73.}
\sub {m-dimensional 47.}
\sub {method of  46$,$ 76.}
\noindent energy expended  67.\par
\noindent entity  17.\par
\sub {extended 23.}
\sub {internal 23.}
\noindent $\epsilon $-infinitesimal microscope 71.\par
\noindent equivalent maps 88.\par
\noindent Euclidean Spaces$,$ definition$,$ 11.\par
\noindent Euclidean n-spaces  14.\par
\noindent Eudoxus [370 BC]  27.\par
\noindent Euler  6.\par
\noindent evident  44.\par
\noindent Example:\par
\sub {3.3.1.  19.}
\sub {3.3.2.  19.}
\sub {3.3.3.  19.}
\sub {3.3.4.  20.}
\sub {3.4.1.  22.}
\sub {3.4.2.  22.}
\sub {3.4.3. An argument  23.}
\sub {4.2.1.  28.}
\sub {4.2.2.  28.}
\sub {4.2.3.  28.}
\sub {4.2.4.  28.}
\sub {4.3.1.  30.}
\sub {4.4.1.A. 33.}
\sub {4.4.1.B. 33.}
\sub {4.4.2. (Fractals) 33.}
\sub {*-transfer  28.}
\noindent extended:\par
\sub {entity  23.}
\sub {relations  17.}
\sub {standard  17.}
\sub {standard constants  23.}
\sub {of constant values  50.}
\sub {the finite  50.}
\noindent extensions of Jordan measures  39.\par
\noindent exterior subrectangles 62.\par
\noindent external constants  23.\par
\noindent {\quad}\par
\noindent {\quad}\par
\noindent {\quad}\par
\noindent fine  35.\par
\sub {partitions 33.}
\noindent finite  30.\par
\sub {extension of 50.}
\sub {summation  31.}
\noindent first-order:\par
\sub {$\epsilon $-resolving power.  71.}
\sub {ideals  69.}
\sub {language  19$,$ 20.}
\noindent flow lines$,$ postulated properties of  44.\par
\noindent flowing heat  98.\par
\noindent fluid motion:\par
\sub {lines of 44.}
\sub {tubes of 44.}
\noindent fluxions  78.\par
\noindent force  76.\par
\noindent formula$,$ bounded 29.\par
\noindent foundations of infinitesimal modeling  78.\par
\noindent fractals and infinite length  33.\par
\noindent free variable  28.\par
\noindent function$,$ limit of 77.\par
\noindent {\quad}\par
\noindent {\quad}\par
\noindent {\quad}\par
\noindent Galilean physics  80.\par
\noindent gauge integral  75.\par
\noindent Gauss 1827; Art 20  46.\par
\noindent geometric elements  59$,$ 73.\par
\sub {for the curve 59.}
\sub {elementary 46.}
\sub {surface element  73.}
\noindent geometry  76.\par
\noindent {\quad}\par
\noindent {\quad}\par
\noindent {\quad}\par
\noindent H-integral  35.\par
\noindent heat conduction  98.\par
\noindent heat:\par
\sub {flowing 98}
\sub {quantity of 98.}
\sub {specific 99.}
\sub {uniformly conducting 98.}
\noindent Henstock [1961]  75.\par
\noindent Herrmann 1980  46.\par
\sub {[1985]  35.}
\sub {[1986]  21.}
\sub {[1989]  33.}
\noindent higher order increments$,$ omission of 81.\par
\noindent Hurd and Loeb [1985]  65.\par
\noindent hyperfinite:\par
\sub {definition 31.}
\sub {partition  35.}
\sub {subsets$,$ basic 30.}
\sub {sum  31.}
\noindent hyperline segments  73.\par
\noindent hyperparallelepiped  74.\par
\noindent hyperpolygonal curves  73.\par
\noindent hyperreals$,$ names for 10.\par
\noindent hyperrectangle$,$ truncated 63.\par
\noindent hypersurface  73.\par
\noindent hypertangent planes  73.\par
\noindent hypertrapezoid element 60.\par
\noindent hypertrapezoids 73.\par
\noindent {\quad}\par
\noindent {\quad}\par
\noindent {\quad}\par
\noindent ideal  6$,$ 11.\par
\sub {maximal 15.}
\noindent idealized  behavior  58$,$ 76.\par
\noindent ideals$,$ order 68.\par
\noindent identity$,$ preserving 17.\par
\noindent {\quad}\par
\noindent {\quad}\par
\noindent {\quad}\par
\noindent {\quad}\par
\noindent inconsistencies  6.\par
\noindent increment:\par
\sub {local 90.}
\sub {{\it n}th order 70.}
\noindent increments$,$ higher order 81.\par
\noindent individual  17.\par
\sub {internal 23.}
\noindent indivisible  18$,$ 36.\par
\noindent induction  20.\par
\noindent inertia$,$ moment of 49.\par
\noindent infinite:\par
\sub {characterized 11.}
\sub {definition 11.}
\sub {length$,$ fractals 39.}
\sub {magnification operator  71.}
\sub {magnification$,$ of a infinitely small portion of the fluid 44.}
\sub {numbers$,$ different from extended real numbers 26.}
\sub {Sum Theorems 37.}
\noindent infinitely close 49.\par
\sub {definition 11. }
\sub {of order n  69$,$ 71.}
\noindent infinitely small 6.\par
\noindent infinitesimal:\par
\sub {approach$,$ pure 58.}
\sub {elements$,$ physical 64.}
\sub {max. and min. rule  48.}
\sub {microscope 71}
\sub {modeling$,$ foundations 78.}
\sub {parallelepiped 74. }
\sub {reasoning  49$,$ 51.}
\noindent infinitesimal:\par
\sub {rectangles  73.}
\noindent infinitesimalizing procedure 49.\par
\noindent infinitesimals:\par
\sub {and local map 88.}
\sub {definition 10.}
\sub {nonnegative 69.}
\noindent inner subrectangles 62.\par
\noindent integral:\par
\sub {Cauchy's definition 35.}
\sub {elementary 48.}
\sub {gauge 75.}
\sub {line 66.}
\sub {M 75.}
\sub {McShane 75.}
\sub {Riemann-complete 75.}
\sub {S 75.}
\sub {surface 73.}
\noindent intermediate partition 36.\par
\noindent internal  22.\par
\sub {constant names  20.}
\sub {constants  23.}
\sub {definition principle  30.}
\sub {entity  23.}
\sub {individual 23.}
\sub {objects and N-world modeling  30.}
\noindent {\quad}\par    
\noindent {\quad}\par
\noindent {\quad}\par
\noindent Jarnik$,$ Kurzweil$,$ Schwabik [1983] 75.\par
\noindent Jordan measures$,$ extended to 39. \par
\noindent Jordan-measurable  42$,$ 47$,$ 62.\par
\noindent {\quad}\par
\noindent {\quad}\par
\noindent {\quad}\par
\noindent $k$-plane  85.\par
\noindent $(k+1)$-consecutive  84.\par
\noindent Keisler  8.\par
\noindent Kepler  6.\par
\noindent {\quad}\par
\noindent {\quad}\par
\noindent {\quad}\par
\noindent laboratory experimentation$,$ basic 46.\par
\noindent lattice structure 95.\par
\noindent Leibniz  6$,$ 8.\par
\noindent Leibniz Principle  21.\par
\sub {restated  23.}
\noindent Leibniz [1701]  18.\par
\noindent length by fine partitions  33.\par
\sub {by limit of a sequence  33.}                                            
\noindent limit of a function  77.\par
\noindent limit theorems  19.\par
\noindent limited numbers definition 11.\par
\noindent line Integral  66.\par
\noindent linear$,$ locally 88.\par
\noindent lines of fluid motion  44.\par
\noindent local:\par
\sub {increment  90.}
\sub {map  48.}
\noindent {\quad}\par
\noindent {\quad}\par
\noindent {\quad}\par
\noindent {\quad}\par
\noindent locally linear  88.\par
\noindent locus  76.\par
\noindent Luxemburg [1962]  15.\par
\noindent {\quad}\par
\noindent {\quad}\par
\noindent {\quad}\par
\noindent m-dimensional elements  47.\par
\noindent M-integral  75.\par
\noindent Machover and Hirschfeld [1969]  15.\par
\noindent magnification:\par
\sub {operator 71}
\sub {infinite 44.}
\noindent map:\par
\sub {equivalent 88.}
\sub {local 88.}
\noindent matter$,$ continuous distribution of 49.\par
\noindent Mawhin [1985]  75.\par
\noindent max. and min. rule 48.\par
\sub {self-evident 49.}
\noindent maximal ideal  15.\par
\noindent max. and min.$,$ differential equation method of 95.\par
\noindent Maxwell [1890]  44.\par
\noindent McShane [1973]  75.\par
\noindent McShane integral 75.\par
\noindent mean value  39.\par
\noindent metamethematically$,$ the necessity of arguing 23.\par
\noindent method of constants  50.\par
\sub {elemental 76.}
\noindent method of elements 76.\par
\noindent microconstruction  68.\par
\noindent microeffects  68.\par
\noindent microguage  68.\par
\noindent micropartition  75.\par
\sub {regular 75.}
\noindent microscope$,$ infinitesimal 71.\par
\noindent modeling:\par
\sub {differential equation 49.}
\sub {infinitesimal 78.}
\sub {internal objects 30.}
\sub {N-world 30.}
\noindent moment 78.\par
\noindent moment of inertia 49.\par
\noindent monad:\par
\sub {definition 12.}
\sub {deleted 77.}
\noindent monadic:\par
\sub {restrictions of physical processes  94.}
\sub {second law of motion  92.}
\sub {second law of motion$,$ alternative  101.}
\noindent Morley's example 52.\par
\sub {[1942] 52.}
\noindent motion 76.\par
\noindent motion of points 76.\par
\noindent moving objects 76.\par
\noindent moving points$,$ tubes of 44.\par
\noindent {\quad}\par
\noindent {\quad}\par
\noindent {\quad}\par
\noindent {\it n}th difference  70.\par
\noindent {\it n}th order ideal 69.\par
\noindent {\it n}th order increment 70.\par
\noindent n-dimensional derivative 90.\par
\noindent N-world  14.\par
\noindent N-world modeling and internal objects 30.\par
\noindent names$,$ internal constant  20.\par
\noindent natural numbers$,$ properties of  20.\par
\noindent natural philosophy  76.\par
\noindent non-Archimedean  31.\par
\noindent non-finite  17.\par
\noindent nonnegative infinitesimals 69.\par
\noindent normal vector 73.\par
\noindent NSP-world  14.\par
\sub {behavior$,$ rules for pure 42.}
\noindent {\quad}\par
\noindent {\quad}\par
\noindent {\quad}\par
\noindent objective reality$,$ both the infinitesimal and infinite 26.\par
\noindent objects$,$ internal constant names  20.\par
\noindent observation$,$ physical 76.\par
\noindent order ideal$,$ {\it n}th 69.\par
\noindent order ideals  68.\par
\noindent order pair  16.\par
\noindent osculating:\par
\sub {$k$-circle  87.}
\sub {$k$-plane  84.}
\noindent {\quad}\par
\noindent {\quad}\par
\noindent {\quad}\par
\noindent parallel$,$ almost 91.\par
\noindent parallelepipeds$,$ 47.\par
\sub {infinitesimal 74.}
\noindent partition  35.\par
\sub {$\delta$ -fine 75.}
\sub {fine 33$,$ 35.}
\sub {hyperfinite 35.}
\noindent partitions$,$ acceptable set of 49.\par
\noindent Peano  20.\par
\noindent phenomenological approach 46.\par
\noindent philosophy$,$ natural 76.\par
\noindent physical:\par
\sub {behavior$,$ simplistic 46.}
\sub {infinitesimal elements 64.}
\sub {observation 76.}
\sub {processes$,$ monadic restrictions 94.}
\noindent Planck  6.\par
\sub {and meaning of physical series 34.}
\noindent plane$,$ tangent 73.\par
\noindent point definable quantities 50.\par
\noindent point-motion 76.\par
\noindent power set operator 17.\par
\noindent predicates in set builder notation 28.\par
\noindent preserves infinitesimals 88.\par
\noindent properties of standard part  operator  14.\par
\noindent prototype$,$ elementary 47.\par
\noindent pure:\par
\sub {infinitesimal approach  58.}
\sub {NSP-world behavior 42.}
\noindent {\quad}\par
\noindent {\quad}\par
\noindent {\quad}\par
\noindent quantifiers used in set builder notation  28.\par
\noindent quantity of heat  98.\par
\noindent quotient ring  15.\par
\noindent {\quad}\par
\noindent {\quad}\par
\noindent {\quad}\par
\noindent radius of curvature  87.\par
\noindent real numbers$,$ positive 10.\par
\noindent rectangle  35.\par
\sub {infinitesimal 73.}
\noindent region$,$ basic 39.\par
\noindent regular micropartition 75.\par
\noindent relation$,$ extended 17.\par
\noindent Riemann:\par  
\sub {integral 35.}
\sub {Stieltjes Integral 65.}
\sub {styled sum 58.}
\sub {sums 35.}
\noindent Riemann-complete integral 75.\par
\noindent ring 11.\par
\noindent ring$,$ quotient 15.\par
\noindent Robinson 8.\par
\sub {[1961]  10$,$ 26.}
\sub {[1966]  15.}
\noindent Robinson and Zakon [1969]  15.\par
\noindent rules  37.\par
\sub {for pure NSP-world behavior 42.}
\sub {of correspondence 76.}
\noindent {\quad}\par
\noindent {\quad}\par
\noindent {\quad}\par
\noindent S-integral 75.\par
\noindent second law of motion:\par
\sub {monadic 92.}
\sub {monadic alternative 101.}
\noindent second order rates of change$,$ comprehending 98.\par
\noindent segments$,$ hyperline 73.\par
\noindent self-evident  49.\par
\sub {max. and min.  49.}
\sub {method of constants 51.}
\noindent sequential convergence  19.\par
\noindent sets$,$ infinite 26.\par
\noindent simple:\par
\sub {behavior 58.}
\sub {fine partition 35.}
\noindent simplistic 76.\par
\sub {physical behavior 46.}
\noindent simply additive 37.\par
\noindent sketching 46.\par
\noindent small$,$ infinitely 6.\par
\noindent sovereign principle 18.\par
\noindent specific heat 99.\par
\noindent Spiegel 8.\par
\noindent Spivak [1965] 35.\par
\noindent * operator$,$ properties of 17.\par
\noindent *-transfer$,$ how two use  22.\par
\noindent {\quad}\par
\noindent {\quad}\par
\noindent *-transfer$,$ proving things by  22.\par
\noindent *-transfer 21.\par
\noindent st operator 14.\par
\noindent standard:\par
\sub {(unstarred) constants  23.}
\sub {constants$,$ extended 23.}
\sub {constants$,$ set of  20.}
\sub {part operator$,$ properties of 14.}
\sub {part operator$,$ using  25.}
\noindent steady state  98.\par
\noindent Stroyan and Luxemburg [1976]  21.\par
\noindent Struik [1961] 84.\par
\noindent subrectangle  37.\par
\sub {boundary 62.}
\noindent subrectangle:\par
\sub {exterior 62.}
\sub {inner 62.}
\noindent sum$,$ hyperfinite 31.\par
\noindent summation$,$ finite 31.\par
\noindent sums$,$ Riemann styled 58.\par
\noindent supernear  38.\par
\noindent supernearness  38.\par
\noindent superstructure  17.\par
\noindent surface:\par
\sub {element 73.}
\sub {integral  73.}
\sub {$\sigma $  73.}
\sub {tangent points 75.}
\noindent Swartz and Thomson [1988]  75.\par
\noindent symbol for a monad  12.\par
\noindent symbol for infinite natural numbers 19.\par
\noindent symbol for infinitely close 11.\par
\noindent symbol for limited  11.\par
\noindent symbol hyperreal 10.\par
\noindent symbol natural numbers 10.\par
\noindent symbol positive real numbers  10.\par
\noindent symbol real numbers 10.\par
\noindent symbols for infinitesimals 10.\par
\noindent Synge and Griffith [1959:\par 173]  49.\par
\noindent {\quad}\par
\noindent {\quad}\par
\noindent {\quad}\par
\noindent tangent:\par
\sub {plane 73.}
\sub {points to the surface 71.}
\sub {to a curve  71.}
\sub {vector to a curve  82.}
\noindent temperature  99.\par
\noindent The Integral  36.\par
\noindent thermal conductivity 99.\par
\noindent truncated hyperrectangular solid 63.\par
\noindent tubes:\par
\sub {of fluid motion 44.}
\sub {of moving points 44.}
\noindent types  15.\par
\noindent {\quad}\par
\noindent {\quad}\par
\noindent {\quad}\par
\noindent ultraproduct  21.\par
\noindent uniform:\par
\sub {rod  49.}
\sub {conducts heat 98.}
\noindent uniformly continuous  32.\par
\noindent unit normal vector to the curve  84.\par
\noindent unstarred constants 23.\par
\noindent {\quad}\par
\noindent {\quad}\par
\noindent {\quad}\par
\noindent variable$,$ free 28.\par
\noindent vector:\par
\sub {normal 73.}
\sub {tangent to a curve 82.}
\sub {unit normal 84.}
\noindent VR1  44.\par
\noindent VR2  44.\par
\noindent {\quad}\par
\noindent {\quad}\par
\noindent {\quad}\par
\noindent Weierstrass 6.\par
\noindent What if...?  25.\par
\noindent work 67.\par
\noindent {\quad}\par
\noindent {\quad}\par
\noindent {\quad}\par
\noindent Zeno  17.\par
\end